\documentclass[10pt]{article}

\usepackage[english]{babel}
\usepackage{inputenc}			
\usepackage[scaled]{helvet}		
\usepackage{courier}			
\usepackage{eulervm}			
\usepackage[bb=boondox,bbscaled=1.05,scr=dutchcal]{mathalfa}
\usepackage{dsfont}
\normalfont

\usepackage{indentfirst}				
\usepackage{tabularx,tabulary,booktabs}	
\usepackage{caption,subcaption}			
\usepackage{fullpage}					
\usepackage[english]{varioref}			
\usepackage[dvipsnames]{xcolor}			
\usepackage{hyperref,bookmark}			
\usepackage{textcomp}
\usepackage{rotating,afterpage,
	lscape}	

\usepackage{graphicx,float}				
\usepackage{tikz}						
\usepackage{tikz-cd}					
\usetikzlibrary{
	arrows,decorations.pathreplacing,decorations.markings,shapes.geometric,positioning,calc,fadings,angles,patterns
	}
\tikzfading[name=inner fade, inner color=transparent!0, outer color=transparent!100]
\graphicspath{{/}}
	
\usepackage{amsmath,amssymb,amsthm,mathtools}
\usepackage{bm,braket,faktor,esint}
\usepackage[export]{adjustbox}

\usepackage{enumerate,enumitem}
\numberwithin{equation}{section}
\usepackage[textsize=footnotesize]{todonotes}
\presetkeys%
		{todonotes}%
		{color=Dandelion}{}%

\newcommand{\de}{\partial}
\newcommand{\dd}{\mathrm{d}}

\newcommand{\NN}{\mathbb{N}}
\newcommand{\ZZ}{\mathbb{Z}}
\newcommand{\QQ}{\mathbb{Q}}
\newcommand{\RR}{\mathbb{R}}

\newcommand{\GG}{\mathbb{G}}

\renewcommand{\emptyset}{\varnothing}

\DeclareMathOperator{\image}{Im}

\newcommand{\Aut}{{\rm Aut}}

\DeclareMathOperator{\Stab}{Stab}

\DeclareMathOperator{\id}{id}

\DeclareMathOperator{\sgn}{sgn}

\DeclareMathOperator*{\Res}{Res}

\DeclareMathOperator{\sh}{sinh}
\DeclareMathOperator{\ch}{cosh}

\DeclareMathOperator{\arch}{arcosh}
\DeclareMathOperator{\cotanh}{coth}

\DeclareMathOperator{\sys}{sys}
\renewcommand{\sp}{{\rm sp}}

\DeclareMathOperator{\Mod}{Mod}

\theoremstyle{plain}
\newtheorem{thm}{Theorem}[section]
\newtheorem{theoremA}{Theorem}
\newtheorem{corollaryA}[theoremA]{Corollary}

\newtheorem{theoremB}{Theorem}

\newtheorem{theoremC}{Theorem}

\newtheorem{theoremD}{Theorem}

\newtheorem{prop}[thm]{Proposition}

\newtheorem{lem}[thm]{Lemma}
\newtheorem{cor}[thm]{Corollary}

\theoremstyle{definition}
\newtheorem{defn}[thm]{Definition}
\newtheorem{rem}[thm]{Remark}
\newtheorem{ex}[thm]{Example}

\makeatletter
\renewenvironment{abstract}{%
	\if@twocolumn
		\section*{\abstractname}%
	\else \normalsize
		\begin{center}%
			{\bfseries \normalsize\abstractname\vspace{\z@}}
		\end{center}%
		\quotation
	\fi}
	{\if@twocolumn\else\endquotation\fi}
\makeatother
	
\definecolor{webbrown}{rgb}{0.65, 0.16, 0.16}	
\hypersetup{
colorlinks=true, linktocpage=true, pdfstartpage=1, pdfstartview=FitV,breaklinks=true, pdfpagemode=UseNone, pageanchor=true, pdfpagemode=UseOutlines,plainpages=false, bookmarksnumbered, bookmarksopen=true, bookmarksopenlevel=1,hypertexnames=true, pdfhighlight=/O,urlcolor=webbrown, linkcolor=RoyalBlue, citecolor=ForestGreen}	

\title{\textsc{
		On the Kontsevich geometry of the combinatorial Teichm\"uller space
}\vspace{.3cm}}
\author{
	\setcounter{footnote}{6}
	J\o{}rgen Ellegaard Andersen\footnote{Center for Quantum Mathematics, Danish Institute for Advanced Study, SDU, Campusvej 55, 5230 Odense, Denmark.}\,\,,
	\setcounter{footnote}{0}
	Ga\"etan Borot \footnote{Max Planck Institut f\"ur Mathematik, Vivatsgasse 7, 53111 Bonn, Germany }\,\,\footnote{Institut f\"ur Mathematik und Institut f\"ur Physik, Humboldt-Universit\"at zu Berlin, Rudower Chaussee 25, 12489 Berlin, Germany.}\,\,,
	S\'everin Charbonnier\footnotemark[1],\\
	Alessandro Giacchetto\footnotemark[1]\,\,,
	Danilo Lewa\'{n}ski \footnotemark[1]\,\,\footnote{Universit\'e Paris-Saclay, CNRS, CEA, Institut de Physique Th\'{e}orique, 91191 Gif-sur-Yvette, 
France.}\, \footnote{Institut des Hautes \'Etudes Scientifiques, le Bois-Marie, 35 route de Chartres, 91440 Bures-sur-Yvette, France.}
	\,\,,
	Campbell Wheeler\footnotemark[1]
}
\date{}

\begin{document}

\maketitle 

\begin{abstract}

For bordered surfaces $\Sigma$, we develop a complete parallel between the geometry of the combinatorial Teichm\"uller space $\mathcal{T}_{\Sigma}^{\textup{comb}}$ equipped with Kontsevich symplectic form $\omega_{\textup{K}}$, and then the usual Weil--Petersson geometry of Teichm\"uller space. The basis for this is an identification of $\mathcal{T}_{\Sigma}^{\textup{comb}}$ with a space of measured foliations with transverse boundary conditions. We equip $\mathcal{T}_{\Sigma}^{\textup{comb}}$ with an analog of the Fenchel--Nielsen coordinates (defined similarly as Dehn--Thurston coordinates) and show they are Darboux for $\omega_{\textup{K}}$ (analog of Wolpert formula).

We then set up the geometric recursion of Andersen--Borot--Orantin to produce mapping class group invariants functions on $\mathcal{T}_{\Sigma}^{\textup{comb}}$ whose integration with respect to Kontsevich volume form satisfy topological recursion. Further we establish an analog of Mirzakhani--McShane identities, and provide applications to the study of the enumeration of multicurves with respect to combinatorial lengths and Masur--Veech volumes. 

The formalism allows us to provide uniform and completely geometric proofs of Witten's conjecture/Kontsevich theorem and Norbury's topological recursion for lattice point count in the combinatorial moduli space, parallel to Mirzakhani's proof of her recursion for Weil--Petersson volumes. We strengthen results of Mondello and Do on the convergence of hyperbolic geometry to combinatorial geometry along the rescaling flow, allowing us to flow systematically natural constructions on the usual Teichm\"uller space to their combinatorial analogue, such as a new derivation of the piecewise linear structure of $\mathcal{T}_{\Sigma}^{\textup{comb}}$ originally obtained in the work of Penner, as the limit under the flow of the smooth structure of $\mathcal{T}_{\Sigma}$. 

\end{abstract}


\newpage
\tableofcontents
\newpage

\thispagestyle{empty}
\bigskip

\section{Introduction}
\label{sec:intro}

\subsection{Overview}

In his proof of Witten's conjecture and following the pioneering works of 
Strebel, Harer, Mumford, Penner and Thurston \cite{Har86,HZ86,Pen88,Str67}, Kontsevich \cite{Kon92} studied a combinatorial model $\mathcal{M}_{g,n}^{\textup{comb}}$ for the moduli space of curves, namely the space of metric ribbon graphs on a surface $\Sigma$ of genus $g$ with $n > 0$ boundaries. He equipped it with a $2$-form $\omega_{\textup{K}}$ whose restriction to the slice with fixed boundary lengths is almost everywhere symplectic, and so that the symplectic volume gives access to the intersection of $\psi$-classes on the Deligne--Mumford moduli space of curves $\overline{\mathcal{M}}_{g,n}$. The corresponding \emph{combinatorial Teichm\"uller space}
\[
	\mathcal{T}^{\textup{comb}}_{\Sigma}(L) =
	\left.
	\Set{\begin{array}{@{}c@{}}
			\text{embedded metric ribbon graphs of genus $g$ with $n$ faces} \\
			\text{of lengths $L = (L_1,\dots,L_n)$ that  are retract of $\Sigma$}
		\end{array}
	} \middle/ \text{isotopy} \right.
\]
was first considered in \cite{Mon04} and its topology studied via the arc complex. Its quotient by the pure mapping class group is $\mathcal{M}_{g,n}^{\textup{comb}}$, and we can equip $\mathcal{T}_{\Sigma}^{\textup{comb}}$ with the pullback of Kontsevich $2$-form, which we still denote $\omega_{\textup{K}}$. 

\medskip

The first purpose of this work is to develop the geometry of $(\mathcal{T}^{\textup{comb}}_{\Sigma},\omega_{\textup{K}})$ in perfect parallel to the Weil--Petersson geometry of the Teichm\"uller space $\mathcal{T}_{\Sigma}$ of hyperbolic metrics on $\Sigma$ with geodesic boundaries. We shall use the name ``combinatorial geometry" exclusively in reference to the former. In particular, we will see that there exists on $\mathcal{T}_{\Sigma}^{\textup{comb}}$ an analogue of Fenchel--Nielsen coordinates, an analogue of Wolpert's formula showing that these coordinates are Darboux for $\omega_{\textup{K}}$, and an analogue of the Mirzakhani--McShane identity. We also set up the geometric recursion of \cite{ABO17} to produce mapping class group invariant functions on $\mathcal{T}_{\Sigma}^{\textup{comb}}$ by a cut-and-paste approach and show that their integration against the Kontsevich volume form on $\mathcal{M}_{g,n}^{\textup{comb}}(L)$ satisfy the Eynard--Orantin topological recursion.

\medskip

In a second part, we combine these results to treat in a uniform way various enumerative problems, parallel to Mirzakhani's proof of her recursion for Weil--Petersson volumes \cite{Mir07simple}. In particular, we describe a fully geometric proof of the Virasoro part of Witten's conjecture \cite{DVV91,Kon92} (bypassing matrix model techniques), obtain a new proof of Norbury's topological recursion \cite{Nor10} for the count of lattice points in $\mathcal{M}_{g,n}^{\textup{comb}}(L)$ (by discrete integration of the combinatorial Mirzakhani--McShane identity), and study the enumeration of multicurves with respect to their combinatorial lengths, their asymptotics and their relation to Masur--Veech volumes of the top stratum of the moduli space of quadratic differentials.

\medskip

In a third part, we study the rescaling flow on the Teichm\"uller space interpolating between hyperbolic geometry and combinatorial geometry. The spine construction originating from the works of Penner \cite{Pen87}, Bowditch--Epstein \cite{BE88} and Luo \cite{Luo07} provides a natural homeomorphism $\sp \colon \mathcal{T}_{\Sigma} \rightarrow \mathcal{T}_{\Sigma}^{\textup{comb}}$ together with a rescaling flow $\sigma \mapsto \sigma^{\beta}$ on $\mathcal{T}_{\Sigma}$ which consists in rescaling the lengths of all edges in the spine by a factor $\beta > 0$, \emph{i.e.} $\sigma^{\beta} = \sp^{-1}(\beta\cdot \sp(\sigma))$. Mondello \cite{Mon09} and Do \cite{Do10} have shown that $(\Sigma,\beta^{-1}\sigma^{\beta})$ converges in the Gromov--Hausdorff sense to the metric graph $\sp(\sigma)$, and that lengths of simple closed curves for $\beta^{-1}\sigma^{\beta}$ converge to combinatorial lengths, pointwise in $\mathcal{T}_{\Sigma}$, and that the suitably rescaled Weil--Petersson Poisson structure converges to the Kontsevich Poisson structure. We strengthen these results by establishing a uniform convergence of lengths and twists for $\sigma^{\beta}$ to their combinatorial analog for $\sp(\sigma)$. This allows us to study the convergence of various natural functions on $\mathcal{T}_{\Sigma}$ to functions on $\mathcal{T}_{\Sigma}^{\textup{comb}}$, in particular the functions obtained by geometric recursion, and Mirzakhani's function describing the constant prefactor in the asymptotic enumeration of large multicurves. We also demonstrate that the Masur--Veech volumes can be approached from the perspective of combinatorial geometry.

\subsection{Context: measured foliations and geometry}

Since Thurston, measured foliations play an important role in the study of the geometry of the Teichm\"uller space, both in the hyperbolic and in the combinatorial case, and we also take this perspective. Before going further, it is perhaps useful to situate our work among the various settings that have already been extensively studied. There are in fact at least three different spaces that have been introduced in relation with Teichm\"uller theory.
\begin{itemize}
	\item
	If $\Sigma$ is a surface which is either closed, punctured, or has boundaries, ${\rm MF}_{\Sigma}$ is Thurston's space of measured foliations in which leaves do not go into the punctures or are included in boundaries (when $\Sigma$ has respectively punctures or boundaries). This space is the topological completion of the set of multicurves with rational coefficients; it has dimension $6g - 6 + 2n$ and carries Thurston's symplectic form. Given a pair of pants decomposition it admits Dehn--Thurston coordinates, which are global coordinates containing lengths and twists parameters \cite{Pen82,Pen84}. For closed surfaces, the twist flow has been related to the Hamiltonian flow of the length functions in \cite{Papado1,Papado2}.

	\item
	If $\Sigma$ is a surface with punctures, the decorated Teichm\"uller space $\widetilde{\mathcal{T}}_{\Sigma}$ introduced by Penner \cite{Pen87,Penbook} is the trivial $\mathbb{R}_{> 0}^n$-bundle over $\mathcal{T}_{\Sigma}$. A tuple $L \in \mathbb{R}_{> 0}^n$ equivalently parametrises horocycles around the punctures, and we denote by $\widetilde{\mathcal{T}}_{\Sigma}(L)$ the fiber over $L$. It can be equipped with the pullback of the Weil--Petersson symplectic form. For bordered surfaces, two equivalent description of the decorated Teichm\"uller have been studied in \cite{Penbord}.

	\item
	If $\Sigma$ is a surface with non-empty boundary, $\mathcal{T}_{\Sigma}(L)$ is the Teichm\"uller space of hyperbolic structures with geodesic boundary components of length $L$. It can be equipped with the Weil--Petersson symplectic form.

	\item
	In this article we will study $\mathcal{T}_{\Sigma}^{\textup{comb}}$ which is defined when $\Sigma$ is a surface with non-empty boundary. Combinatorial structures (\emph{i.e.} elements of this space) can be identified with equivalence classes of measured foliations having no saddle connections and whose leaves are transverse to $\partial \Sigma$. Such a description is already manifest from the arc complex perspective \cite{Luo07,Mon09} which is reviewed in Appendix~\ref{app:topology:comb}. We denote $\mathcal{T}_{\Sigma}^{\textup{comb}}(L)$ the slice of $\mathcal{T}_{\Sigma}^{\textup{comb}}$ in which boundaries have combinatorial lengths $L \in \mathbb{R}_{+}^n$. Although this space differs from ${\rm MF}_{\Sigma}$, Dehn--Thurston coordinates can be defined in an identical way, and coincide with what we call combinatorial Fenchel--Nielsen coordinates. We will keep this name to stress the parallel with hyperbolic geometry that our work reinforces.
\end{itemize}
It should be stressed that $\widetilde{\mathcal{T}}_{\Sigma}(L)$ and $\mathcal{T}_{\Sigma}^{\textup{comb}}(L)$, albeit homeomorphic, are not equipped with the same symplectic form (Weil--Petersson in the first case, Kontsevich in the second case), neither do these symplectic forms agree when identifying the two spaces in a natural way. They can only be compared asymptotically using rescaling flows.
\begin{itemize}
	\item
	If $\Sigma$ is a punctured surface, building on the fact that ${\rm MF}_{\Sigma}$ appears as Thurston's boundary of $\mathcal{T}_{\Sigma}$ and is naturally realised as $\widetilde{\mathcal{T}}_{\Sigma}(\beta L)$ when $\beta \rightarrow 0$, Papadopoulos and Penner \cite{PP93} showed that the Weil--Petersson form on $\widetilde{\mathcal{T}}_{\Sigma}$ approaches Thurston symplectic form in this limit.

	\item
	If $\Sigma$ is a bordered surface, the Weil--Petersson form on $\mathcal{T}_{\Sigma}(\beta L)$ approaches Kontsevich $2$-form on $\mathcal{T}_{\Sigma}^{\textup{comb}}$ when $\beta \rightarrow \infty$ \cite{Mon09,Do10}, see Section~\ref{sec:hyp:comb} for a precise statement.
\end{itemize}

In our study of $\mathcal{T}_{\Sigma}^{\textup{comb}}(L)$, several techniques from measured foliations \cite{FLP12,PH92} are similar to those used in the study of the decorated Teichm\"uller space, but adapted to our setting. For clarity, this paper includes in Section~\ref{sec:comb:spaces} a self-contained description of the topological aspects of $\mathcal{T}_{\Sigma}^{\textup{comb}}$ and the construction of combinatorial Fenchel--Nielsen coordinates, without pretension of originality but in a form that is useful later, for instance when we prove a combinatorial analogue of Mirzakhani--McShane identity. We however prove two properties, which do not seem to be known and are the basis for the second and third part of the paper:
\begin{itemize}
	\item the image of the combinatorial Fenchel--Nielsen coordinates in $(\mathbb{R}_+ \times \mathbb{R})^{3g - 3 + n}$ is (open) dense, with complement of zero measure,
	\item the combinatorial Fenchel--Nielsen are Darboux coordinates for $\omega_{\textup{K}}$.
\end{itemize}

Penner's thesis \cite{Pen82} described the transformation of Dehn--Thurston coordinates on ${\rm MF}_{\Sigma}$ under change of pairs of pants decomposition and accordingly justified it has a piecewise linear integral structure. His results imply the transformation laws for combinatorial Fenchel--Nielsen coordinates, since they are a straightforward adaptation of Dehn--Thurston coordinates to a different type of measured foliations. In the third part of this work, we will give a different proof of this result, by flowing the transformations for hyperbolic Fenchel--Nielsen coordinates with the rescaling flow and using the uniform convergence to combinatorial Fenchel--Nielsen that we establish.

\subsection{Detailed summary of results}

\paragraph{Combinatorial spaces and their geometry.}
Many aspects of the geometry of hyperbolic structures have an analogue for combinatorial structures. First of all, one can define the length $\ell_{\GG}(\gamma)$ of a simple closed curve $\gamma$ with respect to a combinatorial structure $\GG$: realising the latter as a measured foliation $\mathcal{F}_{\GG}$, this is Thurston's intersection pairing between $\mathcal{F}_{\GG}$ and $\gamma$.

\medskip

Thanks to the notion of length, one can parametrise the combinatorial Teichmüller space using a maximal set of simple closed curves, \emph{i.e.} a pants decomposition. On a surface $\Sigma$ of genus $g$ with $n$ boundary components, there are $3g - 3 + n$ length parameters that determine the combinatorial structure on each pair of pants, and there are $3g - 3 + n$ twist parameters that determine how the pairs of pants are glued together. All together, they constitute a combinatorial analogue of Fenchel--Nielsen coordinates, and they are known in the measured foliation and train track settings as Dehn--Thurston coordinates (cf. \cite[Exposé~6]{FLP12} and \cite[Theorem~3.1.1]{PH92}). The main difference with the hyperbolic world lies in the twist $\tau$: for some values of $\tau$, it is not possible to glue combinatorial structures; in the measured foliation description this corresponds to the creation of saddle connections that cannot be removed by Whitehead equivalences, and are not allowed in the combinatorial Teichmüller space. Metrically, they would in fact correspond to nodal surfaces. However, we show that the set of twists for which we cannot perform the gluing is a countable subset of $\RR$ with open dense complement. We give here a concise form of the stronger Theorem~\ref{thm:FN:coordinates}.

\begin{theoremA}\label{thm:A1}~
\begin{enumerate}
	\item (Adaptation of Dehn--Thurston--Penner arguments \cite{Deh22,FLP12,PH92}). Given a seamed pants decomposition for a connected bordered surface $\Sigma$ of genus $g$ with $n$ boundary components, we have an open map
	\begin{align*}
		\mathcal{T}^{\textup{comb}}_{\Sigma}(L) & \longrightarrow \RR_{+}^{3g-3+n} \times \RR^{3g-3+n} \\
		\GG & \longmapsto \bigl( \ell(\GG), \tau(\GG) \bigr)
	\end{align*}
	that is a homeomorphism onto its image.

	\item The image of $\mathcal{T}^{\textup{comb}}_{\Sigma}(L)$ inside $\RR_{+}^{3g-3+n} \times \RR^{3g-3+n}$ has a complement of zero measure.
\end{enumerate}
\end{theoremA}

The second part of this theorem, which to the best of our knowledge is new, will be crucial in the proof of Theorem~\ref{thm:B2}: the integration of the combinatorial Mirzakhani identity produces a recursion relation between volumes of $\mathcal{M}^{\textup{comb}}_{g,n}(L)$.

\medskip

Like in the hyperbolic case, we show that the twist parameters can be recovered from the data of $9g - 9 + 3n$ lengths of simple closed curves, determined by a fixed seamed pants decomposition (Theorem~\ref{thm:9g:minus:9:plus:3n} in the text). Again, this is an adaptation of \cite{FLP12} to 
the different type of measured foliations required in the combinatorial Teichmüller space (Section~\ref{subsec:9g:minus:9:plus:3n}). This result is used in the derivation of Theorems~\ref{thm:C1} and \ref{thm:C2}.

\medskip

In Section~\ref{subsec:Kont:form} we recall the definition of Kontsevich $2$-form $\omega_{\textup{K}}$ \cite{Kon92}, which is symplectic on the top-di\-men\-sio\-nal cells of $\mathcal{M}^{\textup{comb}}_{g,n}(L)$. The volume of $\mathcal{M}^{\textup{comb}}_{g,n}(L)$ with respect to $\dd\mu_{\textup{K}} = \omega_{\textup{K}}^{3g-3+n}/(3g-3+n)!$, denoted
\[
	V^{\textup{K}}_{g,n}(L) = \int_{\mathcal{M}^{\textup{comb}}_{g,n}(L)}\dd\mu_{\textup{K}},
\]
is finite. As functions of $L \in \RR_{+}^n$, they are polynomials whose coefficients compute $\psi$-classes intersections on the Deligne--Mumford 
compactification of the moduli space of curves, see \cite{Kon92} completed by \cite{Zvo02}. After we lift $\omega_{\textup{K}}$ to a mapping class 
group invariant $2$-form on $\mathcal{T}_{\Sigma}^{\textup{comb}}$, the main result of Section~\ref{sec:symplectic:structure} is a combinatorial analogue of Wolpert's formula \cite{Wol85}, expressing Kontsevich's form in terms of combinatorial Fenchel--Nielsen coordinates (Theorem~\ref{thm:Wolpert} in the text).

\begin{theoremA}\label{thm:A2}
	Let $\Sigma$ be a connected bordered surface of genus $g$ with $n$ boundary components, and fix any combinatorial Fenchel--Nielsen coordinates $(\ell_{i},\tau_{i})$ for $\mathcal{T}_{\Sigma}^{\textup{comb}}(L)$. Then
	\[
		\omega_{\textup{K}}
		=
		\sum_{i=1}^{3g-3+n} \dd\ell_{i} \wedge \dd\tau_{i}.
	\]
\end{theoremA}

We note that another, \emph{a priori} different set of Darboux coordinates for $\omega_{\textup{K}}$ have been constructed by Bertola and Korotkin 
\cite{BK18} from periods of quadratic differentials. An advantage of our combinatorial Fenchel--Nielsen coordinates is their compatibility with the cutting and gluing operations. Together with Theorem~\ref{thm:A1}, this 
enables us to deduce from Theorem~\ref{thm:A2} an integration result for mapping class group invariant functions with respect to the measure $\mu_{\textup{K}}$ (Section~\ref{subsec:integration}), which is the analogue of Mirzakhani's integration lemma \cite{Mir07simple} known in the hyperbolic context.

\medskip

We already mentioned that the definition of combinatorial Fenchel---Nielsen coordinates in $\mathcal{T}_{\Sigma}^{\textup{comb}}$ is identical to the definition of Dehn--Thurston on ${\rm MF}_{\Sigma}$. Besides, we observe that Kontsevich $2$-form on $\mathcal{T}_{\Sigma}^{\textup{comb}}$ is defined identically to Thurston symplectic form on ${\rm MF}_{\Sigma}$, compare e.g. with \cite[Section 3]{Bonahon93} in which one should consider the train track dual to the trivalent ribbon graph: switches correspond to corners of the ribbon graph and intersecting transverse arcs correspond to edges meeting at a vertex. Adapting the proof of Theorem~\ref{thm:A2} to ${\rm MF}_{\Sigma}$ therefore leads to the following result, which seems to be unnoticed in the literature --- to the best of our knowledge and also to our surprise.

\begin{corollaryA}
If $\Sigma$ is a punctured surface, the Dehn--Thurston coordinates are almost everywhere Darboux coordinates for Thurston symplectic form on ${\rm MF}_{\Sigma}$. \hfill $\blacksquare$
\end{corollaryA}

\paragraph{Geometric recursion (GR) and topological recursion (TR).} 
The Weil--Petersson volumes of $\mathcal{M}_{g,n}(L)$, the Kontsevich volumes, and the counting of lattice points in $\mathcal{M}_{g,n}^{\textup{comb}}(L)$ all satisfy topological recursion in the sense of \cite{EO07inv}. In particular, they can be computed recursively in $2g - 2 + n$. We demonstrate that this type of recursive relations can be rather systematically obtained from recursions at the level of functions (hence, before integration) on the Teichm\"uller space. More precisely, building on the framework of geometric recursion in the sense of \cite{ABO17} (which in the original paper has been applied to $\mathcal{T}_{\Sigma}$), we set up in Section~\ref{sec:setting:up:gr} the geometric recursion to construct mapping class group invariant functions on $\mathcal{T}_{\Sigma}^{\textup{comb}}$. Let us denote by $P$ a pair of pants and by $T$ a torus with one boundary component. The following result is Theorem~\ref{thm:comb:GR} in the text.
 
\begin{theoremB}[\emph{Combinatorial GR is well-defined}]\label{thm:B1}
	Let $A,B,C$ be measurable functions on $\mathbb{R}_{+}^3$ with $A$ and $C$ symmetric under exchange of their last two variables, $D_{T}$ be a measurable function on $\mathcal{T}_{T}^{\textup{comb}}$ which is mapping class group invariant, and assume they obey the bounds of Definition~\ref{defn:adm:comb}. Then, the following definitions are well-posed, and assign functorially to any bordered surface $\Sigma$ a measurable function $\Xi_{\Sigma}$ on $\mathcal{T}_{\Sigma}^{\textup{comb}}$, called the GR amplitude.
	\begin{itemize}
		\item[$\bullet$] $\Xi_{P}(\GG) = A(\vec{\ell}_{\GG}(\partial P))$, where $\vec{\ell}_{\GG}(\partial P)$ is the triple of combinatorial boundary lengths of $P$.
		\item[$\bullet$] $\Xi_{T} = D_{T}$.
		\item[$\bullet$] If $\Sigma$ is a disjoint union of $\Sigma_1,\ldots,\Sigma_k$, $\Xi_{\Sigma_1 \sqcup \cdots \sqcup \Sigma_k}(\GG_1,\ldots,\GG_k) = \prod_{i = 1}^k \Xi_{\Sigma_i}(\GG_i)$.
		\item[$\bullet$] If $\Sigma$ is connected and has Euler characteristic $\chi_{\Sigma} \leq -2$, 
		\[
			\Xi_{\Sigma}(\GG)
			=
			\sum_{m = 2}^n \sum_{[P] \in \mathcal{P}_{\Sigma}^{m}}
				B(\vec{\ell}_{\GG}(\partial P)) \, \Xi_{\Sigma - P}(\GG|_{\Sigma - P})
			+
			\frac{1}{2} \sum_{[P] \in \mathcal{P}_{\Sigma}^{\varnothing}}
				C(\vec{\ell}_{\GG}(\partial P)) \, \Xi_{\Sigma - P}(\GG|_{\Sigma - P})
		\]
		where $\mathcal{P}_{\Sigma}^{m}$ and $\mathcal{P}_{\Sigma}^{\varnothing}$ (\emph{cf.} Definition~\ref{defn:homotopy:class:pop}) are certain sets of homotopy classes of embedded pairs of pants in $\Sigma$, such that $\Sigma - P$ is stable, and $\GG|_{\Sigma - P}$ is the result of cutting the combinatorial structure $\GG$ and restricting it to $\Sigma - P$.
	\end{itemize}
	Further, the function $\Xi_{\Sigma}$ is invariant under mapping classes of $\Sigma$ preserving $\partial_1\Sigma$.
\end{theoremB}

By means of Theorem~\ref{thm:A2}, we show that integrating GR amplitudes automatically yields functions of boundary lengths that satisfy topological recursion. If $\Sigma$ is a connected bordered surface of type $(g,n)$, let us denote by $\Xi_{g,n}$ the function induced by $\Xi_{\Sigma}$ on $\mathcal{M}_{g,n}^{\textup{comb}}$. The following result is Theorem~\ref{thm:comb:GR:TR} in the text.

\begin{theoremB}[\emph{TR from GR}]\label{thm:B2}
	If $A,B,C$ are measurable functions on $\RR_{+}^3$ and $D_{T}$ is a measurable function on $\mathcal{T}_{T}^{\textup{comb}}$ satisfying the conditions of Definition~\ref{defn:strong:adm}, then the integrals 
	\[
		V\Xi_{g,n}(L) = \int_{\mathcal{M}_{g,n}^{\textup{comb}}(L)} \Xi_{g,n}\,\dd\mu_{\textup{K}}
	\]
	exist, define measurable functions on $L \in \RR_{+}^n$, and we have that
	\[
	\begin{split}
		& \quad V\Xi_{g,n}(L_1,\ldots,L_n) =  \\
		& =
			\sum_{m = 2}^n \int_{\mathbb{R}_{+}} \dd \ell \, \ell \, B(L_1,L_m,\ell) \, V\Xi_{g,n - 1}(\ell,L_2,\ldots,\widehat{L_m},\ldots,L_n) \\
		& \quad + \frac{1}{2} \int_{\mathbb{R}_{+}^2} \dd \ell \dd\ell' \, \ell \ell' \, C(L_1,\ell,\ell')\bigg(
			V\Xi_{g - 1,n + 1}(\ell,\ell',L_2,\ldots,L_n)
			+
			\!\!\!\!\! \sum_{\substack{h + h' = g \\ J \sqcup J' = \{L_2,\ldots,L_n\}}} \!\!\!\!\! V\Xi_{h,1 + \#J}(\ell,J) \, V\Xi_{h',1+\#J'}(\ell',J')
		\bigg).
	\end{split}
	\vspace{-\baselineskip}
	\]
	where by convention $V\Xi_{0,1} = 0$ and $V\Xi_{0,2} = 0$.
\end{theoremB}

A similar result holds if, instead of integrating against $\mu_{\textup{K}}$, we sum over the lattice $\mathcal{M}_{g,n}^{\textup{comb},\ZZ}(L) \subset \mathcal{M}_{g,n}^{\textup{comb}}(L)$, which consists of metric ribbon graphs with integral edge lengths. This has no counterpart in the hyperbolic world. Due to the existence of pathological twists for the gluing 
--- although they are rare in the whole space, they could (and in fact sometimes do) hit the lattice --- this is only possible under extra conditions for the initial data $B$ and $C$.

\begin{theoremB}[\emph{Discrete TR from GR}]\label{thm:B3}
	Let $A,B,C$ be functions on $\mathbb{R}_{+}^3$ and $D_{T}$ be a function on $\mathcal{T}_{T}^{\textup{comb}}$ such that
	\[
		L_1 + L_2 < \ell \quad \Rightarrow \quad B(L_1,L_2,\ell) = 0,
		\qquad \text{and} \qquad
		L_1 < \ell + \ell' \quad \Rightarrow\quad C(L_1,\ell,\ell') = 0.
	\]
	For all $2g - 2 + n > 0$, the lattice sum
	\[
		N\Xi_{g,n}(L) = \sum_{\bm{G} \in \mathcal{M}_{g,n}^{\textup{comb},\ZZ}(L)} \frac{\Xi_{g,n}(\bm{G})}{
	\#\Aut(\bm{G})}
	\]
	defines a function of $L \in \ZZ_{+}^n$ which is zero whenever $\sum_{i = 1}^n L_i$ is odd, and otherwise satisfies the recursion
	\begin{equation*}
	\begin{split}
		& \quad N\Xi_{g,n}(L_1,\ldots,L_n) = \\ 
		& = \sum_{m = 2}^n \sum_{\ell > 0} \ell \, B_{\ZZ}(L_1,L_m,\ell) \, N\Xi_{g,n - 1}(\ell,L_2,\ldots,\widehat{L}_m,\ldots,L_n) \\
		& \quad + \frac{1}{2} \sum_{\ell,\ell' > 0} \ell \ell' \, C_{\ZZ}(L_1,\ell',\ell') \bigg(
			N\Xi_{g-1,n+1}(\ell,\ell',L_2,\ldots,L_n)
			+
			\!\!\!\!\!\sum_{\substack{h + h ' = g \\ J \sqcup J' = \{L_2,\ldots,L_n\}}} \!\!\!\!\! N\Xi_{h,1 + \#J}(\ell,J) \, N\Xi_{h',1+\#J'}(\ell',J')
		\bigg),
	\end{split}
	\vspace{-\baselineskip}
	\end{equation*}
	where $X_{\ZZ}(L_1,L_2,L_3)$ is equal to $X(L_1,L_2,L_3)$ if $L_1 + L_2 + L_3$ is an even integer and $0$ otherwise.
\end{theoremB}

As applications of this more general theory we can re-prove known results 
in a completely geometric and uniform way, as well as obtaining new results. A key role for applications is played by the combinatorial analogue of the Mirzakhani--McShane identity, whose proof transposes the original strategy of Mirzakhani \cite{Mir07simple} to the combinatorial world (where it is much simpler).

\begin{theoremB}\label{thm:B4}
	Denote $[x]_{+} = \max\set{x,0}$ and define the Kontsevich initial data
	\[
	\begin{split}
		A^{\textup{K}}(L_1,L_2,L_3) & = 1, \\
		B^{\textup{K}}(L_1,L_2,\ell) & = \frac{1}{2L_1}\big( [L_1 - L_2 - \ell]_{+} - [-L_1 + L_2 - \ell]_{+} + [L_1 + L_2 - \ell]_{+} \big), \\
		C^{\textup{K}}(L_1,\ell,\ell') & = \frac{1}{L_1} [L_1 - \ell - \ell']_{+}, \\
		D^{\textup{K}}_{T}(\GG) & = \sum_{\gamma \in S_{T}^{\circ}} C^{\textup{K}}\big( \ell_{\GG}(\partial T),\ell_{\GG}(\gamma),\ell_{\GG}(\gamma) \big),
	\end{split}
	\]
	where $S_{T}^{\circ}$ is the set of isotopy classes of simple closed curves in $T$ which are not boundary parallel. The corresponding GR amplitudes are $\Xi_{\Sigma}^{\textup{K}}(\GG) = 1$ for any $\GG \in \mathcal{T}_{\Sigma}^{\textup{comb}}$ and any bordered surface $\Sigma$.
\end{theoremB}

Combining this result with Theorem~\ref{thm:B2} gives a new proof of the topological recursion for Kontsevich volumes, whereas Theorem~\ref{thm:B3} gives a new proof of the topological recursion for the lattice point count. The former is equivalent to a proof of Witten's conjecture, which originally followed from Kontsevich theorem \cite{Kon92} supplemented by \cite{DVV91}. The latter is a result known since Norbury \cite{Nor10}. The enumeration of the numbers $\#\mathcal{M}_{g,n}^{\textup{comb},\ZZ}(L)$ has been connected to matrix integrals in the early works of Chekhov and Makeenko \cite{Che93matrix,CM92multi,CM92hint} and further related to enumeration of chord diagrams in \cite{ACNP1,ACNP2}. At that point, Schwinger--Dyson equations for such models give rise to equations that are eventually but non-obviously equivalent to \cite{Nor10}. The scheme of proofs we 
put forward transcends the algebraic manipulations --- whose geometric meaning is unclear --- pertaining to the realm of matrix integrals and which were necessary in both Kontsevich's original proof and in Chekhov--Makeenko's works.

\medskip

A geometric proof of Witten's conjecture that is in a way similar to ours 
was proposed by Bennett, Cochran, Safnuk and Woskoff in \cite{BCSW12}. In 
this regard, the novel element of our work is firstly to make evident the 
connection of the partition of unity of \cite[Section 4]{BCSW12} with a Mirzakhani--McShane identity. Then the mechanism of integration in \cite{BCSW12} --- that relies on a \emph{local} torus action and was valid only for functions of restricted support such as the Kontsevich initial data ---  gets realised as a special case of the more general Theorem~\ref{thm:B2}, by means of the \emph{global} combinatorial Fenchel--Nielsen coordinates and of Theorem~\ref{thm:A2}.

\paragraph{Flowing from hyperbolic to combinatorial spaces.}

The hyperbolic and combinatorial Teichm\"uller spaces can be identified via the spine homeomorphism $\sp \colon \mathcal{T}_{\Sigma} \to \mathcal{T}_{\Sigma}^{\textup{comb}}$ originating from the work of Penner \cite{Pen87} and Bowditch--Epstein \cite{BE88} and Luo \cite{Luo07}. One can in fact interpolate between their respective geometries with the rescaling flow, defined for $\beta \geq 1$ by
\[
	\sigma \in \mathcal{T}_{\Sigma}(L),
	\qquad
	\sigma^{\beta} =  \sp^{-1}\bigl(\beta \, \sp(\sigma)\bigr) \in \mathcal{T}_{\Sigma}(\beta L),
\]
where the operation $\beta \cdot$ consists in multiplying all edge lengths in $\sp(\sigma)$ by $\beta$.

\medskip

It is known from the works of Mondello~\cite{Mon09} and Do~\cite{Do10} that, for each $\sigma \in \mathcal{T}_{\Sigma}$ as $\beta \to \infty$, the 
metric $\beta^{-1}\sigma^{\beta}$ converges in the Gromov--Hausdorff sense to $\sp(\sigma)$, $\beta^{-1}\ell_{\sigma^{\beta}}(\gamma)$ converges to $\ell_{\sp(\sigma)}(\gamma)$ for each simple closed curve $\gamma$, and 
that the Poisson structure $\beta^2 R_{\beta}^{\ast} \pi_{\textup{WP}}$ converges to $\pi_{\textup{K}}$ on the open cells. For the comparison we have used the map $R_{\beta} \colon \mathcal{T}_{\Sigma}^{\textup{comb}}(L) \rightarrow \mathcal{T}_{\Sigma}(\beta L)$ which is obtained by composing  $\sp^{-1}$ with $\beta \cdot$. In Section~\ref{subsec:convergence:length:twist}, we complete this description by giving uniform bounds for the 
convergence of lengths and twists. Note that to access twists, we rely on 
the aforementioned combinatorial $(9g - 9 + 3n)$-theorem.

\begin{theoremC}\label{thm:C1}
	Let $\Sigma$ be a connected bordered surfaces of type $(g,n)$. For any $\gamma$ simple closed curve in $\Sigma$ and $\sigma \in \mathcal{T}_{\Sigma}$, we have that
	\[
		\lim_{\beta \to \infty} \frac{\ell_{\sigma^{\beta}}(\gamma)}{\beta}
		=
		\ell_{\sp(\sigma)}(\gamma),
	\]
	and the convergence is uniform on the thick parts of $\mathcal{T}_{\Sigma}$. Moreover, given a seamed pants decomposition, we have $3g - 3 + n$ hyperbolic and combinatorial twist parameter functions $\tau_i$. Then for any $\sigma \in \mathcal{T}_{\Sigma}$, we have that
	\[
		\lim_{\beta \to \infty} \frac{\tau_i(\sigma^{\beta})}{\beta}
		=
		\tau_i \big(\sp(\sigma)\big),
	\]
	and the convergence is uniform on every compact of $\mathcal{T}_{\Sigma}$.
\end{theoremC}

Thanks to this uniformity, we can flow quite systematically results in hyperbolic geometry to results in combinatorial geometry, and natural functions on $\mathcal{T}_{\Sigma}$ to natural functions on $\mathcal{T}_{\Sigma}^{\textup{comb}}$. In this regard, $\mathcal{T}_{\Sigma}^{\textup{comb}}$ behaves like a (differential geometer's) tropicalisation of $\mathcal{T}_{\Sigma}$. We use this method in Section~\ref{subsec:PL:structure} to 
rederive (an equivalent expression of) Penner's formulae for the action of the mapping class group on Dehn--Thurston coordinates \cite{Pen82} from 
their hyperbolic analogue \cite{Okai93}. An immediate consequence of Penner's formulae is a piecewise linear structure on the combinatorial Teichmüller space.

\begin{theoremC} \label{thm:C2}  [Penner]
	The combinatorial Fenchel--Nielsen coordinates equip $\mathcal{T}_{\Sigma}^{\textup{comb}}$ with a piecewise linear structure.
\end{theoremC}

Penner \cite{Pen82} proved this result for Dehn--Thurston coordinates by direct combinatorial methods. We note that our proof is different, as we obtain it from the hyperbolic case via the rescaling flow. This piecewise linear structure is what remains of the smooth structure of $\mathcal{T}_\Sigma(\beta L)$ when we let $\beta \rightarrow \infty$. 
\medskip

As a second application, we show in Sections~\ref{subsec:GR:flow}--\ref{subsec:TR:flow} that the flow in the $\beta \to \infty$ limit takes the hyperbolic GR to the combinatorial GR, and does the same for TR after integration.

\begin{theoremC}\label{thm:C3}
	Let $(A_{\beta},B_{\beta},C_{\beta})_{\beta \geq 1}$ be a $1$-parameter family of initial data, satisfying the conditions of Theorem~\ref{thm:rescale:GR} and converging uniformly on any compact to a limit $(\hat{A},\hat{B},\hat{C})$ as $\beta \rightarrow \infty$. Denote by $\Omega_{\Sigma;\beta}$ the GR amplitudes on $\mathcal{T}_{\Sigma}$, which results from the initial data $(A_{\beta},B_{\beta},C_{\beta})$, and by $\widehat{\Xi}_{\Sigma}$ the GR amplitudes on $\mathcal{T}_{\Sigma}^{\textup{comb}}$, which results from the initial data $(\hat{A},\hat{B},\hat{C})$. We have for any bordered surface $\Sigma$ and $\sigma \in \mathcal{T}_{\Sigma}$
	\[
		\lim_{\beta \to \infty} \Omega_{\Sigma;\beta}(\sigma^{\beta})
		=
		\widehat{\Xi}_{\Sigma}\big( \sp(\sigma) \big),
	\]
	and the convergence is uniform on any compact of $\mathcal{T}_{\Sigma}$.
\end{theoremC}

\begin{theoremC}\label{thm:C4}
	Let $(A_{\beta},B_{\beta},C_{\beta})_{\beta \geq 1}$ be a $1$-parameter family of initial data, satisfying the conditions of Theorem~\ref{thm:rescale:TR} and converging uniformly on any compact to a limit $(\hat{A},\hat{B},\hat{C})$ as $\beta \rightarrow \infty$. Denote by $\Omega_{\Sigma;\beta}$ the GR amplitudes on $\mathcal{T}_{\Sigma}$, which results from the initial data $(A_{\beta},B_{\beta},C_{\beta})$, and by $\widehat{\Xi}_{\Sigma}$ the GR amplitudes on $\mathcal{T}_{\Sigma}^{\textup{comb}}$, which results form the initial data $(\hat{A},\hat{B},\hat{C})$. For $2g - 2 + n > 0$ and any $L \in \RR_+^{n}$, we have
	\[
		\lim_{\beta \to \infty} \frac{V\Omega_{\Sigma;\beta}(\beta L)}{\beta^{6g - 6 + 2n}}
		=
		V\widehat{\Xi}_{\Sigma}(L),
	\]
	and the convergence is uniform for $L$ in any set of the form $(0,M]^{n}$ with $M > 0$.
\end{theoremC}

This gives a second proof of the combinatorial Mirzakhani--McShane identity of Theorem~\ref{thm:B4}, by applying the flow to the original hyperbolic Mirzakhani--McShane identity of \cite{Mir07simple}.

\paragraph{Application to statistics of measured foliations with respect to combinatorial structures.}

Along the lines of \cite{ABO17}, we prove in Section~\ref{subsec:twisted:GR} a generalisation of Theorem~\ref{thm:B4} that gives access to statistics of combinatorial lengths of multicurves. Denote by $M_{\Sigma}$ the set of multicurves on $\Sigma$ --- which include the empty multicurve and exclude the boundary components --- and by $M'_{\Sigma}$ the subset of primitive multicurves.

\begin{theoremD}\label{thm:D1}
	Let $f \colon \mathbb{R}_{+} \to \mathbb{R}$ be a mea\-su\-ra\-ble fun\-ct\-ion such that $|f(\ell)|$ decays faster than $\ell^{-s}$, for any $s > 0$. Let $(A,B,C,D)$ be initial data satisfying the conditions of Theorem~\ref{thm:B1}. Denote by $\Xi_{\Sigma}$ the corresponding GR amplitudes, and assume they are invariant under braidings of all boundary components. Then, the GR amplitudes $\Xi_{\Sigma}[f]$ associated with the initial data
	\[
	\begin{split}
		A[f](L_1,L_2,L_3) & = A(L_1,L_2,L_3), \\
		B[f](L_1,L_2,\ell) & = B(L_1,L_2,\ell) + f(\ell)A(L_1,L_2,\ell), \\
		C[f](L_1,\ell,\ell') & = C(L_1,\ell,\ell') + f(\ell) B(L_1,\ell,\ell') + f(\ell') B(L_1,\ell',\ell) + f(\ell)f(\ell')A(L_1,\ell,\ell'), \\
		D_{T}[f](\GG) & = D_{T}(\GG) + \sum_{\gamma \in S_{T}^{\circ}} f(\ell_{\GG}(\gamma))\,A\big(\ell_{\GG}(\partial T),\ell_{\GG}(\gamma),\ell_{\GG}(\gamma)\big),
	\end{split}
	\]
	yield the statistics of combinatorial lengths of primitive multicurves
	\[
		\Xi_{\Sigma}[f](\GG) = \sum_{c \in M'_{\Sigma}} \Xi_{\Sigma - c}(\GG|_{\Sigma - c}) \prod_{\gamma \in \pi_0(c)} f(\ell_{\GG}(\gamma)).
	\]
\end{theoremD}

On this basis, we analyse various questions related to the number of multicurves of length bounded by $t$, as well as their asymptotic growth as $ 
t \rightarrow \infty$. More precisely, for any bordered surface $\Sigma$ of type $(g,n)$ consider the function $\mathcal{N}_{\Sigma}^{\textup{comb}} \colon \mathcal{T}_{\Sigma}^{\textup{comb}} \times \RR_+ \to \RR_+$ counting multicurves on $\Sigma$ with combinatorial length bounded by $t$:
\[
	\mathcal{N}_{\Sigma}^{\textup{comb}}(\GG;t)
	=
	\#\Set{c \in M_{\Sigma} | \ell_{\GG}(c) \le t}.
\]
Thanks to the above formalism, we prove that the counting function --- after Laplace transform in the cut-off parameter $t$ --- is computed by GR, 
and its average over the combinatorial moduli is computed by TR. We also present the analogous results in the hyperbolic setting, giving a complete recursive solution to the enumeration problem of multicurves on both the combinatorial and ordinary Teichmüller spaces.

\begin{theoremD}\label{thm:D2}
	The Laplace transform $\widehat{\mathcal{N}}_{\Sigma}^{\textup{comb}}(\GG;s) = s \int_{\RR_+} \mathcal{N}_{\Sigma}^{\textup{comb}}(\GG;t) e^{-st} \dd t$ of the counting function coincides with the GR amplitudes $\Xi_{\Sigma}^{\textup{K}}[f_s](\GG)$ obtained by twisting the Kontsevich initial data of Theorem~\ref{thm:B4} by $f_s(\ell) = \frac{1}{e^{s\ell} - 1}$, in accordance with Theorem~\ref{thm:D1}. Moreover, they coincides with the statistics of combinatorial lengths of multicurves, weighted by a decaying exponential function:
	\[
		\widehat{\mathcal{N}}_{\Sigma}(\GG;s) = \sum_{c \in M_{\Sigma}} e^{-s\ell_{\GG}(c)}.
	\]
\end{theoremD}

Combining this result with Theorem~\ref{thm:B2} yields a topological recursion for the (Laplace transform of the) average number of multicurves with combinatorial length bounded by $t$. In other words the quantities
\[
	V\mathcal{N}_{g,n}^{\textup{comb}}(L;t)
	=
	\int_{\mathcal{M}_{g,n}^{\textup{comb}}(L)} \mathcal{N}_{g,n}^{\textup{comb}}(\GG;t) \, \dd\mu_{\textup{K}}(\GG)
\]
are well-defined and can be computed recursively on $2g - 2 + n$. Moreover, $V\mathcal{N}_{g,n}^{\textup{comb}}(L;t)$ is a homogeneous polynomial in $\pi^2t^2,L_1^2,\dots,L_n^2$ of degree $3g - 3 + n$ with positive rational coefficients, that is symmetric in $L_1,\dots,L_n$. We provide a list of polynomials $V\mathcal{N}_{g,n}^{\textup{comb}}(L;t)$ for low $2g - 2 + n$ in Table~\ref{tab:multicurve:poly}.

\begin{table}[h]
	\small
	\begin{center}
\begin{tabular}{>{$}c<{$}>{$}c<{$}|>{$}l<{$}}
		\toprule
		g & n &
			V\mathcal{N}_{g,n}^{\textup{comb}}(L;t) \\
		\midrule
		0 & 3 &
			1
		\\[.5ex]
		0 & 4 &
			\tfrac{1}{2} m_{(1)} + \tfrac{\pi^2 t^2}{4}
		\\[.5ex]
		0 & 5 &
			\tfrac{1}{8} m_{(2)} + \tfrac{1}{2} m_{(1,1)}+
			\tfrac{\pi^2 t^2}{4} m_{(1)} +
			\tfrac{\pi^4 t^4}{32}
		\\[.5ex]
		0 & 6 &
			\tfrac{1}{48} m_{(3)} + \tfrac{3}{16} m_{(2,1)} + \tfrac{3}{4} m_{(1,1,1)} +
			\tfrac{5\pi^2 t^2}{48} m_{(2)} + \tfrac{3\pi^2 t^2}{8} m_{(1,1)}
			\tfrac{3\pi^4 t^4}{64} m_{(1)} +
			\tfrac{\pi^6 t^6}{384}
		\\[.5ex]
		0 & 7 &
			\begin{aligned}[t]
				&
				\tfrac{1}{384} m_{(4)} + \tfrac{1}{24} m_{(3,1)} + \tfrac{3}{32} m_{(2,2)} + \tfrac{3}{8} m_{(2,1,1)} + \tfrac{3}{2} m_{(1,1,1,1)} +
				\tfrac{5 \pi^2 t^2}{192} m_{(3)} + \tfrac{5 \pi^2 t^2}{24} m_{(2,1)} \\
				& \quad
				+ \tfrac{3 \pi^2 t^2}{4} m_{(1,1,1)} +
				\tfrac{61 \pi^4 t^4}{2304} m_{(2)} + \tfrac{3 \pi^4 t^4}{32} m_{(1,1)} +
				\tfrac{\pi^6 t^6}{192} m_{(1)} +
				\tfrac{\pi^8 t^8}{6144}
			\end{aligned}
		\\[.5ex]
		\midrule
		1 & 1 &
			\tfrac{1}{48} m_{(1)} + \tfrac{\pi^2 t^2}{24}
		\\[.5ex]
		1 & 2 &
			\tfrac{1}{192} m_{(2)} + \tfrac{1}{96} m_{(1,1)} +
			\tfrac{\pi^2 t^2}{48} m_{(1)} +
			\tfrac{\pi^4 t^4}{384}
		\\
\bottomrule
	\end{tabular}
	\caption{A list of polynomials $V\mathcal{N}_{g,n}^{\textup{comb}}(L;t)$ for low values of $2g - 2 + n$. Here $m_{\lambda}$ is the monomial symmetric polynomial associated to the partition $\lambda$, evaluated at $L_1^2,\dots,L_n^2$.}
	\label{tab:multicurve:poly}
	\end{center}
\end{table}

Since the work of Mirzakhani \cite{Mir08earth}, it is known that the Masur--Veech volumes $MV_{g,n}$ of the principal stratum of the moduli space of quadratic differentials on curves of genus $g$ with $n$ punctures can be accessed through the enumerative geometry of measured foliations with respect to hyperbolic length functions. We show that the role of hyperbolic geometry is not essential. Namely, a similar relation holds with respect to combinatorial length functions. Let ${\rm MF}_{\Sigma}$ be the set of measured foliations, equipped with its Thurston measure $\mu_{{\rm Th}}$ coming from integer point counting.

\begin{theoremD}\label{thm:D3}
	For any bordered surface $\Sigma$ of type $(g,n)$, there exists a continuous function $\mathscr{B}^{\textup{comb}}_{\Sigma} \colon \mathcal{T}_{\Sigma}^{\textup{comb}} \to \mathbb{R}_{+}$ such that
	\[
		\forall \GG \in \mathcal{T}_{\Sigma}^{\textup{comb}},
		\qquad
		\mathscr{B}^{\textup{comb}}_{\Sigma}(\GG)
		=
		\lim_{t \to \infty} \frac{\#\set{c \in M_{\Sigma} | \ell_{\GG}(c) \leq t}}{t^{6g - 6 + 2n}}
		=
		\mu_{\textup{Th}}\big(\set{\lambda \in {\rm MF}_{\Sigma} | \ell_{\GG}(\lambda) \leq 1}\big).
	\]
	Further, for any $L \in \mathbb{R}_{+}^n$, $\mathscr{B}^{\textup{comb}}_{\Sigma}$ is integrable over $\mathcal{M}_{g,n}^{\textup{comb}}(L)$ with respect to $\mu_{\textup{K}}$, and
	\[
		MV_{g,n} = 2^{4g - 2 + n}(4g - 4 + n)! \cdot (6g - 6 + 2n) \cdot \int_{\mathcal{M}_{g,n}^{\textup{comb}}(L)} \mathscr{B}_{\Sigma}^{\textup{comb}}\,\dd\mu_{\textup{K}}.
	\] 
\end{theoremD}

In \cite{ABCDGLW19}, the authors and Delecroix have described the Masur--Veech volumes as the constant term of a family of polynomials $MV_{g,n}(L)$ satisfying the topological recursion. In the proof of Theorem~\ref{thm:D3}, we find a geometric interpretation of the Masur--Veech polynomials in terms of multicurve length statistics on the combinatorial moduli space: they coincide with the Laplace transform of $V\mathcal{N}_{g,n}^{\textup{comb}}(L;t)$ evaluated at Laplace variable $s = 1$.

\medskip

Although the equality between Masur--Veech volumes and the integral of $\mathscr{B}^{\textup{comb}}_{\Sigma}$ goes through an indirect comparison between a certain sum over stable graphs (see \cite{DGZZ19,ABCDGLW19}), we believe it would be interesting to find a direct geometric proof of the 
equality, similarly to the result of \cite{Mir08earth} in the hyperbolic setting. In other words, we expect an almost everywhere defined map
\[
	F \colon \mathcal{T}_{\Sigma}^{\textup{comb}} \times {\rm MF}_{\Sigma}
	\longrightarrow
	\mathcal{QT}_{\Sigma}
	\qquad \text{s.t.} \qquad
	F^{\ast}\mu_{\textup{MV}} = \mu_{\textup{K}} \times \mu_{\textup{Th}}.
\]
Here $\mathcal{QT}_{\Sigma}$ is the space of holomorphic quadratic differentials on $\Sigma$. We also expect the above map to respect the integral structure, and this would give a geometric proof the square-tiled surfaces counting of \cite{ABCDGLW19} in terms of lattice points polynomials.

\medskip

We conclude by giving nine equivalent ways to compute $MV_{g,n}$ and discuss the logical relations between them, adding to the contributions of earlier work the perspectives offered by combinatorial geometry.

\paragraph{Appendices.}

The paper is supplemented with three appendices and an \hyperref[sec:notation]{index of notations}. In Appendix~\ref{app:topology:comb} we recall 
the original definition of the topological space $\mathcal{T}_{\Sigma}^{\textup{comb}}$ from the arc complex, whose equivalence with the one induced by combinatorial length functions is justified in Corollary~\ref{cor:topologies:comb:Teich}. In Appendix~\ref{app:cut:gluing} we give detailed examples of cutting, gluing and Fenchel--Nielsen coordinates. In Appendix~\ref{app:integral:points:factor:2} we present self-contained proofs of facts --- known as folklore but for which we could not identify a clear reference --- related to the appearance of factors of $2$ in the study of the integral structure of $\mathcal{M}_{g,n}^{\textup{comb}}(L)$. 
 
\paragraph{Notation for statements and proofs.}The symbol $\qedsymbol$ is used at the end of a proof as usual. The symbol $\blacksquare$ is used at the end of those statements that do not come with a formal proof in the text: they are either proved in other sources in the literature or the argument before their statement already constitutes a proof for them. Those statements that are not followed immediately by a proof or a $\blacksquare$ are proved later in the text. 

\medskip
\begin{center}
	\noindent\rule{8cm}{0.4pt}	
\end{center}
\medskip

\paragraph{Acknowledgements.}

We thank Vincent Delecroix for numerous discussions, Leonid Chekhov, Ralph Kaufmann and Athanase Papadopoulos for references. We are grateful to Robert Penner for bringing to our attention his work on measured foliations for the decorated Teichm\"uller space, which led us to clarify similarities and differences with our work. J.E.A. was supported in part by the Danish National Sciences Foundation Centre of Excellence grant ``Quantum Geometry of Moduli Spaces" (DNRF95) at which part of this work was carried out. J.E.A. and D.L. are currently funded in part by the ERC Synergy grant ``ReNewQuantum". G.B., S.C., A.G., D.L. and C.W. are supported by the Max-Planck-Gesellschaft. 
D.L. is supported by the Institut de Physique Th\'eorique (IPhT), CEA, CNRS, and by the Institut des Hautes \'{E}tudes Scientifiques (IHES).

\newpage
\section{Topology of combinatorial Teichm\"uller spaces}
\label{sec:comb:spaces}

We begin this section with a basic introduction to the combinatorial moduli space $\mathcal{M}_{g,n}^{\textup{comb}}$, then we study the topology of the combinatorial Teichm\"uller space $\mathcal{T}^{\textup{comb}}_{\Sigma}$, which is the orbifold universal cover of $\mathcal{M}_{g,n}^{\textup{comb}}$. 
We review its relation with the ordinary Teichm\"uller space of hyperbolic structures, we propose a description in terms of measured foliations, and then we introduce the associated length functions. We discuss cutting and gluing after a twist, and finally introduce the combinatorial analogue of Fenchel--Nielsen coordinates, that are shown in Section~\ref{subsec:PL:structure} to endow $\mathcal{T}_{\Sigma}^{\textup{comb}}$ with an atlas having piecewise linear coordinate transformations.

\subsection{Combinatorial moduli spaces}
\label{subsec:comb:moduli}

\begin{defn}\label{defn:graph}
	A \emph{graph} is a triple $(\vec{E},i,\sim)$, where $\vec{E}$ is a finite set (oriented edges), $i \colon \vec{E} \to \vec{E}$ is a fixed-point-free involution (reversing edge orientations) and $\sim$ is an equivalence relation on $\vec{E}$ (classes of which are the vertices). The set of $i$-orbits is denoted by $E$ and describes the unoriented edges; the equivalence class of an oriented edge $\vec{e}$ is simply denoted by $e$. The set $V = \vec{E}/\sim$ describes the vertices. For a fixed $v \in V$, $\#\set{\vec{e} \in \vec{E} | [\vec{e}] = v}$ is called the valency of $v$. Connectedness of graphs is defined in the natural way.
\end{defn}

\begin{defn}\label{defn:MRG}
	A \emph{ribbon graph} is a triple $G = (\vec{E},i,s)$ where $\vec{E}$ is a finite set (oriented edges), $i \colon \vec{E} \to \vec{E}$ is a fixed-point-free involution and $s \colon \vec{E} \to \vec{E}$ is a permutation. The set of $i$-orbits is denoted by E=$E_{G}$ and describes the unoriented edges; the equivalence class of an oriented edge $\vec{e}$ is simply denoted by $e$. The set of $s$-orbits is denoted by $V=V_G$ and describes the vertices. Let $\phi = i \circ s^{-1}$; the set of $\phi$-orbits is denoted by $F=F_G$ and describes the faces (or boundaries) of $G$. Given a ribbon graph $(\vec{E},i,s)$, the underlying graph is $(\vec{E},i,\sim)$, where we declare that $\vec{e} \sim \vec{e}'$ if they are in the same $s$-orbit. A ribbon graph is connected if the underlying graph is.

	\smallskip

	The genus of a connected ribbon graph $G$, denoted by $g$, is defined by
	\begin{equation}
		2 - 2g = \#V_G - \#E_G + \#F_G.
	\end{equation}
	If furthermore it has $n=\#F_G$ boundaries, it is said to be of type $(g,n)$. We call $G$ reduced if all its vertices have valency $\ge 3$, and labelled if its boundaries are labelled $\de_1 G,\dots, \de_n G$. We define $\mathcal{R}_{g,n}$ to be the set of connected, reduced, labelled ribbon graphs of type $(g,n)$.
	
	\begin{rem}
		One can think of a ribbon graph as a graph with a cyclic order on the half-edges incident to each vertex. Using this, one can thicken the edges into a ribbons, which at each vertex glues together to form a smooth surface with boundary (see Figure \ref{fig:geometric:realisations}), hence the name ribbon graph.
	\end{rem}

	A \emph{metric ribbon graph} $\bm{G}$ consists of a ribbon graph $G = (\vec{E},i,s)$, together with the assignment of $\ell_{\bm{G}} \colon E_{G} \to \RR_{+}$. The {\em perimeter} of a boundary component $\de_m G$ consisting of the cyclic sequence of edges $(e_1,\ldots,e_k)$ is defined by
	\begin{equation}
		\ell_{\bm{G}}(\de_m G) = \sum_{i = 1}^k \ell_{\bm{G}}(e_i).
	\end{equation}
\end{defn}

\begin{defn} \label{defn:automorphism}
	An automorphism of a ribbon graph $G = (\vec{E},i,s)$ is a permutation $\varphi \colon \vec{E} \to \vec{E}$ that commutes with $i$ and $s$ and acts trivially on the set of faces $F$. We denote by $\Aut(G)$ the automorphism group of $G$. An automorphism of a metric ribbon graph $\bm{G}$ is an automorphism of the underlying ribbon graph preserving $\ell_{\bm{G}}$. We denote by $\Aut(\bm{G})$ the automorphism group of $\bm{G}$; it is a subgroup of $\Aut(G)$.
\end{defn}

For any ribbon graph $G$, its automorphism group $\Aut(G)$ acts on the space of metrics on $G$, by edge permutation. Given a point in $\RR_{+}^{E_{G}}$, namely a metric ribbon graph $\bm{G}$, its stabiliser under the action of $\Aut(G)$ is precisely $\Aut(\bm{G})$. It is then natural to make the following definition.

\begin{defn}\label{defn:Mgncomb}
	For $2g - 2 + n > 0$, the \emph{combinatorial moduli space} $\mathcal{M}^{\textup{comb}}_{g,n}$ is the orbicell-complex parametrising metric ribbon graphs of type $(g,n)$, \emph{i.e.}
	\begin{equation}
		\mathcal{M}^{\textup{comb}}_{g,n} = \bigcup_{G \in \mathcal{R}_{g,n}} \frac{\RR_{+}^{E_{G}}}{\Aut(G)},
	\end{equation}
	where the cells naturally glue together via edge degeneration: when an edge length goes to zero, the edge contracts to give a metric ribbon graph with fewer edges as discussed in \cite{Kon92}. Define the perimeter map $p \colon \mathcal{M}^{\textup{comb}}_{g,n} \to \RR_{+}^n$ by setting
	\begin{equation}
		p(\bm{G}) = \bigl( \ell_{\bm{G}}(\de_1 G), \ldots, \ell_{\bm{G}}(\de_n G) \bigr).
	\end{equation}
	We also denote $\mathcal{M}^{\textup{comb}}_{g,n}(L) = p^{-1}(L)$ for $L \in \RR_{+}^n$.
\end{defn}

From the above definition, it follows that the combinatorial moduli space $\mathcal{M}^{\textup{comb}}_{g,n}$ has a natural real orbifold structure of dimension $6g - 6 + 3n$. As we remarked before, for a point $\bm{G}$ in the moduli space its orbifold stabiliser is $\Aut(\bm{G})$. Further, we have an orbicell decomposition given by the sets $\RR_{+}^{E_{G}}/\Aut(G)$ consisting of those metric ribbon graphs whose underlying ribbon graph is $G$.  Notice that the dimension of such a cell is $\# E_{G}$, and in particular the top-dimensional cells are the ones associated to trivalent metric ribbon graphs, for which $\#E_{G} = 6g - 6 + 3n$.

\begin{ex}\label{ex:comb:moduli}
	The moduli space $\mathcal{M}_{0,3}^{\textup{comb}}$ is homeomorphic through the perimeter map to the open cone $\RR_{+}^3$, obtained as the union of seven cells corresponding to the seven ribbon graphs of type $(0,3)$ (see Figure~\ref{fig:comb:moduli:03}). In this case there are no orbifold points and $\mathcal{M}_{0,3}^{\textup{comb}}(L_1,L_2,L_3)$ is just a point.

	\smallskip

	A more complicated example is that of $\mathcal{M}_{1,1}^{\textup{comb}}$. The space is obtained by gluing two orbicells, as depicted in Figure~\ref{fig:comb:moduli:11}.
	\begin{itemize}
		\item The first orbicell is $\RR_{+}^3$ quotiented by $\ZZ/6\ZZ$, that is the automorphism group of the unique ribbon graph of type $(1,1)$ with three edges:
		\begin{center}

		\end{center}
	All metric ribbon graphs $\bm{G}_1$ whose underlying graph is $G_1$ have $\Aut(\bm{G}_1) = \ZZ/2\ZZ$ except for the metric ribbon graph where all three edges have the same length which has $\Aut(\bm{G}_1) = \ZZ/6\ZZ$ instead.
		\item
		The second orbicell is $\RR_{+}^2$ quotiented by $\ZZ/4\ZZ$, the automorphism group of the unique ribbon graph of type $(1,1)$ with two edges:
		\begin{center}
		%
		\end{center}
	All metric ribbon graphs $\bm{G}_2$ whose underlying graph is $G_2$ have $\Aut(\bm{G}_2) = \ZZ/2\ZZ$, except for the metric ribbon graph where both edges have the same length which rather has $\Aut(\bm{G}_1) = \ZZ/4\ZZ$.
	\end{itemize}
	\begin{figure}
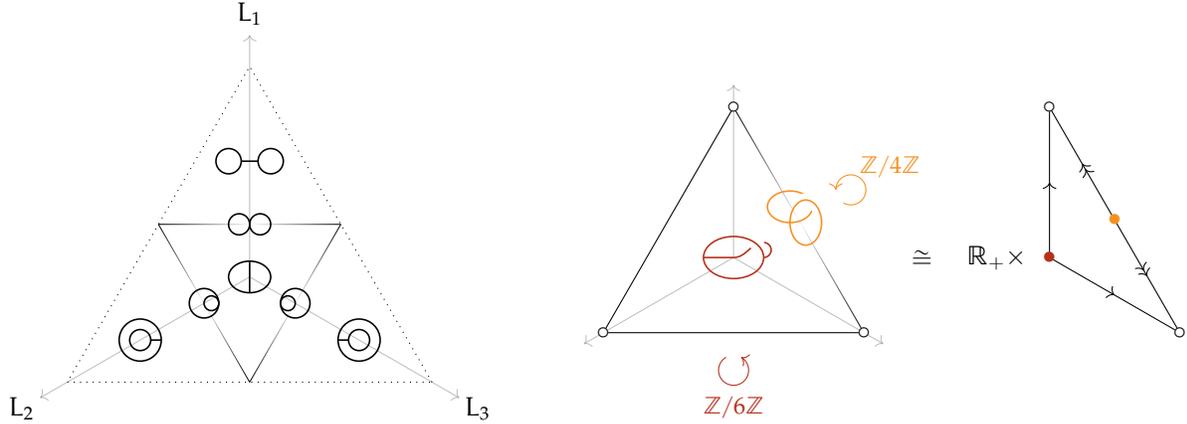

		\centering
		\begin{subfigure}[t]{.41\textwidth}
		\centering
		%
		\subcaption{The $7$ cells composing $\mathcal{M}_{0,3}^{\textup{comb}}$. The pic\-tu\-re re\-pre\-sents the cone $\mathbb{R}_{+}^3$ together with a slice $\Set{L_1+L_2+L_3 = {\rm const}}$. The dotted lines are not part of the space.}
		\label{fig:comb:moduli:03}
		\end{subfigure}
		\hfill%
		\begin{subfigure}[t]{.56\textwidth}
		\centering
		%
		\subcaption{The two orbicells composing $\mathcal{M}_{1,1}^{\textup{comb}}$. On the right, the red point has stabiliser $\ZZ/6\ZZ$, the orange point has stabiliser $\ZZ/4\ZZ$, and all other points have stabiliser $\ZZ/2\ZZ$. The white points are not part of the space.}
		\label{fig:comb:moduli:11}
		\end{subfigure}
		\caption{The combinatorial moduli spaces of type $(0,3)$ and $(1,1)$.}
		\label{fig:comb:moduli}
	\end{figure}
\end{ex}

The main reason for considering $\mathcal{M}^{\textup{comb}}_{g,n}(L)$ is the following result.

\begin{thm}\label{Riemann:combinatorial:homeo}
	The moduli spaces $\mathfrak{M}_{g,n}$ and $\mathcal{M}^{\textup{comb}}_{g,n}(L)$ are orbifold-homeomorphic.
	\hfill $\blacksquare$
\end{thm}

There are at least two different ways to define the above homeomorphism. The first, originating from Jenkins and Strebel \cite{Jen57, Str67}, relies on the geometry of meromorphic quadratic differentials. The second one uses hyperbolic geometry and is due to Penner and Bowditch--Epstein \cite{Pen87, BE88}.  In the next section we discuss the latter and its generalisation to combinatorial Teichm\"uller spaces, due to Luo and Mondello \cite{Luo07,Mon09}.

\subsection{Combinatorial Teichm\"uller spaces}
\label{subsec:comb:Teich}

One can think of points in the combinatorial Teichm\"uller space as points in the combinatorial moduli space together with a marking, the advantage being that there is now a well-defined notion of length of curves. For later purposes, we review its known relation with hyperbolic surfaces via the spine construction, and emphasise its equivalent description via measured foliations.

\subsubsection{Primary definitions}

\begin{defn}\label{defn:bord:MCG}
	A \emph{bordered surface} $\Sigma$ is a non-empty topological, compact, oriented surface with labelled boundary components $\de_1\Sigma,\ldots,\de_n\Sigma$. We assume that each connected component has non-empty boundary and is stable, \emph{i.e.} its Euler characteristic is negative. If furthermore $\Sigma$ is connected of genus $g$, we call $(g,n)$ the type of $\Sigma$. We use $P$ (resp. $T$) to refer to bordered surfaces with the topology of a pair of pants (resp. of a torus with one boundary component).

	\smallskip

	Denote by $\Mod_{\Sigma}$ the \emph{mapping class group} of $\Sigma$:
	\begin{equation}
		\Mod_{\Sigma} = {\rm Diff}^{+}(\Sigma) / {\rm Diff}^{+}_0(\Sigma),
	\end{equation}
	where ${\rm Diff}^{+}(\Sigma)$ is the group of orientation-preserving diffeomorphisms of the surface $\Sigma$ and ${\rm Diff}^{+}_0(\Sigma)$ denotes its subgroup consisting of those diffeomorphisms isotopic to the identity. The \emph{pure mapping class group} $\Mod_{\Sigma}^{\de}$ is the subgroup of $\Mod_{\Sigma}$ consisting of mapping classes which preserve the labellings of boundary component of $\Sigma$.
\end{defn}

To describe the combinatorial Teichm\"uller space it is convenient to introduce some explicit geometric structures associated with a metric ribbon graph.

\begin{defn} \label{defn:geom:realis:G}
	The \emph{geometric realisation of a graph} $\Gamma$ is the $1$-di\-men\-sio\-nal CW-complex
	\begin{equation}
		|\Gamma| = \bigsqcup_{e \in E } [0,1] \Big/ \sim
	\end{equation}
	obtained by identifying endpoints corresponding to the same vertex in $\Gamma$.
\end{defn}

\begin{figure}
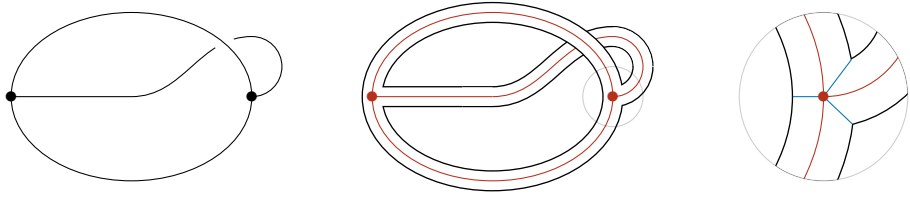

	\centering

	\caption{The geometric realisation of a graph $\Gamma$ (left), of a metric ribbon graph $\bm{G}$ (center), and the local picture around a vertex of $|\bm{G}|$ (right). In red we draw $\bm{G}$ embedded into $|\bm{G}|$ and the segments along which the stripes are glued are depicted in blue.}
	\label{fig:geometric:realisations}
\end{figure}

\begin{defn}\label{defn:geom:realis:MRG}
	The \emph{geometric realisation of a metric ribbon graph} $\bm{G}$ of type $(g,n)$ is the oriented, compact, genus $g$ surface with $n$ boundary components
	\begin{equation}
		|\bm{G}| =
			\Set{
				(u,t,\vec{e}) \in [-1,1] \times \RR \times \vec{E} | t \in [0,\ell_{\bm{G}}(e)]
			} \Big/ \sim
	\end{equation}
	where the equivalence relation $\sim$ is given by
	\begin{equation}
		(u,t, \vec{e}) \sim (u, \ell_{\bm{G}}(e) - t, i(\vec{e}))
		\qquad\qquad
		{\rm for}\,\,u\in [-1,1],
		\;
		\,\,\vec{e}\in \vec{E},
		\;
		\,\, t \in [0, \ell_{\bm{G}}(e)]
	\end{equation}
	and
	\begin{equation}
		(u, \ell_{\bm{G}}(e), \vec{e}) \sim (-u, 0, s(\vec{e}))
		\qquad\qquad
		{\rm for}\,\, u\in [0,1],
		\;
		\,\, \vec{e} \in \vec{E}.
	\end{equation}
The geometric realisation of the underlying graph can be seen as a subset of $|\bm{G}|$, and the inclusion is a deformation retract; see Figure~\ref{fig:geometric:realisations} for an example.
\end{defn}

\begin{defn}\label{defn:comb:struct}
	A \emph{combinatorial marking} of a bordered surface $\Sigma$ is an ordered pair $(\bm{G},f)$ where $\bm{G}$ is a metric ribbon graph and $f \colon \Sigma \to |\bm{G}|$ is a homeomorphism respecting the labellings of boundaries of $\bm{G}$ and $\Sigma$. The \emph{combinatorial Teichm\"uller space} is defined as
	\begin{equation}
		\mathcal{T}_{\Sigma}^{\textup{comb}} = \Set{ (\bm{G},f) | \phantom{\Big|} \text{$(\bm{G},f)$ is a combinatorial marking on $\Sigma$}} \Big/ \sim.
	\end{equation}
	Here, $(\bm{G},f) \sim (\bm{G}',f')$ if and only if there exists a metric ribbon graph isomorphism $\varphi \colon \bm{G} \to \bm{G}'$ such that $\varphi \circ f$ is isotopic to $f'$. We call such equivalence classes \emph{combinatorial structures} and denote by $\GG = [\bm{G},f]$ the points in $\mathcal{T}^{\textup{comb}}_{\Sigma}$.
\end{defn}

Notice that, for a fixed $\GG \in \mathcal{T}^{\textup{comb}}_{\Sigma}$ each boundary component $\de_m \Sigma$ corresponds to a unique face of the embedded graph. Thus, we can talk about the length (with respect to $\GG$) of the boundary $\de_m \Sigma$, denoted by $\ell_{\GG}(\de_m \Sigma)$. In particular, we can define the perimeter map $p \colon \mathcal{T}^{\textup{comb}}_{\Sigma} \to \RR_{+}^n$ by setting
\begin{equation}\label{eqn:comb:perimeter}
	p(\GG) = \bigl( \ell_{\GG}(\de_1 \Sigma), \ldots, \ell_{\GG}(\de_n \Sigma) \bigr).
\end{equation}
We also set $\mathcal{T}^{\textup{comb}}_{\Sigma}(L) = p^{-1}(L)$.

\medskip

We remark that, for a fixed representative $(\bm{G},f)$ of $\GG \in \mathcal{T}^{\textup{comb}}_{\Sigma}$, we have a retraction from the surface to the geometric realisation of the graph underlying $\bm{G}$. Thus, we can picture elements of $\mathcal{T}^{\textup{comb}}_{\Sigma}$ as metric ribbon graphs embedded into $\Sigma$ up to isotopy, such that the embedded graph is a retract of the surface.

\medskip

We recall in Appendix~\ref{app:topology:comb} that the combinatorial Teichm\"uller space can be endowed with a natural topology via the so-called proper arc complex that makes it into a cell complex. The cells, denoted $\mathfrak{Z}_{\Sigma,G}$, are indexed by isotopy classes of embedding of $G$ into $\Sigma$, onto which $\Sigma$ retracts, and they parametrise the metrics of $G$.  The practical consequence of this topological discussion is that the edge lengths form a coordinate system in each cell. In particular, this makes it easier to check whether a function defined on $\mathcal{T}^{\textup{comb}}_{\Sigma}$ is piecewise continuous, once it is expressed in terms of edge lengths. 
\begin{figure}[t]
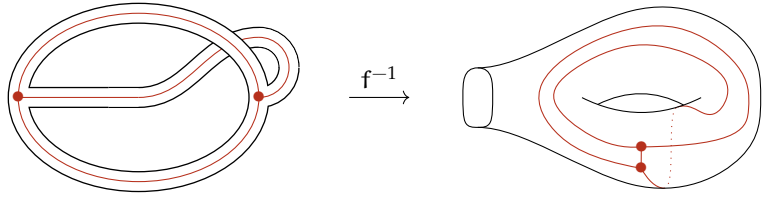

	\centering

	\caption{An example of combinatorial marking on a torus with one boundary component.}
	\label{fig:comb:marking}
\end{figure}

\medskip

The mapping class group of $\Sigma$ naturally acts on $\mathcal{T}^{\textup{comb}}_{\Sigma}$ by setting
\begin{equation}
	[\phi].[\bm{G},f] = [\bm{G},f \circ \phi],
\end{equation}
for $[\phi] \in \Mod_{\Sigma}$ and $[\bm{G},f] \in \mathcal{T}^{\textup{comb}}_{\Sigma}$. When $\Sigma$ has type $(g,n)$, by forgetting the marking we have the isomorphism
\begin{equation}
	\mathcal{M}_{g,n}^{\textup{comb}}(L) \cong \mathcal{T}^{\textup{comb}}_{\Sigma}(L) / \Mod_{\Sigma}^{\de}.
\end{equation}
We denote the quotient $\mathcal{T}^{\textup{comb}}_{\Sigma}(L) / \Mod_{\Sigma}^{\de}$ by $\mathcal{M}_{\Sigma}^{\textup{comb}}(L)$ when we want to refer to the actual surface $\Sigma$.

\begin{ex}
	For a pair of pants $P$, the perimeter map gives the isomorphism $\mathcal{T}_{P}^{\textup{comb}} \cong \mathbb{R}_{+}^3$. The pure mapping class group is trivial, so that $\mathcal{M}_{0,3}^{\textup{comb}} \cong \mathbb{R}_{+}^3$, see Figure~\ref{fig:comb:Teich:P}.

	\medskip

	For a torus $T$ with one boundary component, $\mathcal{T}_{T}^{\textup{comb}}$ is the union of infinitely many cells homeomorphic to $\RR_{+}^3$, glued together through infinitely many cells homeomorphic to $\RR_{+}^2$. In Figure~\ref{fig:comb:Teich:T} we presented some adjacent cells of $\mathcal{T}_{T}^{\textup{comb}}$. In the quotient by $\Mod_{T}^{\de}$, all the $3$- and $2$-cells are identified, and we are left with a further action of $\ZZ/6\ZZ$ for the top-dimensional cell and an action of $\ZZ/4\ZZ$ for the codimension $1$ cell. The result is the combinatorial moduli space $\mathcal{M}_{1,1}^{\textup{comb}}$ described in Example~\ref{ex:comb:moduli}.
	\begin{figure}[t]
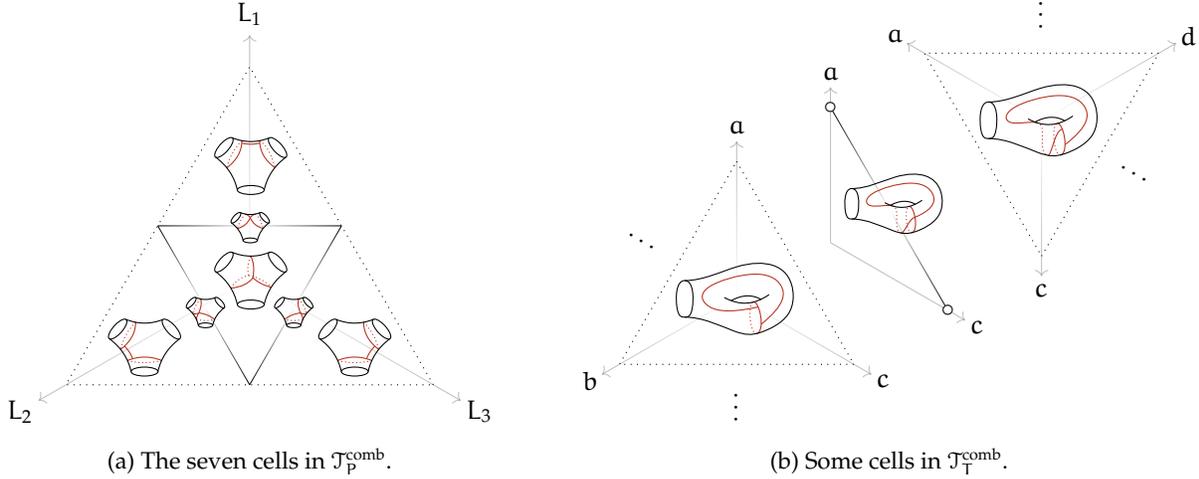

		\centering
		\begin{subfigure}[t]{.41\textwidth}
		\centering

		\subcaption{The seven cells in $\mathcal{T}^{\textup{comb}}_P$.}
		\label{fig:comb:Teich:P}
		\end{subfigure}
		\hfill%
		\begin{subfigure}[t]{.56\textwidth}
		\centering
		%
		\subcaption{Some cells in $\mathcal{T}^{\textup{comb}}_{T}$.}
		\label{fig:comb:Teich:T}
		\end{subfigure}
		\caption{The combinatorial Teichm\"uller spaces of type $(0,3)$ and $(1,1)$.}
		\label{fig:comb:Teich}
	\end{figure}
\end{ex}

\subsubsection{Relation with ordinary Teichm\"uller spaces}
\label{subsubsec:ordinary:combinatorial}

For a fixed connected bordered surface $\Sigma$ of type $(g,n)$, we can consider the ordinary Teichm\"uller space. Recall that a hyperbolic marking on $\Sigma$ is a pair $(X,f)$ where $X$ is a hyperbolic surface with labelled geodesic boundaries and $f \colon \Sigma \to X$ is an orientation-preserving diffeomorphism respecting the labelling. Define the Teichm\"uller space as
\begin{equation}\label{defn:Teichm}
	\mathcal{T}_{\Sigma}
	=
	\Set{ (X,f) | \phantom{\Big|} \text{$(X,f)$ is a hyperbolic marking on $\Sigma$}} \Big/ \sim ,
\end{equation}
where $(X,f) \sim (X',f')$ if and only if there exists an isometry $\varphi \colon X \to X'$ respecting the labelling of the boundaries and such that $\varphi \circ f$ is isotopic to $f'$. We denote points in $\mathcal{T}_{\Sigma}$ by $\sigma = [X,f]$, and we call them hyperbolic structures on $\Sigma$. By considering hyperbolic length of the boundaries, we have a perimeter map $p \colon \mathcal{T}_{\Sigma} \to \RR_{+}^n$ and we set $\mathcal{T}_{\Sigma}(L) = p^{-1}(L)$ for $L \in \mathbb{R}_{+}^n$.

\medskip

It is well-known that the quotient of $\mathcal{T}_{\Sigma}(L)$ by the pure mapping class group $\Mod_{\Sigma}^{\de}$ is orbifold-ho\-meo\-mor\-phic to
\begin{equation}\label{defn:MgnL}
	\mathcal{M}_{g,n}(L) = \left.
			\Set{ X | \begin{array}{@{}c@{}}
				\text{$X$ is a hyperbolic surface of type $(g,n)$ } \\
				\text{ with labelled boundary of lengths $L$}
			\end{array}
			}
		\middle/ \sim \right. ,
\end{equation}
where $X \sim X'$ if and only if there exists an isometry from $X$ to $X'$ preserving the labelling of the boundary components. Further, for each $L \in \mathbb{R}_{+}^n$, $\mathcal{M}_{g,n}(L)$ is orbifold-homeomorphic to the moduli space of complex curves $\mathfrak{M}_{g,n}$.  For later use, we review the description of a mapping class group equivariant homeomorphism between the combinatorial Teichm\"uller space and the ordinary one, due to Luo and Mondello \cite{Luo07,Mon09}. It lifts to the level of Teichm\"uller spaces the construction of Penner and Bowditch--Epstein \cite{Pen87,BE88}.

\medskip

Let $\sigma \in \mathcal{T}_{\Sigma}(L)$ and fix a representative $(X,f)$. Define the valency $\nu_{\sigma}(q)$ of a point $q$ in the interior of $\Sigma$ as the number of shortest geodesics joining $q$ to $\de\Sigma$ that realise the distance ${\rm dist}_{\sigma}(q,\de\Sigma)$. Clearly $\nu_{\sigma}(q) \ge 1$. Define the loci $A_{\sigma} = \set{q \in \Sigma | \nu_{\sigma}(q) = 2}$ and $V_{\sigma} = \set{q \in \Sigma | \nu_{\sigma}(q) \ge 3}$. We have that $A_{\sigma}$ is the disjoint union of simple, open geodesic arcs $\alpha_{e}$ indexed by $E = \pi_0(A_{\sigma})$, the set of edges, and $V_{\sigma}$ is a finite collection of points, the vertices. We define the \emph{spine} $\sp(\sigma)$ of $\sigma$ as the $1$-dimensional CW-complex embedded in $\Sigma$ given by $V_{\sigma} \cup A_{\sigma}$.

\medskip

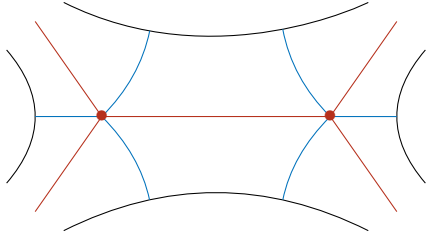
\begin{figure}
	\centering
	\begin{tikzpicture}[x=1pt,y=1pt,xscale=.45,yscale=.45]

		\draw (128, 768) .. controls (256, 704) and (384, 768) .. (384, 768);
		\draw (128, 576) .. controls (213.3333, 618.6667) and (298.6667, 618.6667) .. (384, 576);
		\draw (432, 616) .. controls (400, 653.3333) and (400, 690.6667) .. (432, 728);
		\draw (80, 728) .. controls (112, 690.6667) and (112, 653.3333) .. (80, 616);
		\draw[BrickRed](160, 672) -- (104, 752);
		\draw[shift={(160.001, 671.999)}, xscale=1.1667, yscale=1.1111, BrickRed] (0, 0) -- (-48, -72);
		\draw[shift={(352.001, 671.999)}, xscale=1.1667, yscale=1.1111, BrickRed] (0, 0) -- (48, -72);
		\draw[BrickRed](352, 672) -- (408, 752);
		\draw[BrickRed](160, 672) -- (352, 672) -- (352, 672);
		\draw[RoyalBlue] (160, 672) -- (104, 671.997);
		\draw[RoyalBlue] (352, 672) -- (408, 671.971);
		\draw[RoyalBlue] (160, 672) .. controls (181.3333, 693.3333) and (194.7553, 717.344) .. (200.266, 744.032);
		\draw[RoyalBlue] (160, 672) .. controls (181.3333, 650.6667) and (194.71, 627.301) .. (200.13, 601.903);
		\draw[RoyalBlue] (352, 672) .. controls (330.6667, 650.6667) and (317.3847, 627.2803) .. (312.154, 601.841);
		\draw[RoyalBlue] (352, 672) .. controls (330.6667, 693.3333) and (317.377, 717.819) .. (312.131, 745.457);
		\node[BrickRed] at (160, 672) {$\bullet$};
		\node[BrickRed] at (352, 672) {$\bullet$};

	\end{tikzpicture}
	\caption{Example of the spine construction. In red, the spine $\sp(\sigma)$. In blue, the ribs emanating from two vertices.}
	\label{fig:spine}
\end{figure}

We can naturally assign to $\sp(\sigma)$ a metric $\ell_{\sp(\sigma)} \colon E \to \RR_{+}$ in the following way. For each vertex $q$ of $\sp(\sigma)$, consider the $\nu_{\sigma}(q)$ shortest geodesics from $q$ to the boundary -- which are called \emph{ribs}. Cutting $\Sigma$ along its ribs yields a union of hexagons; the diagonal of each hexagon whose endpoints are the vertices of the spine corresponds to the edges of the spine (see Figure~\ref{fig:spine}). We assign to it the length of the side of the hexagon which lies along the boundary of $\Sigma$; there are two such sides, but they have the same length since the reflection with respect to the edge is a hyperbolic isometry of the hexagon. In this way, $\sp(\sigma)$ induces a combinatorial marking on $\Sigma$ and the perimeters of $\sp(\sigma)$ correspond precisely to the hyperbolic lengths of the boundaries of $\Sigma$. Further, isotopy classes of hyperbolic markings correspond to isotopy classes of combinatorial markings. Thus, we are led to the following definition.

\begin{defn}\label{defn:spine:map}
	There exists a well-defined map
	\begin{equation}
		\sp \colon \mathcal{T}_{\Sigma}(L) \longrightarrow \mathcal{T}_{\Sigma}^{\textup{comb}}(L),
	\end{equation}
	called the \emph{spine map}.
\end{defn}

It is possible (although more difficult) to construct the inverse map and actually show that it is a homeomorphism, equivariant with respect to the action of $\Mod_{\Sigma}^{\de}$.
	
\begin{thm}\label{thm:homeo:ord:comb:Teich}\cite{Luo07,Mon09}
	The spine map $\sp \colon \mathcal{T}_{\Sigma}(L) \to \mathcal{T}^{\textup{comb}}_{\Sigma}(L)$ is a homeomorphism, equivariant under the action of the pure mapping class group.
	\hfill $\blacksquare$
\end{thm}

In \cite{Luo07,Mon09} the theorem is stated in terms of the proper arc complex, rather than the combinatorial Teichm\"uller space, but their relation is recalled in Appendix~\ref{app:topology:comb}. As a direct consequence, we find that for each fixed $L \in \RR_{+}^n$, there is an orbifold homeomorphism $\mathcal{M}_{g,n}(L) \cong \mathcal{M}^{\textup{comb}}_{g,n}(L)$ and thus with $\mathfrak{M}_{g,n}$ (Theorem~\ref{Riemann:combinatorial:homeo}).

\subsubsection{Relation with measured foliations}
\label{subsubsec:measured:foliations}

For a given metric ribbon graph, its geometric realisation is naturally endowed with a measured foliation, as we now explain. We refer to \cite[Section~5.1]{FLP12} for a complete discussion about measured foliations, but to be self-contained we recall here the basic definitions.

\begin{defn}\label{defn:measured:foliations}
	Let $\Sigma$ be a bordered surface, and $\mathcal{F}$ a \emph{foliation} on $\Sigma$ with isolated singularities. A \emph{transverse invariant measure} on $\mathcal{F}$ is a measure $\mu$ defined on  each arc transverse to the foliation, invariant under isotopy of arcs through transverse arcs whose endpoints remain in the same leaf. If the arc passes through a singularity, the transversality pertains to all points of the arc belonging to a regular leaf.
\end{defn}

\begin{figure}[hb]
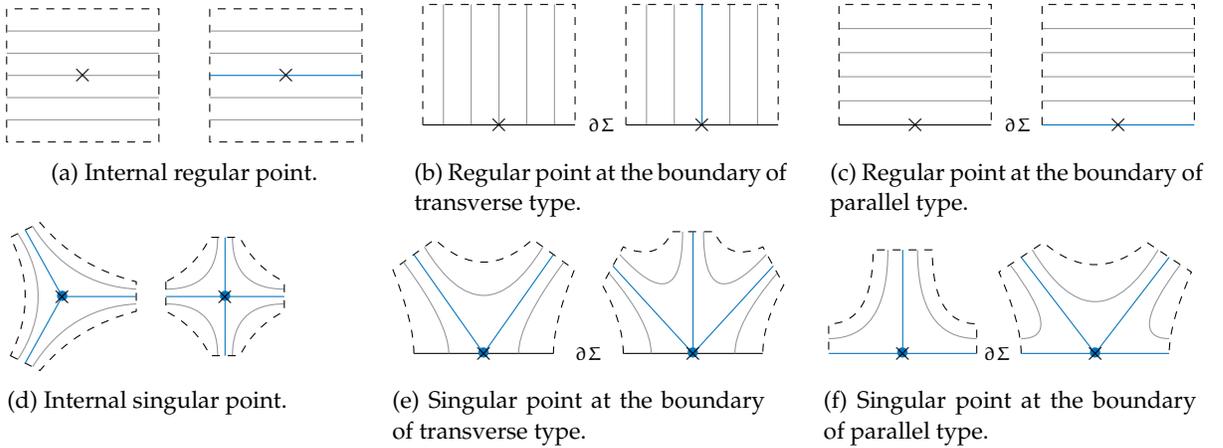

\centering
	\hfill
	\begin{subfigure}[t]{.3\textwidth}
	\centering

	\caption{Internal regular point.}
	\label{fig:points:foliation:a}
	\end{subfigure}
	\hfill
	\begin{subfigure}[t]{.3\textwidth}
	\centering
	%
	\caption{Regular point at the boundary of transverse type.}
	\label{fig:points:foliation:b}
	\end{subfigure}
	\hfill
	\begin{subfigure}[t]{.3\textwidth}
	\centering
	%
	\caption{Regular point at the boundary of parallel type.}
	\label{fig:points:foliation:c}
	\end{subfigure}
	\hfill
	\begin{subfigure}[t]{.3\textwidth}
	\centering
	%
	\caption{Internal singular point.}
	\label{fig:points:foliation:d}
	\end{subfigure}
	\hfill
	\begin{subfigure}[t]{.3\textwidth}
	\centering
	%
	\caption{Singular point at the boundary of transverse type.}
	\label{fig:points:foliation:e}
	\end{subfigure}
	\hfill
	\begin{subfigure}[t]{.3\textwidth}
	\centering
	%
	\caption{Singular point at the boundary of parallel type.}
	\label{fig:points:foliation:f}
	\end{subfigure}
	\hfill
	\caption{Possible models for points in a foliation. The singular leaves are depicted in blue, while smooth leaves in grey. Only singular points of low valency are depicted, but any higher valency is allowed.}
	\label{fig:points:foliation}
\end{figure}

\begin{samepage}
	In what follows, we also suppose that:
	\begin{itemize}
		\item the measure is regular with respect to the Lebesgue one: every regular point of $\Sigma$ admits a smooth chart $U \ni (x,y)$ where the foliation is defined by $\dd y$ and the measure on each transverse arc is induced by $|\dd y|$;
		\item each point of $\Sigma$ has a neighbourhood that is the domain of a chart isomorphically foliated as one of the models of Figure~\ref{fig:points:foliation}.
	\end{itemize}
\end{samepage}

\begin{figure}[ht]
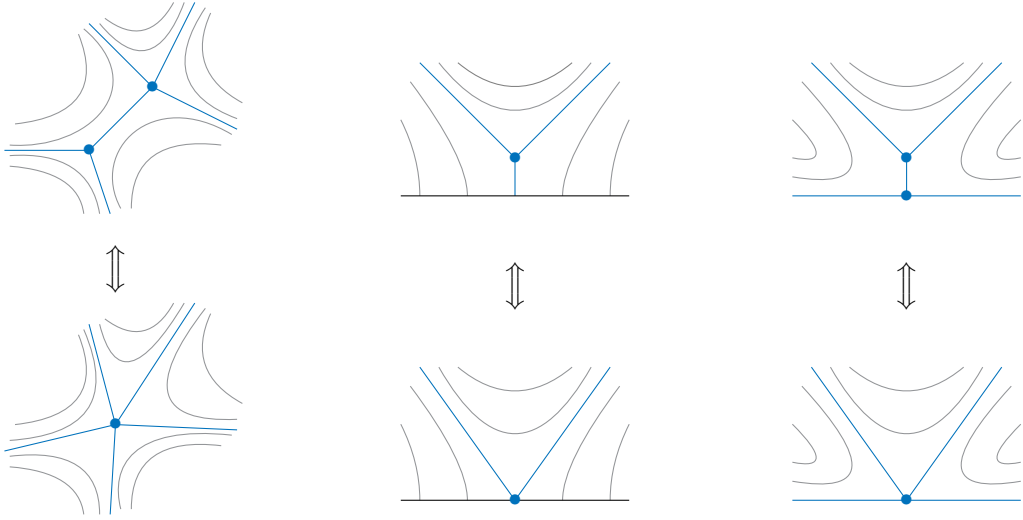

\centering
	\begin{subfigure}[t]{.31\textwidth}
	\centering

	\end{subfigure}
	\caption{Whitehead moves.}
	\label{fig:Whitehead:moves}
\end{figure}

\begin{defn}
	We say that two measured foliations on $\Sigma$ are \emph{Whitehead equivalent} if they differ by isotopy or a Whitehead move, (see Figure~\ref{fig:Whitehead:moves}). We denote by ${\rm MF}_{\Sigma}^{\star}$ the set of Whitehead equivalence classes of measured foliations  on $\Sigma$.
\end{defn}

We can now discuss the relation between foliations and metric ribbon graphs.
	 
\begin{defn}\label{defn:mf:geom:realis:MRG}
	Given a metric ribbon graph $\bm{G}$ of type $(g,n)$, the geometric realisation $|\bm{G}|$ has a unique measured foliation $(\mathcal{F}_{\bm{G}},\mu_{\bm{G}})$ such that:
	\begin{itemize}
		\item the singularities of $(\mathcal{F}_{\bm{G}},\mu_{\bm{G}})$ are the vertices of the embedded graph,
		\item the measured foliation is transverse to the embedded graph,
		\item $(\mathcal{F}_{\bm{G}},\mu_{\bm{G}})$ on the hexagon around each edge $e \in E$ agrees with $|\dd t|$, where $t$ is the natural coordinate on $[0, \ell_{\bm{G}}(e)]$ as in Definition~\ref{defn:geom:realis:MRG}.
	\end{itemize}
\end{defn}

\begin{figure}
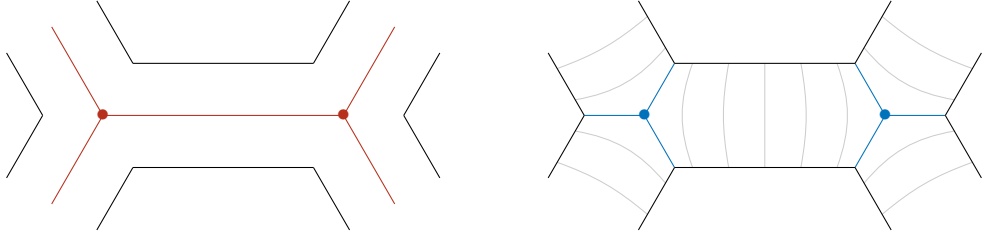

	\centering

	\caption{The geometric realisation of a metric ribbon graph (left) and the associated measured foliation. The edges of the embedded graph are depicted in red, the singular leaves of the associated foliation in blue.}
	\label{fig:geom:realis:foliation}
\end{figure}

The singular leaves of the measured foliation cut $|\bm{G}|$ into hexagons, each with two opposite edges consisting of arcs along the boundary of $|\bm{G}|$, and the remaining four edges consisting of singular leaves (see Figure~\ref{fig:geom:realis:foliation}). The diagonals parallel to the boundary arcs are nothing but the edges of the embedded graph. In the description of \S~\ref{subsubsec:ordinary:combinatorial}, the singular leaves correspond to the ribs, and the singular points to the vertices of the spine. Notice that such a hexagon decomposition of $|\bm{G}|$, together with the assignment of a positive length to each diagonal, is sufficient to reconstruct the metric ribbon graph $\bm{G}$. For an element $\GG = [\bm{G},f]$ in $\mathcal{T}^{\textup{comb}}_{\Sigma}$, we get a natural isotopy class of measured foliations $(\mathcal{F}_{\GG},\mu_{\GG})$ on $\Sigma$ by pushing $(\mathcal{F}_{\bm{G}},\mu_{\bm{G}})$ forward along $f$. In the following, we omit the transverse measure $\mu_{\GG}$ when there is no ambiguity.

\medskip

The above construction defines a map
\begin{equation}\label{eqn:comb:Teich:to:MF}
	\mathcal{F}_{\ast} \colon \mathcal{T}_{\Sigma}^{\textup{comb}} \longrightarrow {\rm MF}_{\Sigma}^{\star},
\end{equation}
whose image is the set of classes of measured foliations admitting a representative, with respect to Whitehead equivalence, whose leaves are all transverse to the boundary of $\Sigma$ (\emph{cf.} Figure~\ref{fig:points:foliation:b} and \ref{fig:points:foliation:e}), and there is no singular leaf connecting two singular points -- such singular leaf is called a \emph{saddle connection}. Moreover, the measured foliation can be used to give a hexagonal decomposition of the surface, and reconstruct the embedded metric ribbon graph. To summarise, we have the following lemma.

\begin{lem}
	The map $\mathcal{F}_{\ast}$ is injective. Its image consists of measured foliations admitting a representative with respect to Whitehead equivalence whose leaves are all transverse to the boundary of the surface and with no saddle connections.
	\hfill $\blacksquare$
\end{lem}

In the remaining part of this section, we establish various topological properties of the combinatorial Teichmüller space. Many of these properties are established in \cite{FLP12} but a different space of measured foliations (namely, the completion of the set of multicurves), and in \cite{PH92} in the language of train tracks. We include proofs and explicitly highlight results that are not simple adaptations from \cite{FLP12,PH92}. The details of these proofs will in fact be used in later parts of the article.

\subsection{Length functions and topology}
\label{subsec:length:top}

We introduce now combinatorial length functions for simple closed curves, and show that combinatorial structures are (locally in $\mathcal{T}_{\Sigma}^{\textup{comb}}$) completely determined by the knowledge of finitely many of these lengths. This gives an alternative description of the topology on $\mathcal{T}_{\Sigma}^{\textup{comb}}$.

\subsubsection{Notations for curves}
\label{subsubsec:notation:curves}

We denote by
\begin{itemize}
	\item
		$S_{\Sigma}$ the set of homotopy classes of simple closed curves on $\Sigma$, and its subset $S_{\Sigma}^{\circ}$ consisting of non boundary parallel (also called \emph{essential}) curves,
	\item
		$M_{\Sigma}$ the set of multicurves, \emph{i.e.} homotopy classes of finite unions of pairwise disjoint essential simple closed curves on $\Sigma$,
	\item
		$M_{\Sigma}'$ the set of primitive multicurves, \emph{i.e.} those multicurves whose components are pairwise non-homotopic.\end{itemize}
From now on, curves are always considered up to homotopy. By convention $M_{\Sigma}$ and $M_{\Sigma}'$ contain the empty multicurve, whereas $S_{\Sigma}$ and $S_{\Sigma}^{\circ}$ do not.

\medskip

If $\Sigma$ is a bordered surface and $\gamma \in S_{\Sigma}^{\circ}$ (or more generally in $M_{\Sigma}'$), we can consider the closed surface $\Sigma_{\gamma}$ defined as the result of cutting $\Sigma$ along a chosen representative of $\gamma$. The assumptions on $\gamma$ imply that every connected component of $\Sigma_{\gamma}$ is stable. Among such connected components, there is one containing $\de_{1}\Sigma$. We label the boundaries of such surface by putting the components of $\de \Sigma$ first (in the order they appear in $\Sigma$), followed by those of $\gamma$ (in some order). For the connected components of $\Sigma_{\gamma}$ that do not contain $\de_{1}\Sigma$, we label the boundaries by putting the components of $\gamma$ first (in some order), followed by those of $\de \Sigma$ (in the order they appear in $\Sigma$). In the following we specify the choice of order only in case it has an actual relevance in the argument.

\subsubsection{Combinatorial length functions}
\label{subsubsec:comb:length}

Consider a combinatorial structure $\GG \in \mathcal{T}^{\textup{comb}}_{\Sigma}$ and a simple closed curve $\gamma \in S_{\Sigma}$. As discussed in the previous section, $\GG$ induces an isotopy class of measured foliations $(\mathcal{F}_{\GG},\mu_{\GG})$ on $\Sigma$, so we have a well-defined length function $\ell_{\GG}$ defined on $S_{\Sigma}$ (\emph{cf.} \cite[Section~5.3]{FLP12})\footnote{In \cite{FLP12} this quantity is denoted by $I(\mathcal{F},\mu;\gamma)$.}
\begin{equation}\label{eqn:length:curve}
	\ell_{\GG}(\gamma) = \inf_{\gamma_0} \; \sup_{\alpha_i} \left(\sum_{i=1}^k \mu_{\GG}(\alpha_i) \right),
\end{equation}
where the infimum is taken over all representatives $\gamma_0$ in the homotopy class $\gamma$, and the supremum is taken over all the possible ways of writting $\gamma_0$ as a sum of arcs $\alpha_1,\dots,\alpha_k$, mutually disjoint and transverse to $\mathcal{F}_{\GG}$.

\medskip

It can be shown \cite[Proposition~5.7]{FLP12} that the infimum is reached by quasitransverse representatives of $\gamma$.

\begin{defn}
	We say that a curve $\gamma_0$ is \emph{quasitransverse} to a foliation $\mathcal{F}$ if each connected component of $\gamma_0$ minus the singularities of $\mathcal{F}$ is either a leaf or is transverse to $\mathcal{F}$. Further, in a neighbourhood of a singularity, we require that no transverse arc lies in a sector adjacent to an arc contained in a leaf, and that consecutive transverse arcs lie in distinct sectors.
\end{defn}

In terms of the embedded metric ribbon graph, we have the following effective way to compute the length $\ell_{\GG}(\gamma)$: consider the (unique) representative $\gamma_0$ of $\gamma$ that has been homotoped to the embedded graph and is non-backtracking. The length of the curve can now be computed as the sum of the lengths of the edges visited by the curve, since such $\gamma_0$ is a quasitransverse representative, perhaps after performing a sequence of Whitehead moves in a small disc neighbourhood of the vertices of $\GG$, making a measured foliation $\mathcal{F}_{\GG}'$ Whitehead equivalent to $\mathcal{F}_{\GG}$ and thus representing the same point in ${\rm MF}_{\Sigma}^{\star}$. One can also conclude from \cite[Proposition 5.9]{FLP12} that the length is always positive.

\begin{figure}[t]
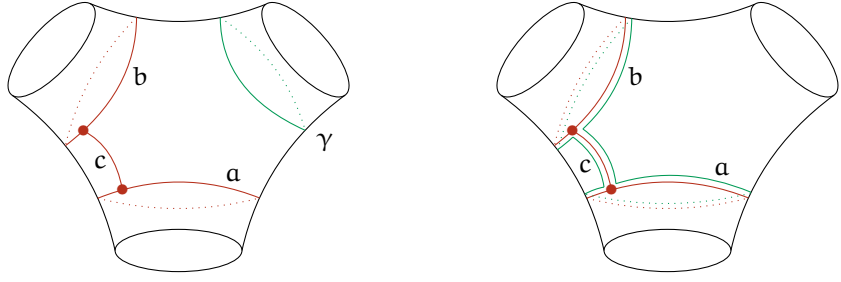

	\centering

	\caption{The length of a simple closed curve $\gamma$ with respect to a combinatorial structure $\GG$ on a pair of pants. We have $\ell_{\GG}(\gamma) = a + c + b + c = a + b + 2c$.}
	\label{fig:example:length}
\end{figure}

\medskip

Notice that the assignment $\GG \mapsto \ell_{\GG}(\gamma)$ is continuous on $\mathcal{T}_{\Sigma}^{\textup{comb}}$, as it is a sum of edge lengths on the closure of each open cell. If the curve is homotopic to one of the boundary components, the notion of boundary length described before (Equation~\eqref{eqn:comb:perimeter}) agrees with this more general definition. The notion of length naturally extends to $M_{\Sigma}$ and $M_{\Sigma}'$ by adding the length of components. Using Equation~\eqref{eqn:length:curve}, we can also define the length of homotopy classes (relative to $\de\Sigma$) of arcs between boundaries. This is again a continuous assignment, but it can now take zero values.

\medskip

In \S~\ref{subsubsection:comb:pairs:of:pants}-\ref{sec:proof:length:functional} we prove the following result regarding the combinatorial length spectrum, which is the analogue of \cite[Theorem~1.3]{FLP12} for measured foliations on a closed surface, and the analogue of \cite[Theorem~1.4]{FLP12} for the usual Teichmüller space of hyperbolic structures, again on a closed surface.

\begin{thm}\label{thm:length:functional}
	Let us equip $\RR_{+}^{S_{\Sigma}}$ with the product topology. The combinatorial length of simple closed curves gives a map $\ell_{\ast} \colon \mathcal{T}_{\Sigma}^{\textup{comb}} \rightarrow \RR_{+}^{S_{\Sigma}}$ which is a homeomorphism onto its image.
\end{thm}

\begin{cor}\label{cor:topologies:comb:Teich}
	The topology on  $\mathcal{T}_{\Sigma}^{\textup{comb}}$ defined via the arc-complex (Appendix~\ref{app:topology:comb}) and the initial to\-po\-lo\-gy induced by the combinatorial length map $\ell_{\ast} \colon \mathcal{T}_{\Sigma}^{\textup{comb}} \rightarrow \RR_{+}^{S_{\Sigma}}$ coincide.
	\hfill $\blacksquare$
\end{cor}

Here we prove the theorem by a direct construction of the inverse map from the image of $\ell_{\ast}$, via the study of lengths of curves and their relation with embedded pairs of pants. These computations are used again in Section~\ref{sec:setting:up:gr}.

\subsubsection{Embedded pairs of pants}
\label{subsubsection:comb:pairs:of:pants}

We are interested in embedded pairs of pants and the way the lengths of their boundary components (with respect to a combinatorial structure) determine the combinatorial structure itself.

\begin{defn}\label{defn:homotopy:class:pop}
	When $\Sigma$ is a bordered surface of type $(g,n)$ such that $2g - 2 + n > 1$, and $m_0 \in \{1,\ldots,n\}$, we let $\mathcal{P}_{\Sigma,m_0}$ be the set of homotopy classes of embeddings $\varphi \colon P \hookrightarrow \Sigma$ of pairs of pants $P$ such that
	\begin{itemize}
		\item[$\bullet$] $\varphi(\de_1 P) = \de_{m_0} \Sigma$,
		\item[$\bullet$] $\overline{\Sigma \setminus \varphi(P)}$ (where the closure is taken in $\Sigma$) is stable,
		\item[$\bullet$] if $\varphi(\de_i P) = \de_{m} \Sigma$ for some $m \neq m_0$, then $i = 2$.
	\end{itemize}
	We denote by $\mathcal{P}_{\Sigma,m_0}^{\varnothing}$ the set of homotopy classes of pair of pants $[P]$'s for which $\de_{i} P \subset \Sigma^{\circ}$ for $i = 2,3$, and by $\mathcal{P}_{\Sigma,m_0}^{m}$ the set of such $[P]$'s for which $\de_{2} P = \de_{m}\Sigma$ for some $m \neq m_0$. We have a partition
	\begin{equation}
		\mathcal{P}_{\Sigma,m_0} = \bigg(\bigsqcup_{m \neq m_0} \mathcal{P}_{\Sigma,m_0}^{m}\bigg) \sqcup \mathcal{P}_{\Sigma,m_0}^{\varnothing}.
	\end{equation}
\end{defn}

These pairs of pants can be described alternatively in terms of arcs between boundaries \cite{ABO17}, which we recall here for the convenience of the reader.

\begin{defn}\label{defn:homotopy:class:arcs}
	Denote by $\mathcal{A}_{\Sigma,m_0}$ the set of non-trivial homotopy classes of proper embeddings\footnote{If $X$ and $Y$ are topological manifolds with boundaries, a continuous map $f \colon X \rightarrow Y$ is called a proper embedding if $f^{-1}(\partial Y) = \partial X$, and we use here the natural notion of homotopies among such.} $a \colon [0,1] \hookrightarrow \Sigma$ with $a(0) \in \de_{m_0}\Sigma$. We denote the class of a proper embedding $a$ with the letter $\alpha$. There is a partition
	\begin{equation}
		\mathcal{A}_{\Sigma,m_0} = \bigg(\bigsqcup_{m \neq m_0} \mathcal{A}_{\Sigma,m_0}^{m}\bigg) \sqcup \mathcal{A}_{\Sigma,m_0}^{\varnothing},
	\end{equation}
	where the elements of $\mathcal{A}_{\Sigma,m_0}^{m}$ are those with a representative $a$ such that $a(1) \in \de_m\Sigma$, and the elements of $\mathcal{A}_{\Sigma,m_0}^{\varnothing}$ are those with a representative $a$ such that $a(1) \in \de_{m_0}\Sigma$.
\end{defn}

\begin{figure}
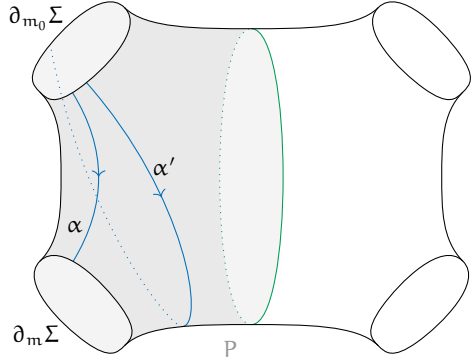

\begin{center}

	\caption{Arcs and pairs of pants: $Q_{m_0}(\alpha) = Q_{m_0}(\alpha') = Q_{m_0}(-\alpha') = [P]$.}
	\label{fig:arc:pants:corresp}
\end{center}
\end{figure}

We obtain a surjective map
\begin{equation}\label{defn:arcs:to:pants}
	Q_{m_0} \colon \mathcal{A}_{\Sigma,m_0} \longrightarrow \mathcal{P}_{\Sigma,m_0}
\end{equation}
by assigning to an arc $\alpha = [a]$ the homotopy class of pairs of pants $[P]$ whose boundaries are formed by the boundary of a closed tubular neighbourhood of $a$ and curves homotopic to the boundaries joined by $a$ (Figure~\ref{fig:arc:pants:corresp}). The boundaries of $P$ are labelled as follows. We always set $\de_1 P = \de_{m_0} \Sigma$. Then if $\varphi(\de_i P) = \de_{m} \Sigma$ for some $m \neq m_0$, then $i = 2$ and $\de_3 P$ is determined; otherwise, we define $\de_2 P$ (resp. $\de_3 P$) to be the boundary component on the left-hand side (resp. right-hand side) of the curve $\alpha$ oriented from $0$ to $1$.

\begin{rem}\label{rem:Q:not:inj}
	The restriction of $Q_{m_0}$ to $\mathcal{A}_{\Sigma,m_0}$ is not injective. More precisely, $Q^{-1}_{m_0}([P])$ contains a single element when $[P] \in \mathcal{P}_{\Sigma,m_0}^{\varnothing}$, while it contains three elements when $[P] \in \mathcal{P}_{\Sigma,m_0}^{m}$ for any $m \ne m_0$. Indeed any given $[P] \in \mathcal{P}_{\Sigma,m_0}^{m}$ can be obtained by an arc $\alpha$ from $\partial_{m_0}\Sigma$ to $\partial_m\Sigma$, but also by an arc $\alpha'$ from $\de_{m_0}\Sigma$ to itself and its inverse $-\alpha'$ (Figure~\ref{fig:arc:pants:corresp}). Notice that in this case $\alpha$ is the only arc in $\mathcal{A}_{\Sigma,m_0}^{m}$, while $\alpha', -\alpha' \in \mathcal{A}_{\Sigma,m_0}^{\varnothing}$.
\end{rem}

We introduce the notion of small pairs of pants, which play a role in the rest of the paper.

\begin{defn}\label{defn:small:pop}
	Let $\GG \in \mathcal{T}_{\Sigma}^{\textup{comb}}$. We say that $[P] \in \mathcal{P}_{\Sigma,m_0}$ is $\GG$-small if
	\begin{equation}
		\ell_{\GG}(\varphi(\partial P) \cap \partial \Sigma)
		\geq
		\ell_{\GG}(\varphi(\partial P) \cap \Sigma^{\circ}).
	\end{equation}
	The definition does not depend on the representative chosen for $[P]$. If the inequalities are strict, then we say the pair of pants is strictly $\GG$-small.
\end{defn}

In other words, dropping the embedding from the notation:
\begin{itemize}
	\item
		when $[P] \in \mathcal{P}_{\Sigma,m_0}^{m}$, we have $\ell_{\GG}(\de_{3}P) \leq \ell_{\GG}(\de_{1}P)+\ell_{\GG}(\de_{2}P)$,
	\item
		when $[P] \in \mathcal{P}_{\Sigma,m_0}^{\varnothing}$, we have $\ell_{\GG}(\de_{2}P)+\ell_{\GG}(\de_{3}P) \leq \ell_{\GG}(\de_{1}P)$.
\end{itemize}
We can characterise small pairs of pants in terms of the corresponding arcs as follows.

\begin{lem}\label{lem:small:pop:char}
	Let $[P]\in\mathcal{P}_{\Sigma,m_0}$.
	\begin{itemize}
		\item
			$[P]\in\mathcal{P}_{\Sigma,m_0}^{m}$ is $\GG$-small if and only if for the unique $\alpha \in Q_{m_0}^{-1}([P])\cap\mathcal{A}_{\Sigma,m_0}^{m}$ we have $\ell_{\GG}(\alpha)=0$,
		\item
			$[P]\in\mathcal{P}_{\Sigma,m_0}^{\varnothing}$ is $\GG$-small if and only if for the unique $\alpha \in Q_{m_0}^{-1}([P])$ we have $\ell_{\GG}(\alpha)=0$.
	\end{itemize}
\end{lem}
\begin{proof}
	Let $\alpha \in \mathcal{A}_{\Sigma,m_0}^{m}$ for $m \neq m_0$ and assume that $\ell_{\GG}(\alpha) > 0$. It is then uniquely represented by some non backtracking edgepath $\alpha$ in $\GG$ with initial and final vertices adjacent to $\de_{m_0}\Sigma$ and $\de_{m}\Sigma$. Then, $\de_{3}Q_{m_0}(\alpha)$ can be homotoped to a non backtracking edgepath consisting in travelling along $\alpha$, then going around $\de_{m}\Sigma$, then travelling backwards along $\alpha$, and going around $\de_{m_0}\Sigma$. This implies that
	\[
		\ell_{\GG}(\partial_{m_0}\Sigma) + \ell_{\GG}(\partial_{m}\Sigma) < \ell_{\GG}\big(\partial_3Q_{m_0}(\alpha)\big),
	\]
	hence $Q_{m_0}(\alpha)$ is not $\GG$-small, and we conclude by the contrapositive. A similar argument works when $\alpha \in \mathcal{A}_{\Sigma,m_0}^{\varnothing}$.
\end{proof}

\begin{rem}\label{rem:small:pop}
	Notice that $\alpha$ can only have length zero with respect to $\GG$ if it has a quasitransverse representative given by an oriented non-singular leaf in the foliation associated to $\GG$. There are finitely many homotopy classes of oriented leaves relative to the boundary, and they are in bijection with the oriented edges of the metric ribbon graph. Since $\mathcal{A}_{\Sigma,m_0}$ surjects onto $\mathcal{P}_{\Sigma,m_0}$ for $\Sigma$ of type $(g,n)$ and $\GG \in \mathcal{T}_{\Sigma}^{\textup{comb}}$, there are at most $2(6g - 6 + 3n)$ $\GG$-small pairs of pants in $\mathcal{P}_{\Sigma,m_0}$.
\end{rem}

Small pairs of pants can be characterised in terms of the support of two functions that play an important role in Section~\ref{subsec:Mirz:McS}. Let $[x]_{+} = \max\set{x,0}$ and consider the functions on $\mathcal{T}^{\textup{comb}}_P\cong\RR_{+}^{3}$ defined by
\begin{equation}\label{eqn:Kon:BC}
\begin{aligned}
	B^{\textup{K}}(L_1,L_2,\ell)	&= \frac{1}{2L_1}\Bigl(
		[L_1 - L_2 - \ell]_+ - [ -L_1 + L_2 - \ell]_+ + [L_1 + L_2 - \ell]_+
	\Bigr), \\
	C^{\textup{K}}(L_1,\ell,\ell')	&= \frac{1}{L_1} [L_1 - \ell - \ell']_+.
\end{aligned}
\end{equation}
It is easy to check that these functions only take non-negative values. They encode aspects of the geometry of combinatorial pairs of pants, as described in the following lemma, and are the combinatorial analogs of the functions $\mathcal{D}$ and $\mathcal{R}$ introduced by Mirzakhani in \cite{Mir07simple} in the hyperbolic context.

\begin{lem}\label{lem:geo:meaning:BC}
	The function $B^{\textup{K}}$ associates to $(\ell_{\GG}(\de_1 P),\ell_{\GG}(\de_2 P),\ell_{\GG}(\de_3 P)) \in \RR_{+}^3 \cong \mathcal{T}_{P}^{\textup{comb}}$ the fraction of the $\de_1 P$ that is not common with $\de_3 P$ (once retracted to the graph). Similarly, $C^{\textup{K}}$ associates to a point $\vec{\ell}_{\GG}(\de P) \in \RR_{+}^3 \cong \mathcal{T}_{P}^{\textup{comb}}$ the fraction of $\de_1 P$ that is not common with $\de_2 P \cup \de_3 P$.
\end{lem}

\begin{proof}
	The result follows from direct computations in  the closure of each of the four open cells of $\mathcal{T}_{P}^{\textup{comb}}$. The various inequalities that define the cells are used to simplify $B^{\textup{K}}$ and $C^{\textup{K}}$. Consider for example Figure~\ref{fig:BC:meaning:a}: we have $\ell_{\GG}(\de_3 P)\geq\ell_{\GG}(\de_1 P)+\ell_{\GG}(\de_2 P)$ so that $B^{\textup{K}}(L_1,L_2,\ell) = 0$, and the portion of $\de_1 P$ that is not common with $\de_3 P$ is zero (all of $\de_1 P$ intersects with $\de_3 P$). Similarly, in the case of Figure~\ref{fig:BC:meaning:b} we find $\ell_{\GG}(\de_2 P)\geq\ell_{\GG}(\de_1 P)+\ell_{\GG}(\de_3 P)$, so that
	\[
		B^{\textup{K}}(L_1,L_2,\ell)
		= \frac{1}{2\ell_{\GG}(\de_1 P)}
			\Bigl(
				\ell_{\GG}(\de_1 P)+\ell_{\GG}(\de_2 P)-\ell_{\GG}(\de_3 P)
				-
				\bigl( -\ell_{\GG}(\de_1 P)+\ell_{\GG}(\de_2 P)-\ell_{\GG}(\de_3 P) \bigr)
			\Bigr)
		= 1,
	\]
	and indeed all of $\de_1 P$ intersects with $\de_2 P$. Finally, in Figure~\ref{fig:BC:meaning:c} we find $\ell_{\GG}(\de_i P)\leq\ell_{\GG}(\de_j P)+\ell_{\GG}(\de_k P)$ for all $i,j,k\in\set{1,2,3}$, which is the condition to be in this cell. We also see that the length of the edge adjacent to $\de_1 P$ and $\de_2 P$ is given by $\frac{1}{2}\left(\ell_{\GG}(\de_1 P)+\ell_{\GG}(\de_2 P)-\ell_{\GG}(\de_3 P)\right)$ and therefore its fraction of the first boundary is given by exactly
	\[
		B^{\textup{K}}(L_1,L_2,\ell)
		= \frac{1}{2\ell_{\GG}(\de_1 P)}
			\bigl(
				\ell_{\GG}(\de_1 P)+\ell_{\GG}(\de_2 P)-\ell_{\GG}(\de_3 P)
			\bigr).
	\]
	The computation in other cells is similar, and analogously for $C^{\textup{K}}$.
	\begin{figure}
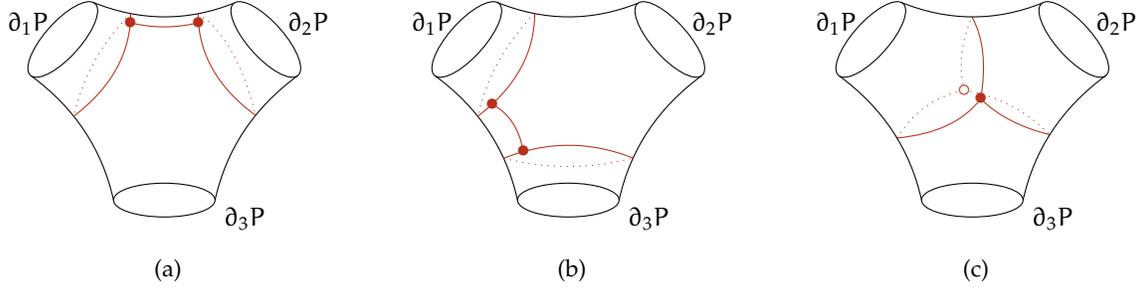

	\centering
		\begin{subfigure}[t]{.32\textwidth}
		\centering

		\caption{}
		\label{fig:BC:meaning:a}
		\end{subfigure}
		\begin{subfigure}[t]{.32\textwidth}
		\centering
		%
		\caption{}
		\label{fig:BC:meaning:b}
		\end{subfigure}
		\begin{subfigure}[t]{.32\textwidth}
		\centering
		%
		\caption{}
		\label{fig:BC:meaning:c}
		\end{subfigure}
		\caption{Three different cells in $\mathcal{T}_P^{\textup{comb}}$.}
		\label{fig:BC:meaning}
	\end{figure}
\end{proof}

Notice that if $B^{\textup{K}}$ or $C^{\textup{K}}$ are non-zero, then there is at least an edge with a corresponding arc, that defines the pair of pants passing through the edge transversely, and therefore of zero length. As an immediate consequence this or simply by considering the support of $B^{\textup{K}}$ and $C^{\textup{K}}$, we deduce the following characterisation of small pairs of pants.

\begin{cor}\label{cor:zero:BC}
	Let $\GG \in \mathcal{T}^{\textup{comb}}_{\Sigma}$:
	\begin{itemize}
		\item $[P] \in \mathcal{P}^{m}_{\Sigma,m_0}$ is $\GG$-small if and only if $B^{\textup{K}}(\vec{\ell}_{\GG}(\de P)) > 0$,
		\item $[P] \in \mathcal{P}^{\varnothing}_{\Sigma,m_0}$ is $\GG$-small if and only if $C^{\textup{K}}(\vec{\ell}_{\GG}(\de P)) > 0$.
	\end{itemize}
	Here $\vec{\ell}_{\GG}(\partial P)$ is the ordered triple of combinatorial lengths of the boundary components of $P$.
	\hfill $\blacksquare$
\end{cor}

\subsubsection{A partial inverse of the combinatorial length spectrum map}
\label{sec:proof:length:functional}

In this paragraph we exhibit an inverse on the image of the combinatorial length map of Theorem~\ref{thm:length:functional}. We first observe that, given $\GG = [\bm{G},f] \in\mathcal{T}_{\Sigma}^{\textup{comb}}$ and an oriented edge $\vec{e}$ of the embedded graph, we can take the dual arc $\alpha_{\vec{e}}$: its starting point (resp. ending point) lies on the component of $\de\Sigma$ adjacent to $f(\vec{e})$ on the right (resp. left) and intersects the embedded graph exactly once through this edge. Note that the two boundary components adjacent to $\vec{e}$ may coincide. We thus obtain a map $\vec{E}_{\bm{G}} \to \mathcal{A}_{\Sigma}^{{\rm all}} := \bigsqcup_{m_0 = 1}^n \mathcal{A}_{\Sigma,m_0}$. Composing with $Q_{m_0}$, we can associate to each oriented edge of $\GG$ a class of embedded pair of pants in $\mathcal{P}^{{\rm all}}_{\Sigma} = \bigsqcup_{m_0 = 1}^n \mathcal{P}_{\Sigma,m_0}$.

\begin{lem}\label{lem:reconstruct:edge:lengths}
	Let $e$ be an edge in $\GG\in\mathcal{T}_{\Sigma}^{\textup{comb}}$ and fix an arbitrary orientation $\vec{e}$.
	\begin{itemize}
		\item If $e$ is adjacent to $\de_{m_1}\Sigma \neq \de_{m_2}\Sigma$, let $[P_{1}]\in\mathcal{P}_{\Sigma,m_1}^{m_2}$ and $[P_{2}]\in\mathcal{P}_{\Sigma,m_2}^{m_1}$ be the pairs of pants corresponding to $\vec{e}$ and $i(\vec{e})$ respectively. Then
		\begin{equation}
		\begin{split}
			\ell_{\GG}(e)
				& =
				\ell_{\GG}(\de_{m_1}\Sigma)
					\Big(
						B^{\textup{K}}(\vec{\ell}_{\GG}(\de P_{1})) - C^{\textup{K}}(\vec{\ell}_{\GG}(\de P_{1}))
					\Big) \\
				& =
				\ell_{\GG}(\de_{m_2}\Sigma)
					\Big(
						B^{\textup{K}}(\vec{\ell}_{\GG}(\de P_{2})) - C^{\textup{K}}(\vec{\ell}_{\GG}(\de P_{2}))
					\Big).
		\end{split}
		\end{equation}
		\item If $e$ is adjacent to $\de_{m_0}\Sigma$ on both sides, let $[P] \in \mathcal{P}_{\Sigma,m_0}$ be the pair of pants corresponding to $\vec{e}$ or $i(\vec{e})$ (the two pairs of pants coincide). Then
		\begin{equation}
			\ell_{\GG}(e) = \frac{1}{2} \ell_{\GG}(\de_{m_0}\Sigma)\,C^{\textup{K}}(\vec{\ell}_{\GG}(\partial P)).
		\end{equation}
	\end{itemize}
\end{lem}
\begin{proof}
	Given a boundary component, we can represent the edges around it by a polygon where some edges and vertices are identified. The sequence of edges around the boundary is non-backtracking.

	\smallskip

	If $e$ is adjacent to $\de_{m_1}\Sigma \neq \de_{m_2}\Sigma$, we have two polygons around $\de_{m_1}\Sigma$ and $\de_{m_2}\Sigma$. When neither of these are $1$-gons, then we can represent $\gamma = \de_{3}P_{1}=\de_{3}P_{2}$ as in Figure~\ref{fig:polygons:around:faces:a}. Due to the absence of bivalent vertices, this is again non backtracking. We therefore have $\ell_{\GG}(\gamma)=\ell_{\GG}(\de_{m_1}\Sigma)+\ell_{\GG}(\de_{m_2}\Sigma)-2\ell_{\GG}(e)$. Notice that we also have $|\ell_{\GG}(\de_{m_1}\Sigma)-\ell_{\GG}(\de_{m_2}\Sigma)| < \ell_{\GG}(\gamma)$. Therefore, comparing the expression for $\ell_{\GG}(e)$ with the expression in the statement, we see they agree.

	\smallskip

	Now suppose without loss of generality that $m_2$ is a $1$-gon. This implies that $m_1$ is not a $1$-gon, as gluing two $1$-gons together would produce a cylinder. It is also clear that $\ell_{\GG}(e)=\ell_{\GG}(\de_{m_2}\Sigma)$. We can represent $\gamma = \de_{3}P_{1}=\de_{3}P_{2}$ by Figure~\ref{fig:polygons:around:faces:b}. This implies that $\ell_{\GG}(\de_{m_1}\Sigma)\geq\ell_{\GG}(\de_{m_2}\Sigma)+\ell_{\GG}(\gamma)$ and therefore, using this to calculate the expression in the statement, we see that they agree.

	\smallskip

	If $e$ is adjacent to $\de_{m_0}\Sigma$ on both sides, we have one polygon around $\de_{m_0}\Sigma$. We can then represent $\de_{2}P$ and $\de_{3}P$ as in Figure~\ref{fig:polygons:around:faces:c}. In absence of bivalent vertices, backtracking cannot occur around either side of $e$. Therefore we have $\ell_{\GG}(\de_{m_0}\Sigma)=\ell_{\GG}(\de_{2}P)+\ell_{\GG}(\de_{3}P)+2\ell_{\GG}(e)$. Comparing the expression for $\ell_{\GG}(e)$ with the expression in the statement, we see they agree.
	\begin{figure}
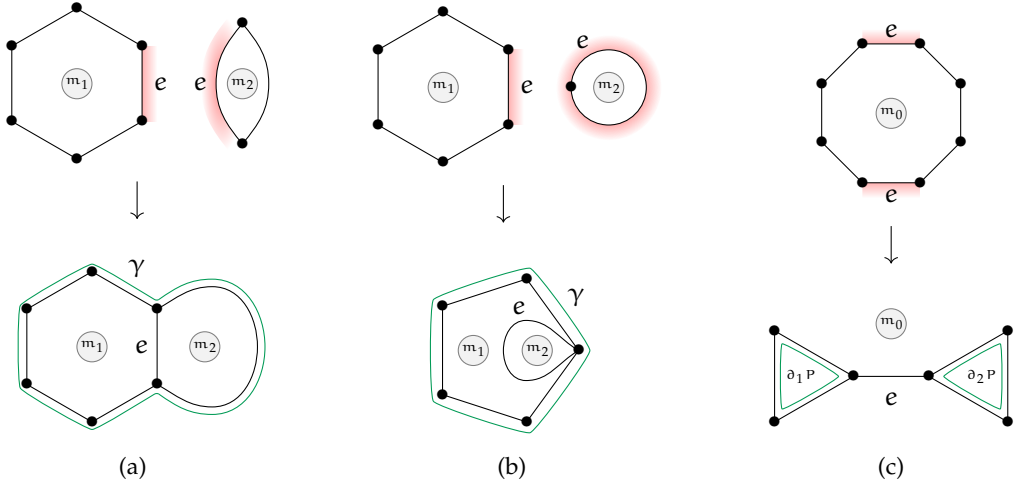

		\centering
		\begin{subfigure}[t]{.3\textwidth}
		\centering

		\caption{}
		\label{fig:polygons:around:faces:a}
		\end{subfigure}
		\begin{subfigure}[t]{.3\textwidth}
		\centering
			%
		\caption{}
		\label{fig:polygons:around:faces:b}
		\end{subfigure}
		\begin{subfigure}[t]{.3\textwidth}
		\centering
			\centering
			%
		\caption{}
		\label{fig:polygons:around:faces:c}
		\end{subfigure}
		\caption{The three cases examined in Lemma~\ref{lem:reconstruct:edge:lengths}.}
		\label{fig:polygons:around:faces}
	\end{figure}
\end{proof}

This lemma shows that we can reconstruct edge lengths from lengths of simple closed curves, which in turn proves that Theorem~\ref{thm:length:functional} holds when restricted to the closure of a cell. To fully reconstruct the ribbon graph (\emph{i.e.} determine in which cell we are) from just the lengths of simple closed curves, we use the characterisations of small pairs of pants and their relation to edges in the embedded graph.

\begin{proof}[Proof of Theorem~\ref{thm:length:functional}]
	We can now exhibit a global inverse on the image of the combinatorial length spectrum map $\ell_{\ast} \colon \mathcal{T}_{\Sigma}^{\textup{comb}} \rightarrow \RR_{+}^{S_{\Sigma}}$. We first consider the following composition:
	\[
	\begin{tikzcd}[column sep=small,row sep=0ex]
		\RR_{+}^{S_{\Sigma}} \arrow[r] & (\RR_{+}^3)^{\mathcal{P}_{\Sigma}^{{\rm all}}} \arrow[r] & \RR_{\geq 0}^{\mathcal{A}_{\Sigma}^{{\rm all}}} \\
		\lambda \arrow[rr,mapsto] && l_{\lambda}.
	\end{tikzcd}
	\]
	The first arrow associates to a length functional $\lambda \in \RR_{+}^{S_{\Sigma}}$ the functional on $\mathcal{P}_{\Sigma}^{{\rm all}}$ defined by $[P] \mapsto \vec{\lambda}(\de P)$, where $\vec{\lambda}(\partial P)$ is the ordered triples of boundary lengths. The second arrow associates to a functional $\vec{\lambda}$ on $\mathcal{P}_{\Sigma}^{{\rm all}}$ a functional on $\mathcal{A}_{\Sigma}^{{\rm all}}$ given by
	\[
		\alpha \longmapsto
			\begin{cases}
				\lambda_1(P_{\alpha})
					\Big(
						B^{\textup{K}}(\vec{\lambda}(P_{\alpha}))
						-
						C^{\textup{K}}(\vec{\lambda}(P_{\alpha}))
					\Big)
					& \text{if } \alpha \in \mathcal{A}_{\Sigma,m_0}^{m} \\[1ex]
				\frac{1}{2}\lambda_1(P_{\alpha}) \, C^{\textup{K}}(\vec{\lambda}(P_{\alpha}))
					& \text{if } \alpha \in \mathcal{A}_{\Sigma,m_0}^{\varnothing}
			\end{cases}
	\]
	where $[P_{\alpha}] = Q_{m_0}(\alpha)$ and $\lambda_1$ is the first component of the triple $\vec{\lambda}$. Note that these are exactly the expressions appearing in Lemma~\ref{lem:reconstruct:edge:lengths}.

	\smallskip

	Consider now the restriction to $\image(\ell_{\ast})$. If $\lambda = \ell_{\ast}(\GG)$ for some $\GG \in \mathcal{T}_{\Sigma}^{\textup{comb}}$, we know from Lemma~\ref{lem:small:pop:char} and Remark~\ref{rem:small:pop} that $l_{\lambda}$, as a functional on $\mathcal{A}_{\Sigma}^{{\rm all}}$, is supported on the complement of those oriented arcs $\alpha$ homotopic to non-singular oriented leaves in the foliation associated to $\GG$. These are in bijection with the oriented edges of $\GG$. Forgetting about the orientation, consider the finite collection $(\alpha_1,\dots,\alpha_k)$ of all such arcs, choosing representatives that are pairwise non-intersecting in $\Sigma$. This defines a proper simplex in the arc complex $\mathcal{A}_{\Sigma}$ (see Appendix~\ref{app:topology:comb}), and taking the dual we obtain a ribbon graph $\bm{G}$ with an embedding into $\Sigma$. We can equip it with a metric, assigning the length $l_{\lambda}(\alpha_i)$ to the edge dual to $\alpha_i$. This represents a point in $\mathcal{T}_{\Sigma}^{\textup{comb}}$, so that we have a map
	\[
		l_{\ast} \colon \image(\ell_{\ast}) \longrightarrow \mathcal{T}_{\Sigma}^{\textup{comb}}.
	\]
	By construction, $l_{\ast} \circ \ell_{\ast} = \id$. Thus, $\ell_{\ast}$ is injective.

	\smallskip

	The map $\ell_{\ast}$ is clearly continuous, since the lengths of simple closed curves are linear combinations of lengths of edges. The inverse map $l_{\ast}$ is also continuous on the $\ell_{\ast}$-image of each cell, since we realised the edge lengths as piecewise linear (and thus continuous) functions of the length of locally finitely many simple closed curves. This completes the proof.
\end{proof}

\subsection{Cutting and gluing}
\label{subsec:cut:glue}

Before describing the Dehn--Thurston coordinates, which in this setting we will call combinatorial Fenchel--Nielsen coordinates, we need the notion of cutting a combinatorial structure along an essential simple closed curve, and the reciprocal notion of gluing combinatorial structures along boundary components of the same length. The main difference with \cite{FLP12} is in the gluing: the combinatorial Teichmüller space does not contain measured foliations with saddle connections, but saddle connections can be created from the gluing process. The main result of this section, Proposition~\ref{prop:open:set:of:tau}, is the analysis of these ``pathological gluings'', which turn out to occur only for a negligible set of twists.

\subsubsection{Basic definitions}

\textsc{Cutting.} Consider a bordered surface $\Sigma$, fix $\GG \in \mathcal{T}^{\textup{comb}}_{\Sigma}$ and $\gamma$ an essential simple closed curve. We want to define a combinatorial structure on the surface $\Sigma_{\gamma}$ obtained by cutting $\Sigma$ along a chosen representative of $\gamma$. To this end, choose a representative $(\bm{G},f)$, so that we have an induced structure of measured foliation on $\Sigma$. If necessary, perform a minimal sequence of local Whitehead moves in small disc neighbourhoods of the vertices, in such a way that $\gamma$ is transversal to the resulting foliation. We then restrict the measured foliation to $\Sigma_{\gamma}$, which is induced from a unique metric ribbon graph $\bm{G}_{\gamma}$ with an embedding which up to isotopy does not depend on the choices made. This defines a combinatorial structure $\GG_{\gamma} \in \mathcal{T}^{\textup{comb}}_{\Sigma_{\gamma}}$.

\medskip

\begin{figure}[ht]
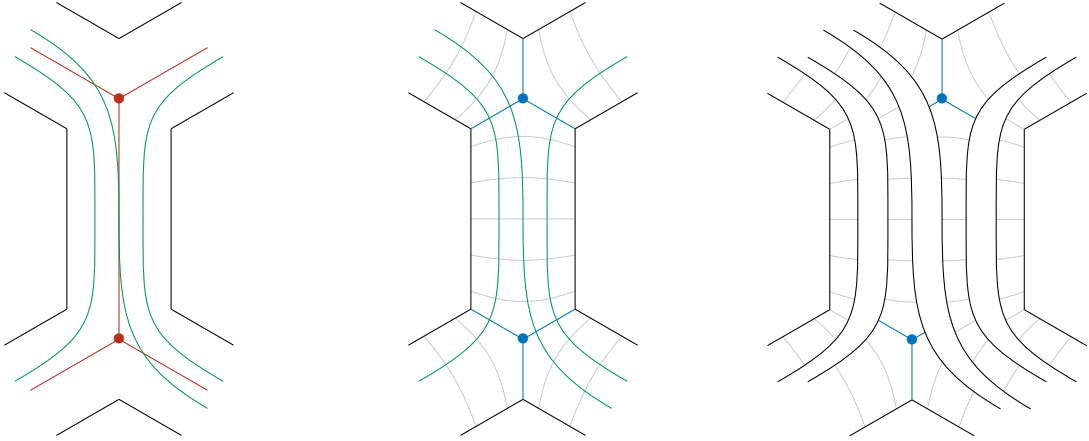

\centering
	\begin{subfigure}[t]{.32\textwidth}
	\centering

	\end{subfigure}
	\caption{Cutting/gluing algorithm: the combinatorial structure $\GG$ (in red) pictured with the sin\-gu\-lar leaves (in blue), and the curve $\gamma$ (in green).
	\label{fig:cut}}
\end{figure}

\begin{figure}[ht]
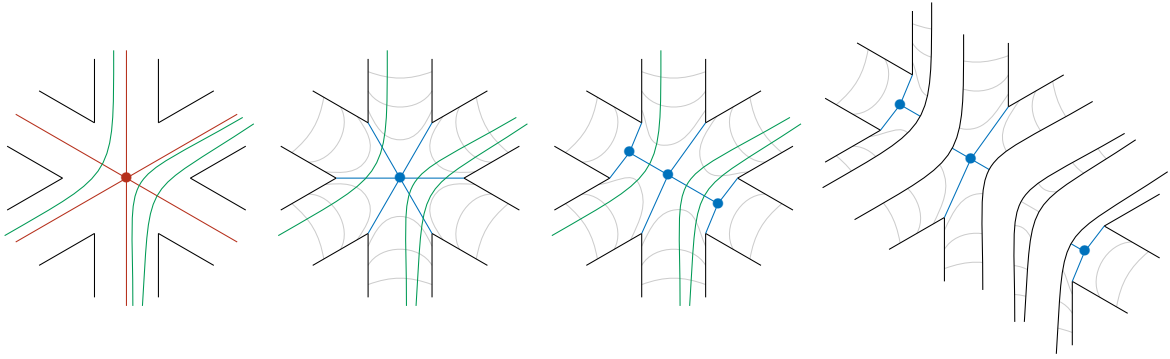

\centering
	\begin{subfigure}{.215\textwidth}
	\centering
	%
	\end{subfigure}
	\caption{Cutting/gluing algorithm for vertices of higher valency. Two Whitehead moves are performed.}
	\label{fig:cut:high:valency}
\end{figure}

Cutting also makes sense when $\gamma$ is a primitive multicurve, and it is equivalent to cutting along each component of $\gamma$ in an arbitrary order. Note that the lengths of edges after cutting are again linear combinations of the edge lengths which agree on the closure of the open cells. This shows that the cutting, viewed as a map $\mathcal{T}^{\textup{comb}}_{\Sigma}\rightarrow\mathcal{T}^{\textup{comb}}_{\Sigma_{\gamma}}$, is continuous. See Figure~\ref{fig:cut} and Figure~\ref{fig:cut:high:valency} for a local illustration of the cutting, and Appendix~\ref{app:cut:gluing} for some global examples.


\medskip

\textsc{Gluing.} Consider a bordered surface $\Sigma$, possibly disconnected, with a choice of two  boundary components $\gamma_{-}$ and $\gamma_{+}$. Let $\GG \in \mathcal{T}^{\textup{comb}}_{\Sigma}$ be such that $\ell_{\GG}(\gamma_{-}) = \ell_{\GG}(\gamma_{+})$. We want to define a combinatorial structure on the surface obtained by topologically gluing $\gamma_{-}$ and $\gamma_{+}$. Fix a representative $(\bm{G},f)$, so that we have an induced structure of a measured foliation $\mathcal{F}$ on $\Sigma$. First, we observe that once we pick a point $p_-$ on $\gamma_-$, there is a unique action of $\RR$ on $\gamma_-$ which preserves the induced measure and orientation on $\gamma_-$. We let $p_-^\tau$ be the image of $p_-$ under the action of $\tau \in \RR$. Pick now a point $p_+$ on $\gamma_+$, and identify $\gamma_-$ with $\gamma_+$ in a measure preserving way, such that $p^\tau_-$ is identified with $p_+$ in an orientation reversing way. This means that we have a unique measured foliation $\mathcal{F}^\tau$ induced on the glued surface, which we denote $\Sigma^\tau$.

\subsubsection{Admissible gluings}

\begin{figure}[ht]
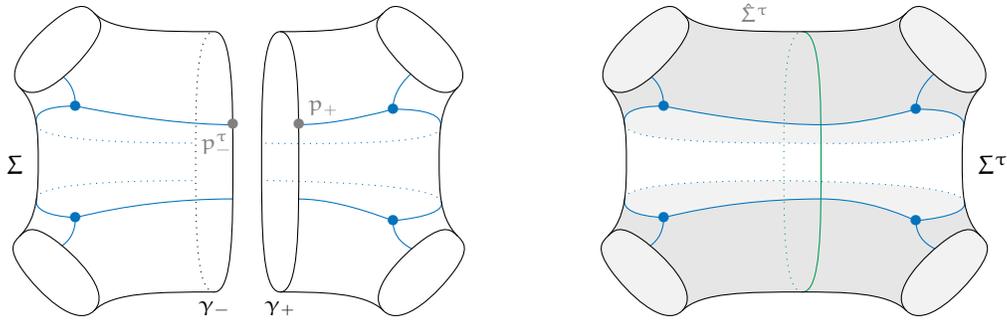

\centering
	\begin{subfigure}[t]{.45\textwidth}
	\centering

	\end{subfigure}
	\caption{A glued measured foliation that is not dual to a combinatorial structure. Notice in grey $\hat{\Sigma}^{\tau}$, that is properly contained in $\Sigma^{\tau}$, and the presence of saddle connections on the boundary of $\hat{\Sigma}^{\tau}$ that cannot be removed by Whitehead moves.}
	\label{fig:bad:twist}
\end{figure}

What is not clear from the above construction is whether the measured foliation $\mathcal{F}^\tau$ is associated to a combinatorial structure on $\Sigma^{\tau}$. If this is true, we call such $\tau$ an \emph{admissible twist}. We refer to Figure~\ref{fig:cut} and Figure~\ref{fig:cut:high:valency} -- read from right to left -- for a local illustration of the gluing, Appendix~\ref{app:cut:gluing} for some global examples, and Figure~\ref{fig:bad:twist} for an example of $\mathcal{F}^{\tau}$ that is not associated to a combinatorial structure. 

\begin{prop}\label{prop:when:gluing:breaks}
	There exists a unique metric ribbon graph $\bm{G}^\tau$ and a unique marking $f \colon \Sigma^\tau \to |\bm{G}^{\tau}|$ up to isotopy, such that the measured foliation induced on $\Sigma^\tau$ agrees with $\mathcal{F}^\tau$ if and only if $\mathcal{F}^\tau$ has a representative without saddle connections, \emph{i.e.} no leaf between two singularities.
\end{prop}

\begin{proof}
	Perform a maximal sequence of Whitehead moves, \emph{i.e.} that reduces the connected components of the the compact singular leaves to a graph with one vertex. Let $\Lambda(\mathcal{F}^\tau)$ be the set of leaves of $\mathcal{F}^\tau$, and define $\hat{\Sigma}^\tau = \set{ \lambda \in \Lambda(\mathcal{F}^\tau) | \lambda \cap \de \Sigma^\tau \neq \varnothing }$ (\emph{cf.} Figure~\ref{fig:bad:twist}). Then from Poincar\'{e} recurrence \cite[Theorem~5.2]{FLP12} this is nothing but the union of all leaves which go from boundary to boundary together with the finitely many leaves which connect the boundary to a singular point of $\mathcal{F}^\tau$ (\emph{i.e.} no leaves starting from the boundary spiral in the surface). If $\hat{\Sigma}^\tau = \Sigma^\tau$, then $\mathcal{F}^\tau$ has no closed singular leaves and we see that the singular leaves of $\mathcal{F}^\tau$ split the surface into hexagons, which in turn determines $\bm{G}^\tau$ uniquely and its marking up to isotopy. If not, choose a good atlas for $\mathcal{F}^\tau$ as defined in \cite[Section~5.2]{FLP12}, and observe that the complement of the singular leaves in $\hat{\Sigma}^\tau$ is a finite disjoint union of squares $S_i$, each with a non-singular foliation transverse to two open arcs of the boundary of $\Sigma^\tau$, running between endpoints of the singular leaves in $\hat{\Sigma}^\tau$, and such that $\bar S_i - S_i \subset \Sigma^\tau$ are made up of a finite number of compact singular leaves of $\mathcal{F}^\tau$. But since we are assuming $\hat{\Sigma}^\tau \neq \Sigma^\tau$, there must be at least one of these which connects two singular points in the interior of $\Sigma^\tau$. As we took a representative of $\mathcal{F}^\tau$ with one singular point for each connected component of the compact singular leaves, we see that this implies there must be a cycle of singular leaves which cannot define a combinatorial structure.
\end{proof}

We observe that, even though the foliation associated to $\GG$ has no saddle connections, they may occur for $\mathcal{F}^\tau$ as a result of the gluing process. However, Proposition~\ref{prop:when:gluing:breaks} together with the next result imply that this is generically not the case, and that it is never the case for $\GG$-strictly small pairs of pants.

\begin{lem}\label{lem:all:twists}
	Let $\GG$ be such that every vertex around $\gamma_-$ has exactly one singular leaf reaching $\gamma_-$ and no singular leaf reaching $\gamma_+$. Then all twists $\tau \in \RR$ are admissible. The same statement holds if we exchange $\gamma_-$ and $\gamma_+$.
\end{lem}
\begin{proof}
	Under such hypothesis, all smooth leaves starting from $\gamma_-$ end at a boundary component of $\Sigma$ that is neither $\gamma_-$ nor $\gamma_+$. Therefore for every twist $\tau$, we glue the singular leaves to leaves of $\gamma_+$ which immediately reach a boundary component of $\de\Sigma^{\tau}$. Gluing the singular leaves of $\gamma_-$ can result in the creation of a leaf which returns back to $\gamma_+$ again. This however corresponds two gluing to leaves of $\gamma_+$ which immediately reach a boundary component of $\de\Sigma^{\tau}$.
\end{proof}

\begin{cor}\label{cor:twist:small:pop}
	Let $\Sigma$ be a bordered surface of Euler characteristic $\leq -2$, take $[P] \in \mathcal{P}_{\Sigma,m_0}$ for some $m_0 \in \{1,\ldots,n\}$, $\GG \in \mathcal{T}_{\Sigma}^{\textup{comb}}$, and consider the operation of cutting along $\partial P \cap \Sigma^{\circ}$, twisting and gluing back. If $P$ is $\GG$-strictly small, then any twist $\tau \in \mathbb{R}^{\pi_0(\partial P \cap \Sigma^{\circ})}$ is admissible. The same is true if $\Sigma = T$ and we self-glue after twisting the pair of pants obtained from $T$ by cutting along $\gamma \in S_{T}^{\circ}$.
		\hfill $\blacksquare$
\end{cor}
\begin{prop}\label{prop:open:set:of:tau}
	For $\GG \in \mathcal{T}^{\textup{comb}}_{\Sigma}$ with $l = \ell_{\GG}(\gamma_{-}) = \ell_{\GG}(\gamma_{+})$ the set of admissible twists is an open dense subset of $\RR$ with countable complement.
\end{prop}

We need the result of the following lemma before proving the proposition. Let $S^- = \set{\ell_1^-,\dots,\ell_M^-}$ and $S^+ = \set{\ell_1^+,\dots,\ell_N^+}$ be the finite sets of lengths of edges in $\GG$ into which $\gamma_-$ and $\gamma_+$ decompose into respectively. Without loss of generality, assume that $p_-$ and $p_+$ are contained in singular leaves of $\GG$.

\begin{lem}\label{lem:count:tau:bad}
	If we choose the points $p_{\pm}$ such that $\tau = 0$ identifies two singular leaves, then for $\tau \notin {\rm span}_{\QQ}(S^- \cup S^+)$, $\mathcal{F}^\tau$ has a Whitehead representative with no singular leaf between two singularities.
\end{lem}

\begin{figure}[ht]
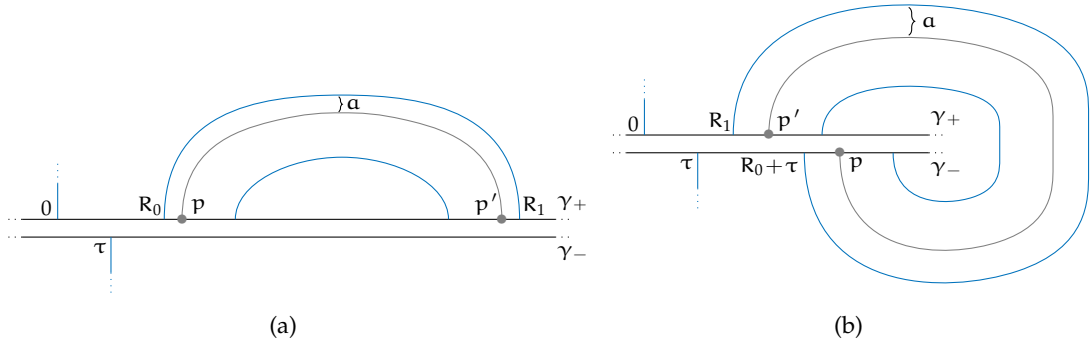

\centering
	\begin{subfigure}[t]{.45\textwidth}
	\centering

	\subcaption{}
	\label{fig:leaves:dynamic:a}
	\end{subfigure}
	\begin{subfigure}[t]{.45\textwidth}
	\centering
	%
	\subcaption{}
	\label{fig:leaves:dynamic:b}
	\end{subfigure}
	\caption{Examples of leaf dynamics, induced on $\gamma$ by the foliation $\mathcal{F}$. The singular leaves before gluing are depicted in blue, the leaf connecting $p$ and $p'$ in grey (here it is depicted as a smooth leaf, \emph{i.e.} $a > 0$). In between the two fully depicted singular leaves there is a strip of smooth leaves all homotopic to each other.}
	\label{fig:leaves:dynamic}
\end{figure}

\begin{proof}
	Denote by $\gamma$ the curve in $\Sigma^{\tau}$ which is the image of $\gamma_{\pm}$, and identify $\gamma \sim \RR/l\ZZ$ where $l = \ell_{\GG}(\gamma_-) = \ell_{\GG}(\gamma_+)$ and $0$ corresponds to a singular leaf on the $\gamma_{+}$ side. Consider a point $p \in \gamma$. If we follow along the leaf passing through $p$ (in either the $\gamma_{+}$ side or the $\gamma_{-}$ side of the glued surface) until it gets back to $\gamma$ at a point $p'$, we find that there exists some $R \in {\rm span}_{\ZZ}(S^- \cup S^+)$ for each of the following cases, such that
	\begin{itemize}
		\item $p' = - p + R$, for a leaf going from $\gamma_{+}$ to $\gamma_{+}$;
		\item $p' = - p + 2\tau + R$, for a leaf going from $\gamma_{-}$ to $\gamma_{-}$;
		\item $p' = p + \tau + R$, for a leaf going from $\gamma_{+}$ to $\gamma_{-}$;
		\item $p' = p - \tau + R$, for a leaf going from $\gamma_{-}$ to $\gamma_{+}$.
	\end{itemize}
	Indeed, we firstly notice that all singular leaves on the $\gamma_{+}$ side are identified as some points in ${\rm span}_{\ZZ}(S^+)$, while on the $\gamma_{-}$ side they are identified as some points in $\tau + {\rm span}_{\ZZ}(S^-)$.

	\smallskip

	Suppose first that $p'$ is obtained from $p$ by following a leaf going from $\gamma_{+}$ to $\gamma_{+}$ (see Figure~\ref{fig:leaves:dynamic:a}). We notice that $p$ is given by $p = R_{0} + a$, where $R_{0} \in {\rm span}_{\ZZ}(S^+)$ is the distance from the chosen singular leaf at $0$ to the singular leaf just before $p$ on the $\gamma_{+}$ side, and $a \ge 0$. Then, following the leaf, we find that the singular leaf just before $p$ becomes the singular leaf just after $p'$, and $p' = R_{1} - a = - p + (R_{0} + R_{1})$ where $R_{1} \in {\rm span}_{\ZZ}(S^+)$ is the distance from the chosen singular leaf at $0$ to the leaf just after $p'$ (following the orientation of $\gamma_{+}$). In particular, we obtain the claim with $R = R_0 + R_1 \in {\rm span}_{\ZZ}(S^+)$.

	\smallskip

	Similarly, suppose now that $p'$ is obtained from $p$ by following a leaf going from $\gamma_{-}$ to $\gamma_{+}$ (see Figure~\ref{fig:leaves:dynamic:b}). Now we have $p = R_{0} + \tau + a$, where $R_0 \in {\rm span}_{\ZZ}(S^-)$ and $a > 0$ (here the singular leaf just before $p$ on the $\gamma_{-}$ side is at distance $R_{0} + \tau$ from the chosen singular leaf at $0$). Then, following the leaf, we find that the singular leaf at $R_{0} + \tau$ is identified with a singular leaf at distance $R_{1} \in {\rm span}_{\ZZ}(S^+)$ from the chosen singular leaf at $0$. Therefore $p' = R_{1} + a = p - \tau + (R_{0} - R_{1})$. Thus, the claim with $R = R_{0} - R_{1} \in {\rm span}_{\ZZ}(S^- \cup S^+)$. 

	\smallskip

	The other cases follow similarly. Now, if $p$ is a point at a singular leaf on the $\gamma_{\pm}$ side, we see by induction that, after gluing, the singular leaf passes through $\gamma$ at some other points of the form $\pm( p \mp n\tau) + R$ for some $n \in \ZZ_{+}$ and $R \in {\rm span}_{\ZZ}(S^- \cup S^+)$. This implies that, if $\GG^\tau$ has two singular points connected by a leaf, then a non-zero integral multiple of $\tau$ is contained in ${\rm span}_{\ZZ}(S^- \cup S^+)$, or equivalently $\tau \in {\rm span}_{\QQ}(S^- \cup S^+)$.
\end{proof}

\begin{figure}[t]
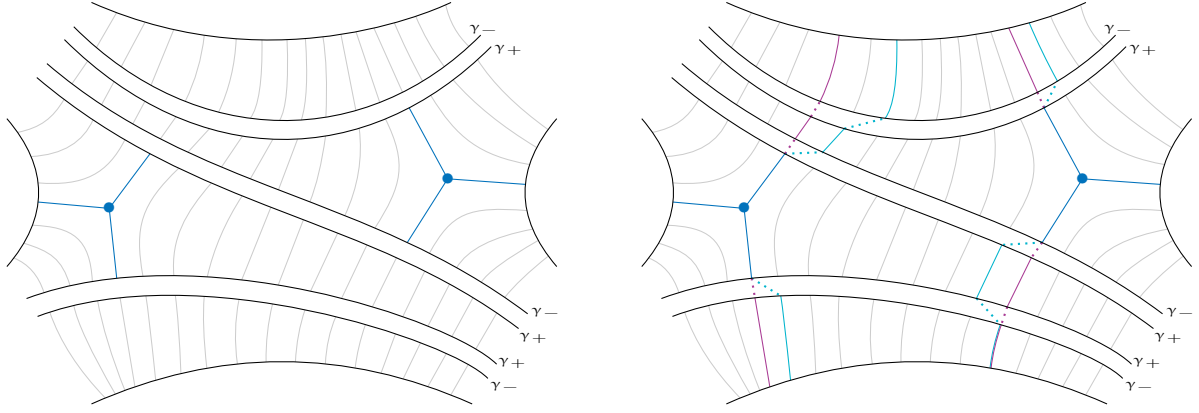

\centering
	\begin{subfigure}[t]{.44\textwidth}
	\centering

	\end{subfigure}
	\caption{A combinatorial structure $\GG$ glued to $\GG^{\tau}$ and $\GG^{\tau'}$ for $|\tau-\tau'|$ small. The singular leaves of $\GG$ are shown in blue. After gluing, they are prolonged with the purple leaves in $\GG^{\tau}$, and with the light blue leaves in $\GG^{\tau'}$. The dotted lines indicate the identification of $\gamma_-$ and $\gamma_+$ for $\tau$ and $\tau'$.}
	\label{fig:twist}
\end{figure}

\begin{proof}[Proof of Proposition~\ref{prop:open:set:of:tau}]\label{proof:open:set:of:tau}
	We first show that the set of admissible twists is an open subset of $\RR$. Consider $\tau \in \RR$ an admissible twist, and denote by $\GG^\tau$ the associated combinatorial structure. For each edge $e$ of $\GG^\tau$, let $n_\gamma(e)$ be the number of times which $\gamma$ travels through the edge $e$. If we take $\tau'$ such that $|\tau'-\tau| < \epsilon$, then we see that the distance between the two singularities of the foliation changes by at most $\epsilon n_{\gamma}(e)$ (\emph{cf.} Figure~\ref{fig:twist}). This also holds at the boundary of the cells, when we have vertices of higher valency whose original distance would be zero. Therefore, if we choose $\epsilon > 0$ smaller than
	\[
		\min_{e}{\frac{\ell_{\GG}(e)}{n_{\gamma}(e)}}
	\]
	where $e$ runs over the edges of $\GG^\tau$ visited by $\gamma$, then all singularities stay at the same or at a positive distance from each other. These lengths are realised by curves homotopic to the original length realising curve. As a consequence, $\mathcal{F}^{\tau'}$ cannot admit a cycle of singular leaves connecting singularities, and thus $\tau'$ is an admissible twist. The countable complement property follows from Lemma~\ref{lem:count:tau:bad}.
\end{proof}

We remark that the set of non admissible twists can have accumulation points, and its set of accumulation points can be non isolated. However what is crucial for the next section is that the non-admissible twists form a measure-zero set in $\RR$.

\subsection{Combinatorial Fenchel--Nielsen coordinates}
\label{subsec:FN:coords}

With the notions of cutting and gluing in the combinatorial spaces defined, we have available the key tools to adapt the definition of Dehn--Thurston coordinates to our framework. The main difference with the hyperbolic setting is that the image of such a coordinate system does not cover the whole codomain, and this is due to the fact that combinatorially certain (rare) values of the twists are forbidden.

\subsubsection{Seams and pants decompositions}

Firstly, we need a technical ingredient, the pants seams, that allows us to define a canonical way of gluing pairs of pants.

\begin{figure}
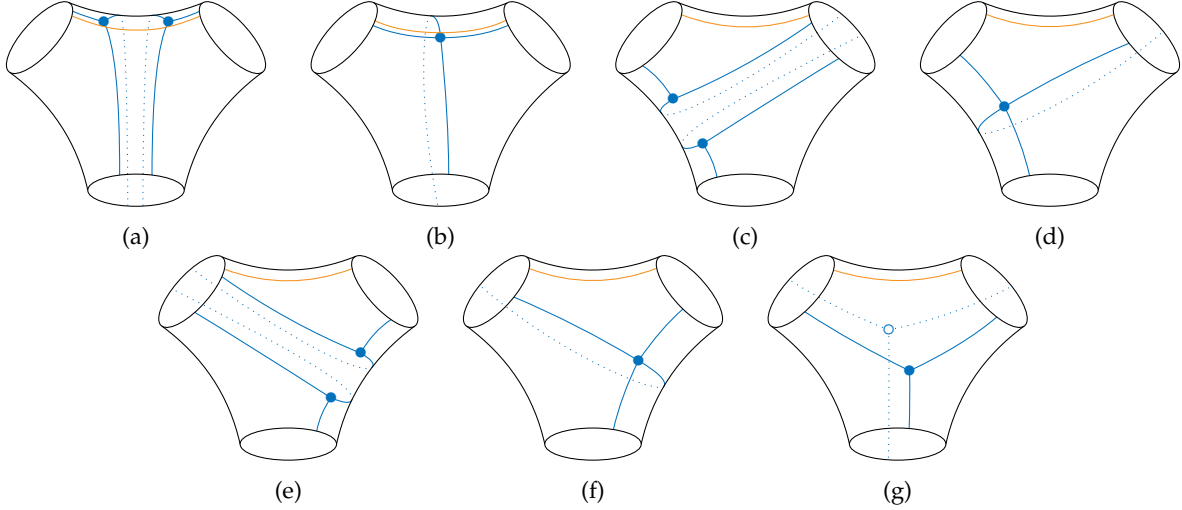

	\centering
	\begin{subfigure}[t]{.24\textwidth}
	\centering

	\caption{}
	\label{fig:seams:pop:a}
	\end{subfigure}
	\begin{subfigure}[t]{.24\textwidth}
	\centering
		%
	\caption{}
	\label{fig:seams:pop:b}
	\end{subfigure}
	\begin{subfigure}[t]{.24\textwidth}
	\centering
		%
	\caption{}
	\label{fig:seams:pop:c}
	\end{subfigure}
	\begin{subfigure}[t]{.24\textwidth}
	\centering
		%
	\caption{}
	\label{fig:seams:pop:d}
	\end{subfigure}
	\begin{subfigure}[t]{.24\textwidth}
	\centering
		%
	\caption{}
	\label{fig:seams:pop:e}
	\end{subfigure}
	\begin{subfigure}[t]{.24\textwidth}
	\centering
		%
	\caption{}
	\label{fig:seams:pop:f}
	\end{subfigure}
	\begin{subfigure}[t]{.24\textwidth}
	\centering
		%
	\caption{}
	\label{fig:seams:pop:g}
	\end{subfigure}
	\caption{Combinatorial seams (in orange) on each cell of $\mathcal{T}^{\textup{comb}}_{P}$. The singular leaves are depicted in blue.}
	\label{fig:seams:pop}
\end{figure}

\begin{defn}\label{defn:seams}
	Consider a combinatorial marking $(\bm{G},f)$ on a pair of pants $P$, with associated foliation $\mathcal{F}$. Define the \emph{combinatorial seam} connecting two distinct boundary components $\gamma$ and $\gamma'$ of $P$ to be the quasitransverse arc connecting $\gamma$ and $\gamma'$, as indicated in Figure~\ref{fig:seams:pop}. In the cases \ref{fig:seams:pop:c}--\ref{fig:seams:pop:g}, the seams are smooth leaves, located at exactly the same distance from the adjacent singular leaves.
\end{defn}

A notion of pants seams connecting a boundary component of $P$ to itself can be found in \cite[Section~6.3]{FLP12}. We remark that the combinatorial seam realises the minimum among lengths of all essential arcs connecting one boundary component to another (the length can be zero). For a point $(L_1,L_2,L_3) \in \RR_{+}^3 \cong \mathcal{T}^{\textup{comb}}_{P}$, the length of the seam connecting $\de_1 P$ and $\de_2 P$ is given by the formula
\begin{equation}\label{eqn:comb:length:seam}
	\ell_{\textup{comb}}(L_1,L_2,L_3) = \left[ \frac{L_3 - L_1 - L_2}{2} \right]_{+},
\end{equation}
while the length of a seam connecting $\de_1 P$ to itself is given by
\begin{equation}
\begin{aligned}\label{eqn:comb:length:seam:same:boundary}
	\ell_{\textup{comb}}(L_1,L_2,L_3)
		& =
			\left[ \frac{L_2 + L_3 - L_1}{2} \right]_{+}
			+
			\left[ \frac{L_2 - L_1 - L_3}{2} \right]_{+}
			+
			\left[ \frac{L_3 - L_1 - L_2}{2} \right]_{+} \\
		& = \max\left\{ \frac{L_2 + L_3 - L_1}{2}, L_2 - L_1, L_3 - L_1, 0 \right\}
\end{aligned}
\end{equation}
Notice that for a hyperbolic marking $(X,\varphi)$ on $P$, there exist a notion of hyperbolic seam connecting $\gamma$ and $\gamma'$, that is the shortest geodesic arc connecting the boundary components $\gamma$ and $\gamma'$. On the other hand, we can consider the combinatorial marking $(\bm{G},f)$ on $P$ associated to $(X,\varphi)$ defined through the spine map of Definition~\ref{defn:spine:map}. The next elementary lemma, of which we omit the proof, shows that the hyperbolic and combinatorial seams are the same arcs.

\begin{lem}\label{lem:hyp:comb:seam}
	Consider a hyperbolic marking $(X,\varphi)$ on $P$, and the associated combinatorial marking $(\bm{G},f)$. Through their aforementioned identification, the hyperbolic and combinatorial seams connecting two boundary components of $P$ coincide.
	\hfill $\blacksquare$
\end{lem}
\begin{rem}
	In the hyperbolic case, for a point $(L_1,L_2,L_3) \in \RR_{+}^3 \cong \mathcal{T}_{P}$, the hyperbolic length of the seam connecting $\de_1 P$ and $\de_2 P$ is given by the formula
	\begin{equation}\label{eqn:hyp:length:seam}
		\cosh\left( \ell_{\textup{hyp}}(L_1,L_2,L_3) \right)
			=
			\frac{\cosh\left(\tfrac{L_3}{2}\right)}{\sinh\left(\tfrac{L_1}{2}\right)\sinh\left(\tfrac{L_2}{2}\right)}
			+
			\cotanh\left(\tfrac{L_1}{2}\right)\cotanh\left(\tfrac{L_2}{2}\right),
	\end{equation}
	while the hyperbolic length of a seam connecting $\de_1 P$ to itself is given by
	\begin{equation}\label{eqn:hyp:length:seam:same:boundary}
		\cosh^2\left( \frac{\ell_{\textup{hyp}}(L_1,L_2,L_3)}{2} \right)
		=
		\frac{
			\cosh^2\left( \tfrac{L_1}{2} \right) + \cosh^2\left( \tfrac{L_2}{2} \right) + \cosh^2\left( \tfrac{L_3}{2} \right)
			+
			2 \cosh\left( \tfrac{L_1}{2} \right) \cosh\left( \tfrac{L_2}{2} \right) \cosh\left( \tfrac{L_3}{2} \right)
			-
			1
		}{\sinh^2\left(\tfrac{L_1}{2}\right)}.
	\end{equation}	
	Formulae \eqref{eqn:comb:length:seam}--\eqref{eqn:comb:length:seam:same:boundary} and  can be recovered from equations \eqref{eqn:hyp:length:seam}--\eqref{eqn:hyp:length:seam:same:boundary} by taking the following limit:
	\begin{equation}
		\lim_{\beta\to\infty} \frac{\ell_{\mathrm{hyp}}(\beta L_1,\beta L_2,\beta L_3)}{\beta} = \ell_{\textup{comb}}(L_1,L_2,L_3).
	\end{equation}
	This fact is revisited and generalised in Section~\ref{sec:hyp:comb} where such a limit is shown to hold for many other expressions.
\end{rem}

We can now define the seamed pants decomposition associated to a bordered surface $\Sigma$, consisting of a pants decomposition of $\Sigma$ together with a collection of curves and arcs.

\begin{samepage}
	\begin{defn}
	Given a bordered surface $\Sigma$ of type $(g,n)$, a \emph{seamed pant decomposition} consists of
	\begin{itemize}
		\item
			a pants decomposition $\mathscr{P} = (\gamma_1,\dots,\gamma_{3g-3+n})$, that is a maximal collection of pairwise non-homotopic, essential, simple closed curves, labelled by $1, \dots, 3g-3+n$,
		\item
			a collection $\mathscr{S}$ of non-homotopic, essential simple closed curves or simple arcs connecting boun\-da\-ry components of $\Sigma$, pairwise non-homotopic relative to the boundary, such that the intersection of $\mathscr{S}$ with any of the pair of pants $P$ in the decomposition specified by $\mathscr{P}$ is a union of three disjoint arcs connecting the boundary components of $P$ pairwise.
	\end{itemize}
\end{defn}
\end{samepage}
Notice that, given $\mathscr{P}$, we can construct an $\mathscr{S}$ by first choosing three disjoint arcs on each pair of pants and then matching up endpoints in any fashion.

\subsubsection{Combinatorial Fenchel--Nielsen coordinates as a.e. global coordinates}

Fix once and for all a seamed pants decomposition $(\mathscr{P},\mathscr{S})$ on $\Sigma$. We define the \emph{length parameters} of a point $\GG \in \mathcal{T}^{\textup{comb}}_{\Sigma}$ to be the tuple of positive real numbers
\begin{equation}
	\ell(\GG) = \bigl( \ell_1(\GG),\dots,\ell_{3g-3+n}(\GG) \bigr),
\end{equation}
where $\ell_i(\GG) = \ell_{\GG}(\gamma_{i})$.

\medskip

As a first step towards the definition of twist parameters, consider a combinatorial marking $(\bm{G},f)$ of a pair of pants $P$ and an arc $\alpha$ connecting two distinct boundary components $\gamma$ and $\gamma'$ of $P$. Let $\delta$ be the combinatorial seam connecting $\gamma$ and $\gamma'$ -- which depends only on $(\bm{G},f)$. Let $\tilde{P}$ be a universal cover of $P$. It contains lifts $\tilde{\gamma}$ and $\tilde{\gamma}'$ of $\gamma$ and $\gamma'$ respectively. Notice that they acquire orientation from $\tilde{P}$. Let $d = \delta \cap \gamma$ and $a = \alpha \cap \gamma$. We choose a lift $\tilde{d}$ of $d$ and call $\tilde{a}$ the first lift of $a$ met by travelling from $\tilde{d}$ along $\tilde{\gamma}$ following its orientation. This determines lifts $\tilde{\delta}$ (resp. $\tilde{\alpha}$) of $\delta$ (resp. $\alpha$) starting from $\tilde{d}$ (resp. $\tilde{a}$). Now let $\tilde{d}' = \tilde{\delta}\cap\tilde{\gamma}'$ and $\tilde{a}' = \tilde{\alpha} \cap \tilde{\gamma}'$. Consider the path $c_{\tilde{d}'\tilde{a}'}$ along $\tilde{\gamma}'$ starting at $\tilde{d}'$ and ending at $\tilde{a}'$. The measured foliation associated to $(\bm{G},f)$ lifts to a measured foliation $\tilde{\mathcal{F}}$ on the universal cover, and we can measure the length of $c_{\tilde{d}'\tilde{a}'}$. We then set $\sgn(c_{\tilde{d}'\tilde{a}'}) = \pm 1$ depending on whether the orientation of $c_{\tilde{d}'\tilde{a}'}$ agrees with the one of $\tilde{\gamma}'$. We define the twisting number of $\alpha$ along $\gamma'$ in $P$ to be
\begin{equation}\label{eqn:twist:boundary}
	t_{\alpha,\gamma'}^{P}(\bm{G},f)
	=
	\sgn(c_{\tilde{d}'\tilde{a}'}) \, \ell_{\tilde{\mathcal{F}}}(c_{\tilde{d}'\tilde{a}'}).
\end{equation}
The definition does not depend on the choice of $\tilde{d}$, since all different choices are related by deck transformations which leave $t_{\alpha,\gamma'}^{P}(\bm{G},f)$ fixed.

\medskip

Given $\GG \in \mathcal{T}^{\textup{comb}}_{\Sigma}$, we define the $i$-th twist parameter $\tau_i(\GG)$ as follows. Fix a marking $(\bm{G},f)$ such that $\gamma_i$ is quasitransverse to the measured foliation induced by the marking. Let $\alpha_i$ be one of the two arcs in $\mathscr{S}$ crossing $\gamma_i$. There are two pairs of pants $Q_i'$ and $Q''_i$ (possibly the same) on each side of $\gamma_i$, and $\alpha_i$ determines two arcs $\alpha_i' = \alpha_i \cap Q'_i$ and $\alpha_i'' = \alpha_i \cap Q''_i$. The $i$-th twist parameter of $\GG$ is defined to be
\begin{equation}\label{eqn:twist:numb}
	\tau_i(\GG) = t_{\alpha'_i,\gamma_i}^{Q_i'}(\bm{G}|_{Q_i'},f|_{Q_i'}) + t_{\alpha''_i,\gamma_i}^{Q_i''}(\bm{G}|_{Q_i''},f|_{Q_i''}).
\end{equation}
This twist parameter is invariant under isotopies, \emph{i.e.} does not depend on the representative of $\GG$. Besides, it does not depend on the choice of the arc in $\mathscr{S}$ crossing $\gamma_i$. This can be seen by passing to the universal cover of a neighbourhood of $\gamma_i$ -- \emph{cf.} \cite[Section~10.6.1]{FM11} for the analogue in the hyperbolic case. Finally, it only depends on the homotopy class of $\alpha_i$, since a different choice of representative would modify both $t'$ and $t''$ by the same quantity, but with different signs. Thus, we have well-defined \emph{twist parameters}
\begin{equation}
	\tau(\GG) = \bigl( \tau_1(\GG),\dots,\tau_{3g-3+n}(\GG) \bigr).
\end{equation}
Notice that we may homotope the representative of $\alpha_i$ such that it is intersecting $\gamma_i$ at a vertex of the combinatorial structure, showing that $t'$ and $t''$ in \eqref{eqn:twist:numb} can be expressed as a sum of edge lengths in $\GG$ with half-integer coefficients. The half is coming from the definition of combinatorial seams, which were required to be equidistant from the adjacent singular leaves in the cells depicted in Figure~\ref{fig:seams:pop:c}--\ref{fig:seams:pop:g}.

\begin{defn}
	Let $\Sigma$ be a bordered surface of type $(g,n)$ equipped with a seamed pants decomposition $(\mathscr{P},\mathscr{S})$. \emph{Combinatorial Fenchel--Nielsen coordinates} relative to $(\mathscr{P},\mathscr{S})$  is by definition the map $\Phi_{L} \colon \mathcal{T}^{\textup{comb}}_{\Sigma}(L) \to \RR_{+}^{3g-3+n} \times \RR^{3g-3+n}$ defined by
	\begin{equation}
		\Phi_{L}(\GG)
		=
		\bigl(\ell(\GG), \tau(\GG) \bigr).
	\end{equation}
\end{defn}

Using the gluing we can establish the following result. The first part is an adaptation of arguments by Dehn \cite{Deh22}, Thurston \cite{FLP12} and Penner \cite{PH92}, who proved similar results for the set of multicurves, measured foliations and train tracks respectively. The second part, \emph{i.e.} the zero-measure statement, will be crucial in Section~\ref{subsec:integration} where we provide a general formula for the integral of mapping class group invariant functions over the combinatorial moduli space.

\begin{thm}\label{thm:FN:coordinates} ~
	\begin{enumerate}
		\item
		For any $L \in \RR_{+}^n$, the map
		\begin{equation}
			\Phi_{L} \colon
			\mathcal{T}^{\textup{comb}}_{\Sigma}(L) \to (\RR_{+} \times \RR)^{3g-3+n}
		\end{equation}
		is a homeomorphism onto its image.

		\item
		The image of $\Phi_{L}$ is an open dense subset whose complement has zero measure. Moreover, if $\Phi_{L}(\GG) = (\ell, \tau)$, then the image of the map $\tau_i$ restricted to
		\begin{equation}
			\Set{ \GG \in \mathcal{T}^{\textup{comb}}_{\Sigma}(L) |  \ell_j(\GG) = \ell_j \; \forall j,\;\tau_{k}(\GG)=\tau_{k}\text{ for }k\neq i}
		\end{equation}
		has a complement of zero measure in $\RR$.
	\end{enumerate}
\end{thm}

\begin{proof}
	To prove the theorem, we use the gluing to construct a partial inverse map. More precisely, define the partial map $\Psi_{L} \colon (\RR_{+} \times \RR)^{3g - 3 + n} \to \mathcal{T}^{\textup{comb}}_{\Sigma}(L)$ by setting
	\[
		\Psi_{L}\bigl(\ell, \tau \bigr)
		=
		\GG,
	\]
	where $\GG$ is defined as follows.
	\begin{itemize}
		\item For each pair of pants bounded by curves in $\mathscr{P}$, we assign boundary lengths defined by the $L$s and $\ell$s. This  determines a unique combinatorial structure on each pair of pants.
		\item We glue the pairs of pants along each $\gamma_i$ after twisting by $\tau_{i}$. The twist zero corresponds to gluing the combinatorial seams of the pairs of pants together.
	\end{itemize}
	By partial map we mean that $\Psi_L$ is not defined on the whole of $(\RR_{+} \times \RR)^{3g-3+n}$, as the gluing does not always define an embedded metric ribbon graph. Notice also that $\Psi_{L}$ does not depend on the order on which we glue the pairs of pants together.

	\smallskip

	We proceed now with the proof. Firstly notice that the definition of the twist parameters implies that gluing with twist zero amounts to gluing all pairs of pants with matching seams. Also, since gluing back a cut combinatorial structure gives back the original one, we can see that $\Psi_{L}$ is defined on the image of $\Phi_{L}$ and that $\Psi_{L} \circ \Phi_{L}$ is the identity on $\mathcal{T}^{\textup{comb}}_{\Sigma}(L)$. Hence, $\Phi_{L}$ is a bijection onto its image.

	\smallskip	
	
	On the closure of each open cell, the length and twists are linear functions of the edge lengths. Therefore, we have bijective linear functions that agree on boundaries of the open cells, and therefore, the inverse has the same properties and is therefore continuous which shows that $\Phi_{L}$ is a homeomorphism onto its image.

	\smallskip
	
	Now to prove openess we note that given any point $\GG\in\mathcal{T}_{\Sigma}^{\textup{comb}}(L)$ there exists a neighbourhood intersecting finitely many cells as there are finitely many ways to expand a singularity using Whitehead equivalence. We can therefore construct a finite simplicial complex containing $\GG$ as a vertex such that the intersection of a $k$-cell of $\mathcal{T}_{\Sigma}^{\textup{comb}}(L)$ is a union of $k$-dimensonal simplices. Then, as $\Phi_{L}$ is linear on each cell and a homeomorphism onto it's image, $\Phi_{L}$ maps the simplicial complex to a simplicial complex in $(\RR_{+} \times \RR)^{3g-3+n}$.
	
	\smallskip	

	A point on a finite simplicial complex in an Euclidean space of the same dimension is on the boundary if and only if it is contained in a codimension one simplex that is on the boundary of only one top-dimensional simplex. Every codimension one simplex containing $\GG$ is contained in two top-dimensional simplices and therefore $\Phi_{L}(\GG)$ is in the interior of the image of the simplex. Thus, $\Phi_{L}$ is open.

	\smallskip

	Finally, notice that by gluing one curve at a time and using Proposition~\ref{prop:when:gluing:breaks} and Lemma~\ref{lem:count:tau:bad} for each $L$ and $\ell$, we can see that the image is dense and its complement has zero measure.
\end{proof}

\subsection{A combinatorial \texorpdfstring{$(9g - 9 + 3n)$-theorem}{nine-g-minus-9-plus-3n theorem}}
\label{subsec:9g:minus:9:plus:3n}

In this paragraph we establish a combinatorial analogue of the hyperbolic $(9g - 9 + 3n)$-theorem (see \emph{e.g.} \cite[Theorem~10.7]{FM11}), that is, any combinatorial structure can be reconstructed from the data of the combinatorial lengths of $(9g - 9 + 3n)$ simple closed curves. Similar computations can be found in \cite[Expos\'e~6]{FLP12}, where only measured foliations on closed surfaces are considered. These results are used in Section~\ref{subsubsec:conv:twist} to compare the hyperbolic and combinatorial twists, and in Section~\ref{subsec:PL:structure} to give a new proof of Penner's formulae \cite{Pen82} for the action of the mapping class group on Dehn--Thurston coordinates.

\medskip

\begin{figure}
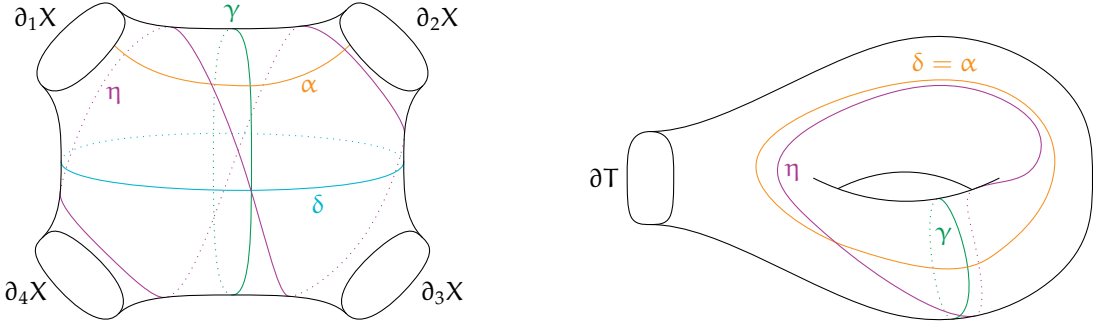

	\centering
	\begin{subfigure}[t]{.48\textwidth}
	\centering

	\end{subfigure}
	\caption{The curves $\delta$ and $\eta$ (we omit the subscript).}
	\label{fig:notation:delta:eta}
\end{figure}

Let $\Sigma$ be of type $(g,n)$ and fix a seamed pants decomposition $(\mathscr{P},\mathscr{S})$, with $\mathscr{P} = (\gamma_1,\dots,\gamma_{3g-3+n})$. The union of the pair of pants in the decomposition that are adjacent to $\gamma_i$ is a surface of type $(0,4)$ or $(1,1)$, and we choose $\alpha_i \in \mathscr{S}$ crossing $\gamma_i$ in $\Sigma_i$. We choose some order on the boundaries such that $\alpha_{i}$ connects $\partial_{1}X$ to $\partial_{2}X$ and $\partial_{4}X$ is in the same pair of pants as $\partial_{1}X$. We now define two other homotopy classes of curves in $\Sigma_i$ (see Figure~\ref{fig:notation:delta:eta}).
\begin{itemize}
	\item If $\Sigma_i$ has type $(0,4)$, we let $\delta_i$ be the curve determined by a tubular neighbourhood of $\alpha_i$ union the boundary component it connects. If $\Sigma_i$ has type $(1,1)$, we let $\delta_i$ be the curve $\alpha_i$.
	\item Let $\eta_i$ be the image of $\delta_i$ after a positive Dehn twist along $\gamma_i$.
\end{itemize}
In the $(0,4)$ case there are two possible choices of $\alpha_i$ as above but both choices give the same $(\delta_i,\eta_i)$.

\begin{thm}\label{thm:9g:minus:9:plus:3n}
	Let $\Sigma$ be a bordered surface of type $(g,n)$ and $(\mathscr{P},\mathscr{S})$ a seamed pants decomposition. The following map is continuous and injective:
	\begin{equation}
	\begin{aligned}
		\mathcal{T}_{\Sigma}^{\textup{comb}}(L) & \longrightarrow \RR_{+}^{9g-9+3n} \\
		\GG & \longmapsto \bigl( \ell_{\GG}(\gamma), \ell_{\GG}(\delta), \ell_{\GG}(\eta) \bigr).
	\end{aligned}
	\end{equation}
\end{thm}

As a preparation to the proof, we present in Lemmata~\ref{lem:delta:eta:04} and \ref{lem:delta:eta:11} closed formulae for $\ell_{\GG}(\delta_i)$ and $\ell_{\GG}(\eta_i)$ in the $(0,4)$ and $(1,1)$ cases respectively. For this purpose we can work locally on $\GG|_{\Sigma_i}$ with a fixed seamed pants decomposition, which we denote by $\ell_i = \ell_{\GG}(\gamma_i)$, $\ell_i' = \ell_{\GG}(\delta_i)$ and $\ell_i'' = \ell_{\GG}(\eta_i)$.

\subsubsection{Four-holed sphere}
\label{subsubsec:four:holed:sphere}

Let $\Sigma_i = X$ be a four-holed sphere. We remove the index $i$ from the notation of $\gamma_i$, $\alpha_i$, $\delta_i$ and $\eta_i$, as well as $\ell_i$, $\ell_i'$, $\ell_i''$ and $\tau_i$. Label by $\de_1 X, \dots, \de_4 X$ the boundary components of $X$, so that $\gamma$ is separating the components $\de_1 X$ and $\de_4 X$ from $\de_2 X$ and $\de_3 X$, and $\alpha$ is connecting the components $\de_1 X$ and $\de_2 X$. Finally, denote $L_i = \ell_{\GG}(\de_i X)$.

\begin{lem}\label{lem:delta:eta:04}
	In the above setting, we have
	\begin{equation}\label{eqn:length:delta:04}
		\ell'(\ell,\tau)
		=
		\max\big\{
			L_1 + L_3 - \ell,L_2 + L_4 - \ell,2|\tau| + M_{1,4}(\ell) + M_{2,3}(\ell)
		\big\},
	\end{equation}
	where $M_{i,j}(\ell) = \max\Big\{ 0,L_i - \ell,L_j - \ell,\frac{L_i + L_j - \ell}{2} \Big\}$. Further, $\ell''(\ell,\tau) = \ell'(\ell,\tau + \ell)$.
\end{lem}

\begin{proof}
	Let us assume that the ribbon graph underlying $\GG$ is trivalent, and fix a marking of it. Let us cut $\GG$ along the curve $\gamma$. There are sixteen possibilities for the cut combinatorial structure: the marked ribbon graph on each pair of pants can belong to each of the four top-dimensional cells of the Teichm\"uller space of a pair of pants. Therefore, in order to check that Equation~\eqref{eqn:length:delta:04} holds for any $\GG$, it is sufficient to check that it is satisfied in each of the sixteen cases. By symmetry considerations, the number of cases can actually be reduced to seven. We show the detailed argument for three particularly representative cases out of those seven, and argue that the other cases can be proven following the same strategy.

	\smallskip

	\afterpage{\clearpage
	\begin{landscape}
		\begin{figure}
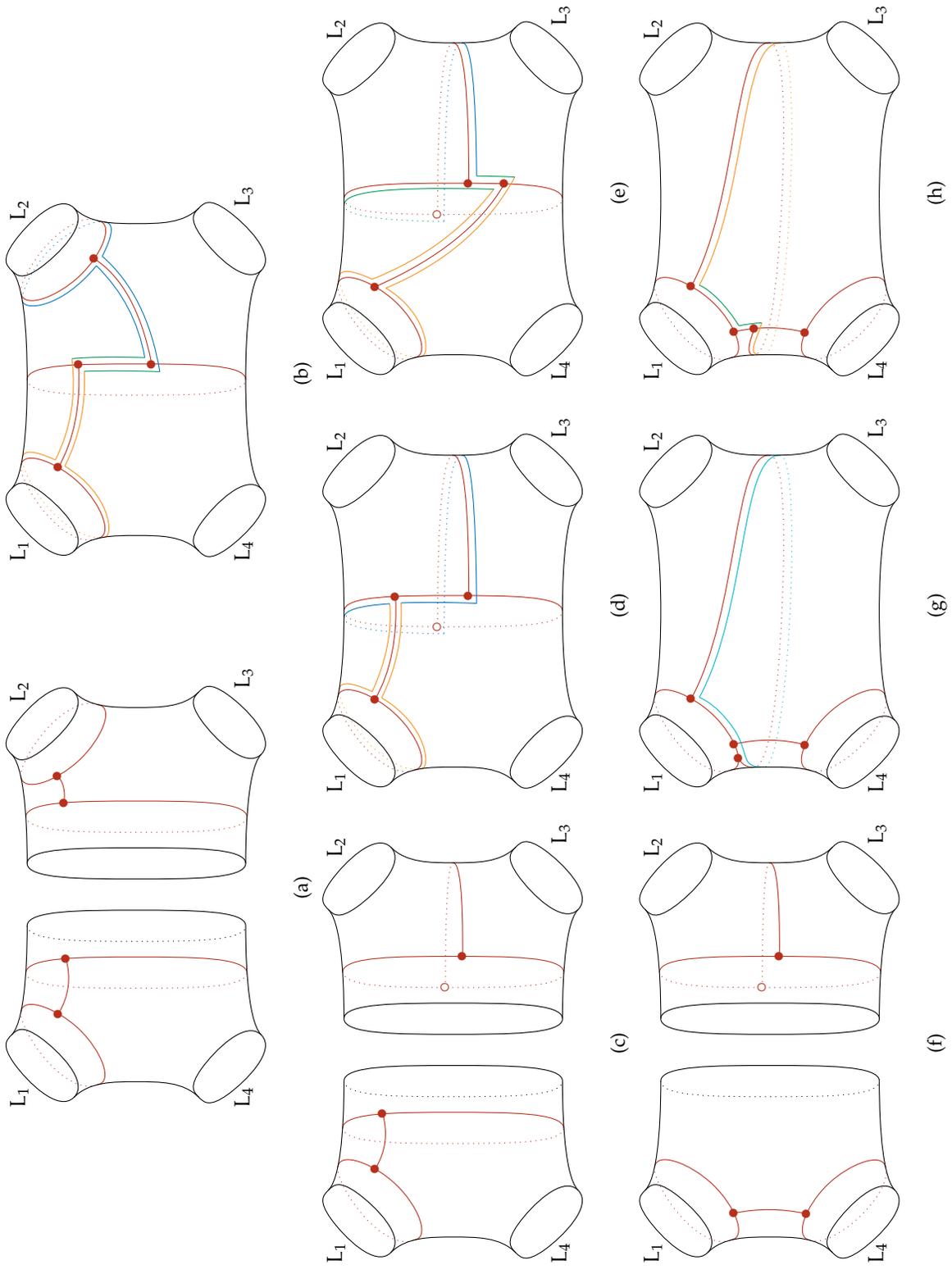

			\centering
			\begin{subfigure}[t]{.5\textwidth}
			\centering

			\subcaption{}
			\label{fig:ell:dash:04:a}
			\end{subfigure}
			\begin{subfigure}[t]{.5\textwidth}
			\centering
			%
			\subcaption{}
			\label{fig:ell:dash:04:b}
			\end{subfigure}
			\begin{subfigure}[t]{.45\textwidth}
			\centering
			%
			\subcaption{}
			\label{fig:ell:dash:04:c}
			\end{subfigure}
			\begin{subfigure}[t]{.4\textwidth}
			\centering
			%
			\subcaption{}
			\label{fig:ell:dash:04:d}
			\end{subfigure}
			\begin{subfigure}[t]{.4\textwidth}
			\centering
			%
			\subcaption{}
			\label{fig:ell:dash:04:e}
			\end{subfigure}
			\begin{subfigure}[t]{.45\textwidth}
			\centering
			%
			\subcaption{}
			\label{fig:ell:dash:04:f}
			\end{subfigure}
			\begin{subfigure}[t]{.4\textwidth}
			\centering
			%
			\subcaption{}
			\label{fig:ell:dash:04:g}
			\end{subfigure}
			\begin{subfigure}[t]{.4\textwidth}
			\centering
			%
			\subcaption{}
			\label{fig:ell:dash:04:h}
			\end{subfigure}
			\caption{The three cases examined in the proof of Lemma~\ref{lem:delta:eta:04}.}
			\label{fig:ell:dash:04}
		\end{figure}
	\end{landscape}
	\clearpage}
	
	Firstly, suppose that $L_4 > L_1 + \ell$ and $L_3 > L_2 + \ell $ (see Figure~\ref{fig:ell:dash:04:a}). Then
	\begin{equation*}
		\max\left\{
				0, L_1 - \ell, L_4 - \ell, \frac{L_1 + L_4 - \ell}{2}
			\right \}
		= L_4 - \ell,
		\qquad
		\max\left\{
				0, L_2 - \ell, L_3 - \ell, \frac{L_2 + L_3 - \ell}{2}
			\right\}
		= L_3 - \ell,
	\end{equation*}
	so the right-hand side of Equation~\eqref{eqn:length:delta:04} reduces to 
	\[
		\max\big\{
				L_1 + L_3 - \ell, L_2 + L_4 - \ell,2 |\tau| + L_3 + L_4 - 2\ell\big\}
		=
		2|\tau| + L_3 + L_4 - 2\ell.
	\]
	In Figure~\ref{fig:ell:dash:04:b}, a quasitransverse representative of $\delta$ is shown. Its orange part has length $L_4 - \ell$, its blue part has length $L_3 - \ell$, and its green part has length $2|\tau|$. In the end, we have $\ell'=  2|\tau| + L_3 + L_4 - 2\ell$, which is consistent with Equation~\eqref{eqn:length:delta:04}.
	 
	\smallskip
	
	Secondly, suppose that $L_4 > L_1 + \ell$ and  $|L_2 - L_3| < \ell < L_2 + L_3$ (see Figure~\ref{fig:ell:dash:04:c}). Then
	\begin{equation*}
		\max\left\{
				0,L_1 - \ell, L_4 - \ell, \frac{L_1 + L_4 - \ell}{2}
			\right\}
		= L_4 - \ell,
		\qquad
		\max\left\{
				0,L_2 - \ell, L_3 - \ell, \frac{L_2 + L_3 - \ell}{2}
			\right\}
		= \frac{L_2 + L_3 - \ell}{2}.
	\end{equation*}
	In this case, we also have $L_1 + L_3 < L_2 + L_4$, so the right-hand side of Equation~\eqref{eqn:length:delta:04} reduces to 
	\[
		\max\left\{
				L_1 + L_3 - \ell, L_2 + L_4 - \ell, 2|\tau| + L_4 - \ell + \frac{L_2 + L_3 - \ell}{2}
			\right\}
		=
		L_2 + L_4 - \ell + \left[ 2|\tau| - \frac{L_2 + \ell - L_3}{2} \right]_{+}.
	\]
	Suppose first that $2|\tau| < \frac{L_2 + \ell - L_3}{2}$, which is depicted in Figure~\ref{fig:ell:dash:04:d} together with a quasitransverse representative of $\delta$. The orange part of $\delta$ has length $L_4 - \ell$, while the blue part of $\delta$ has length $L_2$. Therefore
	\[
		\ell' = L_2 + L_4 - \ell
		=
		L_2 + L_4 - \ell + \left[ 2|\tau| - \frac{L_2 + \ell - L_3}{2} \right]_{+}.
	\]
	Suppose now that $2|\tau| > \frac{L_2 + \ell - L_3}{2}$, see Figure~\ref{fig:ell:dash:04:e}. The orange part of $\delta$ has length $L_4 - \ell$, the blue part of $\delta$ has length $\frac{L_2 + L_3 - \ell}{2}$, and the green part of $\delta$ has length $2|\tau|$. Thus:
	\[
		\ell' = 2|\tau| + L_4 - \ell + \frac{L_2 + L_3 - \ell}{2}
		=
		L_2 + L_4 - \ell + \left[ 2|\tau| - \frac{L_2 + \ell - L_3}{2} \right]_{+}.
	\]
	Again, in both cases Equation~\eqref{eqn:length:delta:04} is satisfied.
	
	\smallskip
	
	Thirdly, suppose that $\ell > L_1 + L_4$ and $|L_2 - L_3| < \ell < L_2 + L_3$ (see Figure~\ref{fig:ell:dash:04:f}). Then
	\begin{equation*}
		\max\left\{
				0,L_1 - \ell, L_4 - \ell, \frac{L_1 + L_4 - \ell}{2}
			\right\}
		= 0,
		\qquad
		\max\left\{
				0,L_2 - \ell, L_3 - \ell, \frac{L_2 + L_3 - \ell}{2}
			\right\}
		= \frac{L_2 + L_3 - \ell}{2}.
	\end{equation*}
	Without loss of generality, we can assume that $L_1 + L_3 > L_2 + L_4$. Then, the right-hand side of Equation~\eqref{eqn:length:delta:04} reduces to 
	\[
		\max\left\{ L_1 + L_3 - \ell, L_2 + L_4 - \ell, 2|\tau| + \frac{L_2 + L_3 - \ell}{2} \right\}
		=
		\max\left\{ L_1 + L_3 - \ell, 2|\tau| + \frac{L_2 + L_3 - \ell}{2} \right\}.
	\]
	The case where $2|\tau| < L_1 + \frac{L_3 - L_2 - \ell}{2}$ is depicted in Figure~\ref{fig:ell:dash:04:g}. The length of $\delta$ is then
	\[
		\ell' =	L_1 + L_3 - \ell
		=
		\max\left\{
			L_1 + L_3 - \ell, 2|\tau| + \frac{L_2 + L_3 - \ell}{2}
		\right\}.
	\]
	In the case where $2|\tau| > L_1 + \frac{L_3 - L_2 - \ell}{2}$, depicted in Figure~\ref{fig:ell:dash:04:h}, the orange part of $\delta$ has length $\frac{L_2 + L_3 - \ell}{2}$, and the green part of $\delta$ has length $2|\tau|$, therefore
	\[
		\ell' =	2|\tau| + \frac{L_2 + L_3 - \ell}{2}
		=
		\max\left\{
			L_1 + L_3 - \ell, 2|\tau| + \frac{L_2 + L_3 - \ell}{2}
		\right\}
	\]
	In both cases, Equation~\eqref{eqn:length:delta:04} is satisfied.

	\smallskip
	
	The case of $\GG$ with higher valencies can be obtained from the trivalent case by continuity of the combinatorial lengths and twists. Lastly, since $\eta$ is obtained from $\delta$ by performing a positive Dehn twist along $\gamma$, its length is given by $\ell'' = \ell'(\ell,\tau + \ell)$.
\end{proof}

The above lemma expresses the lengths $\ell'$ and $\ell''$ as functions of the Fenchel--Nielsen coordinates $(\ell,\tau)$. We can invert the perspective, expressing $\tau$ as a function of $\ell$, $\ell'$ and $\ell''$.

\begin{cor}\label{cor:twist:04}
	In the previous situation, we have
	\begin{equation}\label{eqn:tau:04}
		\tau = \begin{cases}
			\dfrac{1}{2}(\ell'' - M_{1,4} - M_{2,3}) - \ell
			& \text{if }
			\ell' = \max\Set{L_1 + L_3 - \ell,L_2 + L_4 - \ell }, \\[2ex]
			- \dfrac{1}{2}(\ell' - M_{1,4} - M_{2,3})
			& \text{if }
			\ell'' = \max\Set{ L_1 + L_3 - \ell,L_2 + L_4 - \ell }, \\[2ex]
			\dfrac{1}{2\ell}\left(\dfrac{\ell'' - M_{1,4} - M_{2,3}}{2}\right)^2 - \dfrac{1}{2\ell}\left(\dfrac{\ell' - M_{1,4} - M_{2,3}}{2}\right)^2 - \dfrac{\ell}{2}
			&
			\text{otherwise.}
		\end{cases}
	\end{equation} 
\end{cor}

\begin{proof}
	Let us denote $p = 2|\tau| + M_{1,4} + M_{2,3}$ and $q = \max\Set{L_1 + L_3 - \ell,L_2 + L_4 - \ell}$, so that $\ell' = \max\Set{p,q}$. We claim that $\ell' = q$ implies $2|\tau| \leq \ell$. If $L_2 + L_4 \geq L_1 + L_3$, this comes from the observation that $q = L_2 + L_4 - \ell$ and $\ell' = q + [\lambda]_{+}$ with
	\[
	\begin{split}
		\lambda
		& = p - L_2 - L_4 + \ell \\
		& = 2|\tau| - \ell +
			\max\left\{\ell - L_4,L_1 - L_4,0,\frac{L_1 - L_4 + \ell}{2}\right\}
			+
			\max\left\{\ell - L_2,0,L_3 -  L_2,\frac{L_3 - L_2 + \ell}{2}\right\}
		\geq 2|\tau| - \ell .
	\end{split}
	\]
	If $L_1 + L_3 \geq L_2 + L_4$, we rather have $q = L_1 + L_3 - \ell$ and the claim follows by writing $\ell' = q + [\mu]_{+}$ with
	\[
	\begin{split}
		\mu
		& = p - L_1 - L_3 + \ell \\
		& = 2|\tau| - \ell +
			\max\left\{\ell - L_1,0,L_4 - L_1,\frac{L_4 - L_1 + \ell}{2}\right\}
			+
			\max\left\{\ell - L_3,L_2 - L_3,0,\frac{L_2 - L_3 + \ell}{2}\right\}
		\geq 2|\tau| - \ell .
	\end{split}
	\]
	Therefore, if $\ell' = \max\set{L_1 + L_3 - \ell,L_2 + L_4 - \ell}$, we must have $|\tau| \leq \ell/2$, hence $|\tau + \ell| = \tau + \ell$. From Equation~\eqref{eqn:length:delta:04} we then find $\ell'' = 2|\tau + \ell| + M_{1,4} + M_{2,3}$, and solving for $\tau$ we get the first case of Equation~\eqref{eqn:tau:04}. The case $\ell'' = \max\Set{L_1 + L_4 - \ell,L_2 + L_3 - \ell}$ is similar. Finally, if none of those conditions are satisfied, then
	\[
		\ell' = 2|\tau| + M_{1,4} + M_{2,3},
		\qquad
		\ell'' = 2|\tau + \ell| + M_{1,4} + M_{2,3}.
	\]
This covers the last case in Equation~\eqref{eqn:tau:04}.
\end{proof}

\subsubsection{One-holed torus}
\label{subsubsec:one:holed:torus}

Let $\Sigma_i = T$ be a one-holed torus. We remove the index $i$ from the notation of $\gamma_i$, $\alpha_i$, $\delta_i$ and $\eta_i$, as well as $\ell_i$, $\ell_i'$, $\ell_i''$ and $\tau_i$, and denote $L = \ell_{\GG}(\de T)$.

\begin{lem}\label{lem:delta:eta:11}
	In the above setting, we have
	\begin{equation} \label{eqn:length:delta:11}
		\ell'(\ell,\tau)
		=
		|\tau| + \left[ \frac{L - 2\ell}{2} \right]_{+}.
	\end{equation}
	Further, $\ell''(\ell,\tau) = \ell'(\ell,\tau + \ell)$.
\end{lem}

\begin{proof}
	As before, we assume $\GG$ to be trivalent and we fix a marking. There are $4$ cases, corresponding to the $4$ open cells of the Teichm\"uller space of the pair of pants we obtain after cutting along $\gamma$.

	\begin{figure}[t]
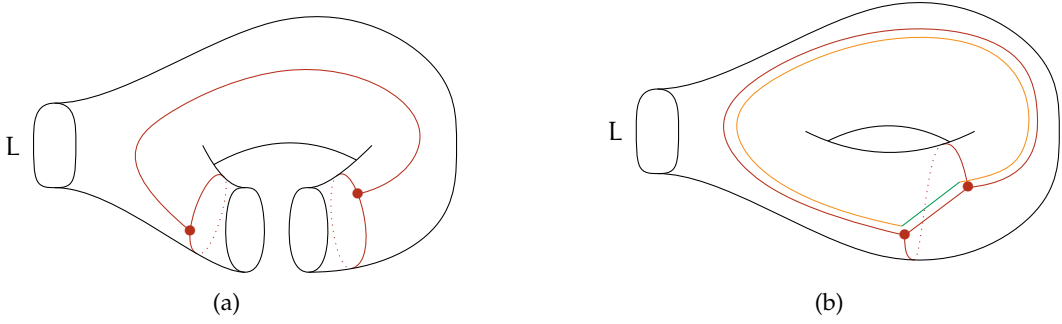

		\centering
		\begin{subfigure}[t]{.48\textwidth}
		\centering

		\subcaption{}
		\label{fig:ell:dash:11:a}
		\end{subfigure}
		\begin{subfigure}[t]{.48\textwidth}
		\centering
		%
		\subcaption{}
		\label{fig:ell:dash:11:b}
		\end{subfigure}
		\caption{The case examined in the proof of Lemma~\ref{lem:delta:eta:11}.}
		\label{fig:ell:dash:11}
	\end{figure}
	
	\smallskip

	We detail the case of Figure~\ref{fig:ell:dash:11:a}, where $L > 2\ell$. In Figure~\ref{fig:ell:dash:11:b}, a quasitransverse representative of $\delta$ is shown. Its orange part has length $L/2-\ell$, its green part has length $|\tau|$. Thus, we find $\ell' = |\tau| + L/2-\ell$, which is consistent with Equation~\eqref{eqn:length:delta:11} under the assumption $L > 2\ell$.

	\smallskip

	The other cases are analogous. Again, the formula extends to higher valency by continuity, and since the curve $\eta$ is obtained from the curve $\delta$ by performing a positive Dehn twist along $\gamma$, its length is given by $\ell'' = \ell'(\ell,\tau + \ell)$.
\end{proof}

Again, we can recover from the above Lemma an expression for the twist parameter $\tau$ as a function of the lengths $\ell$, $\ell'$ and $\ell''$. The proof is similar to the $(1,1)$ case.

\begin{cor}\label{cor:twist:11}
	In a one-holed torus $T$ with the above setting, the twist parameter is given as a function of $(\ell,\ell',\ell'')$ by
	\begin{equation}
		\tau
		=
		\frac{1}{2 \ell}
			\left( \ell'' - \left[\frac{L - 2\ell}{2}\right]_{+} \right)^2
			-
			\frac{1}{2\ell} \left( \ell' - \left[\frac{L - 2\ell}{2}\right]_{+} \right)^2
	 - \frac{\ell}{2}.
	\end{equation}
	\hfill $\blacksquare$
\end{cor}

\subsubsection{Proof of the combinatorial \texorpdfstring{$(9g - 9 + 3n)$-theorem}{nine-g-minus-9-plus-3n theorem}}

\begin{proof}[Proof of Theorem~\ref{thm:9g:minus:9:plus:3n}]
	The map is clearly continuous. Further, if $\GG, \GG' \in \mathcal{T}_{\Sigma}^{\textup{comb}}(L)$ are mapped to the same vector of lengths, then Corollaries~\ref{cor:twist:04} and \ref{cor:twist:11} would give the same length and twist parameters. As the combinatorial Fenchel--Nielsen map is a homeomorphism into the image, we deduce that $\GG = \GG'$. This justifies the injectivity.
\end{proof}

\newpage
\section{The symplectic structure}
\label{sec:symplectic:structure}

The Kontsevich symplectic form was originally introduced by Kontsevich in \cite{Kon92} as an ingredient in the proof of Witten's conjecture. Its main feature is the connection with intersection theory of $\overline{\mathfrak{M}}_{g,n}$. In this section, we prove a combinatorial analogue of Wolpert's formula (Theorem \ref{thm:Wolpert}) for the Kontsevich symplectic form, saying that Fenchel--Nielsen coordinates are Darboux. As a direct consequence, we achieve an explicit integration formula for a certain class of measurable functions on the combinatorial moduli space with respect to Kontsevich's volume form (Proposition~\ref{prop:integration}). This integration is the key operation that connects topological and geometric recursion, and is the combinatorial analogue of Mirzakhani's integration lemma on the moduli space of bordered Riemann surfaces \cite[Theorem~7.1]{Mir07simple}.

\subsection{Kontsevich's form}
\label{subsec:Kont:form}

Consider a ribbon graph $G$ of type $(g,n)$. For each index $i \in \{1,\ldots,n\}$, we make the choice of a first edge $e^{[i]}_1$ on the $i$-th face of $G$. We label the edges around the $i$-th face by $e^{[i]}_1, \ldots , e^{[i]}_{N_i}$ following the orientation of the face, which is opposite to the orientation of the boundary. Notice that every edge has a double label, as it bounds two faces or it appears twice in the cycle of a single face.

\medskip

Let now $\Sigma$ be a connected bordered surface of type $(g,n)$ and fix a combinatorial marking $(G,f)$. On each cell  of $\mathcal{T}_{\Sigma}^{\textup{comb}}$, we have functions $\ell^{[i]}_j \colon \mathfrak{Z}_{\Sigma,G}(L) \rightarrow \RR_{+}$, which associates to each $\GG \in \mathfrak{Z}_{\Sigma,G}(L)$ the length $\ell_{\GG}(e^{[i]}_j)$.

\begin{defn}\label{defn:Psi:class:comb}
	For each $i \in \{1,\ldots,n\}$, consider the differential $2$-form $\Psi_i$ on the cell complex $\mathcal{T}^{\textup{comb}}_{\Sigma}(L)$, defined on each cell $\mathfrak{Z}_{\Sigma,G}(L)$ by
	\begin{equation}\label{eqn:comb:psi}
		\Psi_i = \sum_{1 \le k < m \le N_i} \frac{\dd\ell^{[i]}_k}{L_i} \wedge \frac{\dd\ell^{[i]}_m}{L_i}.
	\end{equation}
	The form is $\Mod_{\Sigma}^{\de}$-invariant (it depends only on the ribbon graph underlying the marking) and we denote the induced form on the quotient $\mathcal{M}_{g,n}^{\textup{comb}}(L)$ with the same symbol\footnote{See \cite[Section~5.2]{Zvo02} for a discussion on the differential geometry of cell complexes. What we need here is that the combinatorial Teichm\"uller spaces and the combinatorial moduli spaces have a well-defined notion of polytopal differential forms, and that the associated polytopal de~Rham cohomology groups coincide with the usual cohomology groups over $\RR$. In particular, we can consider the cohomology class $[\Psi_i] \in H^{2}(\mathcal{M}_{g,n}^{\textup{comb}}(L))$.}.
\end{defn}

It can be shown that the definition of $\Psi_i$ does not depend on the choice of the first edge: the difference between two possible choices is of the form $L_i^{-2}\dd L_i \wedge \vartheta_i$ for some differential $1$-forms $\vartheta_i$, and thus is zero along the fibres $\mathcal{M}_{g,n}^{\textup{comb}}(L)$ of the perimeter map.

\medskip

In his original work \cite{Kon92}, Kontsevich related the above differential forms to the geometry of a certain circle bundle over $\mathcal{M}_{g,n}^{\textup{comb}}(L)$.

\begin{defn}\label{defn:comb:psi:classes}
	For each $i \in \{1,\ldots,n\}$, define $\mathcal{S}_i^{\textup{comb}}$ as the space of ordered pairs $(\bm{G},q)$ where $\bm{G} \in \mathcal{M}_{g,n}^{\textup{comb}}(L)$ and $q$ is a point belonging to an edge that borders the $i$-th face of $|\bm{G}|$. Its topology is the one induced by the natural cell structure. This defines a topological circle bundle $\mathcal{S}_i^{\textup{comb}} \to \mathcal{M}_{g,n}^{\textup{comb}}(L)$.
\end{defn}

\begin{thm}\cite{Kon92,Zvo02}\label{thm:comb:cotangent:classes}
	The class $[\Psi_i] \in H^2(\mathcal{M}_{g,n}^{\textup{comb}}(L))$ equals $-c_1(\mathcal{S}_i^{\textup{comb}})$. Further, under the identification $\mathfrak{M}_{g,n} \cong \mathcal{M}^{\textup{comb}}_{g,n}(L)$, the pullback of $\Psi_i$ extends continuously to the Deligne--Mumford compactification $\overline{\mathfrak{M}}_{g,n}$ and the associated cohomology class equals $\psi_i = c_1(\mathcal{L}_i) \in H^2(\overline{\mathfrak{M}}_{g,n})$, where $\mathcal{L}_i$ is the relative cotangent bundle at the $i$-th marked point.
	\hfill $\blacksquare$
\end{thm}

\begin{defn}\label{defn:Kform}
	Define the \emph{Kontsevich $2$-form} on $\mathcal{T}^{\textup{comb}}_{\Sigma}(L)$ as
	\begin{equation}\label{eqn:Kon:form}
		\omega_{\textup{K}} = \frac{1}{2} \sum_{i=1}^n L_i^2 \, \Psi_i.
	\end{equation}
\end{defn}

\begin{thm}\cite{Kon92}
	The differential form $\omega_{\textup{K}}$ is non-degenerate when restricted to strata corresponding to graphs with no vertices of even valency.
	\hfill $\blacksquare$
\end{thm}

\emph{A fortiori}, $\omega_{\textup{K}}$ descends to a symplectic form on the top-dimensional stratum of $\mathcal{M}^{\textup{comb}}_{g,n}(L)$, that is denoted with the same symbol.

\begin{defn}\label{defn:Kvolume}
	Define the \emph{Kontsevich measure}
	\begin{equation}\label{eqn:Kon:volume:form}
		\dd\mu_{\textup{K}} = \frac{\omega_{\textup{K}}^{d_{g,n}}}{d_{g,n}!},
		\qquad\quad
		d_{g,n} = 3g - 3 + n.
	\end{equation}
\end{defn}

Strictly speaking, $\dd\mu_{\textup{K}}$ is not a volume form on the whole $\mathcal{M}^{\textup{comb}}_{g,n}(L)$, although it is a volume form on the top-dimensional stratum. In any case, we have a notion of volume
\begin{equation}\label{eqn:Kon:volume}
	V^{\textup{K}}_{g,n}(L) = \int_{\mathcal{M}_{g,n}^{\textup{comb}}(L)} \dd\mu_{\textup{K}}
\end{equation}
where the integral is taken over the top-dimensional stratum. Equivalently, if $G$ is a trivalent ribbon graph of type $(g,n)$ and $P_{G}(L) \subseteq \RR_{+}^{E_{G}}$ is the polytope corresponding to those metrics on $G$ with fixed perimeter $L \in \RR_{+}^n$, then we set
\begin{equation}\label{eqn:Kon:volume:sum:ribbon:graphs}
	V^{\textup{K}}_{g,n}(L) = \sum_{\substack{G \in \mathcal{R}_{g,n} \\ \text{trivalent}}}
		\frac{1}{\#\Aut(G)} \int_{P_{G}(L)} \dd\mu_{\textup{K}}.
\end{equation}
We remark that, by abuse of notation, we are using the same symbol to denote the measure on $P_G(L)$ and its quotient $P_G(L)/\Aut(G)$.

\medskip

Notice that the volumes are finite, because the pullback of $\omega_{\textup{K}}$ extends continuously to $\overline{\mathfrak{M}}_{g,n}$. We moreover remark that, from Definition~\ref{defn:Kform} and Theorem~\ref{thm:comb:cotangent:classes}, the volumes are polynomials in the variables $L_i$, with coefficients given by $\psi$-classes intersections.

\begin{cor} \cite{Kon92} \label{cor:Kont:volume:psi:classes}
	The volume $V^{\textup{K}}_{g,n}(L)$ is a polynomial in $L_1^2,\ldots,L_n^2$ and satisfies
	\begin{equation}
		V^{\textup{K}}_{g,n}(L)
		=
		\sum_{\substack{d_1,\dots,d_n \geq 0 \\ d_1 + \cdots+d_n = d_{g,n}} } \bigg(\int_{\overline{\mathfrak{M}}_{g,n}} \prod_{i = 1}^n \psi_i^{d_i}\bigg) \prod_{i=1}^n \frac{L_i^{2d_i}}{2^{d_i}d_i!}.
	\end{equation}
	\hfill $\blacksquare$
\end{cor}
It is useful to record the expression of the Kontsevich measure in terms of edge lengths.

\begin{lem}\label{lem:Kon:Lebesgue}
	If $G$ is a trivalent ribbon graph of type $(g,n)$, let $(l_e)_{e \in E_{G}}$ be the edge lengths. We have the equality of measures in $\RR_+^{E_{G}}$
	\begin{equation}
		\dd\mu_{\textup{K}} \cdot \prod_{i = 1}^{n} \dd L_{i} = 2^{2g - 2 + n} \prod_{e \in E_{G}} \dd l_{e}.
	\end{equation}
\end{lem}

\begin{proof}
	From \cite[Appendix C]{Kon92} or \cite{CMS11}, we get
	\[
		\frac{1}{d_{g,n}!}
			\left(
				\sum_{i=1}^n L_i^2 \, \Psi_i
			\right)^{d_{g,n}}
		\prod_{i = 1}^{n} \dd L_{i}
		=
		2^{5g-5+2n} \prod_{e \in E_{G}} \dd l_{e}.
	\]
	Dividing on both sides by $2^{d_{g,n}}$ yields the result.
\end{proof}

\subsection{A combinatorial analogue of Wolpert's formula}

The purpose of this section is to show that combinatorial Fenchel--Nielsen coordinates are Darboux for $\omega_{\textup{K}}$.
\begin{thm}\label{thm:Wolpert}
	Let $\Sigma$ be a bordered surface of type $(g,n)$ and fix any combinatorial Fenchel--Nielsen coordinates $(\ell_{i},\tau_{i})$ for $\mathcal{T}_{\Sigma}^{\textup{comb}}(L)$. Denote by $\iota_L \colon \mathfrak{Z}_{\Sigma,G}(L) \hookrightarrow \mathcal{T}^{\textup{comb}}_{\Sigma}(L)$ the inclusion of a cell $\mathfrak{Z}_{\Sigma,G}(L)$. We have
	\begin{equation}\label{eqn:Wolpert}
		\omega_{\textup{K}}
		=
		\iota_{L}^{\ast}\left( \sum_{i=1}^{3g-3+n} \dd\ell_{i} \wedge \dd\tau_{i} \right).
	\end{equation}
\end{thm}

The advantage of this formula is that, while the left-hand side is clearly pure mapping class group invariant and it does not depend on the pants decomposition, the right-hand side has a simple expression in terms of the global coordinates and does not rely on the cells decomposition.

\subsubsection{Symplectic properties of the twist}

The main technical ingredient for the proof is Proposition~\ref{thm:Hamiltonianity}: the vector field associated to the twist along a simple closed curve is the Hamiltonian vector field of the length function of the curve. This is analogous to the situation in the hyperbolic case explored by Wolpert~\cite[Theorem~1.3]{Wol83}, and has the same flavour of a result in the space of measured foliations on closed surfaces proved by Papadopoulos \cite{Papado1,Papado2}. In the latter case, the symplectic structure is given by Thurston's from, and the length function coincide with the one considered here, \emph{i.e.} the intersection pairing between a foliation and a simple closed curve. Moreover, Proposition~\ref{thm:Hamiltonianity} generalises a result previously proved locally in~\cite[Lemma~3.2]{BCSW12} and only for very special curves cutting out small pairs of pants.

\medskip

To prove such a result, we need to understand how small changes in the twist parameter affect the metric on the embedded ribbon graph. More precisely, fix $\GG \in \mathcal{T}_{\Sigma}^{\textup{comb}}$ in an open cell, and $\gamma$ an essential, simple closed curve. Notice that, from Proposition~\ref{prop:open:set:of:tau}, we can cut $\GG$ along $\gamma$ and glue it back after twisting by a small amount $\tau \in \RR$. In particular, since we are in an open cell, it makes sense to talk about the vector field $\de_{\tau}$ generated by infinitesimal changes in the twist.

\medskip

\begin{figure}
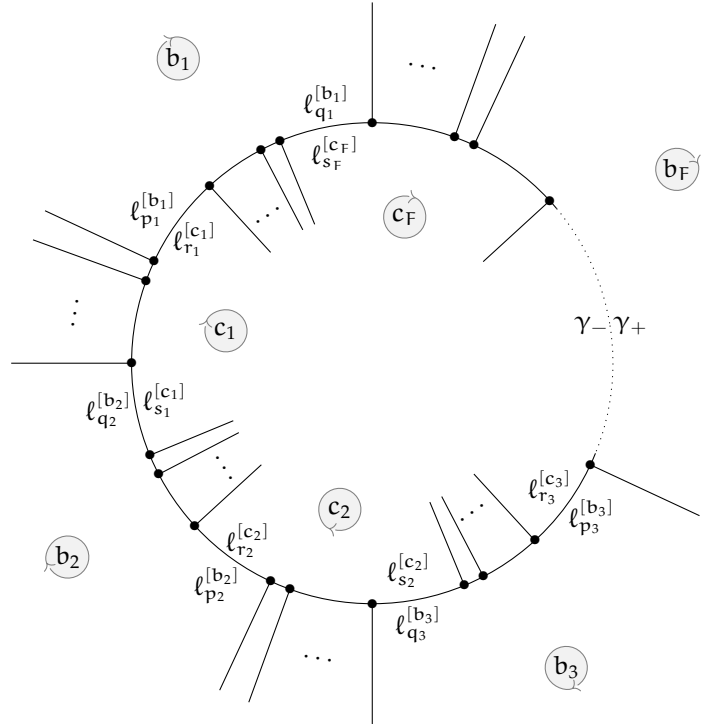

	\centering

\caption{A schematic picture of $\gamma$ used to calculate the vector field associated to the twist. The labels can be redundant if $\gamma$ visits an edge multiple times.}
\label{fig:local:description:vector:field}
\end{figure}

To get an expression for $\de_{\tau}$, we observe that each time $\gamma$ passes along an edge, a twist by small $\tau$ has the effect of either adding $\tau$ to the edge length, subtracting $\tau$ to the edge length, or leaving the edge length invariant (see the proof of Proposition~\ref{prop:open:set:of:tau} and Figure~\ref{fig:twist}). This depends on the direction taken by $\gamma$ at two consecutive vertices. Then one simply sums the changes in the length of the edges visited by $\gamma$. In the notation of Figure~\ref{fig:local:description:vector:field} (edges may appear twice along $\gamma$), the vector field describing the twisting along $\gamma$ is given by
\begin{equation}
	\de_{\tau}
	=
	\sum_{i=1}^{F}
		\left(
			\frac{\de}{\de \ell^{[b_i]}_{p_i}} - \frac{\de}{\de \ell^{[b_i]}_{q_i}}
		\right)
	=
	\sum_{i=1}^{F}
		\left(
			\frac{\de}{\de \ell^{[c_i]}_{r_i}} - \frac{\de}{\de \ell^{[c_i]}_{s_i}}
		\right),
\end{equation}
or in a more symmetric expression,
\begin{equation}\label{eqn:twist:vector:field}
	\de_{\tau}
	=
	\frac{1}{2} \sum_{i=1}^{F}
		\left(
			\frac{\de}{\de \ell^{[b_i]}_{p_i}} - \frac{\de}{\de \ell^{[b_i]}_{q_i}}
			+
			\frac{\de}{\de \ell^{[c_i]}_{r_i}} - \frac{\de}{\de \ell^{[c_i]}_{s_i}}
		\right).
\end{equation}

\medskip

To compute the contraction $\iota_{\de_{\tau}}\omega_{\textup{K}}$, we need the following technical lemma. Note that an edge $e_p^{[i]}$ is either adjacent to two different faces \emph{i.e.} $e_p^{[i]} = e_{r}^{[j]}$ for $i \neq j$, or adjacent to the same face on both sides \emph{i.e.} $e_p^{[i]} = e_{r}^{[i]}$ for $p \neq r$.

\begin{lem}\label{lem:contraction}
	In the interior of a top-dimensional cell, we have with the above notations
	\begin{equation}
	\label{eqn:contraction} \iota_{\de^{[i]}_p + \de^{[j]}_r} \omega_{\textup{K}}
				=
				\sum_{k = p+1}^{N_i} \dd\ell^{[i]}_{k} - \sum_{k = 1}^{p - 1} \dd\ell^{[i]}_{k}
				+
				\sum_{u = r+1}^{N_j} \dd\ell^{[j]}_{u} - \sum_{u = 1}^{r - 1} \dd\ell^{[j]}_{u}.
	\end{equation}
\end{lem}

\begin{proof}
	We recall that $\omega_{\textup{K}} = \frac{1}{2} \sum_{i = 1}^n L_i^2 \Psi_i$ and $\Psi_i$ defined in \eqref{defn:Psi:class:comb} only involves edges around the $i$-th face. Consider first the case $e_p^{[i]} = e_r^{[j]}$ for $i \neq j$. The interior product only receives contributions from $\Psi_i$ and $\Psi_j$. The interior product with $\Psi_i$ gives
	\[
		\frac{2}{L_i^2}
		\left(
			\sum_{k = p+1}^{N_i} \dd\ell^{[i]}_{k}
			-
			\sum_{k = 1}^{p - 1} \dd\ell^{[i]}_{k}
		\right),
	\]
	whereas the insertion into $\Psi_j$ gives
	\[
		\frac{2}{L_j^2}
		\left(
			\sum_{u = r+1}^{N_j} \dd\ell^{[i]}_{u}
			-
			\sum_{u = 1}^{r - 1} \dd\ell^{[j]}_{u}
		\right).
	\]
	Therefore we obtain Equation~\eqref{eqn:contraction}. In the second case, we have $e_p^{[i]} = e_r^{[i]}$ for $p \neq r$, and the interior product only receives a contribution from $\Psi_i$. Assuming $p < r$, the interior product with $\Psi_i$ is
	\begin{align*}
		& \frac{2}{L_i^2}
		\left(
			\sum_{k = p+1}^{r-1} \dd\ell^{[i]}_{k}
			+
			2 \sum_{k = r+1}^{N_i} \dd\ell^{[i]}_{k}
			-
			2 \sum_{k = 1}^{p-1} \dd\ell^{[i]}_{k}
			-
			\sum_{k = p+1}^{r-1} \dd\ell^{[i]}_{k}
		\right) \\
		& \qquad\qquad
		=
		\frac{2}{L_i^2}
		\left(
			\sum_{k = p+1}^{N_i} \dd\ell^{[i]}_{k} - \sum_{k = 1}^{p - 1} \dd\ell^{[i]}_{k}
			+
			\sum_{u = r+1}^{N_i} \dd\ell^{[i]}_{u} - \sum_{u = 1}^{r - 1} \dd\ell^{[i]}_{u}
		\right).
	\end{align*}
	Therefore, once again, we obtain Equation~\eqref{eqn:contraction}. The case $p > r$ is similar.
\end{proof}

We are ready now to state the main property of the twist vector field: it is the Hamiltonian vector field associated to the combinatorial length function $\ell \colon \GG \mapsto \ell_{\GG}(\gamma)$.

\begin{prop}\label{thm:Hamiltonianity}
	On the top-dimensional cells of $\mathcal{T}^{\textup{comb}}_{\Sigma}(L)$, we have
	\begin{equation}
		\dd\ell
		=
		\iota_{\de_{\tau}}\omega_{\textup{K}}.
	\end{equation}
\end{prop}
\begin{proof}
	Fix a top-dimensional cell and suppose that $\gamma$ is given by the schematic of Figure~\ref{fig:local:description:vector:field}. Then we have
	\[
		\dd\ell
		=
		\sum_{i=1}^{F} \Biggl(
			\frac{1}{2} \dd\ell_{p_i}^{[b_{i}]}
			+
			\frac{1}{2} \dd\ell_{q_{i}}^{[b_{i}]}
			+
			\sum_{p_i \prec k \prec q_{i}} \dd\ell_{k}^{[b_{i}]}
			+
			\frac{1}{2} \dd\ell_{r_i}^{[c_{i}]}
			+
			\frac{1}{2} \dd\ell_{s_{i}}^{[c_{i}]}
			+
			\sum_{r_i \prec u \prec s_{i}} \dd\ell_{u}^{[c_{i}]}
		\Biggr),
	\]
	where the symbol $\sum_{\mu \prec \lambda \prec \nu}$ indicates the sum over all edges of a certain face indexed by $\lambda$, that are between the edges indexed by $\mu$ and $\nu$, following the orientation of the face and excluding the extremes $\mu$ and $\nu$. Notice that the orientation of the face is opposite to the orientation of the corresponding boundary. On the other hand, we can reduce the computation of $\iota_{\de_{\tau}}\omega_{\textup{K}}$ to the insertion of the coordinate vector fields appearing in Equation~\eqref{eqn:twist:vector:field}.

	\smallskip 

	Let us explain why this calculation leads to a well defined answer. Notice first that Definition~\ref{defn:comb:psi:classes} of $\Psi_i$ makes perfect sense on the whole $\mathcal{T}_{\Sigma}^{\textup{comb}}$, where the perimeter is not fixed, but the definition depends on the choice of a first edge in the $i$-th face. However, its pullback to $\mathcal{T}_{\Sigma}^{\textup{comb}}(L)$ is independent of such a choice. Secondly,  observe that the vector fields $\de_e$ are defined on $\mathcal{T}^{\textup{comb}}_{\Sigma}$ but do not have a meaning on $\mathcal{T}^{\textup{comb}}_{\Sigma}(L)$, as they do not preserve the boundary lengths $L$. However, particular linear combinations of them, such as $\de_{\tau}$, do. Therefore, it is legitimate to compute each contribution $\iota_{\de_e} \Psi$ separately on $\mathcal{T}^{\textup{comb}}_{\Sigma}$ (\emph{i.e.} we can safely use Lemma~\ref{lem:contraction}), then sum them up to obtain $\iota_{\de_{\tau}} \omega_{\textup{K}}$, and eventually take the pullback to $\mathcal{T}_{\Sigma}^{\textup{comb}}(L)$.

	\smallskip 

	This said, a repeated use of Lemma~\ref{lem:contraction} results in the following computation.
	\[
	\begin{split}
		\iota_{\de_{\tau}} \omega_{\textup{K}}
		&=
		\frac{1}{2} \sum_{i=1}^{F}
		\Biggl[
			\Bigl(
				\iota_{\de^{[b_i]}_{p_i} + \de^{[c_i]}_{r_i}}
			\Bigr) \omega_{\textup{K}}
			-
			\Bigl(
				\iota_{\de^{[b_i]}_{q_i} + \de^{[c_i]}_{s_i}}
			\Bigr) \omega_{\textup{K}}
		\Biggr] \\
		&=
		\frac{1}{2} \sum_{i=1}^{F}
		\Biggl[
			\Biggl(
				\sum_{k = p_i+1}^{N_{b_i}} \dd\ell^{[b_i]}_{k} - \sum_{k = 1}^{p_i - 1} \dd\ell^{[b_i]}_{k}
				+
				\sum_{u = r_i+1}^{N_{c_i}} \dd\ell^{[c_i]}_{u} - \sum_{u = 1}^{r_i - 1} \dd\ell^{[c_i]}_{u}
			\Biggr) \\
		& \qquad\qquad
			-
			\Biggl(
				\sum_{m = q_i+1}^{N_{b_i}} \dd\ell^{[b_i]}_{m} - \sum_{m = 1}^{q_i - 1} \dd\ell^{[b_i]}_{m}
				+
				\sum_{v = s_i+1}^{N_{c_i}} \dd\ell^{[c_i]}_{v} - \sum_{v = 1}^{s_i - 1} \dd\ell^{[c_i]}_{v}
			\Biggr)
		\Biggr] \\
		& =
		\sum_{i=1}^{F} \Biggl(
			\frac{1}{2} \dd\ell_{p_i}^{[b_{i}]}
			+
			\frac{1}{2} \dd\ell_{q_{i}}^{[b_{i}]}
			+
			\sum_{p_i \prec k \prec q_{i}} \dd\ell_{k}^{[b_{i}]}
			+
			\frac{1}{2} \dd\ell_{r_i}^{[c_{i}]}
			+
			\frac{1}{2} \dd\ell_{s_{i}}^{[c_{i}]}
			+
			\sum_{r_i \prec u \prec s_{i}} \dd\ell_{u}^{[c_{i}]}
		\Biggr).
	\end{split}
	\]
This indeed coincides with $\dd \ell$.
\end{proof}

\subsubsection{Proof of the combinatorial Wolpert's formula}

For a fixed oriented surface $\Sigma$, denote by $\bar{\Sigma}$ the surface with opposite orientation. The following lemma is based on the work of Wolpert~\cite[Lemma~1.1]{Wol85}.

\begin{lem}\label{Wol:Lem}
	Let $\Sigma$ be a bordered surface, $\rho \colon \Sigma\rightarrow\bar{\Sigma}$ be an isotopy class of orientation-reversing diffeomorphism that restricts to the identity on the boundary. Fix $\gamma\in S_{\Sigma}^{\circ}$. Then $\rho$ induces a homeomorphism $\mathcal{T}_{\Sigma}^{\textup{comb}}\rightarrow\mathcal{T}_{\bar{\Sigma}}^{\textup{comb}}$ and
	\begin{itemize}
		\item $\rho^{*} \dd\ell(\gamma) = \dd\ell(\rho \circ \gamma)$,
		\item $\rho^{*}\omega_{\textup{K}} = -\omega_{\textup{K}}$,
		\item if $\rho$ fixes $\gamma$, then $\rho^{*} \dd\tau(\gamma) = -\dd\tau(\gamma) + \frac{m}{2} \dd\ell(\gamma)$ for some $m \in \ZZ$.
	\end{itemize}
	Here $d\tau$ is the differential form dual with respect to $\omega_{\textup{K}}$ to the vector field $\de_{\tau}$ of Equation\eqref{eqn:twist:vector:field}.
\end{lem}
\begin{proof}
	The aforementioned induced map $\mathcal{T}_{\Sigma}^{\textup{comb}} \rightarrow \mathcal{T}_{\bar{\Sigma}}^{\textup{comb}}$ is the composition of $\rho$ with the marking. It inverts the orientations of all curves, but it fixes the length functions, hence the first point. Further, the $\Psi$-classes are going to be calculated using the opposite orientation, which yields the sign for the second point. The last point follows from the fact that the elements of $\Stab(\gamma)$ are generated by (half-) Dehn twists along curves that do not intersect $\gamma$. Then, as $\rho$ reverses the orientation of the surface, $\dd\tau(\gamma)$ acquires a sign from the orientation reversal and an ambiguity of $\frac{1}{2}\ZZ d\ell(\gamma)$ from a potential (half-) Dehn twists along $\gamma$.
\end{proof}

We are now ready to give a proof of the combinatorial Wolpert's formula.

\begin{proof}[Proof of Theorem \ref{thm:Wolpert}]
	Fix a seamed pants decomposition. We know that $(\ell_{i},\tau_{i})$ give global coordinates on $\mathcal{T}_{\Sigma}^{\textup{comb}}(L)$. Therefore, on the top-dimensional cells
	\[
		\omega_{\textup{K}}
		=
		\sum_{i<j} a_{ij} \, \dd\ell_{i}\wedge \dd\ell_{j}
		+
		\sum_{i<j} b_{ij} \, \dd\tau_{i}\wedge \dd\tau_{j}
		+
		\sum_{i,j} c_{ij} \, \dd\ell_{i}\wedge \dd\tau_{j}
	\]
	for some functions $a_{ij}$, $b_{ij}$ and $c_{ij}$. Notice that from Proposition~\ref{thm:Hamiltonianity} we have
	\[
		\iota_{\de_{\tau_{i}}} \omega_{\textup{K}}
		=
		\sum_{i<j} b_{ij} \, \dd\tau_{j} - \sum_{j<i} b_{ji} \, \dd\tau_{j} + \sum_{j}c_{ji} \, \dd\ell_{j}
		=
		\dd\ell_{i},
	\]
	and hence $b_{ij} = 0$, $c_{ij} = \delta_{ij}$. Finally, if $\rho$ is the isotopy class of an orientation-reversing diffeomorphism fixing $\gamma_{i}$, we have
	\[
		\rho_{*}\Big(\frac{\de}{\de\ell_{i}}\Big)
		=
		\frac{\de}{\de\ell_{i}} + \frac{m_{i}}{2} \frac{\de}{\de\tau_{i}}.
	\]
	Therefore from Lemma~\ref{Wol:Lem}, for $i<j$, we find	
	\begin{equation*}
	\begin{split}
		a_{ij} = & \omega_{\textup{K}} \left(\frac{\de}{\de\ell_{i}},\frac{\de}{\de\ell_{j}}\right) =	\omega_{\textup{K}} \left(
			\frac{\de}{\de\ell_{i}} + \frac{m_{i}}{2} \frac{\de}{\de\tau_{i}},
			\frac{\de}{\de\ell_{j}} + \frac{m_{j}}{2} \frac{\de}{\de\tau_{j}}
		\right) =
		\omega_{\textup{K}} \left(
			\rho_{*}\frac{\de}{\de\ell_{i}},
			\rho_{*}\frac{\de}{\de\ell_{j}}
		\right) =
		\rho^{*} \omega_{\textup{K}} \left(
			\frac{\de}{\de\ell_{i}},
			\frac{\de}{\de\ell_{j}}
		\right) \\
		 = &
		- \omega_{\textup{K}} \left(
			\frac{\de}{\de\ell_{i}},
			\frac{\de}{\de\ell_{j}}
		\right).
	\end{split}
	\end{equation*}
	and thus $a_{ij} = 0$. This proves the result on the top-dimensional cells.

	\smallskip

	To extend it to cells $\iota_{L} \colon \mathfrak{Z}_{\Sigma,G}(L) \hookrightarrow \mathcal{T}_{\Sigma}^{\textup{comb}}(L)$ of positive codimension, we can consider it at the boundary of a top-dimensional cell. Then $\iota^{*}_{L}$ simply sets $\dd\ell_{e} = 0$ for each edge $e$ of zero length on the boundary of the top-dimensional cell. This has exactly the same affect as excluding such edges from the sum in Equation~\eqref{defn:Psi:class:comb} of $\Psi_{i}$, which coincides with the definition of $\Psi_{i}$ on cells of positive codimension.
\end{proof}

\subsection{Integration over the combinatorial moduli space}
\label{subsec:integration}

In this section we establish an integration result, analogous to \cite[Theorem~7.1]{Mir07simple} for $(\mathcal{M}_{g,n}(L),\omega_{\textup{WP}})$, exploiting the symplectic structure of $(\mathcal{M}_{g,n}^{\textup{comb}}(L),\omega_{\textup{K}})$ for which the key is our combinatorial Wolpert's formula~\eqref{eqn:Wolpert}. This improves the results of \cite[Theorem~1.1]{BCSW12}, which in fact can be extended from their original use to the integration of functions with support restricted\footnote{This restriction on the support is related to the one appearing, \textit{e.g.}, in \eqref{eqn:vanishing:small:pop}.} to ``small pairs of pants''.

\medskip

Let us introduce some notation. Consider a bordered surface $\Sigma$ of type $(g,n)$ and let $\gamma$ be a primitive multicurve with ordered components $(\gamma_j)_{j = 1}^k$. We denote by $\Gamma$ the orbit $\Mod_{\Sigma}^{\de}.\gamma$ (although it is not important for what follows, such an object can be seen as a stable graph with ordered edges, see Section~\ref{subsec:twisted:GR}). Furthermore, consider an assignment
\begin{equation}
	\Sigma' \longmapsto \Xi_{\Sigma'} \in {\rm Mes}(\mathcal{T}_{\Sigma'}^{\textup{comb}},\mu_{\textup{K}})
\end{equation}
of a measurable function on the combinatorial Teichm\"uller space to each bordered surface $\Sigma'$ diffeomorphic to the cut surface $\Sigma_{\gamma}$. We assume that for any diffeomorphism  $\phi \colon \Sigma' \rightarrow \Sigma''$ which preserves the labelling of the boundary components, we have $\phi_*\Xi_{\Sigma'} = \Xi_{\Sigma''}$ where $\phi_*$ is the map induced between the combinatorial Teichm\"uller spaces. In particular, $\Xi_{\Sigma'}$ is invariant under the action of $\Mod_{\Sigma'}^{\de}$ and descends to a function $\Xi_{\Gamma}$ on the moduli space
\begin{equation}
	\mathcal{M}_{\Gamma}^{\textup{comb}} = \prod_{v \in \pi_0(\Sigma_{\gamma})} \mathcal{M}_{g(v),n(v)}^{\textup{comb}},
\end{equation}
which depends on $\Gamma$ only. We can decompose $\pi_0(\de\Sigma') = \pi_0(\gamma) \sqcup \pi_0(\gamma) \sqcup \pi_0(\de \Sigma)$, so it makes sense to consider\footnote{As explained in \S~\ref{subsubsec:notation:curves}, the boundary components of the cut surface are labeled in a specific way. Thus, with the symbol $\mathcal{M}_{\Gamma}^{\textup{comb}}(\ell,\ell,L)$, we mean the product of moduli spaces with fixed boundary lengths from $\ell_1,\dots,\ell_k,L_1,\dots,L_n$ ordered in such a way that they match the labeled boundary components of the cut surface.} $\mathcal{M}_{\Gamma}^{\textup{comb}}(\ell,\ell,L)$. We further assume that, for almost every $(L,\ell) \in \RR_{+}^{n} \times \RR_{+}^{k}$, $\Xi_{\Gamma}$ is integrable  with respect to the Kontsevich measure on $\mathcal{M}_{\Gamma}^{\textup{comb}}(\ell,\ell,L)$, and we denote
\begin{equation}
	V\Xi_{\Gamma}(\ell,\ell,L)
	=
	\int_{\mathcal{M}^{\textup{comb}}_{\Gamma}(\ell,\ell,L)} \Xi_{\Gamma} \, \dd\mu_{\textup{K}}.
\end{equation}
Finally, consider a measurable function $f \colon \RR_{+}^{n} \times \RR_{+}^{k} \rightarrow \RR$ and define a new function $\Xi_{\Sigma}^{f,\Gamma}$ on $\mathcal{T}_{\Sigma}^{\textup{comb}}$ by setting
\begin{equation}\label{eqn:Xi:f:gamma}
	\Xi_{\Sigma}^{f,\Gamma}(\GG)
	=
	\sum_{\alpha \in \Gamma}
		f\big(\vec{\ell}_{\GG}(\de \Sigma), \vec{\ell}_{\GG}(\alpha)\big) \, \Xi_{\Sigma_{\alpha}}(\GG|_{\Sigma_{\alpha}}),
\end{equation}
where $\vec{\ell}_{\GG}(\de \Sigma) = (\ell_{\GG}(\de_i \Sigma))_{i = 1}^n$ and  $\vec{\ell}_{\GG}(\alpha) = (\ell_{\GG}(\alpha_j))_{j = 1}^k$. When the series \eqref{eqn:Xi:f:gamma} is absolutely convergent, it defines a $\Mod_{\Sigma}^{\de}$-invariant function, which descends to a function $\Xi_{g,n}^{f,\Gamma}$ on the moduli space $\mathcal{M}_{g,n}^{\textup{comb}}$.
  
\begin{prop}\label{prop:integration}
	Assume that the series \eqref{eqn:Xi:f:gamma} is absolutely convergent, and that for almost every $L \in \RR_{+}^n$ its limit is integrable with respect to $\mu_{\textup{K}}$ on $\mathcal{M}_{\Sigma}^{\textup{comb}}(L)$. Assume as well that for almost every $(L,\ell) \in \RR_+^{n} \times \RR_+^{k}$ the function $\Xi_{\Gamma}$ is integrable on $\mathcal{M}_{\Gamma}^{\textup{comb}}(\ell,\ell,L)$ with respect to the Kontsevich measure. Then
	\begin{equation}\label{eqn:integration}
		\int_{\mathcal{M}_{g,n}^{\textup{comb}}(L)}
			\Xi_{g,n}^{f,\Gamma} \, \dd\mu_{\textup{K}}
		=
		\int_{\RR_{+}^{k}} f(L,\ell) \, V\Xi_{\Gamma}(\ell,\ell,L)
			\prod_{j = 1}^k \ell_{j} \, \dd\ell_{j}.
 	\end{equation}
\end{prop} 

\begin{proof}
	We adapt Mirzakhani's proof of \cite[Theorem~7.1]{Mir07simple}, which concerned the hyperbolic setting with $\Xi_{\Sigma'} = 1$, functions $f(L,\ell) = F(\ell_1 + \cdots + \ell_k)$ and $\gamma$ a primitive multicurve with unordered components. Because the components are ordered, our formula does not contain  automorphism factors. The main difference is instead that, in the combinatorial setting, we have to remove the zero measure set of ill-defined twists.

	\smallskip

	Consider the space
	\[
		\mathcal{M}_{g,n}^{\textup{comb},\Gamma}(L) = \mathcal{T}_{\Sigma}^{\textup{comb}}(L)\Big/\bigcap_{j = 1}^k \Stab(\gamma_j),
	\]
	where $\Stab(\gamma_j)$ is the stabiliser of $\gamma_j$ in $\Mod_{\Sigma}^{\de}$. We denote $\Pi^{\Gamma} \colon \mathcal{M}_{g,n}^{\textup{comb},\Gamma}(L) \to \mathcal{M}_{g,n}^{\textup{comb}}(L)$ the natural projection. Notice that
	\[
		\mathcal{M}_{g,n}^{\textup{comb},\Gamma}(L)
		\cong
		\Set{ (\bm{G},\alpha) | \bm{G} \in \mathcal{M}_{g,n}^{\textup{comb}}(L), \quad \alpha \in \Gamma}.
	\]
	Since the symplectic structure on $\mathcal{T}_{\Sigma}^{\textup{comb}}$ is invariant under the action of the pure mapping class group, it induces a symplectic structure on $\mathcal{M}_{g,n}^{\textup{comb},\Gamma}(L)$, which is the same as the pullback $\left(\Pi^{\Gamma}\right)^\ast\omega_{\textup{K}}$. We denote the associated measure by $\mu_{\textup{K}}^{\Gamma}$.

	\smallskip

	Consider now the map $\mathcal{T}^{\textup{comb}}_{\Sigma}(L) \to \RR_{+}^{k}$ given by the tuple of combinatorial lengths of the components of $\gamma$. It descends to a map $\mathcal{L}^{\Gamma} \colon \mathcal{M}^{\textup{comb},\Gamma}_{g,n}(L) \to \RR_{+}^{k}$. We denote by $\mathcal{M}^{\textup{comb},\Gamma}_{g,n}(L)[\ell] = (\mathcal{L}^{\Gamma})^{-1}(\ell)$ the level sets for $\ell \in \mathbb{R}_{+}^{k}$. We have a map
	\begin{equation}\label{eqn:proj:gamma}
		\Pi \colon \mathcal{M}^{\textup{comb},\Gamma}_{g,n}(L)[\ell] \longrightarrow \mathcal{M}^{\textup{comb}}_{\Gamma}(\ell,\ell,L)
	\end{equation}
	defined in the natural way: given an element $(\bm{G},\alpha) \in \mathcal{M}^{\textup{comb},\Gamma}_{g,n}(L)[\ell]$, we take a lift $\GG \in \mathcal{T}^{\textup{comb}}_{\Sigma}(L)$ of $\bm{G}$, we restrict $\GG$ to the cut surface $\Sigma_{\alpha}$ as explained in Section~\ref{subsec:cut:glue} and we project the restriction to the moduli space $\mathcal{M}^{\textup{comb}}_{\Sigma_{\alpha}}(\ell,\ell,L) \cong \mathcal{M}^{\textup{comb}}_{\Gamma}(\ell,\ell,L)$. The result does not depend on the choice of the lift $\GG$ since we are projecting to the combinatorial moduli space of $\Sigma_{\alpha}$ after restriction.

	\smallskip

	Notice that the spaces on both sides of \eqref{eqn:proj:gamma} have a natural measure: $\mathcal{M}^{\textup{comb},\Gamma}_{g,n}(L)[\ell]$ is equipped with the disintegration of $\mu_{\textup{K}}^{\Gamma}$ along $\mathcal{L}^{\Gamma}$, and $\mathcal{M}^{\textup{comb}}_{\Gamma}(\ell,\ell,L)$ with its Kontsevich measure. By construction of \eqref{eqn:Xi:f:gamma} and the property of disintegration,
	\begin{equation}\label{eqn:disintegration:Xi}
	\begin{split}
		\int_{\mathcal{M}_{g,n}^{\textup{comb}}(L)} \Xi^{f,\Gamma}_{g,n}\dd\mu_{\textup{K}}
		& =
			\int_{\mathcal{M}_{g,n}^{\textup{comb},\Gamma}(L)} (f \circ \mathcal{L}^{\Gamma}) \cdot (\Xi_{\Gamma} \circ \Pi)\,\dd\mu_{\textup{K}}^{\Gamma} \\
		& =
			\int_{\RR_{+}^{k}} f(L,\ell) \left(
				\int_{\mathcal{M}^{\textup{comb},\Gamma}_{g,n}(L)[\ell]} (\Xi_{\Gamma} \circ \Pi)\, \dd\mu_{\textup{K}}^{\Gamma}
			\right) \prod_{j = 1}^k \dd \ell_j.
	\end{split}
	\end{equation}
	We can complete $\gamma$ into a seamed pants decomposition and use the Fenchel--Nielsen coordinates to describe the space $\mathcal{M}^{\textup{comb},\Gamma}_{g,n}(L)[\ell]$. The combinatorial Wolpert's formula~\eqref{eqn:Wolpert} implies that the measure $\mu_{\textup{K}}^{\Gamma}$ has a product structure with respect to the fibration $\Pi$. Besides, the fibre of $\bm{G}' \in \mathcal{M}^{\textup{comb}}_{\Gamma}(\ell,\ell,L)$ is identified with an open subset of full measure in $\prod_{j = 1}^k \RR/ 2^{-\mathfrak{t}_{j}}\ell_{j}\ZZ$, where
	\[
		\mathfrak{t}_{j} =
		\begin{cases}
			1 & \text{if $\gamma_j$ separates off a torus with one boundary,} \\
			0 & \text{otherwise.}
		\end{cases}
	\]
	This follows from the description of the image of the Fenchel--Nielsen coordinates in Theorem~\ref{thm:FN:coordinates}. The factor of $2^{-\mathfrak{t}_{j}}$ in the case when $\gamma_j$ separates off a torus with one boundary is due to the fact that any element in $\mathcal{M}^{\textup{comb}}_{1,1}(L)$ comes with an elliptic involution, so $\Stab(\gamma_i)$ contains the half-twist along $\gamma_j$ and the fundamental region of the twist coordinate in the combinatorial moduli space becomes $[0,\ell_{j}/2]$ minus a measure-zero set. So, for any open set $U \subseteq \mathcal{M}^{\textup{comb}}_{\Gamma}(\ell,\ell,L)$, we have
	\[
		\prod_{j = 1}^k 2^{-\mathfrak{t}_{j}} \int_{\Pi^{-1}(U)}  (\Xi_{\Gamma} \circ \Pi) \, \dd\mu_{\textup{K}}^{\Gamma}
		=
		\prod_{j = 1}^k 2^{-\mathfrak{t}_{j}} \ell_j \; \int_{U}\Xi_{\Gamma} \, \dd\mu_{\textup{K}}
	\]
	whenever the functions we wish to integrate are integrable. We noticed that the integrals over $\mathcal{M}_{1,1}^{\textup{comb}}$ have an extra factor of $\tfrac{1}{2}$ due to the presence of the elliptic involution, while such a factor is not present on $\mathcal{M}^{\textup{comb},\Gamma}_{g,n}(L)[\ell]$. We must therefore include an extra factor $2^{-\mathfrak{t}_{j}}$ in the right-hand side, and this cancels the factor of $2^{-\mathfrak{t}_{j}}$ coming from the half-twist. By a partition of unity argument, we obtain
	\[
		\int_{\mathcal{M}^{\textup{comb},\Gamma}_{g,n}(L)[\ell]} (\Xi_{\Gamma} \circ \Pi)\,\dd\mu_{\textup{K}}^{\Gamma}
			=
			V\Xi_{\Gamma}(\ell,\ell,L) \prod_{j = 1}^k \ell_{j},
	\]
	which we insert in \eqref{eqn:disintegration:Xi} to complete the proof. 
\end{proof}

\begin{rem}
	In Mirzakhani's work \cite{Mir07simple} there is an unnatural convention for the integral over $\mathcal{M}_{1,1}(L)$ which does not include the factor of $\tfrac{1}{2}$ coming from the elliptic involution. For this reason, she finds an extra factor $\prod_{j=1}^k 2^{-\mathfrak{t}_j}$.
\end{rem}

\begin{rem}\label{rem:integration:lemma:automorphism}
	As briefly mentioned before, an equivalent way to state Proposition~\ref{prop:integration}, which used in Section~\ref{subsec:unit:balls:comb}, relies on the notion of stable graphs (introduced in Section~\ref{subsec:twisted:GR}). In this language, the set of connected components of $\Sigma_{\gamma}$ is the set $V_{\Gamma}$ of vertices of the stable graph with ordered edges $\Gamma = \Mod_{\Sigma}^{\de}.\gamma$, the set of components of $\gamma$ descends to the set $E_{\Gamma}$ of edges of $\Gamma$, and the set of boundary component of $\Sigma$ is the set $\Lambda_{\Gamma}$ of leaves of $\Gamma$. With this notation, Equation~\eqref{eqn:integration} becomes
	\begin{equation}\label{eqn:integration:stable:graph}
		\int_{\mathcal{M}_{g,n}^{\textup{comb}}(L)} \Xi_{g,n}^{f,\Gamma} \dd\mu_{\textup{K}}
		=
		\int_{\mathbb{R}_{+}^{E_{\Gamma}}}
			f\big((L_{\lambda})_{\lambda\in\Lambda_{\Gamma}},(\ell_e)_{e\in E_{\Gamma}}\big) \,
			\prod_{v\in V_{\Gamma}} V\Xi_{g(v),n(v)}\big((\ell_e)_{e\in E_{v}},(L_{\lambda})_{\lambda\in\Lambda_v}\big)
			\prod_{e \in E_{\Gamma}} \ell_e \, \dd \ell_e.
	\end{equation}
	If in the above arguments we do not suppose that the components of the primitive multicurve are ordered, then we have to include the automorphism factor $\frac{1}{\#\Aut(\Gamma)}$ multiplying the right-hand side of Equation~\eqref{eqn:integration:stable:graph}. Moreover, if the multicurve is not primitive but comes with a weight $a \in \ZZ_+^{E_{\Gamma}}$, we have to consider the automorphism factor $\frac{1}{\#\Aut(\Gamma,a)}$ of the weighted multicurve.
\end{rem}

\begin{rem}\label{rem:integration:lemma:nonnegative}
	We can get rid of the integrability assumption for $\Xi_{g,n}^{f,\Gamma}$, if we suppose that each term in the series \eqref{eqn:Xi:f:gamma} is non-negative. In this case, the disintegration property still holds in \eqref{eqn:disintegration:Xi}, and the result of Proposition~\ref{prop:integration} becomes an equality between (possibly infinite) integrals.
\end{rem}

\newpage
\section{Functions from geometric recursion}
\label{sec:setting:up:gr}

In this section, we review the construction of mapping class group invariant functions on the combinatorial and hyperbolic Teich\-m\"uller spaces via the general framework of geometric recursion (GR) proposed in \cite{ABO17}. We moreover show that the integral on the moduli space of these functions against the Kontsevich measure satisfy topological recursion (TR). We shall establish the combinatorial version of the Mirzakhani--McShane identity: the constant function $1$ on the combinatorial Teichm\"uller space is an output of the geometric recursion for some particular input data. As a corollary, we can present purely geometric proofs of Witten conjecture/Kontsevich theorem \cite{Wit91,Kon92} and of Norbury's lattice point count by topological recursion \cite{Nor10}, which are completely parallel to Mirzakhani's proof of the topological recursion for the Weil--Petersson volumes in \cite{Mir07simple}.
\medskip

In Section~\ref{sec:hyp:comb} we study a particular flow related to the boundary lengths in the classical sense, which interpolates between the hyperbolic and the combinatorial geometries, allowing us to transfer several of the results discussed from $\mathcal{T}_{\Sigma}$ to $\mathcal{T}_{\Sigma}^{\textup{comb}}$ directly.

\subsection{Definition of the geometric recursion}

We begin by reviewing in a simplified form the framework of geometric recursion introduced in \cite{ABO17}.

\medskip

Let us introduce the category $\mathcal{B}_1$ whose objects are bordered surfaces (see Definition~\ref{defn:bord:MCG}) and morphisms are isotopy classes of diffeomorphisms relatively to the boundary which preserve $\de_1\Sigma$ but are allowed to permute the labelling of the other boundary components.

\medskip

Let us discuss the excision of pairs of pants. We have introduced in Definition~\ref{defn:homotopy:class:pop} various sets of homotopy classes of embeddings $\varphi \colon P \hookrightarrow \Sigma$. From now on, we only consider the ones corresponding to $m_0 = 1$, \emph{i.e.} $\varphi(\de_1 P) = \de_1 \Sigma$, and drop the index $m_0$ from the notations. We therefore have a set partition
\begin{equation}
	\mathcal{P}_{\Sigma} = \bigg(\bigsqcup_{m = 2}^{n} \mathcal{P}_{\Sigma}^{m}\bigg) \sqcup \mathcal{P}_{\Sigma}^{\emptyset}.
\end{equation}
If $[P] \in \mathcal{P}_{\Sigma}$, we denote by $\Sigma - P$ the bordered surface $\overline{\Sigma \setminus \varphi (P)}$ obtained from some representative $\varphi \colon P \hookrightarrow \Sigma$ of $[P]$. To define the labelling of the boundary components of $\Sigma - P$, we say that the (labelled) boundary components of $P$ that appear in $\Sigma - P$ are put first, followed by the (labelled) boundary components of $\Sigma$ that appear in $\Sigma - P$. Choosing another representative of $[P]$ gives a canonically isomorphic object.

\medskip

The geometric recursion (abbreviated GR) starts with a functor $E$ from $\mathcal{B}_{1}$ to a category of topological vector spaces and aims at constructing $E$-valued functorial assignments $\Sigma \mapsto \Omega_{\Sigma} \in E(\Sigma)$. For this purpose, $E$ must come with extra functorial data that satisfy a number of axioms, subsumed in the notion of \emph{target theory}. Instead of repeating the fully general definition of target theories and associated geo\-me\-tric recursion \cite{ABO17}, we shall describe in concrete terms the geometric recursion for the two examples of $E$ used in the present work, namely the spaces of measurable functions on the Teichm\"uller space of $\Sigma$, either seen as a space of hyperbolic or combinatorial structures on $\Sigma$. The first one was developed in \cite[Sections 7--10]{ABO17} while the second one is new. In general, if $\mathcal{X}$ is a topological space, we denote by ${\rm Mes}(\mathcal{X})$ the space of real-valued measurable functions on $\mathcal{X}$. Most results of this section still hold true after replacing everywhere  ``measurable'' by ``continuous''. 

\subsubsection{The hyperbolic case}
\label{subsubsec:GR}

For any bordered surface, we let $E(\Sigma) = {\rm Mes}(\mathcal{T}_{\Sigma})$. GR initial data consists of a quadruple $(A,B,C,D)$ where $A,B,C$ are measurable functions on $\mathcal{T}_{P} \cong \mathbb{R}_{+}^3$ such that $X(L_1,L_2,L_3) = X(L_1,L_3,L_2)$ for $X \in \{A,C\}$ and $D_T$ is a measurable function on $\mathcal{T}_{T}$. We require $T \mapsto D_T$ to be functorial. We let
\begin{equation}
	\Omega_{\varnothing} = 1,
	\qquad
	\Omega_{P}(\sigma) = A(\vec{\ell}_{\sigma}(\de P)),
	\qquad
	\Omega_{T} = D_T,
\end{equation}
where $\vec{\ell}_{\sigma}(\de P)$ is the ordered triple of hyperbolic lengths of the boundary components of $P$. For disconnected surfaces we set
\begin{equation}
	\Omega_{\Sigma_1 \sqcup \cdots \sqcup \Sigma_k}(\sigma_1,\ldots,\sigma_k) = \prod_{i = 1}^k \Omega_{\Sigma_i}(\sigma_i).
\end{equation}
It remains to define $\Omega_{\Sigma}$ for connected surfaces with Euler characteristic $\chi_{\Sigma} \leq -2$. This is done inductively $\chi_{\Sigma}$ by geometric recursion:
\begin{equation}\label{eqn:GRdef}
	\Omega_{\Sigma}(\sigma)
	=
	\sum_{m = 2}^n \sum_{[P] \in \mathcal{P}_{\Sigma}^{m}}
		B(\vec{\ell}_{\sigma}(\de P))\,\Omega_{\Sigma - P}(\sigma|_{\Sigma - P})
	+
	\frac{1}{2} \sum_{[P] \in \mathcal{P}_{\Sigma}^{\varnothing}}
		C(\vec{\ell}_{\sigma}(\de P))\,\Omega_{\Sigma - P}(\sigma|_{\Sigma - P}).
\end{equation}
Here, we choose as representative of $[P]$ the embedding as a pair of pants with geodesic boundaries. In this case, the restriction of $\sigma$ to $\Sigma - P$ makes it a hyperbolic surface with geodesic boundaries. The sum \eqref{eqn:GRdef} has countably many terms and therefore its convergence should be discussed.

\medskip

\begin{defn}\label{defn:adm:hyp}
	Let $\mathcal{T}_{\Sigma}^{(\epsilon)} \subset \mathcal{T}_{\Sigma}$ be the $\epsilon$-thick part of the Teichm\"uller space, \emph{i.e.} the set of $\sigma \in \mathcal{T}_{\Sigma}$ such that $\ell_{\sigma}(\gamma) \geq \epsilon$ for any $\gamma \in S_{\Sigma}$. We say that initial data $(A,B,C,D)$ are \emph{admissible} if for any $\epsilon > 0$ there exists $t \geq 0$  such that, for all $s \geq 0$, there exists $M_{\epsilon,s} > 0$ for which
	\begin{equation}\label{eqn:decay}
	\begin{aligned}
		& \sup_{L_1,L_2,L_3 \geq \epsilon} \frac{|A(L_1,L_2,L_3)|}{\big((1 + L_1)(1 + L_2)(1 + L_3)\big)^{t}} \leq  M_{\epsilon,0}
		&& \qquad
		\sup_{\sigma \in \mathcal{T}_{T}^{(\epsilon)}} \frac{|D_T(\sigma)|}{\bigl(1 + \ell_{\sigma}(\de T)\bigr)^{t}} \leq M_{\epsilon,0}
		\\
		&
		\sup_{L_1,L_2,\ell \geq \epsilon} \frac{|B(L_1,L_2,\ell)|\,(1 + [\ell - L_1 - L_2]_{+})^s}{\big((1 + L_1)(1 + L_2)\big)^{t}} \leq M_{\epsilon,s}
		&& \qquad
		\sup_{L_1,\ell,\ell' \geq \epsilon} \frac{|C(L_1,\ell,\ell')|\,(1 + [\ell + \ell' - L_1]_{+})^s}{(1 + L_1)^{t}} \leq M_{\epsilon,s}
	\end{aligned}
	\end{equation}
\end{defn}

\begin{thm} \cite[Corollary~8.3]{ABO17} \label{thm:hyp:GR}
	If $(A,B,C,D)$ are admissible initial data, then for any bordered surface $\Sigma$ 
	\begin{itemize}
		\item[$\bullet$] the series \eqref{eqn:GRdef} converges absolutely and uniformly on any compact of $\mathcal{T}_{\Sigma}$,
		\item[$\bullet$] $\Sigma \mapsto \Omega_{\Sigma} \in {\rm Mes}(\mathcal{T}_{\Sigma})$ is a well-defined functorial assignment (in particular, $\Omega_{\Sigma}$ is $\Mod_{\Sigma}^{\de}$-invariant),
		\item[$\bullet$] there exists $u \geq 0$ depending only on the topological type of $\Sigma$, such that for any $\epsilon > 0$ we have
		\begin{equation}\label{eqn:bound:hyp:GR}
			\sup_{\sigma \in \mathcal{T}_{\Sigma}^{(\epsilon)}} |\Omega_{\Sigma}(\sigma)| \leq K_{\epsilon} \prod_{b \in \pi_0(\de \Sigma)} \bigl( 1 + \ell_{\sigma}(b) \bigr)^{u}
		\end{equation}
		for some constant $K_{\epsilon}$ depending only on $\epsilon$ and the topological type of $\Sigma$.
	\end{itemize} 
	\hfill $\blacksquare$
\end{thm}
\subsubsection{The combinatorial case}\label{subsubsec:GRcomb}

For any bordered surface, we let $E(\Sigma) = {\rm Mes}(\mathcal{T}_{\Sigma}^{\textup{comb}})$. GR initial data consist of a quadruple $(A,B,C,D)$ such that $A,B,C$ are measurable functions on $\mathcal{T}_{P}^{\textup{comb}} \cong \mathbb{R}_{+}^3$ such that $X(L_1,L_2,L_3) = X(L_1,L_3,L_2)$ for $X \in \{A,C\}$ and $D$ is a measurable function on $\mathcal{T}_{T}^{\textup{comb}}$. We require $T \mapsto D_T$ to be functorial. We let
\begin{equation}
	\Xi_{\varnothing} = 1,
	\qquad
	\Xi_{P}(\GG) = A(\vec{\ell}_{\GG}(\de P)),
	\qquad
	\Xi_{T} = D_T,
\end{equation}
For disconnected surfaces we set
\begin{equation}
	\Xi_{\Sigma_1 \sqcup \cdots \sqcup \Sigma_k}(\GG_1,\ldots,\GG_k) = \prod_{i = 1}^k \Xi_{\Sigma_i}(\GG_i).
\end{equation}
For connected surfaces $\Sigma$ with $\chi_{\Sigma} \leq -2$, we define inductively
\begin{equation}\label{eqn:comb:GRdef}
	\Xi_{\Sigma}(\GG) = \sum_{m = 2}^n \sum_{[P] \in \mathcal{P}_{\Sigma}^{m}}
		B(\vec{\ell}_{\GG}(\de P))\,\Xi_{\Sigma - P}(\GG|_{\Sigma - P})
	+
	\frac{1}{2} \sum_{[P] \in \mathcal{P}_{\Sigma}^{\varnothing}}
		C(\vec{\ell}_{\GG}(\de P))\,\Xi_{\Sigma - P}(\GG|_{\Sigma - P}),
\end{equation}
where $\GG|_{\Sigma - P}$ have been defined by the cutting procedure in Section~\ref{subsec:cut:glue}. We remark that, given a homotopy class of embedded pair of pants, we can restrict the combinatorial structure $\GG$ to a representative of $P$, so that the triple of lengths $\vec{\ell}_{\GG}(\de P)$ makes sense.

\medskip

\begin{defn}\label{defn:adm:comb}
	Let $\mathcal{T}_{\Sigma}^{\textup{comb},(\epsilon)} \subset \mathcal{T}_{\Sigma}^{\textup{comb}}$ be the $\epsilon$-thick part of the combinatorial Teichm\"uller space, \emph{i.e.} the set of $\GG \in \mathcal{T}_{\Sigma}^{\textup{comb}}$ such that $\ell_{\GG}(\gamma) \geq \epsilon$ for any $\gamma \in S_{\Sigma}$. We say that initial data $(A,B,C,D)$ are \emph{admissible} if they satisfy the same conditions appearing in Definition~\ref{defn:adm:hyp}, except we use $\mathcal{T}_{T}^{\textup{comb},(\epsilon)}$ in the condition for $D_T$.
\end{defn}

\begin{thm}\label{thm:comb:GR}
	If $(A,B,C,D)$ are admissible initial data, then for any bordered surface $\Sigma$
	\begin{itemize}
		\item[$\bullet$] the series \eqref{eqn:comb:GRdef} converges absolutely and uniformly on any compact of $\mathcal{T}_{\Sigma}^{\textup{comb}}$,
		\item[$\bullet$] $\Sigma \mapsto \Xi_{\Sigma} \in {\rm Mes}(\mathcal{T}_{\Sigma}^{\textup{comb}})$ is a well-defined functorial assignment (in particular, $\Xi_{\Sigma}$ is $\Mod_{\Sigma}^{\de}$-invariant),
		\item[$\bullet$] there exists $u \geq 0$ depending only on the topological type of $\Sigma$, such that for any $\epsilon > 0$ we have
		\begin{equation}\label{eqn:bound:comb:GR}
			\sup_{\GG \in \mathcal{T}_{\Sigma}^{\textup{comb},(\epsilon)}} |\Xi_{\Sigma}(\GG)| \leq K_{\epsilon} \prod_{b \in \pi_0(\de \Sigma)} \bigl( 1 + \ell_{\GG}(b) \bigr)^{u}
		\end{equation}
		for some constant $K_{\epsilon}$ depending only on $\epsilon$ and the topological type of $\Sigma$.
	\end{itemize} 
\end{thm}

Although the spaces ${\rm Mes}(\mathcal{T}_{\Sigma})$ and ${\rm Mes}(\mathcal{T}_{\Sigma}^{\textup{comb}})$ can be identified via the spine homeomorphism of Theorem~\ref{thm:homeo:ord:comb:Teich}, the way we measure lengths and we restrict to $P$ and $\Sigma - P$ is different and as a result the hyperbolic/combinatorial structure in \eqref{eqn:GRdef} and \eqref{eqn:comb:GRdef} are completely different. So, for identical initial data, the hyperbolic and the combinatorial GR do not produce the same functions (even after identification of their domains). The relation between the hyperbolic and combinatorial GR is elucidated in Section~\ref{sec:hyp:comb}.

\begin{proof}
	The result follows from the general theory of \cite{ABO17} after proving that ${\rm Mes}(\mathcal{T}_{\Sigma}^{\textup{comb}})$ is a target theory. We present a self-contained proof which does not rely on these general notions, by specialising the strategy of \cite{ABO17} to this simpler setting.
	
	\smallskip
	
	It is enough to prove the result for connected surfaces. By definition of the initial data, the result holds for connected surfaces of Euler characteristic $-1$. Let us assume it holds for all surfaces of Euler characteristic strictly greater than $\chi$. Let $\Sigma$ be a bordered surface of type $(g,n)$ with $2 - 2g - n = \chi$, take $\epsilon > 0$ and fix $\GG \in \mathcal{T}_{\Sigma}^{\textup{comb},(\epsilon)}$. For any $[P] \in \mathcal{P}_{\Sigma}$, we have as well $\GG|_{\Sigma - P} \in \mathcal{T}_{\Sigma - P}^{\textup{comb},(\epsilon)}$. Therefore, by induction hypothesis, there exist $u \geq 0$ and $K_{\epsilon} > 0$ which we can choose to depend only on $\epsilon$ and the topological type of $\Sigma$, such that
	\[
		\big|\Xi_{\Sigma - P}(\GG_{\Sigma - P})\big| \leq K_{\epsilon} \, \prod_{b \in \pi_0(\de(\Sigma - P))} \bigl( 1 + \ell_{\GG}(b) \bigr)^{u}.
	\]
	We now study the absolute convergence of the GR series \eqref{eqn:comb:GRdef}. We use the notation $X_{P}$ for the function $B$ when $[P] \in \mathcal{P}_{\Sigma}^{m}$ and for the function $\frac{1}{2} C$ when $[P] \in \mathcal{P}_{\Sigma}^{\emptyset}$. We first isolate the sum over $\GG$-small pairs of pants. Using the fact that there are at most $2(6g - 6 + 3n)$ $\GG$-small pairs of pants (Remark~\ref{rem:small:pop}), together with the admissibility conditions on $X_{P}$ (Definition~\ref{defn:adm:comb}) and the inequality $(1 + L_1 + L_2)^{t} \leq (1 + L_1)^{t}(1 + L_2)^{t}$ for any $L_1,L_2 > 0$, we get
	\begin{equation}\label{eqn:bound:small:pop}
		\sum_{\substack{[P] \in \mathcal{P}_{\Sigma} \\ \GG\textrm{-small}}} \big|X_{P}(\vec{\ell}_{\GG}(\de P))\,\Xi_{\Sigma - P}(\GG|_{\Sigma - P})\big|
		\leq
		2(6g - 6 + 3n) \, M_{\epsilon,0} K_{\epsilon} \prod_{i = 1}^n \bigl( 1 + \ell_{\GG}(\de_i\Sigma)\bigr)^{\max\set{u,t}}.
	\end{equation}
	We now turn to the contributions of the $\GG$-big pairs of pants in $\mathcal{P}_{\Sigma}^{m}$. We have for any $s > 0$
	\[
	\begin{split}
		& \sum_{\substack{[P] \in \mathcal{P}_{\Sigma}^{m} \\ \GG \textrm{-big}}} \big|B(\vec{\ell}_{\GG}(\de P)) \, \Xi_{\Sigma - P}(\GG|_{\Sigma - P})\big| \\
		\leq &
		M_{\epsilon,s} K_{\epsilon}
		\Biggl(
			\prod_{i \neq 1,m} \bigl( 1 + \ell_{\GG}(\de_i\Sigma) \bigr)^{u}
		\Biggr)
		\Biggl(
			\sum_{L \in \ell_{\GG}(\de_1\Sigma) + \ell_{\GG}(\de_m\Sigma) + \NN}
				\frac{(2 + L)^{t}\,\#\set{\gamma \in S_{\Sigma}^{\circ} | L \leq \ell_{\GG}(\gamma) < L + 1}}{\big(1 + L - \ell_{\GG}(\de_1\Sigma) - \ell_{\GG}(\de_m\Sigma)\big)^{s}}
		\Biggr) \\
		\leq &
		M_{\epsilon,s} K_{\epsilon} m_{\epsilon}
		\Biggl(
			\prod_{i \neq 1,m} (1 + \ell_{\GG}(\de_i\Sigma))^{u}
		\Biggr)
		\Biggl(
			\sum_{L \geq 1}  \bigl( 1 + L + \ell_{\GG}(\de_1\Sigma) + \ell_{\GG}(\de_m\Sigma) \bigr)^{t + 6g - 6 + 2n}L^{-s}
		\Biggr).
	\end{split}
	\]
	In the last line, we invoked the polynomial growth of the number of multicurves with respect to combinatorial length, justified later in Proposition~\ref{prop:poly:growth}. Specialising to $s = (6g - 6 + 2n + t) + 2$ makes the sum in brackets converging to a polynomial of degree $t' = t + 6g - 6 + 2n$ in the variable $\ell_{\GG}(\de_1\Sigma) + \ell_{\GG}(\de_m\Sigma)$ and, together with \eqref{eqn:bound:small:pop}, it implies the existence of a constant $K'_{\epsilon} > 0$ such that
	\[
		\sum_{\substack{[P] \in \mathcal{P}_{\Sigma}^{m}}} \big|B(\vec{\ell}_{\GG}(\de P))\,\Xi_{\Sigma - P}(\GG|_{\Sigma - P})\big| \leq K'_{\epsilon} \prod_{i = 1}^n \bigl( 1 + \ell_{\GG}(\de_i\Sigma) \bigr)^{\max\set{u,t'}}.
	\]
	A similar argument shows that
	\[
		\sum_{\substack{[P] \in \mathcal{P}_{\Sigma}^{\emptyset}}} \big|C(\vec{\ell}_{\GG}(\de P))\,\Xi_{\Sigma - P}(\GG|_{\Sigma - P})\big| \leq K'_{\epsilon} \prod_{i = 1}^n \bigl( 1 + \ell_{\GG}(\de_i\Sigma) \bigr)^{\max\set{u,t'}}
	\]
	for a perhaps larger constant $K'_{\epsilon}$. Consequently, the series
	\[
		\sum_{[P] \in \mathcal{P}_{\Sigma}} X_{P}(\vec{\ell}_{\GG}(\de P))\,\Xi_{\Sigma - P}(\GG|_{\Sigma - P})
	\]
	converges absolutely and uniformly on any compact of $\mathcal{T}_{\Sigma}^{\textup{comb},(\epsilon)}$, to a limit that we denote $\Xi_{\Sigma}$. Further, the bounds that we just proved imply that this limit satisfies
	\[
		\forall \GG \in \mathcal{T}_{\Sigma}^{\textup{comb},(\epsilon)},\qquad |\Xi_{\Sigma}(\GG)| \leq K''_{\epsilon} \prod_{i = 1}^n \bigl( 1 + \ell_{\GG}(\de_i\Sigma) \bigr)^{u'}
	\]
	for some constant $K_{\epsilon}'' > 0$ and $u' = \max\set{u,t'}$. The proof is then completed by induction.
\end{proof}

\subsection{(Discrete) integration and topological recursion}
\label{subsec:(disc):integr:TR}

Since the functions produced by the hyperbolic or combinatorial GR (generically called GR amplitudes) are pure mapping class group invariant, they descend to functions on the corresponding moduli spaces. For a connected surface $\Sigma$ of type $(g,n)$ and hyperbolic (resp. combinatorial) GR amplitudes $\Omega_{\Sigma}$ (resp. $\Xi_{\Sigma}$), we denote by $\Omega_{g,n}$ (resp. $\Xi_{g,n}$) the functions induced on the associated moduli spaces. For the initial datum $T \mapsto D_T$, we denote by $D$ the induced function on $\mathcal{M}_{1,1}$ or $\mathcal{M}_{1,1}^{\textup{comb}}$.

\medskip

In the first part of this section we discuss how to integrate hyperbolic and combinatorial GR amplitudes respectively against the Weil--Petersson measure and the Kontsevich measure. In the second part, which belongs exclusively to the combinatorial setting, we discuss how to define discrete integration on the combinatorial moduli space via sums over integral metric ribbon graphs.

\subsubsection{Integration and topological recursion}
\label{subsubsec:integr:TR}

For $L \in \mathbb{R}_{+}^n$, $\mathcal{M}_{g,n}(L)$ is naturally endowed with the Weil--Petersson measure $\mu_{\textup{WP}}$, while $\mathcal{M}_{g,n}^{\textup{comb}}(L)$ is equipped with the Kontsevich measure $\mu_{\textup{K}}$. It is important to note that these measures are compatible with cutting, as expressed in the hyperbolic setting by \cite[Theorem~7.1]{Mir07simple} and in the combinatorial setting by Proposition~\ref{prop:integration}. Thus, integration of functions obtained by GR against these measures over the moduli space with fixed boundary lengths produce functions on $\mathbb{R}_{+}^n$ that also satisfy a recursion on the Euler characteristic, called topological recursion (TR). In order to guarantee integrability we are going to introduce stronger assumptions on the initial data. 

\begin{defn}\label{defn:strong:adm}
	In the hyperbolic setting, the initial data $(A,B,C,D)$ are cal\-led \emph{strong\-ly admissible} if there exists $\eta \in [0,2)$ and $t \geq 0$ such that, for any $s \geq 0$, there exists $M_{s} > 0$ such that for any $L_1,L_2,L_3,\ell,\ell' > 0$
	\begin{equation}
	\begin{split}
		|A(L_1,L_2,L_3)| & \leq M_0\,\big((1 + L_1)(1 + L_2)(1 + L_3)\big)^{t}, \\
		|B(L_1,L_2,\ell)| & \leq \frac{M_{s}\,\big((1 + L_1)(1 + L_2)\big)^{t}}{\ell^{\eta}\big(1 + [\ell - (L_1 + L_2)]_{+}\big)^{s}}, \\
		|C(L_1,\ell,\ell')| & \leq \frac{M_s\,(1 + L_1)^{t}}{(\ell\ell')^{\eta}\big(1 + [\ell + \ell' - L_1]_{+}\big)^{s}},
	\end{split} 
	\end{equation}
	and $D$ is integrable on $\mathcal{M}_{1,1}(L_1)$ and satisfies
	\begin{equation}\label{eqn:D:cond:hyp}
		\bigg|\int_{\mathcal{M}_{1,1}(L_1)} D \, \dd\mu_{\textup{WP}}\bigg| \leq M_0\,(1 + L_1)^{t}.
	\end{equation}
	We define strong admissibility of the initial data $(A,B,C,D)$ in the combinatorial setting by the same conditions, except for \eqref{eqn:D:cond:hyp} which gets substituted by
	\begin{equation}\label{eqn:D:cond:comb} 
		\bigg|\int_{\mathcal{M}_{1,1}^{\textup{comb}}(L_1)} D \, \dd\mu_{\textup{K}}\bigg| \leq M_0\,(1 + L_1)^{t}.
	\end{equation}
\end{defn}

Whenever these integrals are defined, we use the generic notations
\begin{equation}\label{defn:Vamplitudes}
	V\Omega_{g,n}(L) = \int_{\mathcal{M}_{g,n}(L)} \Omega_{g,n} \,\dd\mu_{\textup{WP}},
	\qquad \qquad
	V\Xi_{g,n}(L) = \int_{\mathcal{M}_{g,n}^{\textup{comb}}(L)} \Xi_{g,n} \,\dd\mu_{\textup{K}}.
\end{equation}

\begin{thm} \cite[Theorem 8.8]{ABO17} \label{thm:hyp:GR:TR}
	Let $(A,B,C,D)$ be strongly admissible initial data for the hyperbolic GR, and $\Omega_{\Sigma}$ be the resulting functions. Then, $\Omega_{\Sigma}$ is integrable against $\mu_{\textup{WP}}$ on $\mathcal{M}_{\Sigma}(L)$ for any $L \in \mathbb{R}_{+}^n$, and the integrals satisfy the following recursion on $2g - 2 + n > 1$.
	\begin{equation}\label{eqn:hyp:TR}
	\begin{split}
		&
		V\Omega_{g,n}(L_1,\ldots,L_n) = \\
		= &
		\sum_{m = 2}^n \int_{\mathbb{R}_{+}} \!\!\!\! \, B(L_1,L_m,\ell) \, V\Omega_{g,n - 1}(\ell,L_2,\ldots,\widehat{L_m},\ldots,L_n) \, \ell \dd\ell \\
		&
		+ \frac{1}{2} \int_{\mathbb{R}_{+}^2} \!\!\!\! C(L_1,\ell,\ell')\bigg(V\Omega_{g - 1,n + 1}(\ell,\ell',L_2,\ldots,L_n) + \!\!\!\!\! \sum_{\substack{h + h' = g \\ J \sqcup J' = \{L_2,\ldots,L_n\}}} \!\!\!\!\! V\Omega_{h,1 + \#J}(\ell,J) \, V\Omega_{h',1 + \#J'}(\ell',J')\bigg)  \ell \ell' \dd\ell \, \dd\ell' 
	\end{split}
	\end{equation}
	with the conventions $V\Omega_{0,1} = 0$ and $V\Omega_{0,2} = 0$, and the base cases
	\begin{equation}
		V\Omega_{0,3}(L_1,L_2,L_3) = A(L_1,L_2,L_3)
		\qquad \text{and}\qquad
		V\Omega_{1,1}(L_1) = \int_{\mathcal{M}_{1,1}(L_1)} D \, \dd\mu_{\textup{WP}}.
	\end{equation}
	\hfill $\blacksquare$
\end{thm}

\begin{thm} \label{thm:comb:GR:TR}
	Let $(A,B,C,D)$ be strongly admissible initial data for the combinatorial GR, and $\Xi_{\Sigma}$ be the resulting functions. Then, $\Xi_{\Sigma}$ is integrable against $\mu_{\textup{K}}$ on $\mathcal{M}_{\Sigma}^{\textup{comb}}(L)$ for any $L \in \mathbb{R}_{+}^n$, and the integrals $V\Xi_{g,n}(L)$ satisfy the same recursion as in \eqref{eqn:hyp:TR}, \emph{i.e.} 
	\begin{equation}\label{eqn:comb:TR}
	\begin{split}
		&
		V\Xi_{g,n}(L_1,\ldots,L_n) = \\
		= &
		\sum_{m = 2}^n \int_{\mathbb{R}_{+}} \!\!\!\!\,  \, B(L_1,L_m,\ell) \, V\Xi_{g,n - 1}(\ell,L_2,\ldots,\widehat{L_m},\ldots,L_n)  \ell \dd \ell \\
		&
		+ \frac{1}{2} \int_{\mathbb{R}_{+}^2} \!\!\!\! C(L_1,\ell,\ell')\bigg(V\Xi_{g - 1,n + 1}(\ell,\ell',L_2,\ldots,L_n) + \!\!\!\!\! \sum_{\substack{h + h' = g \\ J \sqcup J' = \{L_2,\ldots,L_n\}}} \!\!\!\!\! V\Xi_{h,1 + \#J}(\ell,J) \, V\Xi_{h',1 + \#J'}(\ell',J')\bigg) \ell\ell'   \dd \ell \, \dd \ell' 
	\end{split}
	\end{equation}
	with the conventions $V\Xi_{0,1} = 0$ and $V\Xi_{0,2} = 0$, and the base cases
	\begin{equation}
		V\Xi_{0,3}(L_1,L_2,L_3) = A(L_1,L_2,L_3)
		\qquad \text{and} \qquad
		V\Xi_{1,1}(L_1) = \int_{\mathcal{M}_{1,1}^{\textup{comb}}(L_1)} D \, \dd\mu_{\textup{K}}.
	\end{equation}
\end{thm}

\begin{proof}
	We first note that the initial data $V\Xi_{0,3}$ is well-defined as $A$ is, and that $V\Xi_{1,1}$ is well-defined by strong admissibility.

	\smallskip

	Now, for a connected surface $\Sigma$ of type $(g,n)$ with $2g - 2 + n > 1$, apply the integration over $\mathcal{M}_{g,n}^{\textup{comb}}(L)$ with respect to $\mu_K$ to both sides of the combinatorial GR \eqref{eqn:comb:GRdef}. We analyse the integration of the sum over $[P] \in \mathcal{P}_{\Sigma}^{m}$. Let $\Gamma$ be the $\Mod_{\Sigma}^{\de}$-orbit of a simple closed curve bounding a pair of pants together with $\de_1 \Sigma$ and $\de_m \Sigma$. We have
	\[
		\sum_{[P] \in \mathcal{P}_{\Sigma}^{m}} B(\vec{\ell}_{\GG}(\de P))\,\Xi_{\Sigma - P}(\GG|_{\Sigma - P})
		=
		\sum_{\alpha \in \Gamma} B(L_1,L_m,\ell_{\GG}(\alpha)) \, \Xi_{\Sigma - P_{\alpha}}(\GG|_{\Sigma-P_{\alpha}}),
	\]
	where $P_{\alpha}$ is the pair of pants bounded by $\de_1 \Sigma$, $\de_m \Sigma$ and $\alpha$. Now applying Proposition~\ref{prop:integration} to $f(L,\ell) = B(L_{1},L_{m},x)$ and the assignment $\Sigma_{\alpha} \mapsto \Xi_{\Sigma -P_{\alpha}}$, we find
	\[
	\begin{split}
		& \int_{\mathcal{M}_{g,n}^{\textup{comb}}(L)} \sum_{\alpha \in \Gamma} B(L_1,L_m,\ell_{\ast}(\alpha)) \, \Xi_{g,n-1} \, \dd \mu_K \\
		= & 
		\int_{\mathbb{R}_{+}} \!\!\!\! \,B(L_1,L_m,\ell) \left(\int_{\mathcal{M}_{g,n-1}^{\textup{comb}}(\ell, L_2,\ldots,\widehat{L_m},\ldots,L_n)} \Xi_{g,n-1} \,\dd \mu_K\right)  \ell \dd \ell\\
		= &
		\int_{\mathbb{R}_{+}} \!\!\!\! \,B(L_1,L_m,\ell) \, V\Xi_{g,n-1}(\ell,L_2,\ldots,\widehat{L_m},\ldots,L_n) \ell \dd \ell.
	\end{split}
	\]
	The treatment of the $C$ summands is similar, the main difference being that the excised pair of pants has two simple closed curves in $\Sigma^{\circ}$, whose lengths are part of the combinatorial Fenchel--Nielsen coordinates over which we need to integrate.
\end{proof}

\subsubsection{Discrete integration and topological recursion}
\label{subsubsec:disc:integr:TR}

As the combinatorial moduli spaces have an integral structure, we can also study discrete integration of GR amplitudes. Again, the geometric recursion here is the key property that guarantees a topological recursion for the discrete integrals, \emph{i.e.} the sum over lattice points. 

\begin{defn} \label{defn:MgncombZ}
	For $2g-2+n>0$, let $\mathcal{M}^{\textup{comb},\ZZ}_{g,n} \subset \mathcal{M}^{\textup{comb}}_{g,n}$ be the set of classes of metric ribbon graphs whose edge lengths are positive integers, and  $\mathcal{M}^{\textup{comb},\ZZ}_{g,n}(L)$ the one with fixed perimeter $L \in \ZZ_{+}^n$. We denote likewise $\mathcal{T}^{\textup{comb},\ZZ}_{\Sigma}$ and $\mathcal{T}^{\textup{comb}}_{\Sigma}(L)$ their preimages in the combinatorial Teichm\"uller space of a bordered surface $\Sigma$ of type $(g,n)$.
\end{defn}

Since for any $\bm{G} \in \mathcal{M}^{\textup{comb},\ZZ}_{g,n}(L)$, we have $\sum_{i = 1}^n L_i = \sum_{e \in E_{G}} 2\ell_{e}$, the set $\mathcal{M}^{\textup{comb},\ZZ}_{g,n}(L)$ is finite for any fixed $L$, and it is empty if $\sum_{i = 1}^{n} L_i$ is odd. For instance,
\begin{equation}
	\mathcal{T}_{P}^{\textup{comb},\ZZ}
	\cong
	\mathcal{M}_{0,3}^{\textup{comb},\ZZ}
	\cong
	\Set{ (L_1,L_2,L_3) \in \ZZ_+^3 | L_1 + L_2 + L_3 \text{ is even} }.
\end{equation}
To handle this situation, for $X \colon \RR_{+}^n \to \RR$ we introduce the notation
\begin{equation}
	X_{\ZZ}(L)
	=
	\begin{cases}
		X(L) & \text{if $L \in \ZZ_{+}^n$ and $\sum_{i = 1}^n L_i$ is even,} \\
		0 & \text{otherwise.}
	\end{cases}
\end{equation}
Further, for any function $\Xi_{g,n}$ on $\mathcal{M}^{\textup{comb}}_{g,n}$, set
\begin{equation}\label{eqn:lattice:pnt:count}
	N\Xi_{g,n}(L) = \sum_{\bm{G} \in \mathcal{M}^{\textup{comb},\ZZ}_{g,n}(L)} \frac{\Xi_{g,n}(\bm{G})}{\#\Aut(\bm{G})} = \sum_{G \in \mathcal{R}_{g,n}} \frac{1}{\#\Aut(G)} \sum_{x \in P_{G}(L)\cap \ZZ_{+}^n} \Xi_{g,n}(x),
\end{equation}
where we recall that $P_G(L)$ is the set of metrics on $G$ with perimeter $L$.

\medskip

\begin{thm}\label{thm:disc:TR}
	Let $A,B,C$ be three functions on $\RR_{+}^3$ such that $A$ and $C$ are symmetric under exchange of their last two variables, and satisfying
	\begin{equation}\label{eqn:vanishing:small:pop}
		\ell > L_1 + L_2 \quad \Longrightarrow \quad B(L_1,L_2,\ell) = 0,
		\qquad \text{and} \qquad
		\ell + \ell' > L_1 \quad \Longrightarrow \quad C(L_1,\ell,\ell') = 0.
	\end{equation}
	Let $D$ be $\Mod_{T}$-invariant function on $\mathcal{T}_{T}^{\textup{comb}}$. Denote by $\Xi_{\Sigma}$ the corresponding combinatorial GR amplitudes. We have the following recursion on $2g - 2 + n > 1$.
	\begin{equation}\label{eqn:discr:TR}
		\begin{split}
		& N\Xi_{g,n}(L_1,\ldots,L_n) \\
		= &
		\sum_{m = 2}^n \sum_{\ell \geq 1} \ell \, B_{\ZZ}(L_1,L_m,\ell) \, N\Xi_{g,n - 1}(\ell,L_2,\ldots,\widehat{L_m},\ldots,L_n) \\
		&
		+ \frac{1}{2} \sum_{\ell,\ell' \geq 1} \ell\ell'\,C_{\ZZ}(L_1,\ell,\ell')\bigg(N\Xi_{g-1,n+1}(\ell,\ell',L_2,\ldots,L_n) +  \!\!\!\!\! \sum_{\substack{h + h' = g \\ J \sqcup J' = \{L_2,\ldots,L_n\}}} \!\!\!\!\! N\Xi_{h,1+\#J}(\ell,J) \, N\Xi_{h',1+\#J'}(\ell',J')\bigg)
	\end{split}
	\end{equation}
	with conventions $N\Xi_{0,1} = 0$ and $N\Xi_{0,2} = 0$, and base cases
	\begin{equation}
		N\Xi_{0,3}(L_1,L_2,L_3) = A_{\ZZ}(L_1,L_2,L_3)
		\qquad \text{and} \qquad
		N\Xi_{1,1}(L_1) = \sum_{\bm{G} \in \mathcal{M}_{1,1}^{\textup{comb},\ZZ}(L_1)} \frac{D(\bm{G})}{\#\Aut(\bm{G})}.
	\end{equation}
\end{thm}

Since cutting combinatorial structures preserve their integrability, in order to obtain the GR amplitudes in Theorem~\ref{thm:disc:TR}, it is sufficient to have $(A,B,C)$ defined on $\mathcal{T}_{P}^{\textup{comb},\ZZ}$ and $D$ defined on $\mathcal{T}_{T}^{\textup{comb},\ZZ}$. When the vanishing condition \eqref{eqn:vanishing:small:pop} hold, the GR sums have only finitely many non-zero terms: they are always well-defined, without the need of admissibility conditions for $(A,B,C,D)$.

\medskip

We can recover the continuous integration of Theorem~\ref{thm:comb:GR:TR} as a limit where we rescale the mesh of the lattice down to $0$. If $k > 0$, we let $\mathcal{M}_{g,n}^{\textup{comb},\ZZ/k}$ be the set of metric ribbon graphs in $\mathcal{M}_{g,n}^{\textup{comb}}$ whose edge lengths all become integral after dilation by $k$.

\begin{prop}\label{prop:limit:lattice:TR}
	Assume{} that $(A,B,C,D)$ are continuous functions on their respective combinatorial Teichm\"uller spaces, satisfying the vanishing conditions \eqref{eqn:vanishing:small:pop}, and such that for any fixed $L_1,L_2 > 0$, the functions $\ell \mapsto B(L_1,L_2,\ell)$ and $(\ell,\ell') \mapsto C(L_1,L_2,\ell')$ are bounded, and the function $D_{T}$ is bounded on $\mathcal{T}_{T}^{\textup{comb}}(L_1)$. Then $(A,B,C,D)$ is strongly admissible. If $d$ be a positive integer and $L \in (\ZZ_{+}/d)^{n}$, we have for $2g - 2 + n > 0$
	\begin{equation}\label{eqn:limit:lattice}
		\lim_{\substack{k \to \infty \\ k \in d\ZZ_+}}
		\bigg(\frac{1}{k^{6g-  6 + 2n}} \sum_{\bm{G} \in \mathcal{M}_{g,n}^{\textup{comb},\ZZ/k}(L)} \frac{\Xi_{g,n}(\bm{G})}{\#\Aut(\bm{G})}\bigg)
		=
		\begin{cases}
			2^{- 2g + 3 - n} \, V\Xi_{g,n}(L)
			&
			\text{if $\sum_{i = 1}^n d \cdot L_i$ is even,} \\
			0 & \text{otherwise.}
		\end{cases}
	\end{equation}
\end{prop}

The continuous integration also appears in the asymptotics of the lattice count for large boundary lengths.

\begin{cor}\label{cor:limit:lattice:TR}
	Assume that $(A,B,C,D)$ are continuous functions on their respective combinatorial Teichm\"uller spaces and satisfy the vanishing conditions \eqref{eqn:vanishing:small:pop}. We further assume the existence of $b,c \in \RR$ for which, for any $M > 0$ and $L_1,L_2,\ell,\ell' \in (0,M]$
	\begin{itemize}
		\item[$\bullet$] $k^{-2b -c} \, A(k L_1,k L_2,k L_3)$ converges uniformly to $\hat{A}(L_1,L_2,L_3)$ as $k \to \infty$,
		\item[$\bullet$] $k^{-b} \, B(k L_1,k L_2,k \ell)$ converges uniformly to $\hat{B}(L_1,L_2,\ell)$ as $k \to \infty$,
		\item[$\bullet$] $k^{-c} \, C(k L_1,k \ell,k \ell')$ converges uniformly to $\hat{C}(L_1,\ell,\ell')$ as $k \to \infty$,
		\item[$\bullet$] $k^{-b} \, D_{T}(k\GG)$ converges uniformly to $\hat{D}_{T}(\GG)$ as $k \to \infty$, for $\GG \in \bigcup_{\ell \in (0,M]}\mathcal{T}_{T}^{\textup{comb}}(\ell)$, where $k \GG$ is obtained from $\GG$ by dilation of the metric by a factor $k$,
		\item[$\bullet$] $(\hat{A},\hat{B},\hat{C},\hat{D})$ satisfy the assumptions of Proposition~\ref{prop:limit:lattice:TR}.
	\end{itemize}
	Denoting by $\hat{\Xi}_{\Sigma}$ the GR amplitudes associated to the initial data $(\hat{A},\hat{B},\hat{C},\hat{D})$, we have for any $d \in \ZZ_+$, $L \in (\ZZ_{+}/d)^{n}$ and $2g - 2 + n > 0$ 
	\begin{equation}
		\lim_{\substack{k \to \infty \\ k \in d\ZZ_+}}
		\frac{N\Xi_{g,n}(k L)}{k^{(g - 1)(b + c) + nb} \cdot k^{6g - 6 + 2n}}
		=
		\begin{cases}
			2^{- 2g + 3 - n} \, V\hat{\Xi}_{g,n}(L)
			& \text{if $\sum_{i = 1}^n d\cdot L_i$ even,} \\
			0
			& \text{otherwise.}
		\end{cases}
	\end{equation}
\end{cor}

\begin{proof}[Proof of Theorem~\ref{thm:disc:TR}]
	When $\sum_{i = 1}^n L_i$ is odd, both sides vanish, so we only need to prove the result when $\sum_{i = 1}^n L_i$ is even, which we now assume. The vanishing conditions for $B$ and $C$ imply that for each $L \in \ZZ_{+}^n$, the sums in the right-hand side have finitely many terms.

	\smallskip

	Let us substitute the GR sum for $\Xi_{\Sigma}$ in \eqref{eqn:lattice:pnt:count}. For $m \in \{2,\ldots,n\}$, we examine in detail how to handle the term
	\[
		\Xi_{\Sigma}^{B,m}(\GG) = \sum_{[P] \in \mathcal{P}_{\Sigma}^{m}} B(\vec{\ell}_{\GG}(\de P)) \,\Xi_{g,n-1}(\GG|_{\Sigma - P}).
	\] 
	Let $\Gamma$ be the $\Mod_{\Sigma}^{\de}$-orbit of a simple closed curve that bounds some $[P] \in \mathcal{P}_{\Sigma}^{m}$. Adapting the notation of Section~\ref{subsec:integration}, we denote by $\mathcal{M}^{\textup{comb},\Gamma,\ZZ}_{g,n}(L)$ the integral points in $\mathcal{M}^{\textup{comb},\Gamma}_{g,n}(L)$. Then
	\begin{equation*} 
	\begin{split}
		N\Xi_{g,n}^{B,m}(L)
		& = \sum_{\bm{G} \in \mathcal{M}^{\textup{comb},\ZZ}_{g,n}(L)}
			\frac{\Xi_{g,n}^{B,m}(\bm{G})}{\# \Aut(\bm{G})} \\
		& = \sum_{(\bm{G},\alpha) \in \mathcal{M}^{\textup{comb},\Gamma,\ZZ}_{g,n}(L)}
			\frac{1}{\# \Aut(\bm{G})} \, B(L_1,L_m,\mathcal{L}^{\Gamma}(\bm{G},\alpha)) \, \Xi_{g,n - 1}(\Pi(\bm{G})) \\
		& = \sum_{\ell\geq 1} \sum_{(\bm{G},\alpha) \in \mathcal{M}^{\textup{comb},\Gamma,\ZZ}_{g,n}(L)[\ell]}
			\frac{1}{\# \Aut(\bm{G})} \, B(L_1,L_m,\ell) \,\Xi_{g,n - 1}(\Pi(\bm{G})),
	\end{split}   
	\end{equation*}
	where we recall from Section~\ref{subsec:integration} the map
	\[
		\mathcal{L}^{\Gamma} \colon \mathcal{M}^{\textup{comb},\Gamma,\ZZ}_{g,n}(L) \longrightarrow \ZZ_{+}
	\]
	assigning to $(\bm{G},\alpha)$ the combinatorial length with respect to $\bm{G}$ of $\alpha$. It has fibers $(\mathcal{L}^{\Gamma})^{-1}(\ell) = \mathcal{M}^{\textup{comb},\Gamma,\ZZ}_{g,n}(L)[\ell]$. Moreover, we have the projection map
	\[
		\Pi \colon \mathcal{M}^{\textup{comb},\Gamma,\ZZ}_{g,n}(L)[\ell] \longrightarrow \mathcal{M}^{\textup{comb},\ZZ}_{g,n-1}(\ell,L_2,\dots,\widehat{L_m},\dots,L_n).
	\]
	We want to cluster this sum according to the fibres of the map $\Pi$. We first notice that, as $\mathcal{M}^{\textup{comb},\Gamma,\ZZ}_{g,n}(L)[\ell]$ is empty when $\ell$ and $L_1 + L_m$ have different parity, we can replace $B$ by $B_{\ZZ}$, and now only consider $\ell$ that has same parity as $L_1 + L_m$. In this case, for any $\bm{G}' \in \mathcal{M}^{\textup{comb},\ZZ}_{g,n - 1}(\ell,L_2,\dots,\widehat{L_m},\dots,L_n)$ and any $\bm{G} \in \Pi^{-1}(\bm{G}')$, we remark that $\#\Aut(\bm{G'}) = \#\Aut(\bm{G})$. Thus
	\[
		\frac{1}{\#\Aut(\bm{G})} B_{\ZZ}(L_1,L_m,\ell) \, \Xi_{g,n-1}(\Pi(\bm{G}))
	\]
	is constant on the fibres of $\Pi$. Due to the vanishing conditions \eqref{eqn:vanishing:small:pop}, the points $(\bm{G},\alpha)$ with non-trivial contribution have associated pair of pants $[P_{\alpha}]$ that is small. Therefore from Corollary~\ref{cor:twist:small:pop}, $\Pi^{-1}(\bm{G}')$ is in bijection with the set of $[\tau] \in \RR/\ell \ZZ$ such that the gluing of $\GG'$ to the combinatorial structure of the pair of pants with boundary lengths $(L_1,L_m,\ell)$ after a twist $\tau$ produces a combinatorial structure with integral edge lengths. In order to satisfy the latter condition, the twists must belong to the set $\left(\tau_{\gamma}(\GG)+\ZZ\right)/\ell \ZZ\simeq\ZZ/\ell\ZZ \simeq \Pi^{-1}(\bm{G}')$, whose cardinality is $\ell$. Therefore, our sum becomes
	\begin{equation}\label{eqn:lattice:B}
	\begin{split}
		\Xi_{g,n}^{B}(L)
		& =\sum_{\ell \geq 1} \sum_{\bm{G}' \in \mathcal{M}_{g,n - 1}(\ell,L_2,\ldots,\widehat{L_m},\ldots,L_n)} \frac{1}{\#\Aut(\bm{G}')}\,B_{\ZZ}(L_1,L_m,\ell) \, \Xi_{g,n - 1}(\bm{G}')  \\
		& = \sum_{\ell \geq 1} \ell\,B_{\ZZ}(L_1,L_m,\ell) \, N\Xi_{g,n - 1}(\ell,L_2,\ldots,\widehat{L_m},\ldots,L_n).
	\end{split}
	\end{equation}
	The treatment of the $C$ summand is similar, except that we should be cautious about automorphism factors. For each of the finitely many $\Mod_{\Sigma}^{\partial}$-orbit $\Gamma$ of a simple closed curve that bounds some $[P] \in \mathcal{P}_{\Sigma}^{\varnothing}$, and we observe that
	\[
		\#\Aut(\Gamma) =
		\begin{cases}
			2 & \text{if $\Sigma - P$ is connected,} \\
			1 & \text{otherwise.}
		\end{cases}
	\]
	Since in the analogue of \eqref{eqn:lattice:B} we have for any $\bm{G} \in \Pi^{-1}(\bm{G}')$ 
	\[
		\#\Aut(\Gamma) \cdot \#\Aut(\bm{G}) = \#\Aut(\bm{G}'),
	\]
	the automorphism factors are again naturally included in $N\Xi_{\Gamma}$, and we get the $C$-terms in \eqref{eqn:discr:TR} without extra automorphism factors (as the $\tfrac{1}{2}$ is already present in $C$).
\end{proof}
\begin{rem}
	The vanishing assumptions \eqref{eqn:vanishing:small:pop} are essential to allow the use of Corollary~\ref{cor:twist:small:pop}. If they did not hold, the fibres $\Pi^{-1}(\bm{G}')$ of the gluing fibration could, and do, meet integral non-admissible twists, hence their cardinality could be smaller than $\ell$ or $\ell \ell'$ and depend on $\bm{G}'$. It would then not be possible to derive a recursion for the weighted sum over lattice points. This problem did not arise for the integration against $\mu_{\textup{K}}$ as the set of non-admissible twists has zero measure with respect to $\mu_{\textup{K}}$.
\end{rem}

Before turning to the proof of Proposition~\ref{prop:limit:lattice:TR}, we need two preliminary results.

\begin{proof}[Proof of Proposition~\ref{prop:limit:lattice:TR}]
	The case $(g,n) = (0,3)$ is obvious, as there is equality before taking the limit. This initial case is special with respect to the other topologies since the moduli space is reduced to a point: in the rest of the proof, we suppose $(g,n)\neq(0,3)$. In general, when $\sum_{i = 1}^n d\cdot L_i$ is not even, the set $\mathcal{M}_{g,n}^{\textup{comb},\ZZ/k}(L)$ with $k \in d\ZZ_+$ is empty, so the left-hand side of Equation~\eqref{eqn:limit:lattice} vanishes, which proves half of the result. Hereafter we assume that $k \in d\ZZ_{+}$, and fix $L \in (\ZZ_{+}/d)^n$ such that $\sum_{i = 1}^n d\cdot L_i$ is even.

	\smallskip

	The thesis follows now from a general discussion. Consider a bounded function $f$ defined on $\mathcal{M}_{g,n}^{\textup{comb}}(L)$. The sum over rescaled lattice points is, by definition,
	\[
		\sum_{\bm{G} \in \mathcal{M}_{g,n}^{\textup{comb},\ZZ/k}(L)}
			\frac{f(\bm{G})}{\#\Aut(\bm{G})}
		=
		\sum_{G \in \mathcal{R}_{g,n}}
			\frac{1}{\#\Aut(G)} \sum_{x \in P_{G}(L) \cap \ZZ_{+}^n/k} f(x),
	\]
	where we recall that $P_{G}(L) \subset \RR_{+}^{E_{G}}$ is the set of metrics on $G$ with perimeters $L$. We first estimate the sum over non-trivalent graphs by
	\[
		\Biggl|
		\sum_{\substack{G \in \mathcal{R}_{g,n} \\ \text{non-trivalent}}}
			\frac{1}{\#\Aut(G)} \sum_{x \in P_{G}(L) \cap \ZZ_{+}^n/k} f(x)
		\Biggr|
		\le
		\Big( \sup_{\mathcal{M}_{g,n}^{\textup{comb}}(L)} |f| \Big)
		\sum_{\substack{G \in \mathcal{R}_{g,n} \\ \text{non-trivalent}}}
			\frac{\#\big( P_{G}(L) \cap \ZZ_{+}^n/k \bigr)}{\#\Aut(G)}.
	\]
	By dimensional reasons, the right-hand side is $O(k^{6g-7+2n})$ as $k \to +\infty$. Hence
	\[
		\lim_{\substack{k \to \infty \\ k \in d\ZZ_+}} \bigg(\frac{1}{k^{6g - 6 + 2n}}
			\sum_{\bm{G} \in \mathcal{M}_{g,n}^{\textup{comb},\ZZ/k}(L)}
			\frac{f(\bm{G})}{\#\Aut(\bm{G})} \bigg)
		=
		\lim_{\substack{k \to \infty \\ k \in d\ZZ_+}} \bigg(\frac{1}{k^{6g - 6 + 2n}}
			\sum_{\substack{G \in \mathcal{R}_{g,n} \\ \text{trivalent}}} \frac{1}{\#\Aut(G)} \sum_{x \in P_{G}(L) \cap \ZZ_{+}^n/k} f(x) \bigg).
	\]
	By Lemma~\ref{lem:index:2} and by definition of the Riemann integral, we have 
	\[
		\lim_{\substack{k \to \infty \\ k \in d\ZZ_+}} \bigg(\frac{1}{k^{6g - 6 + 2n}} \sum_{x \in P_{G}(L) \cap \ZZ_{+}^n/k} f(x) \bigg)
		=
		 \int_{P_{G}(L)} f \, \dd \ell_{n+1}\dots \dd \ell_{6g-6+3n}
	\]
	provided $f$ is continuous. Let $G$ be a trivalent ribbon graph, and let $\iota_L \colon P_{G}(L)\hookrightarrow \RR_+^{E_G}$ be the inclusion map. The cell $P_G(L)$ is naturally equipped with the measure $\mu_{\textup{Leb}}$ defined as $\iota_L^*(\prod_{e \in E_{G}} \dd\ell_{e})$. Using Lemma \ref{lem:dav:eyn}, we see that $\prod_{j=n+1}^{6g-6+3n} \dd \ell_j = 2 \mu_{\textup{Leb}}$. Besides, we know from Lemma~\ref{lem:Kon:Lebesgue} that $\mu_{{\rm Leb}} = 2^{2 - 2g - n}\mu_{\textup{K}}$. As a consequence, we find
	\[
		2^{3 - 2g - n} \int_{\mathcal{M}_{g,n}^{\textup{comb}}(L)} f \, \dd\mu_{\textup{K}}
		=
		\lim_{\substack{k \to \infty \\ k \in d\ZZ_+}} \bigg(\frac{1}{k^{6g - 6 + 2n}}
			\sum_{\substack{G \in \mathcal{R}_{g,n} \\ \text{trivalent}}} \frac{1}{\#\Aut(G)} \sum_{x \in P_{G}(L) \cap \ZZ_{+}^n/k} f(x) \bigg).
	\]
\end{proof}
\begin{proof}[Proof of Corollary~\ref{cor:limit:lattice:TR}]
	Let $\Xi_{\Sigma}^{(k)}$ be the GR amplitudes for the initial data
	\[
	\begin{split}
		A^{(k)}(L_1,L_2,L_3)	& = k^{-2b - c} \, A(k L_1,k L_2,k L_3), \\
		B^{(k)}(L_1,L_2,L_3)	& = k^{-b} \, B(k L_1,k L_2,k \ell), \\
		C^{(k)}(L_1,L_2,L_3)	& = k^{-c} \, C(k L_1,k \ell,k \ell'), \\ 
		D^{(k)}_{T}(\GG)		& = k^{-b} \, D_{T}(k \GG).
	\end{split} 
	\]
	Tracking the powers of $k^{-1}$ in the GR sum \eqref{eqn:comb:GRdef}, one can show by induction that
	\[
		\Xi_{\Sigma}(k \GG) = k^{(g - 1)(b + c) + bn} \, \Xi_{\Sigma}^{(k)}(\GG)
	\]
	and thus
	\begin{equation}\label{eqn:before:limit:lattice}
		\frac{N\Xi_{g,n}(kL_1,\ldots,kL_n)}{k^{6g - 6 + 2n}\cdot k^{(g - 1)(b + c) + bn}}
		= 
		\frac{1}{k^{6g - 6 + 2n}}
		\Bigg(
			\sum_{\bm{G} \in \mathcal{M}_{g,n}^{\textup{comb},\ZZ/k}(L)} \frac{\Xi_{\Sigma}^{(k)}(\bm{G})}{\#\Aut(\bm{G})}
		\Bigg).
	\end{equation}
	The vanishing conditions \eqref{eqn:vanishing:small:pop} have two consequences for us.
	\begin{itemize}
		\item[(i)]
		The amplitude $\Xi_{\Sigma}^{(k)}(\GG)$ is given by a finite sum of products with $2g - 2 + n$ factors that can be either $A^{(k)},B^{(k)},C^{(k)},D^{(k)}$.
		\item[(ii)]
		If we fix $L \in \RR_{+}^n$, for any $\GG \in \mathcal{T}_{\Sigma}^{\textup{comb}}(L)$, the factors of $A^{(k)},B^{(k)},C^{(k)}$ are evaluated on triples of lengths that are smaller than $M = \sum_{i = 1}^n L_i$, and $D$ is evaluated on elements of $\bigcup_{\ell \in (0,M]}\mathcal{T}_{T}^{\textup{comb}}(\ell)$.
	\end{itemize}
	Hence $\Xi_{\Sigma}^{(k)}$ converges uniformly to $\hat{\Xi}_{\Sigma}$ on $\mathcal{T}_{\Sigma}^{\textup{comb}}(L)$. We can then replace $\Xi_{\Sigma}^{(k)}$ with $\hat{\Xi}_{\Sigma}$ in \eqref{eqn:before:limit:lattice} up to an error that tends to $0$ when $k \to \infty$, and we conclude by using Proposition~\ref{prop:limit:lattice:TR} for $\hat{\Xi}_{\Sigma}$.
\end{proof} 

\subsection{Remark: inducing \texorpdfstring{$D$}{D} from \texorpdfstring{$C$}{C}}
\label{subsec:D:from:C}

There is a natural way to complete $(A,B,C)$ into an initial data $(A,B,C,D)$, satisfying all the assumptions that we may desire to impose.

\begin{lem}\label{lem:D:from:C}
	If we are only given $(A,B,C)$ satisfying the conditions in Definition~\ref{defn:adm:hyp}, the series
	\begin{equation}\label{eqn:D:choice}
		D_T(\sigma)
		=
		\sum_{\gamma \in S_{T}^{\circ}}
			C\big(\ell_{\sigma}(\de T),\ell_{\sigma}(\gamma),\ell_{\sigma}(\gamma)\big)
	\end{equation}
	converges absolutely on any compact of $\mathcal{T}_{T}$ to a $\Mod_{T}$-invariant function, and $(A,B,C,D)$ are admissible initial data. Furthermore, if $(A,B,C)$ satisfy the conditions in Definition~\ref{defn:strong:adm}, where in the bound for $C$ one assumes $0 \leq \eta < 1$, then $(A,B,C,D)$ are strongly admissible and
	\begin{equation}
		\int_{\mathcal{M}_{1,1}(L)} D\, \dd\mu_{\textup{WP}}
		=
		\int_{\mathbb{R}_{+}} \,C(L,\ell,\ell)  \ell \dd \ell.
	\end{equation}
	\hfill $\blacksquare$
\end{lem}

The same lemma holds in the combinatorial setting -- replacing hyperbolic with combinatorial lengths and $\mu_{\textup{WP}}$ with $\mu_{\textup{K}}$. In the last statement, the stronger condition $\eta < 1$ for $C$ (instead of $\eta < 2$) guarantees that $\ell \mapsto \ell \, C(L,\ell,\ell)$ is integrable near $0$. We omit the proof of Lemma~\ref{lem:D:from:C} and its combinatorial analog, as it is very similar to arguments already used in Theorems~\ref{thm:hyp:GR}, \ref{thm:comb:GR}, \ref{thm:hyp:GR:TR} and \ref{thm:comb:GR:TR}.

\medskip

When we say that $(A,B,C)$ are admissible initial data, we implicitly assume that they should be completed by the choice \eqref{eqn:D:choice} of $D_T$. With this convention, notice that (strong) admissibility of $(A,B,C)$ is the same condition in the hyperbolic and in the combinatorial setting.

\begin{rem}
	When $C(L_1,\ell,\ell')$ vanishes for $L_1 < \ell + \ell'$, the sum \eqref{eqn:D:choice} is finite. It is thus well-defined without  admissibility conditions and satisfies
	\begin{equation}
		\sum_{\bm{G} \in \mathcal{M}_{1,1}^{\textup{comb},\ZZ}(L)} \frac{D(\bm{G})}{\#\Aut(\bm{G})}
		=
		\sum_{\ell = 1}^{\lfloor L/2 \rfloor} \ell\,C(L,\ell,\ell)
	\end{equation}
	This is in particular relevant for the discrete setting, \textit{e.g.} for Theorem~\ref{thm:disc:TR}.
\end{rem}

\subsection{Mirzakhani--McShane identities}
\label{subsec:Mirz:McS}

\subsubsection{In the hyperbolic case}

Mirzakhani proved in \cite[Theorem 1.3]{Mir07simple} a generalisation for all bordered surfaces of the identity discovered by McShane \cite{McS98} for the punctured torus, and used it to prove a topological recursion in the form of \eqref{eqn:hyp:TR} for the Weil--Petersson volumes. It can be reformulated by saying that the constant function $1$ can be obtained from geometric recursion. We state it below and we are going to use it later on.

\begin{thm} \cite{Mir07simple} \label{thm:Mirz:McS:hyp}
	Let $F(x) = 2\ln(1 + e^{x/2})$. The following initial data
	\begin{equation}\label{eqn:Mirz:init:data}
	\begin{split}
		A^{\textup{M}}(L_1,L_2,L_3) & = 1, \\
		B^{\textup{M}}(L_1,L_2,\ell) & = \frac{1}{2L_1}\big(F(L_1 + L_2 - \ell) + F(L_1 - L_2 - \ell) - F(-L_1 + L_2 - \ell) - F(-L_1 - L_2 - \ell)\big), \\
		C^{\textup{M}}(L_1,\ell,\ell') & = \frac{1}{L_1}\big(F(L_1 - \ell - \ell') - F(-L_1 - \ell - \ell')\big),
	\end{split}
	\end{equation}
	are admissible, and lead by geometric recursion to $\Omega_{\Sigma}^{\textup{M}}(\sigma) = 1$ for any $\Sigma$ and $\sigma \in \mathcal{T}_{\Sigma}$. In other words, for any connected bordered surface $\Sigma$ with $\chi_{\Sigma} < - 1$ and any $\sigma \in \mathcal{T}_{\Sigma}$, we have
	\begin{equation}\label{eqn:Mirz:McS:hyp}
		1
		=
		\sum_{m = 2}^n \sum_{[P] \in \mathcal{P}_{\Sigma}^{m}}
			B^{\textup{M}}(\vec{\ell}_{\sigma}(\de P))
		+
		\frac{1}{2} \sum_{[P] \in \mathcal{P}_{\Sigma}^{\varnothing}}
			C^{\textup{M}}(\vec{\ell}_{\sigma}(\de P)),
	\end{equation}
	and for a torus $T$ with one boundary component and any $\sigma \in \mathcal{T}_{T}$, we have
	\begin{equation}
		1 = \sum_{\gamma \in S_{T}^{\circ}} C^{\textup{M}}\big(\ell_{\sigma}(\de T),\ell_{\sigma}(\gamma),\ell_{\sigma}(\gamma)\big).
	\end{equation}
	\hfill $\blacksquare$
\end{thm} 

\subsubsection{In the combinatorial setting}

We are now going to prove a recursion for the constant function $1$ on the combinatorial Teichm\"uller space. The initial data are the functions $(A^{\textup{K}},B^{\textup{K}},C^{\textup{K}})$ 
\begin{equation}\label{eqn:Kon:init:data}
\begin{split}
	A^{\textup{K}}(L_1,L_2,\ell)	&= 1, \\
	B^{\textup{K}}(L_1,L_2,\ell)	&= \frac{1}{2L_1}\Bigl(
		[L_1 - L_2 - \ell]_+ - [ -L_1 + L_2 - \ell]_+ + [L_1 + L_2 - \ell]_+
	\Bigr), \\
	C^{\textup{K}}(L_1,\ell,\ell')	&= \frac{1}{L_1} [L_1 - \ell - \ell']_+.
\end{split}
\end{equation}
\emph{i.e.} the same as in the Theorem~\ref{thm:Mirz:McS:hyp} above, but for $F(x) = [x]_+ = \max\set{x,0}$. Notice that $B^{\textup{K}}$ and $C^{\textup{K}}$ were already introduced in \eqref{eqn:Kon:BC} while discussing homotopy classes of embedded pairs of pants. The proof adapts to the combinatorial setting the strategy used by Mirzakhani to obtain Theorem~\ref{thm:Mirz:McS:hyp}.

\begin{thm} \label{thm:Mirz:McS:comb}
	For any connected bordered surface $\Sigma$ such that $\chi_{\Sigma} < - 1$ and any $\GG \in \mathcal{T}_{\Sigma}^{\textup{comb}}$, we have
	\begin{equation}\label{eqn:Mirz:McS:comb}
		1
		=
		\sum_{m = 2}^n \sum_{[P] \in \mathcal{P}_{\Sigma}^{m}}
			B^{\textup{K}}(\vec{\ell}_{\GG}(\de P))
		+
		\frac{1}{2} \sum_{[P] \in \mathcal{P}_{\Sigma}^{\varnothing}}
			C^{\textup{K}}(\vec{\ell}_{\GG}(\de P)),
	\end{equation}
	and for a torus $T$ with one boundary component and any $\GG \in \mathcal{T}_{T}^{\textup{comb}}$, we have
	\begin{equation}\label{eqn:Mirz:McS:comb:torus}
		1 = \sum_{\gamma \in S_{T}^{\circ}} C^{\textup{K}}\big(\ell_{\GG}(\de T),\ell_{\GG}(\gamma),\ell_{\GG}(\gamma)\big).
	\end{equation}
\end{thm}

\begin{proof}
	Consider the case of connected $\Sigma$ with $\chi_{\Sigma} < - 1$. For simplicity, set $X_{P}^{\textup{K}}$ equal to $B^{\textup{K}}$ or $\frac{1}{2}C^{\textup{K}}$ depending on the type of $[P] \in \mathcal{P}_{\Sigma}$. The basic idea to prove such an identity is to write $\ell_{\GG}(\de_{1}\Sigma)$ as a sum of lengths of the edges around $\de_1 \Sigma$. Recall that in the proof of Theorem~\ref{thm:length:functional}, we introduced a map\footnote{In the proof of Theorem~\ref{thm:length:functional}, the map was assigning to a point in $\mathcal{T}_{\Sigma}^{\textup{comb}}$ a functional on $\mathcal{A}_{\Sigma}^{\textup{all}}$. Here we are restricting the functionals to $\mathcal{A}_{\Sigma} \subset \mathcal{A}_{\Sigma}^{\textup{all}}$, \emph{i.e.} to arcs with initial point in $\de_1\Sigma$, where we recall the convention $\mathcal{A}_{\Sigma} = \mathcal{A}_{\Sigma,1}$ adopted at the beginning of this section.} $\mathcal{T}_{\Sigma}^{\textup{comb}} \to \RR_{\ge 0}^{\mathcal{A}_{\Sigma}}$ that assigns to a combinatorial structure $\GG$ the functional on $\mathcal{A}_{\Sigma}$ given by
	\[
		\alpha \longmapsto
			\begin{cases}
				\ell_{\GG}(\de_1 \Sigma) \Bigl(
					B^{\textup{K}}(\vec{\ell}_{\GG}(\de P_{\alpha}))
					-
					C^{\textup{K}}(\vec{\ell}_{\GG}(\de P_{\alpha}))
				\Bigr)
				& \text{if } \alpha \in \mathcal{A}_{\Sigma}^{m} \\[1ex]
				\frac{1}{2}\ell_{\GG}(\de_1 \Sigma) \, C^{\textup{K}}(\vec{\ell}_{\GG}(\de P_{\alpha}))
				& \text{if } \alpha \in \mathcal{A}_{\Sigma}^{\varnothing}
			\end{cases}
	\]
	where $[P_{\alpha}] = Q(\alpha)$ is the homotopy class of pair of pants determined by the arc $\alpha$. Fix once and for all $\GG \in \mathcal{T}_{\Sigma}^{\textup{comb}}$, and denote by $l_{\GG} \colon \mathcal{A}_{\Sigma} \to \RR_{\ge 0}$ the value of above map at $\GG$. This function has finite support, and the arcs with non-zero contribution are in bijection with the edges of $\GG$ around $\de_1\Sigma$. Furthermore, for $\alpha$ dual to an edge $e$ around $\de_1 \Sigma$, the value $l_{\GG}(\alpha)$ is the combinatorial length of $e$ (\emph{cf.} Lemma~\ref{lem:reconstruct:edge:lengths}). As a consequence,
	\[
		\ell_{\GG}(\de_{1}\Sigma)
		=
		\sum_{\alpha\in\mathcal{A}_{\Sigma}}l_{\GG}(\alpha).
	\]
	Now from Remark~\ref{rem:Q:not:inj}, we know that the map $Q \colon \mathcal{A}_{\Sigma} \to \mathcal{P}_{\Sigma}$ is not injective, but has finite fibers. More precisely, we have the following situation.
	\begin{itemize}
		\item If $[P] \in \mathcal{P}_{\Sigma}^{\varnothing}$, then $Q^{-1}([P])$ consists of a single arc $\alpha_0 \in \mathcal{A}_{\Sigma}^{\varnothing}$. Then
		\[
			\sum_{\alpha \in Q^{-1}([P])} l_{\GG}(\alpha)
			=
			l_{\GG}(\alpha_0)
			=
			\frac{1}{2}\ell_{\GG}(\de_{1}\Sigma)\,C^{\textup{K}}(\vec{\ell}_{\GG}(\de P)).
		\]
		\item If $[P] \in \mathcal{P}_{\Sigma}^{m}$, then $Q^{-1}([P])$ consists of three arcs: $\alpha_0 \in \mathcal{A}_{\Sigma}^{m}$, $\alpha'$ and its inverse $-\alpha'$ in $\mathcal{A}_{\Sigma}^{\varnothing}$. Then
		\[
		\begin{split}
			\sum_{\alpha\in Q^{-1}([P])} l_{\GG}(\alpha)
			& =
				l_{\GG}(\alpha_0) + l_{\GG}(\alpha') + l_{\GG}(-\alpha') \\
			& = 
				\ell_{\GG}(\de_{1}\Sigma) \Big(
						B^{\textup{K}}(\vec{\ell}_{\GG}(\de P)) - C^{\textup{K}}(\vec{\ell}_{\GG}(\de P)) + \tfrac{1}{2}C^{\textup{K}}(\vec{\ell}_{\GG}(\de P)) + \tfrac{1}{2}C^{\textup{K}}(\vec{\ell}_{\GG}(\de P))
					\Big) \\
			& =
				\ell_{\GG}(\de_{1}\Sigma) \, B^{\textup{K}}(\vec{\ell}_{\GG}(\de P)).
		\end{split}
		\]
	\end{itemize}
	Finally, as $Q$ is surjective
	\[
		\ell_{\GG}(\de_{1} \Sigma)
		=
		\sum_{[P] \in \mathcal{P}_{\Sigma}} \sum_{\alpha \in Q^{-1}([P])} l_{\GG}(\alpha)
		=
		\ell_{\GG}(\de_{1}\Sigma)\sum_{[P] \in \mathcal{P}_{\Sigma}} X_P^{\textup{K}}(\vec{\ell}_{\GG}(P)).
	\]
	Dividing by $\ell_{\GG}(\de_{1}\Sigma)$ concludes the proof of the identity \eqref{eqn:Mirz:McS:comb}. The case of the torus with one boundary component can be handled in a similar way and leads to \eqref{eqn:Mirz:McS:comb:torus}.
\end{proof}

A notable difference with the original Mirzakhani--McShane identity is that in \eqref{eqn:Mirz:McS:comb} there is a finite number of non-zero terms; this reflects the much simpler dynamics of geodesics in combinatorial surfaces compared to the hyperbolic ones. As in the hyperbolic case, the identity can be reformulated in terms of geometric recursion.

\begin{cor}\label{GR:comb:WK}
	The initial data \eqref{eqn:Kon:init:data} are admissible, and lead by geometric recursion to $\Xi_{\Sigma}^{\textup{K}}(\GG) = 1$ for any $\Sigma$ and $\GG \in \mathcal{T}_{\Sigma}^{\textup{comb}}$.
\end{cor}

\begin{proof}
	The only thing left to check is the admissibility condition for $(A^{\textup{K}},B^{\textup{K}},C^{\textup{K}},D^{\textup{K}})$, which follows from the fact that the functions $B^{\textup{K}}$ and $C^{\textup{K}}$ are supported on small pairs of pants.
\end{proof}

We can now apply the integration results of Section~\ref{subsec:(disc):integr:TR} to the combinatorial Mirzakhani--McShane identity.
The continuous integration computes the Kontsevich volumes and Theorem~\ref{thm:hyp:GR:TR} establishes again Witten's conjecture \cite{Wit91}, which was first proved by the combination of \cite{Kon92} and \cite{DVV91}: the generating function of $\psi$-class intersection numbers satisfies Virasoro constraints/topological recursion. The discrete integration calculates the number of integral points in $\mathcal{M}_{g,n}^{\textup{comb}}$ and Theorem~\ref{thm:disc:TR} together with Corollary~\ref{cor:limit:lattice:TR} with $b = c = 0$, gives a new proof of Norbury's result \cite{Nor10}, \emph{i.e.} the discrete topological recursion for these counts and their connection with Kontsevich volumes.

\begin{cor}[\emph{Witten's conjecture/Kontsevich theorem}] The Kontsevich volumes $V_{g,n}^{\textup{K}}(L)$ equal $V\Xi_{g,n}^{\textup{K}}(L)$ and satisfy the topological recursion of Theorem~\eqref{thm:comb:GR:TR} with initial data \eqref{eqn:Kon:init:data}.
\end{cor}

\begin{cor}[\emph{Norbury's theorem}]
	The lattice point counting functions
	\begin{equation}
		N_{g,n}(L) = \sum_{\bm{G} \in \mathcal{M}_{g,n}^{\textup{comb},\ZZ}(L)} \frac{1}{\#\Aut(\bm{G})}
	\end{equation}
	equal $N\Xi^{\textup{K}}_{g,n}(L)$ and satisfy the discrete topological recursion \eqref{eqn:discr:TR} with initial data \eqref{eqn:Kon:init:data}. Further,
	\begin{equation}
		\lim_{\substack{k \to \infty \\ k \in d\ZZ_+}}
		\frac{N_{g,n}(k L)}{k^{6g - 6 + 2n}}
		=
		\begin{cases}
			2^{3-2g -n} \, V_{g,n}^{\textup{K}}(L)
			& \text{if $\sum_{i = 1}^n d\cdot L_i$ is even,} \\
			0
			& \text{otherwise.}
		\end{cases}
	\end{equation}
	\hfill $\blacksquare$
\end{cor}

The recursive formula for Kontsevich volumes in this integral form firstly appeared in \cite{BCSW12}, where it is derived by constructing a partition of unity similar to the combinatorial Mirzakhani--McShane identity and exploiting the local torus symmetries on the combinatorial moduli space, whose symplectic quotients are also combinatorial moduli space of higher Euler characteristics. Compared to \cite{BCSW12}, the  new element of the proof that we propose is to make it a \emph{complete} analogue to Mirzakhani's proof of the recursion for Weil--Petersson volumes, by means of Theorem~\ref{thm:Mirz:McS:comb}. Our perspective stresses that, at a general level, recursions between volumes (or more generally, between integrals over the moduli spaces) arise from finer recursions that hold at the geometric level (here between functions on Teichm\"uller spaces).

\subsection{Combinatorial length statistics of multicurves}
\label{subsec:twisted:GR}

Following \cite[Theorem~10.1]{ABO17} in the hyperbolic world, we can generalise the combinatorial Mirzakhani--McShane identity \eqref{eqn:Mirz:McS:comb} to obtain statistics of combinatorial lengths of primitive multicurves via GR. 

\begin{thm} \label{thm:GR:statistics}
	Let $(A,B,C,D)$ be admissible initial data and denote by $\Xi_{\Sigma}$ the associated GR amplitudes. Let $f \colon \mathbb{R}_{+}\to\mathbb{R}$ be a measurable function such that for any $\epsilon > 0$ and $s \geq 0$, there exists $M_{s,\epsilon}$ such that
	\begin{equation}
		\sup_{\ell \ge \epsilon} \; |f(\ell)| \, \ell^s \leq M_{s,\epsilon}.
	\end{equation}
	Then, the following initial data are admissible:
	\begin{equation}\label{eqn:twisted:initial:data}
	\begin{aligned}
		A[f](L_{1},L_{2},L_{3}) & =  A(L_{1},L_{2},L_{3}), \\
		B[f](L_1,L_2,\ell) & = B(L_1,L_2,\ell) + A(L_1,L_2,\ell) \, f(\ell), \\
		C[f](L_1,\ell,\ell') & = C(L_1,\ell,\ell')
			+ B(L_1,\ell,\ell') \, f(\ell) + B(L_1,\ell',\ell) \, f(\ell')
			+ A(L_1,\ell,\ell') \, f(\ell) \, f(\ell'), \\
		D_T[f](\GG) & = D_T(\GG)
			+ \sum_{\gamma\in S_{T}^{\circ}} A\bigl( \ell_{\GG}(\de T),\ell_{\GG}(\gamma),\ell_{\GG}(\gamma)\bigr) \,
			f(\ell_{\GG}(\gamma)).
	\end{aligned}
	\end{equation}
	Denote by $\Xi_{\Sigma}[f]$ the corresponding GR amplitudes. If for all $\Sigma$, $\Xi_{\Sigma}$ is invariant under all braidings of boundary components of $\Sigma$, we have
	\begin{equation}\label{eqn:twisted:stat}
		\Xi_{\Sigma}[f](\GG)
		=
		\sum_{c \in M_{\Sigma}'} \Xi_{\Sigma_c}(\GG|_{\Sigma_c})\prod_{\gamma \in \pi_{0}(c)}f(\ell_{\GG}(\gamma)),
	\end{equation}
	where $\Sigma_c$ is the bordered surface obtained by cutting $\Sigma$ along $c$ (the choice of the first boundary component is irrelevant due to the assumed invariance).
\end{thm}

\begin{rem}\label{rem:pure:statistics:multicurves}
	A useful version of the above result can be stated for combinatorial length statistics of multicurves (not only primitive ones). Consider $F \colon \mathbb{R}_{>0} \to (-1,1)$ a measurable function such that
	\begin{equation}
		f(x) = \sum_{k \ge 1} F(x)^k = \frac{F(x)}{1 - F(x)}
	\end{equation}
	satisfies the conditions of Theorem~\ref{thm:GR:statistics}. Then the GR amplitudes $\Xi_{\Sigma}[f]$ associated to the initial data \eqref{eqn:twisted:initial:data} are given by the combinatorial length statistic of multicurves weighted by $F$:
	\begin{equation}
		\Xi_{\Sigma}[f](\GG) = \sum_{c \in M_{\Sigma}} \Xi_{\Sigma_c}(\GG|_{\Sigma_c}) \prod_{\gamma \in \pi_0(c)} F(\ell_{\GG}(\gamma)).	
	\end{equation}
\end{rem}

\begin{proof}[Sketch of the proof of Theorem~\ref{thm:GR:statistics}.]
	The proof repeats the one of \cite[Theorem~10.1]{ABO17}, with combinatorial leng\-ths instead of hyperbolic ones. We only sketch the induction step here. Notice that, as all series are absolutely convergent, we can apply Fubini's theorem and interchange the summations:
	\[
	\begin{aligned}
		\sum_{c \in M_{\Sigma}'}
			\Xi_{\Sigma_c}(\GG|_{\Sigma_{c}}) \prod_{\gamma \in \pi_{0}(c)} f(\ell_{\GG}(\gamma))
		& =
		\sum_{c \in M_{\Sigma}'}
			\sum_{[P] \in \mathcal{P}_{\Sigma_{c}}}
				X_{P}(\GG) \, \Xi_{\Sigma_{c} - P}(\GG|_{\Sigma_{c} - P})
			\prod_{\gamma \in \pi_{0}(c)} f(\ell_{\GG}(\gamma)) \\
		& =
		\sum_{[P] \in \mathcal{P}_{\Sigma}}
			X_{P}[f](\GG)
			\Biggl(
				\sum_{\gamma\in M_{\Sigma - P}'}
				\Xi_{(\Sigma - P)_{\gamma}}(\GG|_{(\Sigma - P)_{\gamma}})
				\prod_{c\in\pi_{0}(\gamma)}f(\ell_{\GG}(\gamma_{c}))
			\Biggr) \\
		& =
		\sum_{[P] \in \mathcal{P}_{\Sigma}}
			X_{P}[f](\GG) \, \Xi_{\Sigma - P}[f](\GG|_{\Sigma - P}),
	\end{aligned}
	\]
	where $X_{P}$ is either $B$ or $\tfrac{1}{2} C$ depending on the type of $[P]$. The second to last equality follows from a case discussion: to interchange the summands we must also sum over the possible ways in which $\de P$ and $c$ can have homotopic components, weighted by their contributions. For example, for $[P] \in \mathcal{P}_{\Sigma}^{\varnothing}$ with $\de P$ and $c$ sharing one component, say $\de_2 P = \gamma_0$, then $[P] \in \mathcal{P}_{\Sigma_{c}}^{\gamma_0}$. This means that we get a contribution of
	\[
		B\bigl(
			L_{1},\ell_{\GG}(\gamma_0),\ell_{\GG}(\de_{3}P)
		\bigr)
		\,
		f\bigl(\ell_{\GG}(\gamma_0)\bigr)
		\,
		\Xi_{(\Sigma - P)_{c}}\bigl(\GG|_{(\Sigma - P)_{c}}\bigr)
		\prod_{\gamma \in \pi_{0}(c - \gamma_{0})}f\bigl(\ell_{\GG}(\gamma)\bigr),
	\]
	which matches one of the terms in the expression for $C[f]$. We refer to \cite{ABO17} for the complete proof.
\end{proof}

We can then integrate this identity over the combinatorial moduli space with respect to $\mu_{\textup{K}}$, and the result can be calculated in two ways: by the topological recursion, and by direct integration. To express the later, we introduce the set of stable graphs
\begin{equation}
	\mathcal{G}_{g,n} = M'_{\Sigma}/\Mod_{\Sigma}^{\de},
\end{equation}
where $\Sigma$ is a bordered surface of type $(g,n)$. We refer to \cite{ABCDGLW19} for a full definition and for the notation, but let us say that a stable graph $\Gamma$ encodes the class of a primitive multicurve $c$ in the following way: vertices $v \in V_{\Gamma}$ correspond to connected components of $\Sigma_c$ and the genera of the components is recorded in a function $h \colon V_{\Gamma} \to \NN$ which is part of the data of $\Gamma$; edges $e \in E_{\Gamma}$ correspond to the components of $c$, and the set of edges incident to $v$ is denoted by $E(v)$; leaves $\lambda \in \Lambda(v)$ correspond to boundary components of the surface attached to $v \in V_{\Gamma}$; the valency of the vertex $v$ is denoted by $k(v)$. The automorphism group of $\Gamma$ is identified with that of the multicurve $c$.

\begin{cor}\label{cor:sum:stable:graphs}
	Assume that $(A,B,C,D)$ are strongly admissible and $f$ is a measurable function for which there exists $\eta' \in [0,2)$ such that
	\begin{equation}
		\sup_{\ell > 0} \; |f(\ell)| \, \ell^{\eta'} < +\infty.
	\end{equation}
	Then $(A[f],B[f],C[f],D[f])$ are strongly admissible and $V\Xi_{g,n}[f](L) = \int_{\mathcal{M}^{\textup{comb}}_{g,n}(L)} \Xi_{g,n}[f]\,\dd\mu_{\textup{K}}$ satisfy the topological recursion \eqref{eqn:comb:TR} with this initial data. Besides, we have
	\begin{equation}\label{eqn:sum:stable:graphs}
		V\Xi_{g,n}[f](L_1,\ldots,L_n)
		=
		\sum_{\Gamma \in \mathcal{G}_{g,n}} \frac{1}{\#\Aut(\Gamma)}
			\int_{\mathbb{R}_{+}^{E_{\Gamma}}}
			\prod_{v \in V_{\Gamma}} V\Xi_{h(v),k(v)}\big((\ell_e)_{e \in E(v)},(L_{\lambda})_{\lambda \in \Lambda(v)}\big)
			\prod_{e \in E_{\Gamma}} \ell_e \, f(\ell_e) \, \dd\ell_e.
	\end{equation}
\end{cor}

\begin{proof}
	The topological recursion follows from Theorem~\ref{thm:comb:GR:TR} and from the second part from the integration formula of Proposition~\ref{prop:integration}.
\end{proof}

We can also obtain an analogous result for the discrete integration, substituting $V\Xi_{h(v),k(v)}$ with $N\Xi_{h(v),k(v)}$ and the integral over $\RR_{+}^{E_{\Gamma}}$ with the sum over edge decorations of the form $\ell \colon E_{\Gamma} \to \ZZ_{+}$. However, we cannot directly apply Theorem~\ref{thm:disc:TR}, as the twisted initial data are in general not supported on small pairs of pants. Instead, we introduce an (\textit{a priori}) different lattice count, for $L \in \ZZ_{+}^n$
\begin{equation}\label{eqn:tilde:twisted:count}
		N\Xi_{g,n}[f](L)
		:=
			\sum_{\Gamma \in \mathcal{G}_{g,n}} \frac{1}{\#\Aut(\Gamma)}
			\sum_{\ell : E_{\Gamma} \to \ZZ_{+}}
				\prod_{v \in V_{\Gamma}} N\Xi_{h(v),k(v)}\big((\ell_e)_{e \in E(v)},(L_{\lambda})_{\lambda \in \Lambda(v)}\big)
				\prod_{e \in E_{\Gamma}} \ell_e \, f(\ell_e).
	\end{equation}
We omit the proof of the next proposition: it follows the same scheme as the proof of Theorem~\ref{thm:GR:statistics}, in the simpler situation where multicurves are replaced by their mapping class group orbits (stable graphs), and with integrals replaced by discrete sums. The key principle is that topological recursion is preserved under the twisting operation.

\begin{prop}
	Let $(A,B,C,D)$ and $N\Xi_{g,n}$ as in Theorem~\ref{thm:disc:TR}. Let $f \colon \ZZ_{+} \to \mathbb{R}$ be a continuous function such that for any $s > 0$, 
	\[
		\sup_{\ell \in \ZZ_{+}} |f(\ell)|\,\ell^{s} < +\infty.
	\]
		Then $N\Xi_{g,n}[f](L)$ is finite and, for $2g - 2 + n \geq 2$, it is calculated by the discrete topological recursion formula 
	\begin{equation}\label{eqn:disc:twisted:TR}
	\begin{split}
		N\Xi_{g,n}[f](L_1,\ldots,L_n)
		= 
		\sum_{m=2}^n \sum_{\ell\geq 1} \ell \, B_{\ZZ}[f](L_1,L_m,\ell) \, N\Xi_{g,n-1}[f](\ell,L_2,\dots,\widehat{L_m},\dots,L_n)
		& \\
		+ \frac{1}{2} \sum_{\ell,\ell'\geq 1} \ell \ell' \, C_{\ZZ}[f](L_1,\ell,\ell') \Bigg(
			N\Xi_{g-1,n+1}[f](\ell,\ell',L_2,\dots,L_n)
		& \\
		+ \!\!\! \sum_{\substack{h + h' = g \\ J\sqcup J'=\{L_2\dots,L_n\}}} \!\!\! N\Xi_{h,1+\#J}[f](\ell,J) \, N\Xi_{h',1+\#J'}[f](\ell',J') \Bigg), &
	\end{split}
	\end{equation}
with conventions $N\Xi_{0,1}[f] = 0$ and $N\Xi_{0,2}[f] = 0$ and base cases
	\begin{equation}\label{eqn:disc:twisted:TR:init:cond}
	\begin{aligned}
		N\Xi_{0,3}[f](L_1,L_2,L_3)
		& = N\Xi_{0,3}(L_1,L_2,L_3)
		= A_{\ZZ}(L_1,L_2,L_3),
		\\
		N\Xi_{1,1}[f](L_1)
		& = N\Xi_{1,1}(L_1) + \frac{1}{2} \sum_{\ell \geq 1} \ell \,A_{\ZZ}(L_1,\ell,\ell)\,f(\ell)\,.
	\end{aligned}
	\end{equation}
	\hfill $\blacksquare$
\end{prop}

\subsection{Kontsevich amplitudes for hyperbolic GR}

On the one hand, by Corollary~\ref{GR:comb:WK}, the combinatorial GR amplitudes $\Xi_{\Sigma}^{\textup{K}}$ for the initial data $(A^{\textup{K}},B^{\textup{K}},C^{\textup{K}})$ coincide with the constant function $1$ on $\mathcal{T}_{\Sigma}^{\textup{comb}}$. On the other hand, we can consider the hyperbolic GR amplitudes $\Omega^{\textup{K}}_{\Sigma}$ associated to the same initial data. They are rather non-trivial functions of $\sigma \in \mathcal{T}_{\Sigma}$, but since the topological recursion formulae are the same in the combinatorial and hyperbolic setting, we have
\begin{equation}\label{eqn:Kon:poly}
	V^{\textup{K}}_{g,n}(L)
	=
	\int_{\mathcal{M}_{g,n}^{\textup{comb}}(L)} \Xi^{\textup{K}}_{g,n} \, \dd\mu_{\textup{K}}
	=
	\int_{\mathcal{M}_{g,n}(L)} \Omega^{\textup{K}}_{g,n} \, \dd\mu_{\textup{WP}}
	=
	\int_{\overline{\mathfrak{M}}_{g,n}} \exp\bigg(\sum_{i = 1}^n \frac{L_i^2}{2}\,\psi_i\bigg)
\end{equation}
for any $L = (L_1,\ldots,L_n)$. Here we propose a geometric interpretation of $\Omega^{\textup{K}}_{\Sigma}$ and give some of its basic properties.

\subsubsection{A geometric interpretation}

In the following we discuss the combinatorial analogue of the geometric reasoning
 at the core of Mirzakhani's proof of the Mirzakhani--McShane identities.
Consider a connected bordered surface $\Sigma$. Recall from \S~\ref{subsubsec:ordinary:combinatorial} the construction of the spine as a subset $\sp_{\sigma}(\Sigma) \subset \Sigma$ that depends on a hyperbolic marking $(X,f)$ representing $\sigma \in \mathcal{T}_{\Sigma}$. We also denote by $\sp_{\sigma}'(\Sigma) \subset \sp_{\sigma}(\Sigma)$ the complement of the set of vertices of the spine.

\medskip

For a given $\sigma \in \mathcal{T}_{\Sigma}$, we equip $\de_1\Sigma$ with the curvilinear measure induced by $\sigma$. 

\begin{lem}\label{lem:shot:geodesics:spine}
	For all but finitely many $x \in \de_1\Sigma$, the geodesic shot from $x$ orthogonally to $\de_1\Sigma$ intersects the spine for the first time at a point $s_{\Sigma}(x) \in \sp_{\sigma}'(\Sigma)$.
\end{lem}

\begin{proof}
	Cutting out the spine, we have a cylinder around each boundary component of the surface. Consider the one around $\de_1 \Sigma$ and take a geodesic that realises the distance between the two boundaries of the cylinder. Cutting along this, we obtain a hyperbolic polygon. As there are no hyperbolic triangles with two right angles, every geodesic shot orthogonally from $x \in \de_1 \Sigma$ must reach the boundary corresponding to the spine at a certain point $s_{\Sigma}(x)$, that is not a vertex for all but finitely many $x$ that are rib-ends on $\de_1 \Sigma$.
\end{proof}

Let $x \in \de_1\Sigma$ be such that $s_{\Sigma}(x)$ exists and is not a vertex of the spine. By definition of the spine and the vertices, there exists a unique second geodesic joining $s_{\Sigma}(x)$ to $\partial_i \Sigma$ for some $i \in \{1,\ldots,n\}$. The union of these two geodesics form a piecewise geodesic arc, denoted by $\gamma_{x}$. This arc determines a unique homotopy class of embedded pair of pants, and we denote by $P_{x}$ its representative that has geodesic boundaries in $\Sigma$. This pair of pants has a spine $\sp_{\sigma}(P_x)$ of its own, and we can ask whether $s_{\Sigma}(x)$ is still part of $\sp_{\sigma}(P_x)$. With these notations, we define the following process.

\begin{defn}
	Choose a random point $x \in \de_1 \Sigma$ uniformly on $\de_1\Sigma$ -- for the measure coming from the curvilinear abscissa induced by $\sigma$. If the geodesic shot from $x$ orthogonally to $\de_1\Sigma$ hits for the first time $\sp_{\sigma}(\Sigma)$ at a vertex, or if $s_{\Sigma}(x) \notin \sp_{\sigma}(P_x)$, quit the process. Otherwise, consider the bordered hyperbolic surface $\Sigma - P_x$. If it is empty, we have finished the process successfully; if not, we repeat it with $\Sigma - P_x$\footnote{The arcs give an order to the boundary components of $P_{x}$, and by a similar process used in \cite[Section 2.3]{ABO17} we can define the first boundary on $\Sigma - P_{x}$.}, each connected component of $\Sigma - P_{x}$ being treated independently. Denote by $\Pi_{\Sigma}(\sigma)$ the probability that the process ends successfully, \emph{i.e.} by giving a pants decomposition. It is clear that $\Pi_{\Sigma}(\sigma)$ only depends on the projection $[\sigma] \in \mathcal{M}_{\Sigma}$ and the process makes no reference to a marking.
\end{defn}

\begin{figure}[ht]
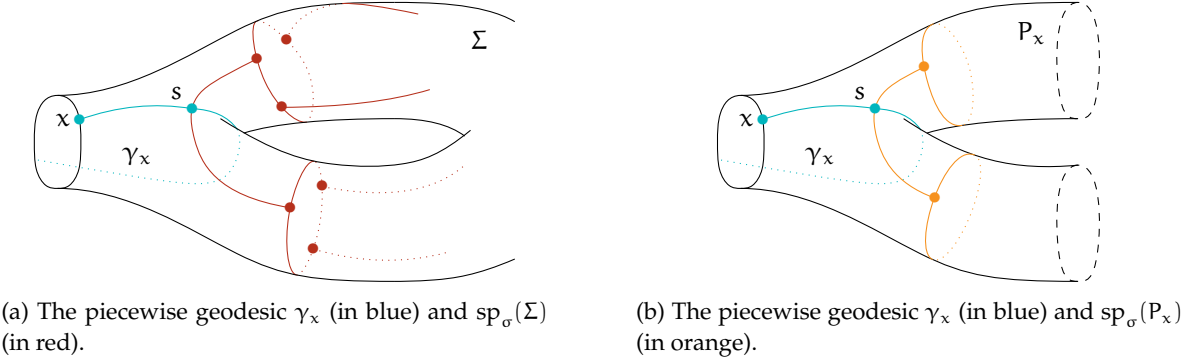

\centering
	\begin{subfigure}[t]{.44\textwidth}
	\centering

	\caption{The piecewise geodesic $\gamma_x$ (in blue) and $\sp_{\sigma}(\Sigma)$ (in red).}
	\label{fig:process:a}
	\end{subfigure}
	\hspace{1cm}
	\begin{subfigure}[t]{.44\textwidth}
	\centering
		%
	\caption{The piecewise geodesic $\gamma_x$ (in blue) and $\sp_{\sigma}(P_x)$ (in orange).}
	\label{fig:process:b}
	\end{subfigure}
	\caption{Example of a successful first step of the process.}
\end{figure}

\begin{prop}
	We have $\Pi_{\Sigma}(\sigma) = \Omega_{\Sigma}^{\textup{K}}(\sigma)$. In particular, $\Omega_{\Sigma}^{\textup{K}}(\sigma) \in [0,1]$, and if we consider the above process for a random hyperbolic surface of type $(g,n)$ with fixed boundary lengths $L \in \mathbb{R}_{+}^n$ (with respect to the Weil--Petersson measure on $\mathcal{M}_{\Sigma}(L)$), we have
	\begin{equation}
		\mathbb{E}\big[\Pi_{\Sigma}(\sigma) \; \big| \; \ell_{\sigma}(\partial \Sigma) = L\big] = \frac{V^{\textup{K}}_{g,n}(L)}{V^{\textup{WP}}_{g,n}(L)},
	\end{equation}
	where $V^{\textup{WP}}_{g,n}(L)$ are the Weil--Petersson volumes of $\mathcal{M}_{g,n}(L)$.
\end{prop}

\begin{proof}
	We prove the result by induction on $2g - 2 + n > 0$. The base case $(g,n) = (0,3)$ is trivial: for any pair of pants $P$, we have $\Pi_{P} \equiv 1$. Furthermore, the following argument can be adjusted to prove that for a one-holed torus $T$,
	\[
		\Pi_T(\sigma)
			=
			\sum_{\gamma \in S_T^{\circ}} C^{\textup{K}}\bigl(
				\ell_{\sigma}(\de T),\ell_{\sigma}(\gamma),\ell_{\sigma}(\gamma)
			\bigr)
			= D_T(\sigma).
	\]
	Consider $\Sigma$ of type $(g,n)$ and	suppose now by induction that the proposition holds for surfaces of Euler characteristic $\chi > -(2g-2+n)$. By definition of the process, we can write $\Pi_{\Sigma}(\sigma)$ as a sum over $\cup_{x\in\de_{1}\Sigma}P_{x}$ (where $P_{x} = \varnothing$ if $x$ does not define a pair of pants):
	\[
		\Pi_{\Sigma}(\sigma) = \sum_{[P] \in \cup_{x\in\de_{1}\Sigma}P_{x}} Y_P(\sigma) \, \Pi_{\Sigma - P}\bigl(\sigma|_{\Sigma - P}\bigr).
	\]
	Here, we identify $[P]$ with its representative $P \subseteq \Sigma$ having geodesic boundaries, and we set
	\begin{equation}\label{eqn:set:Y}
		Y_{P}(\sigma) = \frac{\ell_{\sigma}(\mathcal{Y}_{P}(\sigma))}{\ell_{\sigma}(\partial_1 \Sigma)},
		\qquad
		\mathcal{Y}_{P}(\sigma) = \Set{ x \in \de_1 \Sigma | s_{\Sigma}(x) \in \sp_{\sigma}(P)}.
	\end{equation}
	By induction hypothesis, we have $\Pi_{\Sigma - P}(\sigma|_{\Sigma - P}) = \Omega_{\Sigma - P}^{\textup{K}}(\sigma|_{\Sigma - P})$. So, we only need to show that $Y_P(\sigma) = X_P^{\textup{K}}(\vec{\ell}_{\sigma}(\partial P))$ and that the $[P]$'s outside $\cup_{x\in\de_{1}\Sigma}P_{x}$ do not contribute to the GR sum defining $\Omega^{\textup{K}}_{\Sigma}$.

	\smallskip

	For the latter, we observe that if $[P] \notin \cup_{x\in\de_{1}\Sigma}P_{x}$, then all points of $\sp_{\sigma}(P)$ incident to $\de_{1}\Sigma$ are equidistant to an internal boundary of $P$. Therefore we see that the spine of $\sigma|_P$ has the whole of the first boundary adjacent to an internal boundary, and therefore from Lemma~\ref{lem:geo:meaning:BC} we see that $X^{\textup{K}}_{P}(\vec{\ell}_{\sigma}(\partial P)) = 0$. So, we can indeed restrict the range of the GR sum defining $\Omega^{\textup{K}}_{\Sigma}$ to $[P] \in \cup_{x\in\de_{1}\Sigma}P_{x}$.

	\smallskip

	For the former, we first give a description of $X_P^\textup{K}$ similar to \eqref{eqn:set:Y}. Let $s_P(x)$ be the first intersection point (if it exists) between the geodesic shot from $x$ orthogonally to $\de_1\Sigma$ and the spine of $P$ -- considered as a subset of $\Sigma$. We introduce the sets
	\[
	\begin{aligned}
		\mathfrak{sp}_{\sigma}(P)
		& =
			\begin{cases}
				\Set{s \in \sp_{\sigma}'(P) | \text{$s$ is incident to $\de_m\Sigma$ and $\de_1 \Sigma$ or only to $\de_1\Sigma$}}
					& \text{if $[P] \in \mathcal{P}_{\Sigma}^{m}$} \\
				\Set{s \in \sp_{\sigma}'(P) | \text{$s$ is incident to $\de_1 \Sigma$ on both sides}}
				& \text{if $[P] \in \mathcal{P}_{\Sigma}^{\emptyset}$}
			\end{cases} \\
		\mathcal{X}_{P}(\sigma)
		& = \Set{ x \in \de_1 \Sigma | s_{P}(x) \in \mathfrak{sp}_{\sigma}(P)}
	\end{aligned} 
	\]
	where it is understood that if $s_{P}(x)$ does not exist, $x$ is not in $\mathcal{X}_{P}(\sigma)$. The properties of $B^{\textup{K}}$ and $C^{\textup{K}}$ stressed in Lemma~\ref{lem:geo:meaning:BC} show that
	\[
		X_P^{\textup{K}}(\vec{\ell}_{\sigma}(\partial P)) = \frac{\ell_{\sigma}(\mathcal{X}_{P}(\sigma))}{\ell_{\sigma}(\partial_1\Sigma)}.
	\]
	It remains to justify that $\mathcal{X}_{P}(\sigma) = \mathcal{Y}_{P}(\sigma)$.
	\begin{description}
		\item[($\subseteq$)]
			Consider $x \in \mathcal{X}_{P}(\sigma)$.  From the definition of $\mathfrak{sp}_{\sigma}(P)$, it is clear that $\mathfrak{sp}_{\sigma}(P) \subseteq \sp_{\sigma}'(\Sigma)$ and $s_{P}(x)$ is equidistant to $\de\Sigma$ in exactly two ways. In particular,
			\[
				{\rm dist}_{\sigma}\big(s_{P}(x),\de_1 \Sigma\big) = {\rm dist}_{\sigma}\big(s_P(x),\de_i\Sigma\big) < {\rm dist}_{\sigma}\big(s_P(x),\de\Sigma-(\de_1 \Sigma\cup\de_i\Sigma)\big)
			\]
			for some $i \in \{1,\ldots,n\}$. Thus, $s_{\Sigma}(x) = s_{P}(x)$ and clearly $s_{P}(x) \in \sp_{\sigma}(P)$.
		\item[($\supseteq$)]
			Consider $x \in \mathcal{Y}_{P}(\sigma)$.  Then we have $s_{P}(x) \in \sp_{\sigma}(\Sigma)$, so that $s_{P}(x) = s_{\Sigma}(x)$ and clearly $s_{P}(x) \in \mathfrak{sp}_{\sigma}(P)$.
	\end{description}
\end{proof}

\subsubsection{Properties of the Kontsevich amplitudes}

The geometric interpretation of $\Omega_{\Sigma}^{\textup{K}}$ already shows the non-trivial fact that such functions take values in $[0,1] \subset \RR$. In the remaining part of this section we are going to show two other properties of $\Omega_{\Sigma}^{\textup{K}}$.

\paragraph{Non-invariance under all braidings.}
\label{par:Kont:not:sym}

From their definition, the GR amplitudes are \emph{a priori} only invariant under mapping classes that preserve the first boundary; the invariance under braidings of $\partial_1\Sigma$ with some $\partial_m\Sigma$ for $m \neq 1$ is not guaranteed. This full invariance turns out to hold for the hyperbolic GR amplitudes $\Omega^{{\rm M}}_{\Sigma}$ and the combinatorial GR amplitudes $\Xi^{\textup{K}}_{\Sigma}$, for the obvious reason that they are identically $1$. One could wonder if the full invariance also holds for $\Omega^{\textup{K}}_{\Sigma}$, but we show this is already not the case for four-holed spheres.

\begin{prop}
	Let $X$ be a four-holed sphere, and take $\gamma \in S_{X}^{\circ}$ separating $\de_{1}X$ and $\de_{2}X$ from $\de_{3}X$ and $\de_{4}X$. Consider $\rho = [\phi \colon X \to X] \in \Mod_{X}$ the involution that fixes $\gamma$ and such that $\phi(\de_{1}X) = \de_{3}X$ and $\phi (\de_{2}X) = \de_{4}X$. Then $\rho.\Omega_{X}^{\textup{K}} \neq \Omega_{X}^{\textup{K}}$. 
\end{prop}

\begin{proof}
	The curve $\gamma$ together with an arc from $\de_1X$ to $\de_3X$, determine a seamed pants decomposition of $X$, and we denote by $(L,\ell,\tau) \in \RR_{+}^4 \times \RR_{+} \times \RR$ the corresponding combinatorial Fenchel--Nielsen coordinates. If we choose $\sigma \in \mathcal{T}_{X}$ such that $\ell_{\sigma}(\gamma)$ is small enough, then by the collar lemma we can make the length of any simple closed geodesic intersecting $\gamma$ greater than $\max\{L_{1}+L_{2},L_{3}+L_{4}\}$. As $B^{\textup{K}}(L_{i},L_{j},\ell)$ vanishes for $\ell \geq L_{i}+L_{j}$, we have
	\[
		\Omega^{\textup{K}}_{X}(L_{1},L_{2},L_{3},L_{4},\ell,\tau)
		=
		B^{\textup{K}}(L_{1},L_{2},\ell)
	\]
	and
	\[
		\rho.\Omega^{\textup{K}}_{X}(L_{1},L_{2},L_{3},L_{4},\ell,\tau)
		=
		\Omega^{\textup{K}}_{X}(L_{3},L_{4},L_{1},L_{2},\ell,\tau)
		=
		B^{\textup{K}}(L_{3},L_{4},\ell).
	\]
	Choosing $(L_1,\dots,L_4) \in \RR_{+}^4$ such that $B^{\textup{K}}(L_{1},L_{2},\ell) \neq B^{\textup{K}}(L_{3},L_{4},\ell)$, we obtain the thesis.
\end{proof}

This is the first example where the non-invariance of some GR amplitudes can be established. Nevertheless, we know that after integration over the moduli space against $\mu_{\textup{WP}}$, we obtain the symmetric polynomials \eqref{eqn:Kon:poly}. This is clear by the definition of the Kontsevich volumes, but it can also be directly proved from the theory of the topological recursion, as $(A^{\textup{K}},B^{\textup{K}},C^{\textup{K}})$ obey a set of quadratic relations that guarantee the invariance of the corresponding TR amplitudes under permutations of $L_1,\ldots,L_n$ -- see \textit{e.g.} \cite{Bor20}.

\paragraph{Support with non-empty complement.}

We prove that the Kontsevich amplitudes are actially zero on some open subset. In order to achieve this, we use geometric recursion to construct an auxiliary function which has the same support as the Kontsevich amplitude and takes integer values.
\medskip

Consider the following geometric recursion initial data
\begin{equation}
\begin{aligned}
	A(L_{1},L_{2},L_{3}) & = 1, \\
	B(L_1,L_2,\ell) & = \theta(L_1 + L_2 - \ell), \\
	C(L_1,\ell,\ell') & = \theta(L_1 - \ell - \ell'),
\end{aligned}
\end{equation}
where $\theta(x) = \mathbf{1}_{\mathbb{R}_{+}}(x)$ is the Heaviside theta function. Let $\Omega_{\Sigma}^{\theta} \in {\rm Mes}(\mathcal{T}_{\Sigma})$ be the geometric recursion amplitude associated to $\Sigma$ computed with respect to hyperbolic lengths. The value $\Omega_{\Sigma}^{\theta}(\sigma) \in \NN$ is counting the number of $\sigma$-small pair of pants decompositions\footnote{Such pants decomposition are ``rooted'', in the sense that they always bound the first boundary component of the surface at each step (\emph{cf}. \cite[Section~3.6]{ABO17})} of $\Sigma$. In particular, these are piecewise constant functions on $\mathcal{T}_{\Sigma}$ with values in non-negative integers and, for $\Sigma$ connected of type $(g,n)$, we see that $\Omega_{\Sigma}^{\theta}$ is bounded by $\prod_{k = 1}^{2g - 2 + n} 6k$. The following lemma easily follows by induction on $2g-2+n$.

\begin{lem}
	The GR amplitudes $\Omega_{\Sigma}^{\theta}$ have the same support in $\mathcal{T}_{\Sigma}^{\textup{comb}}$ as $\Omega_{\Sigma}^{\textup{K}}$, and the TR amplitudes $V\Omega_{g,n}^{\theta}(L) = \int_{\mathcal{M}_{g,n}(L)} \Omega_{\Sigma}^{\theta} \, \dd\mu_{\textup{WP}}$ are homogeneous polynomials of degree $6g - 6 + 2n$.
	\hfill $\blacksquare$
\end{lem}

\begin{cor}
	For $\Sigma$ not of type $(0,3)$, the support of $\Omega_{\Sigma}^{\textup{K}}$ has non-empty complement.
\end{cor}

\begin{proof}
	The Weil--Petersson volumes $V^{\textup{WP}}_{g,n}(L)$ are polynomials in $L = (L_1,\ldots,L_n)$ with non-zero constant term while $V\Omega_{g,n}^{\theta}(L)$ are homogeneous of positive degree for $(g,n) \neq (0,3)$. Hence
	\[
		\lim_{L\to 0} \frac{V\Omega_{g,n}^{\theta}(L)}{V^{\textup{WP}}_{g,n}(L)} = 0.
	\]
	It shows that $V\Omega_{g,n}^{\theta}(L)$ is strictly less than $V^{\textup{WP}}_{g,n}(L)$ for small $L$ and we deduce there exists an open set of $\mathcal{T}_{\Sigma}^{\textup{comb}}(L)$ for small $L$ on which $\Omega_{\Sigma}^{\theta} < 1$ and therefore on which $\Omega_{\Sigma}^{\theta} = 0$, hence $\Omega_{\Sigma}^{\textup{K}}$ vanishes.
\end{proof}

\begin{rem}
	To conclude, we remark that the Laplace transform of the TR amplitudes $V\Omega_{g,n}^{\theta}(L)$, namely
	\begin{equation}
		\omega_{g,n}(z_1,\ldots,z_n) = \left(
			\int_{\RR_{+}^n} \prod_{i=1}^n \dd L_i \, L_i\,e^{-z_iL_i}\, V\Omega^{\theta}_{g,n}(L_1,\ldots,L_n)  \,\right)
			\dd z_1 \otimes \cdots \otimes \dd z_n
	\end{equation}
	satisfy a topological recursion \emph{\`{a} la} Eynard--Orantin \cite{EO07inv}:
	\begin{equation}
		\omega_{g,n}(z_1,\ldots,z_n)
		=
		\Res_{z \to 0} \; K(z_1,z) \bigg(\omega_{g - 1,n + 1}(z,-z,z_2,\ldots,z_n)
		\, +
		\!\!\!\!\! \sum_{\substack{h + h ' = g \\ J \sqcup J' = \{z_2,\ldots,z_n\}}}^{\text{no }(0,1)} \!\!\!\!\! \omega_{h,1 + \#J}(z,J) \otimes \omega_{h',1 + \#J'}(-z,J') \bigg),
	\end{equation}
	where
	\begin{equation}
		K(z_{1},z) = \frac{1}{2z(z-z_{1})^{2}}\frac{\dd z_{1}}{\dd z} 
		\qquad \text{and} \qquad 
		\omega_{0,2}(z_{1},z_{2}) =  \frac{\dd z_{1}\otimes \dd z_{2}}{(z_{1}-z_{2})^{2}}.
	\end{equation}  
	However, in this case the recursion kernel $K(z_1,z)$ does not have the usual structure of \cite{EO07inv} and the mul\-ti\-dif\-fe\-ren\-tials $\omega_{g,n}(z_1,\ldots,z_n)$ are not symmetric under permutation of their $n$ variables. It seems to be the first known example where a non-symmetric topological recursion yields a geometrically meaningful quantity.
\end{rem}

\newpage
\section{Rescaling flow: from hyperbolic to combinatorial geometry}
\label{sec:hyp:comb}

The hyperbolic and combinatorial Teichm\"uller spaces can be identified via the spine homeomorphism $\sp \colon \mathcal{T}_{\Sigma} \to \mathcal{T}_{\Sigma}^{\textup{comb}}$ defined in \S~\ref{subsubsec:ordinary:combinatorial}. We introduce a flow that interpolates between their respective geometries, coming from the work of Mondello and Do \cite{Mon09,Do10}.

\begin{defn}\label{defn:rescaling:flow}
	Let $\GG\in\mathcal{T}^{\textup{comb}}_{\Sigma}(L)$ and $\beta\in\mathbb{R}_{+}$. Define $\beta\GG\in\mathcal{T}^{\textup{comb}}_{\Sigma}(\beta L)$ to be the combinatorial structure represented by the same marked ribbon graph, but with all lengths multiplied by $\beta$. This define a flow, called the \emph{rescaling flow}, on $\mathcal{T}^{\textup{comb}}_{\Sigma}$ which preserves the strata and is continuous. It can be lifted to $\Phi_{\beta} \colon \mathcal{T}_{\Sigma}(L) \to \mathcal{T}_{\Sigma}(\beta L)$ by the spine map: for $\sigma \in \mathcal{T}_{\Sigma}(L)$, set
	\begin{equation}\label{eqn:flow}
		\Phi_{\beta}(\sigma) = \sigma^{\beta} = \sp^{-1}\big( \beta \, \sp(\sigma) \big) \in \mathcal{T}_{\Sigma}(\beta L).
	\end{equation}
	The maps $\Phi_{\beta}$ and $\beta \cdot$ are $\Mod_{\Sigma}$-equivariant, and thus descend to the moduli spaces $\mathcal{M}_{g,n}$ and $\mathcal{M}^{\textup{comb}}_{g,n}$. We also define the map $R_{\beta} \colon \mathcal{T}^{\textup{comb}}_{\Sigma}(L) \to \mathcal{T}_{\Sigma}(\beta L)$ by postcomposing $\beta \cdot$ with $\sp^{-1}$.
\end{defn}

\begin{center}
\begin{tikzcd}
	\mathcal{T}_{\Sigma}(L) \arrow[r,"\sp"] \arrow[d,dashrightarrow,"\Phi_{\beta}" left] &
	\mathcal{T}_{\Sigma}^{\textup{comb}}(L) \arrow[d,"\beta \cdot"] \arrow[dl,dashrightarrow,"R_{\beta}" above] \\
	\mathcal{T}_{\Sigma}(\beta L) \arrow[r,leftarrow,"\sp^{-1}"] &
	\mathcal{T}_{\Sigma}^{\textup{comb}}(\beta L) 
\end{tikzcd}
\end{center}

The large $\beta$ asymptotic of the rescaling flow has been previously studied pointwise on $\mathcal{T}_{\Sigma}$.

\begin{thm}\cite{Mon09,Do10}\label{thm:Mon:Do}
	For any bordered surface $\Sigma$, the metric space $(\Sigma,\beta^{-1}\sigma^{\beta})$ converges in the Gromov--Hausdorff topology to the metric ribbon graph $\sp(\sigma)$ when $\beta \to \infty$. Moreover, the Poisson structure $\beta^2 R_{\beta}^* \pi_{\textup{WP}}$ converges pointwise to $\pi_{\textup{K}}$ on the top-dimensional strata of $\mathcal{T}_{\Sigma}^{\textup{comb}}$.
	\hfill $\blacksquare$
\end{thm}

Here we shall complete this description by giving effective bounds on the thick part of the Teichm\"uller space for lengths, and on compacts for twist parameters, showing convergence of the Fenchel--Nielsen coordinates. Under the flow, combinatorial geometry appears as a kind of tropicalisation of hyperbolic geometry. We apply these results to prove that Fenchel--Nielsen coordinates equip $\mathcal{T}_{\Sigma}^{\textup{comb}}$ with a natural piecewise linear structure, and that the flow carries the hyperbolic GR/TR to the combinatorial GR/TR.

\subsection{Convergence of lengths and twists}
\label{subsec:convergence:length:twist}

\subsubsection{Convergence of lengths}

We first obtain an effective comparison between lengths of simple closed curves along the flow, by refining the arguments of \cite{Do10}.

\begin{defn}
	Let $\Sigma$ be bordered surface. The hyperbolic (resp. combinatorial) \emph{systole} with respect to $\sigma \in \mathcal{T}_{\Sigma}$ (resp. $\GG \in \mathcal{T}_{\Sigma}^{\textup{comb}}$) as
	\begin{equation}
		\sys_{\sigma} = \inf_{\gamma \in S_{\Sigma}} \ell_{\sigma}(\gamma),
		\qquad
		\sys_{\GG} = \inf_{\gamma \in S_{\Sigma}} \ell_{\GG}(\gamma).
	\end{equation}
\end{defn}

\begin{prop}\label{lem:length:bound}
	Let $\Sigma$ be a connected bordered surface of topology $(g,n)$. For any $\GG \in \mathcal{T}^{\textup{comb}}_{\Sigma}$, we denote $\sigma = \sp^{-1}(\GG)$ and $\sigma^{\beta} = \sp^{-1}(\beta \GG)$. Then for all $\beta \geq 1$, we have
	\begin{equation}\label{eqn:length:bound:below}
		\forall \gamma \in S_{\Sigma},
		\qquad
		\ell_{\GG}(\gamma) \leq \frac{\ell_{\sigma^{\beta}}(\gamma)}{\beta}.
	\end{equation}
	Moreover, for any $\epsilon > 0$ there exists $\beta_{\epsilon} \geq 1$ and $\kappa_{\epsilon} > 0$ depending only on $\epsilon$ and the topological type of $\Sigma$, such that for any $\beta \geq \beta_{\epsilon}$ and $\GG \in \mathcal{T}^{\textup{comb}}_{\Sigma}$ for which the systole in $\GG$ is bounded from below by $\epsilon$, we have
	\begin{equation}\label{eqn:length:bound:above}
		\forall \gamma \in S_{\Sigma},
		\qquad
		\frac{\ell_{\sigma^{\beta}}(\gamma)}{\beta} \leq \bigg(1 + \frac{\kappa_{\epsilon}}{\beta}\bigg)\ell_{\GG}(\gamma).
	\end{equation}
	We can take $\kappa_{\epsilon}$ and $\beta_{\epsilon}$ increasing with $2g - 2 + n$.
\end{prop}

\begin{cor} \label{prop:rescale:length}
	For any $\gamma \in S_{\Sigma}$ and $\sigma \in \mathcal{T}_{\Sigma}$, we have
	\begin{equation}
		\lim_{\beta \to \infty} \frac{\ell_{\sigma^{\beta}}(\gamma)}{\beta}
		=
		\ell_{\sp(\sigma)}(\gamma),
	\end{equation}
	and the limit is uniform when $\sp(\sigma)$ belongs to the thick part of $\mathcal{T}_{\Sigma}^{\textup{comb}}$.
	\hfill $\blacksquare$
\end{cor}

\begin{proof}[Proof of Proposition~\ref{lem:length:bound}]
	If $\gamma$ is a boundary curve, there is nothing to prove since $\ell_{\sigma^{\beta}}(\gamma) = \beta\ell_{\sigma}(\gamma)$. We now assume $\gamma \in S_{\Sigma}^{\circ}$ and start from the last inequality in the proof of \cite[Lemma 11]{Do10}:
	\begin{equation}\label{eqn:Do:bound}
		\forall \gamma \in S_{\Sigma}^{\circ} \qquad \ell_{\GG}(\gamma) \leq \frac{\ell_{\sigma^{\beta}}(\gamma)}{\beta} \leq \ell_{\GG}(\gamma) + \frac{2E_{\GG}(\gamma)\,r_{\beta}}{\beta},
	\end{equation}
	where $E_{\GG}(\gamma)$ is the number (with multiplicity) of edges along which the non-backtracking representative of $\gamma$ in $\GG$ travels, and $r_{\beta}$ is the maximal of rib lengths with respect to $\sigma^{\beta}$. This already gives the lower bound.
		
	\smallskip

	Recall that $\GG$ contains at most $6g - 6 + 3n$ edges. We say that an edge is $\epsilon$-big if its length in $\GG$ is larger or equal to $\epsilon/(6g - 6 + 3n)$, and that it is $\epsilon$-short otherwise. Denote by $E^{(\epsilon)}_{\GG}(\gamma)$ the number (with multiplicity) of $\epsilon$-big edges along which $\gamma$ travels. If we assume that $\sys_{\GG} \geq \epsilon$, the union of the $\epsilon$-short edges appearing in $\gamma$ must be a forest. As $\gamma$ is a closed loop, it has to exit each tree it passes through via an $\epsilon$-big edge. Hence $E_{\GG}(\gamma) \leq (6g - 6 + 3n)E_{\GG}^{(\epsilon)}(\gamma)$.  Observing that $\frac{\epsilon}{6g - 6 + 3n}\,E^{(\epsilon)}_{\GG}(\gamma) \leq \ell_{\GG}(\gamma)$, we obtain  
	\[
		E_{\GG}(\gamma) \leq \frac{(6g-6+3n)^{2}}{\epsilon} \, \ell_{\GG}(\gamma).
	\]
	We must now bound uniformly the maximum of rib lengths for $\sigma^{\beta}$. For this purpose, we bound the distance between any point of the surface to the boundary for the metric $\sigma^{\beta}$. Recall that the injectivity radius at a point $q \in \Sigma$ for the hyperbolic structure $\sigma$, here denoted $\mathfrak{r}_{\sigma}(q)$, is the supremum of all $\rho > 0$ such that there is a locally isometric embedding of an open hyperbolic disk of radius $\rho$. Since the area of $(\Sigma,\sigma^{\beta})$ is $2\pi(2g - 2 + n)$, which must be greater or equal to the area of such disks, we have
	\begin{equation}\label{eqn:inj:rad:bound} 
		\mathfrak{r}_{\sigma}(q) \leq \sqrt{2(2g - 2 + n)}.
	\end{equation}
	Besides, it is clear that
	\begin{equation}\label{eqn:inj:rad:equal}
		\mathfrak{r}_{\sigma}(q)
		=
		\min\big\{\tfrac{1}{2} \sys_{\sigma}(q), {\rm dist}_{\sigma}(q,\de\Sigma)\big\},
	\end{equation}
	where $\sys_{\sigma}(q)$ is defined to be the infimum over the lengths of non-constant geodesic loops based at $q$. We apply this to the hyperbolic structure $\sigma^{\beta}$. Using the lower bound \eqref{eqn:Do:bound}, we remark that
	\[
		\beta \epsilon \leq \sys_{\beta \GG} \leq \sys_{\sigma^{\beta}} \leq \sys_{\sigma^{\beta}}(q).
	\]
	If $\beta \geq \beta_{\epsilon} = \frac{2\sqrt{2(2g - 2 + n)}}{\epsilon}$, we deduce from \eqref{eqn:inj:rad:bound} and \eqref{eqn:inj:rad:equal} that
	\[
		\forall q \in \Sigma,
		\qquad
		{\rm dist}_{\sigma^{\beta}}(q,\de \Sigma) \leq \sqrt{2(2g - 2 + n)}.
	\]
	In particular, choosing for $q$ the vertices of $\GG$, we find $r_{\beta} \leq \sqrt{2(2g - 2 + n)}$. Together with \eqref{eqn:Do:bound}, it shows that
	\[
		\forall \beta \geq \beta_{\epsilon}
		\quad
		\forall \gamma \in S_{\Sigma}^{\circ},
		\qquad
		\ell_{\GG}(\gamma)
		\leq
		\frac{\ell_{\sigma^{\beta}}(\gamma)}{\beta}
		\leq
		\bigg(1 + \frac{18\sqrt{2}(2g - 2 + n)^{5/2}}{\beta \epsilon}\bigg)\ell_{\GG}(\gamma),
	\]
	which gives the thesis.
\end{proof}

\subsubsection{Convergence of twists}
\label{subsubsec:conv:twist}
Studying the convergence of twists along the flow by a direct geometric method requires bounds on the distance between the hyperbolic and combinatorial seams. Instead of following this direction, we use the hyperbolic $(9g - 9 + 3n)$-theorem to write the hyperbolic twists in terms of hyperbolic lengths of certain curves, to show that these formulae converge with help of Proposition~\ref{lem:length:bound}, and compare the limit to the expressions for the combinatorial twists underlying the combinatorial $(9g - 9 + 3n)$-theorem established in Section~\ref{subsec:9g:minus:9:plus:3n}.

\begin{prop}\label{lem:twist:bound}
	Let $\Sigma$ be a bordered surface of type $(g,n)$, fix a seamed pants decomposition and let $(\tau_i)_{i = 1}^{3g - 3 + n}$ be the associated combinatorial twist parameters. For any compact $K \subset \mathcal{T}_{\Sigma}^{\text{comb}}$ and $\GG \in K$, we denote $\sigma = \sp^{-1}(\GG)$ and $\sigma^{\beta} = \sp^{-1}(\beta \GG)$. There exists constants $\beta_{K} \ge 1$ and $c_{K} > 0$ depending only on $K$ such that, for any $\beta \geq \beta_{K}$ and $i \in \{1,\ldots,3g - 3 + n\}$, we have
	\begin{equation}\label{eqn:tau:bound}
		\bigg|\frac{\tau_i(\sigma^{\beta})}{\beta} - \tau_i(\GG)\bigg|
		\leq
		\frac{c_{K}}{\beta}.
	\end{equation}
\end{prop}
\begin{cor} \label{prop:rescale:twist}
	In any fixed seamed pants decomposition, for $\sigma \in \mathcal{T}_{\Sigma}$, we have
	\begin{equation}
		\lim_{\beta \rightarrow \infty} \frac{\tau_i(\sigma^{\beta})}{\beta}
		=
		\tau_i \big(\sp(\sigma)\big),
	\end{equation}
	and the limit is uniform when $\sigma$ belongs to an arbitrary compact of $\mathcal{T}_{\Sigma}$.
	\hfill $\blacksquare$
\end{cor}

Before starting the proof of Proposition~\ref{lem:twist:bound}, we recall the formulae which allows us for four-holed spheres and one-holed tori to express the change of Fenchel--Nielsen coordinates under a flip of the pair of  pants decomposition. They can be found in \cite[Theorems~1.i and 2.i]{Okai93}, or deduced from \cite[Sections~3.3 and 3.4]{Bus10}.

\medskip

Let $X$ be a four-holed sphere and $\sigma \in \mathcal{T}_{X}(L_1,L_2,L_3,L_4)$. We place ourselves in the situation described in \S~\ref{subsubsec:four:holed:sphere}. Namely, we fix a coordinate system, which in turn defines a simple closed curve $\gamma$ separating $X$ into a pair of pants having boundary components $(\de_1 X,\de_4 X,\gamma)$ and another pair of pants having boundary components $(\de_2X,\de_3X,\gamma)$; we have $\delta$ a simple closed curve intersecting $\gamma$ exactly twice and separating $X$ into a pair of pants having boundary components $(\de_1 X,\de_2 X,\delta)$ and another pair of pants having boundary components $(\de_3 X,\de_4 X,\delta)$; finally, let $\eta$ be the curve obtained from $\gamma$ by applying a Dehn twist along $\delta$. Let $\sigma \in \mathcal{T}_{X}(L_1,L_2,L_3,L_4)$ and denote $\ell = \ell_{\sigma}(\gamma)$, $\ell' = \ell_{\sigma}(\delta)$ and $\ell'' = \ell_{\sigma}(\eta)$, and define $\tau$ to be the hyperbolic twist determined by the seamed pants decomposition. If we define
\begin{equation}\label{eqn:constant:Cij}
	C_{i,j}(\ell)
	=
	\ch^{2}(\tfrac{\ell}{2}) + \ch^{2}(\tfrac{L_i}{2}) + \ch^{2}(\tfrac{L_j}{2})
	+ 2 \ch(\tfrac{L_i}{2}) \ch(\tfrac{L_j}{2}) \ch(\tfrac{\ell}{2})
	- 1,
\end{equation}
the length of $\delta$ is then given in terms of Fenchel--Nielsen coordinates by
\begin{equation}\label{eqn:Okai:elldash:04}
\begin{split}
	& \ch\left(\tfrac{\ell'(\ell,\tau)}{2}\right)
		\sh^{2}\left(\tfrac{\ell}{2}\right)
		\\
		& \quad =
		\ch(\tfrac{L_1}{2}) \ch(\tfrac{L_2}{2}) + \ch(\tfrac{L_3}{2}) \ch(\tfrac{L_4}{2}) 
		+
		\ch(\tfrac{\ell}{2})
		\Big(
			\ch(\tfrac{L_1}{2}) \ch(\tfrac{L_3}{2}) + \ch(\tfrac{L_2}{2}) \ch(\tfrac{L_4}{2})
		\Big) \\
		& \quad\quad +
		\ch(\tau) \sqrt{C_{1,4}(\ell)C_{2,3}(\ell)},
\end{split}
\end{equation}
while the length of $\eta$ is $\ell''(\ell,\tau) = \ell'(\ell,\tau + \ell)$.

\medskip

Likewise, if $T$ is a one-holed torus, $\sigma \in \mathcal{T}_{T}(L)$, and we are in the situation described in \S~\ref{subsubsec:one:holed:torus}, the length of $\delta$ is given by
\begin{equation}\label{eqn:Okai:elldash:11}
	\ch\left(\tfrac{\ell'(\ell,\tau)}{2}\right)
		=
		\frac{\ch(\tfrac{\tau}{2})}{\sh(\tfrac{\ell}{2})}\sqrt{\frac{\ch(\tfrac{L}{2}) + \ch(\ell)}{2}},
\end{equation}
while the length of $\eta$ is $\ell''(\ell,\tau) = \ell'(\ell,\tau + \ell)$.

\begin{proof}[{Proof of Proposition~\ref{lem:twist:bound}}]
	Fix $\sigma$ in a compact $K$ of $\mathcal{T}_{\Sigma}$, and denote $\GG = \sp(\sigma)$. We use repeatedly Proposition~\ref{lem:length:bound}, which implies that for any simple closed curve $\nu$ chosen in a fixed finite subset of $S_{\Sigma}$, we have
	\[
		\ell_{\GG}(\nu) \leq \frac{\ell_{\sigma^{\beta}}(\nu)}{\beta} \leq \ell_{\GG}(\nu) + \frac{c_{K}}{\beta}
	\]
	for any $\beta \ge \beta_{K}$ and some constant $c_{K} > 0$ depending only on the compact $K$ and this finite set. In what follows, $c_{K},c_K',\ldots$ denote positive constants depending on $K$ and whose value may change from line to line. We denote lengths and twists with a superscript $\beta$ to refer to the hyperbolic quantities measured with respect to $\sigma^{\beta}$, while lengths and twists without superscripts denote the combinatorial quantities measured with respect to $\GG$.

	\smallskip

	As the twist parameters $\tau_i^{\beta}$ and $\tau_i$ are computed locally in each piece $\Sigma_i$ of type $(0,4)$ or $(1,1)$ defined by the seamed pants decomposition as in Section~\ref{subsec:9g:minus:9:plus:3n}, we can restrict our attention to each piece separately. On $\Sigma_i$ we have the curve $\delta_i$ defined by the seamed pants decomposition, and for every $k \in \ZZ$ we consider the curve obtained as the image of $\delta_i$ after $k$ Dehn twists along $\gamma_i$. Denote by $\ell_i^{(k),\beta}$ its hyperbolic length with respect to $\sigma^{\beta}$, and by $\ell_i^{(k)}$ its combinatorial length with respect to $\GG$. Compared to the notation of Section~\ref{subsec:9g:minus:9:plus:3n},
	\[
		\ell_i'^{\beta} = \ell_i^{(0),\beta},
		\quad
		\ell_i''^{\beta} = \ell_i^{(1),\beta},
		\qquad
		\ell_i' = \ell_i^{(0)},
		\quad
		\ell_i'' = \ell_i^{(1)}.
	\]
	As we work with a single piece at a time, we in fact often omit the subscript $i$. Suppose $\Sigma_i = X$ is a four-holed sphere. With the labelling matching the one described in \S~\ref{subsubsec:four:holed:sphere}, we denote $L_i^{\beta} = \ell_{\sigma^{\beta}}(\de_iX)$ for $i \in \{1,2,3,4\}$. According to Lemma~\ref{lem:delta:eta:04}, we have
	\begin{equation}\label{eqn:ell:k}
		\ell^{(k)} = \max\big\{L_1 + L_3 - \ell,L_2 + L_4 - \ell,2|\tau + k\ell| + M_{1,4}(\ell) + M_{2,3}(\ell)\big\}
	\end{equation}
	where $M_{i,j}(\ell) = \max\big\{0,L_i - \ell,L_j - \ell,\tfrac{L_i + L_j - \ell}{2}\big\}$. The open sets
	\begin{equation*}
		U^{(k)}
		=
		\Set{ \GG \in \mathcal{T}_{\Sigma}^{\text{comb}} | 2\ell < \tau + k\ell < 4\ell },
		\qquad
		k \in \ZZ.
	\end{equation*}
	cover the compact $K$, so we can select finitely many indices $k_1,\dots,k_N \in \ZZ$ to cover $K$. If $\GG \in U^{(k)} \cap K$ for some $k \in \set{k_1,\dots,k_N}$, the maximum in \eqref{eqn:ell:k} is given by the last argument, that is
	\[
		\ell^{(k)} = 2(\tau + k\ell) + M_{1,4}(\ell) + M_{2,3}(\ell)
	\]
	which we see as an expression of $\tau + k\ell$ in terms of other lengths. We would like to compare it, when $\beta$ is large, to the expression of $\tau^{\beta} + k\ell^{\beta}$ which we can access via Equation~\eqref{eqn:Okai:elldash:04}, namely
	\begin{equation}\label{eqn:cosh:tau:k}
	\begin{split}
		\ch\left( \tau^{\beta} + k\ell^{\beta} \right) 
		=
		\bigg[ &
			\ch\big(\tfrac{\ell^{(k),\beta}}{2}\big)
			\sh^{2}\!\big(\tfrac{\ell^{\beta}}{2}\big)
			-
			\ch\big(\tfrac{L_1^{\beta}}{2}\big) \ch\big(\tfrac{L_2^{\beta}}{2}\big)
			-
			\ch\big(\tfrac{L_3^{\beta}}{2}\big) \ch\big(\tfrac{L_4^{\beta}}{2}\big)
			\\
			& -
			\ch\big(\tfrac{\ell^{\beta}}{2}\big)
			\Big(
				\ch\big(\tfrac{L_1^{\beta}}{2}\big) \ch\big(\tfrac{L_3^{\beta}}{2}\big)
				+
				\ch\big(\tfrac{L_2^{\beta}}{2}\big) \ch\big(\tfrac{L_4^{\beta}}{2}\big)
			\Big)
		\bigg]
		C_{1,4}^{-1/2}(\ell^{\beta}) C_{2,3}^{-1/2}(\ell^{\beta}).
	\end{split}
	\end{equation}
	To obtain an upper bound for the right-hand side of \eqref{eqn:cosh:tau:k}, we can ignore the negative terms and use as upper bound for the numerator 
	\[
		\ch\big(\tfrac{\ell^{(k),\beta}}{2}\big)\sh^{2}\!\big(\tfrac{\ell^{\beta}}{2}\big)
		\le
		\tfrac{1}{4}\ e^{\tfrac{\ell^{(k),\beta}}{2} + \ell^{\beta}},
	\]
	and as lower bound for the factors in the denominator
	\[
		C_{i,j}(\ell^{\beta})
		\geq
		\tfrac{1}{4}\,e^{\max\big\{\ell^{\beta},L_i^{\beta},L_j^{\beta},\frac{L_i^{\beta} + L_j^{\beta} + \ell^{\beta}}{2}\big\}}
		\geq
		\tfrac{1}{4}\,e^{\max\big\{L_i^{\beta},L_j^{\beta},\frac{L_i^{\beta} + L_j^{\beta} + \ell^{\beta}}{2}\big\}}.
	\]
	Combined with $\arch(x) \leq \ln(2x)$ for $x\geq 1$, this results in
	\[
	\begin{split}
		\frac{|\tau^{\beta} + k\ell^{\beta}|}{\beta}
		& \leq
		\frac{1}{\beta}\bigg(\ln 2 + \frac{\ell^{(k),\beta}}{2} + \ell^{\beta} - \tfrac{1}{2}\max\big\{L_1^{\beta},L_4^{\beta},\tfrac{L_1^{\beta} + L_4^{\beta} + \ell^{\beta}}{2}\big\} - \tfrac{1}{2}\max\big\{L_2^{\beta},L_3^{\beta},\tfrac{L_2^{\beta} + L_3^{\beta} + \ell^{\beta}}{2}\big\}\bigg) \\
		& \leq 
		\frac{\ell^{(k)}}{2} + \ell - \tfrac{1}{2}\max\big\{L_1,L_4,\tfrac{L_1 + L_4 + \ell}{2}\big\} - \tfrac{1}{2}\max\big\{L_2,L_3,\tfrac{L_2 + L_3 + \ell}{2}\big\} + \frac{c_{K}'}{\beta} \\
		& \leq
		\frac{\ell^{(k)}}{2} - \frac{1}{2}\big(M_{1,4}(\ell) + M_{2,3}(\ell)\big) + \frac{c_{K}'}{\beta}
	\end{split}
	\]
	and thus
	\[
		\frac{|\tau^{\beta} + k\ell^{\beta}|}{\beta} \leq \tau + k\ell + \frac{c_{K}'}{\beta}.
	\]
	We now look for a bound from below for \eqref{eqn:cosh:tau:k}. We first observe that by Equation~\eqref{eqn:length:bound:below}, we have $\ell^{\beta} \geq \ell$ for $\beta \geq 1$. Since on the compact $K$, $\ell$ is bounded from below by $\epsilon_K > 0$, we deduce that
	\[  
		\sh^2\!\big(\tfrac{\ell^{\beta}}{2}\big) \geq m_{K} e^{\ell^{\beta}}
		\qquad \text{with} \qquad
		m_{K} = \frac{(1 - e^{-\epsilon_K})^2}{4} > 0.
	\]  
	This leads to a (rather crude) lower bound for the numerator of the right-hand side of \eqref{eqn:cosh:tau:k}
	\begin{equation}\label{eqn:bound:num}
	\begin{split}
		&
		\ch\big(\tfrac{\ell^{(k),\beta}}{2}\big)
		\sh^{2}\!\big(\tfrac{\ell^{\beta}}{2}\big)
		-
		\ch\big(\tfrac{L_1^{\beta}}{2}\big) \ch\big(\tfrac{L_2^{\beta}}{2}\big)
		-
		\ch\big(\tfrac{L_3^{\beta}}{2}\big) \ch\big(\tfrac{L_4^{\beta}}{2}\big)
		\\
		& \qquad\qquad\quad
		-
		\ch\big(\tfrac{\ell^{\beta}}{2}\big)
		\Big(
			\ch\big(\tfrac{L_1^{\beta}}{2}\big) \ch\big(\tfrac{L_3^{\beta}}{2}\big)
			+
			\ch\big(\tfrac{L_2^{\beta}}{2}\big) \ch\big(\tfrac{L_4^{\beta}}{2}\big)
		\Big) \\
		& \ge
		\tfrac{m_K}{2} e^{\tfrac{\ell^{(k),\beta}}{2} + \ell^{\beta}}
		-
		e^{\tfrac{L_1^{\beta} + L_2^{\beta}}{2}}
		-
		e^{\tfrac{L_3^{\beta} + L_4^{\beta}}{2}}
		-
		e^{\tfrac{\ell^{\beta} + L_1^{\beta} + L_3^{\beta}}{2}}
		-
		e^{\tfrac{\ell^{\beta} + L_2^{\beta} + L_4^{\beta}}{2}}
		\\
		& \ge
		\tfrac{m_K}{4} e^{\tfrac{\ell^{(k),\beta}}{2} + \ell^{\beta}} + \bigg(\tfrac{m_{K}}{4}e^{\tfrac{\ell^{(k),\beta}}{2} + \ell^{\beta}} - 4\,e^{\frac{1}{2}\max\{L_1^{\beta} + L_2^{\beta},L_3^{\beta} + L_4^{\beta},L_1^{\beta} + L_3^{\beta} + \ell^{\beta},L_2^{\beta} + L_4^{\beta} + \ell^{\beta}\}}\bigg).
	\end{split}
	\end{equation} 
	We split the first term in order to exhibit positivity of our lower bound -- which is therefore not useless. Indeed, using Proposition~\ref{lem:length:bound} and Equation~\eqref{eqn:ell:k} on $U^{(k)} \cap K$ we get
	\[
		\frac{\ell^{(k),\beta}}{2} + \ell^{\beta}  \geq \max\big\{L_1^{\beta} + L_2^{\beta},L_3^{\beta} + L_4^{\beta},L_1^{\beta} + L_3^{\beta} + \ell^{\beta},L_2^{\beta} + L_4^{\beta} + \ell^{\beta}\big\} + \beta c''_{K}
	\]
	for some constant $c''_{K} > 0$. Consequently, there exists $\beta_{K} \geq 1$ such that for any $\beta \geq \beta_{K}$ the expression inside the bracket in \eqref{eqn:bound:num} is positive, and can be ignored in the lower bound. To obtain an upper bound for the denominator of \eqref{eqn:cosh:tau:k}, we write
	\[
		C_{i,j}(\ell^{\beta})
		\leq
		5 \, e^{\max\big\{\ell^{\beta},L_i^{\beta},L_j^{\beta},\tfrac{L_i^{\beta} + L_j^{\beta} + \ell^{\beta}}{2}\big\}}.
	\]
	Combined with ${\rm arccosh}(x) \geq \ln x$, this implies for $\beta \geq \beta_{K}$
	\begin{equation*}
	\begin{split}
		\frac{|\tau^{\beta} + k\ell^{\beta}|}{\beta}
		& \geq
		\frac{1}{\beta} \bigg(
			\ln\left(\tfrac{m_K}{20}\right) + \tfrac{1}{2}\ell^{(k),\beta}
			-
			\tfrac{1}{2} \max\big\{
				\ell^{\beta},L_1^{\beta},L_4^{\beta},\tfrac{L_1^{\beta} + L_4^{\beta} + \ell^{\beta}}{2}
			\big\}
			-
			\tfrac{1}{2} \max\big\{
				\ell^{\beta},L_2^{\beta},L_3^{\beta},\tfrac{L_2^{\beta} + L_3^{\beta} + \ell^{\beta}}{2}
			\big\}
			\bigg) \\
		& \geq \tau + k\ell - \frac{c'''_{K}}{\beta}
	\end{split}
	\end{equation*}
using the arguments we already used for the upper bound. We deduce that on $U^{(k)} \cap K$ and for $\beta \geq \beta_{K}$
	\begin{equation}\label{eqn:bound:tau:k}
		\bigg|\frac{|\tau^{\beta} + k\ell^{\beta}|}{\beta} - (\tau + k\ell)\bigg|
		\leq
		\frac{c_{K}}{\beta}.
	\end{equation}
	A similar argument shows that on $U^{(k)} \cap K$ and for $\beta \geq \beta_K$, we have
	\begin{equation}\label{eqn:bound:tau:k:plus:1}
		\bigg|\frac{|\tau^{\beta} + (k+1) \ell^{\beta}|}{\beta} - \big( \tau + (k+1) \ell \big)\bigg|
		\leq
		\frac{c_{K}}{\beta},
	\end{equation}
	for perhaps larger constants $\beta_{K},c_{K} > 0$. We can now conclude by using \eqref{eqn:bound:tau:k}--\eqref{eqn:bound:tau:k:plus:1} to estimate
	\[
		\frac{\tau^{\beta}}{\beta} = \frac{\beta}{2\ell^{\beta}}\bigg(\bigg|\frac{\tau^{\beta} + (k + 1)\ell^{\beta}}{\beta}\bigg|^2 - \bigg|\frac{\tau^{\beta} + k\ell^{\beta}}{\beta}\bigg|^2\bigg) - \frac{(2k + 1)\ell^{\beta}}{2\beta}.
	\]
	Since $\ell^{\beta} \geq \epsilon_K > 0$ on the compact $K$, we arrive on $U^{(k)} \cap K$ and for $\beta \geq \beta_{K}$ at the inequality
	\begin{equation}\label{eqn:tau:bound:proof}
		\bigg|\frac{\tau^{\beta}}{\beta} - \tau\bigg| \leq \frac{c_K}{\beta},
	\end{equation}
	for perhaps larger constants $\beta_{K},c_{K} > 0$. These arguments were done for a fixed $k \in \{k_1,\ldots,k_N\}$, and as $(U^{(k_i)})_{i = 1}^N$ cover $K$, we can find constants $\beta_{K}$ and $c_{K}$ such that the same estimate \eqref{eqn:tau:bound:proof} holds uniformly over $K$ for $\beta \geq \beta_{K}$.

	\smallskip

	A similar argument can be carried out for $\Sigma_i$ being a one-holed torus, in the situation described in \S~\ref{subsubsec:one:holed:torus}. Instead of \eqref{eqn:cosh:tau:k} we should estimate $\tau^{\beta} + k\ell^{\beta}$ via the formula
	\[ 
		\ch\big(\tfrac{\tau^{\beta} + k\ell^{\beta}}{2}\big) = \ch\big(\tfrac{\ell^{(k),\beta}}{2}\big)\sh\big(\tfrac{\ell^{\beta}}{2}\big)\,\sqrt{\frac{2}{\ch\big(\frac{L^{\beta}}{2}\big) + \ch(\ell^{\beta})}}
	\]
	coming from \eqref{eqn:Okai:elldash:11}, and compare it to the expression of $\tau + k\ell$ deduced from Lemma~\ref{lem:delta:eta:11}. As the estimates in this case are much easier and without surprise, we omit them.
	 
	\smallskip

	We conclude that the desired estimate \eqref{eqn:tau:bound} is valid for any fixed $i \in \{1,\ldots,3g - 3 + n\}$, and since this set is finite we can choose constants $\beta_{K}$ and $c_{K}$ independently of $i$: in the end, they only depend on the fixed seamed pants decomposition and on the compact $K$.
\end{proof}

\subsection{Penner's formulae}
\label{subsec:PL:structure}

The combinatorial Fenchel--Nielsen coordinates on $\mathcal{T}_{\Sigma}^{\text{comb}}$ described in Theorem~\ref{thm:FN:coordinates} depend on a seamed pants decomposition $(\mathscr{P},\mathscr{S})$. Notice that changing $\mathscr{S}$ amounts to changing $\tau_{i}(\GG) \mapsto \tau_{i}(\GG) + k_{i} \ell_{i}(\GG)$ by some $k_{i} \in \ZZ$, while $\ell_{i}(\GG)$ remains unchanged. The most interesting case occurs when we change the pair of pants decomposition $\mathscr{P}$. From~\cite{HT80}, these changes are generated by local changes in four-holed spheres and one-holed tori and were first described by Penner in \cite{Pen82,Pen84} for the case of multicurves. A generalisation of this action to measured laminations can be found in \cite{PH92}.

\medskip

Notice that the same formulae apply to $\mathcal{T}_{\Sigma}^{\text{comb}}$, as the different behaviour of the foliations at the boundary do not affect the coordinates themselves. In this section we give a new proof of Penner's result, by flowing their hyperbolic analogue found by Okai\footnote{In \cite[Theorem~2.ii]{Okai93} there is a misprint in the denominator of the formula giving $\tau'$ for the one-holed torus. More precisely, the first $\ch(\tfrac{L}{2})$ should not be squared, and the second $\ch(\tfrac{L}{2})$ should be replaced with $\ch(\tfrac{L}{2}) - 1$. We report the correct formula \eqref{eqn:taudash:11} here.} \cite{Okai93}. We now report these hyperbolic formulae.

\medskip

\textsc{Four-holed sphere.} For a four-holed sphere $X$ , consider the hyperbolic Fenchel--Nielsen coordinates $(\ell,\tau) \in \mathcal{T}_{X}(L_1,L_2,L_3,L_4)$ relative to the system of curves $(\mathscr{P},\mathscr{S})$ of Figure~\ref{fig:coord:system:04}. Then the change of seamed pants decomposition to  $(\mathscr{P}',\mathscr{S}')$ is given by Equation~\eqref{eqn:Okai:elldash:04} and
\begin{equation}\label{eqn:taudash:04}
\begin{split} 
	\ch\big(\tau'(\ell,\tau)\big)
	& =
	\bigg[\sh^{2}\!\big(\tfrac{\ell'(\ell,\tau)}{2}\big) \ch\big(\tfrac{\ell}{2}\big) - \ch\big(\tfrac{L_1}{2}\big) \ch\big(\tfrac{L_4}{2}\big) - \ch\big(\tfrac{L_2}{2}\big) \ch\big(\tfrac{L_3}{2}\big) \\
	& \quad - \ch\big(\tfrac{\ell'(\ell,\tau)}{2}\big)\Big(\ch\big(\tfrac{L_1}{2}\big) \ch\big(\tfrac{L_3}{2}\big) + \ch\big(\tfrac{L_2}{2}\big) \ch\big(\tfrac{L_4}{2}\big)\Big)\bigg]C_{1,2}(\ell)^{-1/2}C_{3,4}(\ell)^{-1/2}	
\end{split} 
\end{equation} 
with $\sgn(\tau') = - \sgn(\tau)$ and $C_{i,j}(\ell)$ has been defined in \eqref{eqn:constant:Cij}.
\smallskip

\begin{figure}
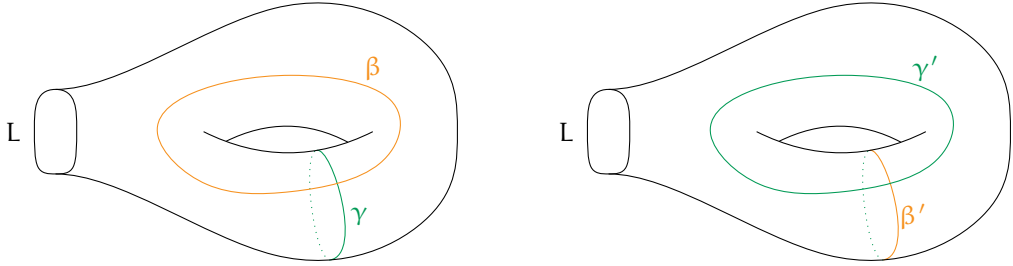

	\centering
	\begin{subfigure}[t]{.44\textwidth}
	\centering
		%
	\end{subfigure}
\caption{Change in the coordinate system of $T$.}
\label{fig:coord:system:11}
\end{figure}

\textsc{One-holed torus.} For a one-holed torus $T$, consider the global coordinates $(\ell,\tau) \in \mathcal{T}_{T}(L)$ relative to $(\mathscr{P},\mathscr{S})$ of Figure~\ref{fig:coord:system:11}. Then the change of coordinate system to $(\mathscr{P}',\mathscr{S}')$ is given by by Equation~\eqref{eqn:Okai:elldash:11} and
\begin{equation}\label{eqn:taudash:11}
	\ch\big(\tfrac{\tau'(\ell,\tau)}{2}\big)
	=
	\ch\big(\tfrac{\ell}{2}\big)
	\sqrt{\frac{
		\ch^2(\tfrac{\tau}{2}) \bigl( \ch(\tfrac{L}{2}) + \ch(\ell) \bigr) - 2 \sh^2(\tfrac{\ell}{2})
		}{
		\ch^2(\tfrac{\tau}{2}) \bigl( \ch(\tfrac{L}{2}) + \ch(\ell) \bigr)
		+ \sh^2(\tfrac{\ell}{2})(\ch(\tfrac{L}{2}) - 1)}}
\end{equation}
with $\sgn(\tau') = - \sgn(\tau)$.
\begin{figure}
	\centering
	\begin{subfigure}[t]{.44\textwidth}
	\centering

		\end{subfigure}
\caption{Change in the coordinate system of $X$.}
\label{fig:coord:system:04}
\end{figure}

\medskip

Using the convergence under the rescaling flow of hyperbolic length and twist parameters to the combinatorial analogues, we can give a new proof of Penner's formulae.

\begin{prop}\label{prop:Penner:formulae} ~

	\smallskip

	\textsc{Sphere with four boundary components.} For a four-holed sphere $X$, consider the global coordinates $(\ell,\tau) \in \mathcal{T}_{X}^{\textup{comb}}(L_1,L_2,L_3,L_4)$ relative to the system of curves $(\mathscr{P},\mathscr{S})$ of Figure~\ref{fig:coord:system:04}. Then the change of coordinate system to $(\mathscr{P}',\mathscr{S}')$ is given by Equation~\eqref{eqn:length:delta:04} and
	\begin{equation}\label{eqn:comb:tau:04}
		|\tau'(\ell,\tau)|
		=
		\frac{1}{2} \bigg|
			2 |\tau|
			+ \ell + M_{1,4}(\ell) + M_{2,3}(\ell)
			- \ell'(\ell,\tau) - M_{1,2}(\ell') - M_{3,4}(\ell')
		\bigg|
	\end{equation}
	with $\sgn(\tau')= - \sgn(\tau)$. We recall here that $M_{i,j}(\ell) = \max\big\{0,L_i - \ell,L_j - \ell,\frac{L_i + L_j - \ell}{2}\big\}$.

	\smallskip
	
	\textsc{Torus with one boundary component.} For a one-holed torus $T$, consider the global coordinates $(\ell,\tau) \in \mathcal{T}_{T}^{\textup{comb}}(L)$ relative to the system of curves $(\mathscr{P},\mathscr{S})$ of Figure~\ref{fig:coord:system:11}. Then the change of coordinate system to $(\mathscr{P}',\mathscr{S}')$ is given by Equation~\eqref{eqn:length:delta:11} and
	\begin{equation}\label{eqn:comb:tau:11}
		|\tau'(\ell,\tau)|
			=
			\left| \ell - \big[ \tfrac{L}{2} - \ell'(\ell,\tau)\big]_{+} \right|
				\end{equation}
	with $\sgn(\tau') = - \sgn(\tau)$.
\end{prop}

\begin{proof}
	The formulae follow from Equations~\eqref{eqn:taudash:04}--\eqref{eqn:taudash:11} by direct computation, together with the convergence of length and twist parameters and using relations of the form
	\[
		\forall a,b \in \RR_{+},
		\qquad
		\lim_{\beta\to\infty} \frac{1}{\beta} \ln{\left( e^{\beta a} + e^{\beta b} \right)}
		=
		\max{\Set{a,b}}.
	\]
The reader can check that the limit of Equations~\eqref{eqn:Okai:elldash:04} and \eqref{eqn:Okai:elldash:11} coincide with Equations~\eqref{eqn:length:delta:04} and \eqref{eqn:length:delta:11}.
\end{proof}

\begin{rem}
	To give a complete description of the action of the mapping class group on combinatorial Fenchel--Nielsen coordinates, one should consider the effect of the change of the coordinate systems described in Figures \ref{fig:coord:system:04} and \ref{fig:coord:system:11} to the twisting numbers at the boundary components (\emph{cf.} Equation~\ref{eqn:twist:boundary}). Although Okai does not consider this, one can easily compute the change on the twisting numbers at the boundary components using hyperbolic trigonometry, and flow such formulae to obtain the analogous results in the combinatorial setting.

	\medskip

	We present here the argument for a one-holed torus $T$, and remark that the same reasoning applies verbatim to a four-holed sphere. In order to talk about twisting number at the boundary, we fix an isotopy class of a hyperbolic metric on $T$, as well as geodesics representatives of curve as in Figure~\ref{fig:one:holed:twist:boundary}. In particular, we have well-defined Fenchel--Nielsen coordinates $(\ell,\tau)$ and $(\ell',\tau')$ relative to $(\mathscr{P},\mathscr{S})$ and $(\mathscr{P}',\mathscr{S}')$ respectively, as in Proposition~\ref{prop:Penner:formulae}, and twisting numbers at the boundary $t$ and $t'$.
	\begin{figure}
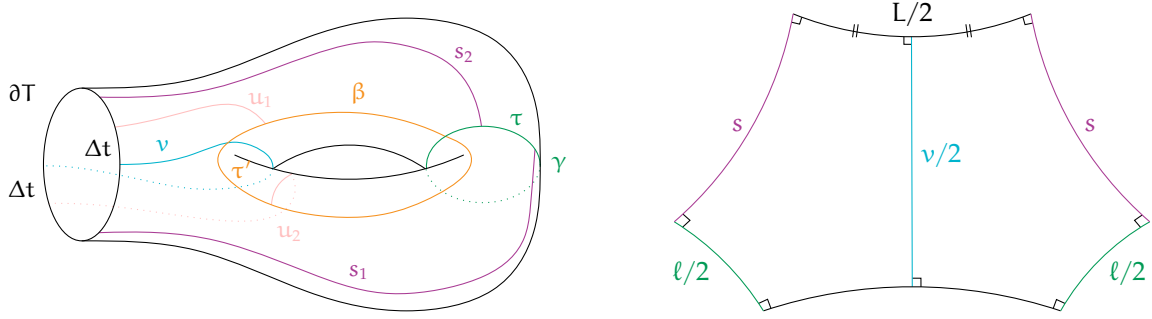

		\centering

		\caption{Twist change at the boundary for $T$.}
		\label{fig:one:holed:twist:boundary}
	\end{figure}

	\medskip

	In order to find the change in the twisting number at the boundary, we need to compute the distance between the seam $s_1$ and the seam $u_1$, which is $t'-t = \frac{L}{4}+\Delta t$. In particular $\Delta t$ is the signed distance between $v$ and $u_1$. First, notice that $\sgn(\Delta t) = - \sgn(\tau')$. On the other hand, from the orthogonal non-convex hyperbolic hexagon $\Delta t \to u_1 \to \tau' \to u_2 \to \Delta t \to v$, we get
	\begin{equation}
		\cosh(\tau') = \sinh^2(\Delta t) \cosh(v)+ \cosh^2(\Delta t).
	\end{equation}
	The orthogonal pentagons on the right-hand side of Figure~\ref{fig:one:holed:twist:boundary} give $\sinh(\tfrac{v}{2}) \sinh(\tfrac{L}{4}) = \cosh(\tfrac{\ell}{2})$, so that
	\begin{equation}
		\cosh(\tau')
		=
		2 \sinh^2(\Delta t) \left(
			\frac{\cosh^2(\tfrac{\ell}{2})}{\sinh^2(\tfrac{L}{4})} + 1
		\right) + 1,
		\qquad
		\sgn(\Delta t) = - \sgn(\tau').
	\end{equation}
	Using the rescaling flow, we obtain the following formula for the combinatorial quantities:
	\begin{equation}
		|\Delta t|
		=
		\frac{|\tau'|}{2} - \frac{1}{2}\left[ \ell-\frac{L}{2} \right]_+ ,
		\qquad
		\sgn(\Delta t) = - \sgn(\tau').
	\end{equation}
\end{rem}

As in the case of measure foliations, completion of the set of multicurves, we can conclude a piecewise linear structure on the combinatorial Teichmüller space.

\begin{cor}\label{cor:PL}
	Consider a bordered surface $\Sigma$ of type $(g,n)$ and fix two seamed pants decompositions $(\mathscr{P},\mathscr{S})$ and $(\mathscr{P}',\mathscr{S}')$, defining combinatorial Fenchel--Nielsen coordinates $\Phi,\Phi' \colon \mathcal{T}_{\Sigma}^{\textup{comb}} \to (\RR_{+} \times \RR)^{3g -3 + n}$. Then the change of coordinate $\Phi' \circ \Phi^{-1}$ between open subsets of $(\RR_{+} \times \RR)^{3g - 3 + n}$ is piecewise linear. In particular, the combinatorial Teichm\"uller space $\mathcal{T}_{\Sigma}^{\textup{comb}}$ is endowed with a canonical piecewise linear structure.
	\hfill $\blacksquare$
\end{cor}

\begin{rem}
	The above transformations are not continuous on the whole of $\RR_{+}\times\RR$. The locus of discontinuity actually identifies a subset of the non-admissible twists, \emph{i.e.} the creation of saddle connections in the measured foliation perspective. For instance, the plot of $\ell'$ and $\tau'$ on $\RR_{+} \times \RR \cong \mathcal{T}^{\textup{comb}}_{T}(L)$ is illustrated in Figure~\ref{fig:change:pop:11}. Notice that along the line
	\begin{equation}
		\mathfrak{l} = \left[\tfrac{L}{2},+\infty\right) \times \set{0},
	\end{equation}
	the function $\ell'$ is identically zero, while $\tau'$ has a discontinuity. But $\mathfrak{l}$ is not in the image of $\mathcal{T}_T^{\textup{comb}}(L)$ under the map $\Phi$ of Theorem~\ref{thm:FN:coordinates}. Thus, having $\ell' = 0$ and $\tau'$ discontinuous along $\mathfrak{l}$ is not contradictory.
	\begin{figure}
		\centering
		\begin{subfigure}[t]{.48\textwidth}
			\centering
			\includegraphics[width=.9\textwidth]{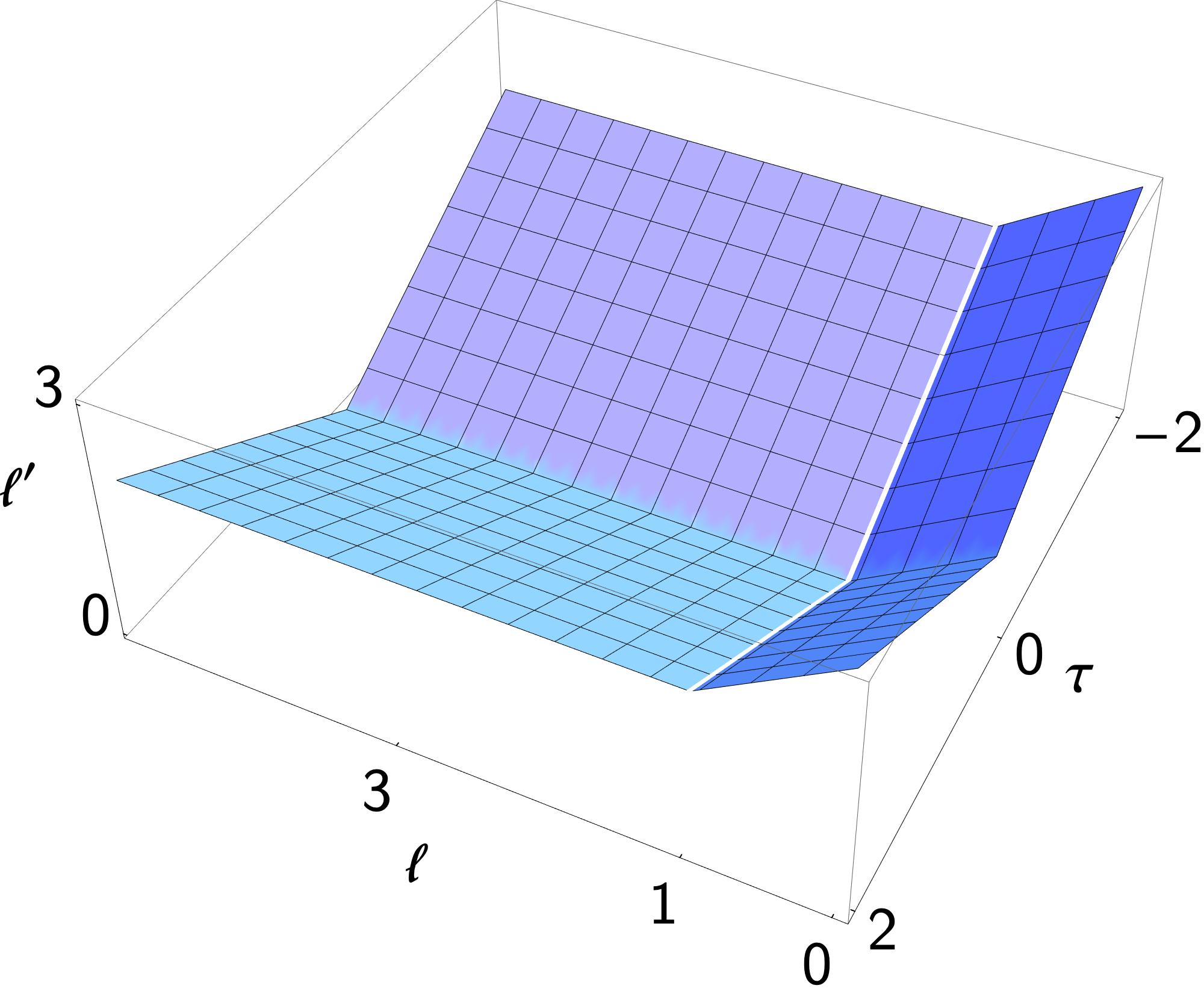}
		\end{subfigure}
		\begin{subfigure}[t]{.48\textwidth}
			\centering
			\includegraphics[width=.9\textwidth]{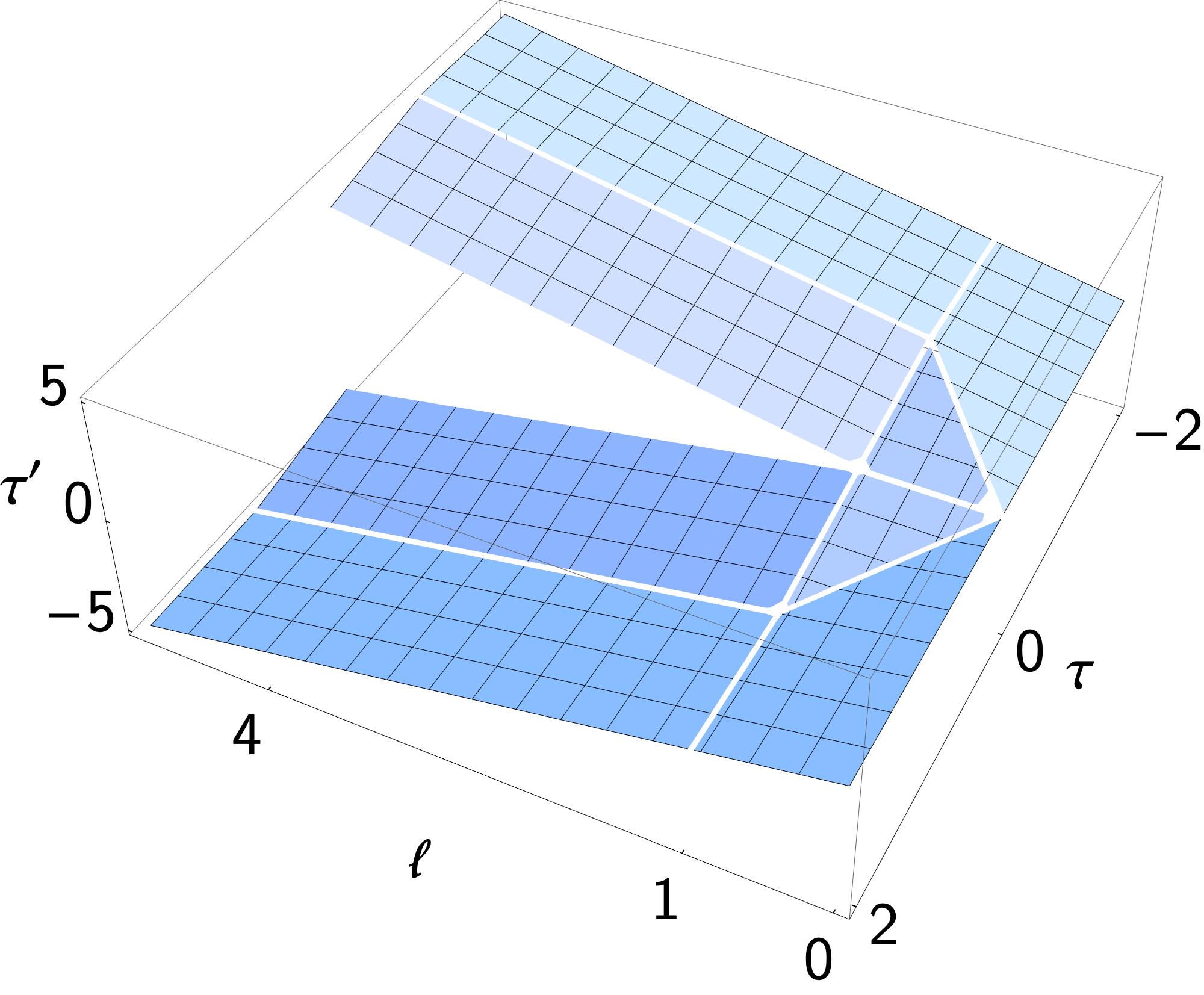}
		\end{subfigure}
		\caption{The graphs of $\ell'(\ell,\tau)$ and $\tau'(\ell,\tau)$, with $L = 2$.}
		\label{fig:change:pop:11}
	\end{figure}
\end{rem}

\subsection{Geometric recursion in the flow}
\label{subsec:GR:flow}

The geometric recursion in the hyperbolic and in the combinatorial setting produces functions respectively on the moduli spaces $\mathcal{M}_{\Sigma}$ and $\mathcal{M}_{\Sigma}^{\textup{comb}}$. The main result of this section is that the flow in the $\beta \to \infty$ limit takes the hyperbolic GR to the combinatorial GR.

\subsubsection{Rescaling initial data}

Before considering the behaviour of geometric recursion amplitudes in the flow, we discuss the rescaling of initial data.

\begin{defn} \label{defn:unif:strong:adm}
	Let $(A_{\beta},B_{\beta},C_{\beta})_{\beta \geq 1}$ be a family of triples of mea\-su\-ra\-ble func\-tions on $\RR_{+}^3$. We say it is \emph{uniformly (strongly) admissible} if the constants in the (strong) admissibility can be chosen to be independent of $\beta \geq 1$.
\end{defn}

We remark that here, following the remark of Section~\ref{subsec:D:from:C}, the triple $(A,B,C)$ is completed by a natural choice of $D$. We introduce the rescaling operator, acting on a function $\phi \colon \RR_{+}^k \to \RR$ by
\begin{equation}
	\rho_{\beta *}\phi(L) = \phi(\beta L).
\end{equation}
We notice two basic properties of this rescaling. Firstly, limits of uniformly admissible rescaled initial data (if they exist) are automatically admissible.

\begin{lem}\label{lem:family:ABC:beta:uniform}
	Let $(A_{\beta},B_{\beta},C_{\beta})_{\beta \geq 1}$ be initial data such that $\rho_{\beta*}(A_{\beta},B_{\beta},C_{\beta})$ is uniformly admissible and converges to a pointwise limit $(\hat{A},\hat{B},\hat{C})$. Then $(\hat{A},\hat{B},\hat{C})$ is admissible.
\end{lem}

\begin{proof}
	This is clear by taking the $\beta \to \infty$ limit in the inequalities specified by the uniform admissibility.
\end{proof}

Secondly, if we rescale initial data that do not depend on $\beta$, the limit (if it exists) must be supported on small pairs of pants.

\begin{lem}\label{lem:family:ABC:uniform}
	If $(A,B,C)$ is admissible, then $\rho_{\beta*}(A,B,C)$ is uniformly admissible, with same bounding constants. Besides, if $\rho_{\beta*}(B,C)$ has a pointwise limit $(\hat{B},\hat{C})$ when $\beta \to \infty$, then 
	\begin{equation}
		\ell > L_1 + L_2
		\quad \Longrightarrow \quad
		\hat{B}(L_1,L_2,\ell) = 0,
		\qquad \text{and} \qquad
		\ell + \ell' > L_1
		\quad \Longrightarrow \quad
		\hat{C}(L_1,\ell,\ell') = 0.
	\end{equation}
\end{lem}

\begin{proof}
	The bound on $\rho_{\beta*}A$ is clear. For the bound on $B$, we write
	\[
	\begin{aligned}
		&
		\sup_{L_1,L_2,\ell \geq \epsilon} \big|B(\beta L_1,\beta L_2,\beta \ell)\big|\,\big(1 + [\ell - L_1 - L_2]_{+}\big)^{s} \\
		\leq &
		\sup_{L_1,L_2,\ell \geq \beta\epsilon} \big|B(L_1,L_2,\ell)\big|\, \big(1 + [\ell - L_1 - L_2]_{+}\big)^{s} \bigg(\frac{1 + [(\ell - L_1 - L_2)/\beta]_{+}}{1 + [\ell - L_1 - L_2]_{+}}\bigg)^{s}.
	\end{aligned}
	\]
	Since $\beta \geq 1$, we can bound the last expression by the supremum over $L_1,L_2,\ell \geq \epsilon$. We also observe that $t \mapsto \frac{1 + [t/\beta]_{+}}{1 + [t]_{+}}$ is equal to $1$ for $t \leq 0$ and is decreasing for $t > 0$, hence it is uniformly bounded by $1$. Using the initial bound for $B$ we get the desired bound for $\rho_{\beta*}B$. The bound for $\rho_{\beta*}C$ is proved similarly.

	\medskip

	To establish the vanishing property for $\hat{B}$, we fix $(L_1,L_2,\ell)$ such that $\ell - (L_1 + L_2) \geq \epsilon$ and $\min\{L_1,L_2,\ell\} \geq \epsilon> 0$ for some $\epsilon > 0$. Specialising the admissibility condition for $B$ at $s = 1$ and rescaled lengths, we have
	\[
		\big|B(\beta L_1,\beta L_2,\beta \ell)\big| \leq \frac{M_{\epsilon,1}}{1 + \beta \epsilon}.
	\]
	Taking the limit $\beta \to \infty$ yields the claim. The vanishing property of $\hat{C}$ is proved similarly.
\end{proof}

\begin{rem}
	If $(A,B,C)$ are continuous functions on $\RR_{+}^3$, one easily checks that $(\rho_{\beta*}A,\rho_{\beta*}B,\rho_{\beta*}C)_{\beta \geq 1}$ forms an equicontinuous family. By Arzela--Ascoli theorem, for any fixed compact, it must admit uniformly converging subsequences, and the vanishing properties in Lemma~\ref{lem:family:ABC:uniform} must hold for any limit point. Therefore, the assumption that $\rho_{\beta*}(A,B,C)$ converges is rather weak.
\end{rem}

\subsubsection{Rescaling geometric recursion amplitudes}

\begin{thm} \label{thm:rescale:GR}
	Let $(A_{\beta},B_{\beta},C_{\beta})_{\beta \geq 1}$ be initial data such that $\rho_{\beta*}(A_{\beta},B_{\beta},C_{\beta})$ is uniformly admissible and converges uniformly on any compact to a limit $(\hat{A},\hat{B},\hat{C})$. Let us denote by $\Omega_{\Sigma;\beta}$ the result of the ${\rm Mes}(\mathcal{T}_{\Sigma})$-valued GR with initial data $(A_{\beta},B_{\beta},C_{\beta})$, and by $\widehat{\Xi}_{\Sigma}$ the result of the ${\rm Mes}(\mathcal{T}_{\Sigma}^{\textup{comb}})$-valued GR with initial data $(\hat{A},\hat{B},\hat{C})$. We have for any bordered surface $\Sigma$ and $\sigma \in \mathcal{T}_{\Sigma}$
	\begin{equation}
		\lim_{\beta \to \infty} \Omega_{\Sigma;\beta}(\sigma^{\beta})
		=
		\widehat{\Xi}_{\Sigma}\big(\sp(\sigma)\big)
	\end{equation}
	and the convergence is uniform for $\sigma$ in any compact of $\mathcal{T}_{\Sigma}$. Besides, there exists $t \geq 0$ depending only on the topology of $\Sigma$, such that, for any $\epsilon > 0$ there exists $M_{\epsilon} > 0$ for which we have, for any $\beta \geq 1$, any $\sigma \in \mathcal{T}_{\Sigma}$ such that $\sys_{\sp(\sigma)} \geq \epsilon$, 
	\begin{equation}
		\big|\Omega_{\Sigma;\beta}(\sigma^{\beta})\big| \leq M_{\epsilon} \prod_{b \in \pi_0(\partial \Sigma)} \big( 1 + \ell_{\sp(\sigma)}(b) \big)^{t}
	\end{equation}
	and the same inequality holds for the limit $\widehat{\Xi}_{\Sigma}(\sp(\sigma))$.
\end{thm}

\begin{rem}
	If we specialise $(A_{\beta},B_{\beta},C_{\beta})$ to $(A^{\textup{M}},B^{\textup{M}},C^{\textup{M}})$ appearing in Theorem~\ref{thm:Mirz:McS:hyp}, this gives another proof of the combinatorial Mirzakhani--McShane of Theorem~\ref{thm:Mirz:McS:comb} or its equivalent form Corollary~\ref{GR:comb:WK}. Indeed, the hyperbolic Mirzakhani--McShane identities immediately show that $\Omega_{\Sigma}^{\textup{M}}(\sigma^{\beta}) = 1$ for any $\sigma \in \mathcal{T}_{\Sigma}$ and $\beta \geq 1$, and the convergence
	\begin{equation}
		\lim_{\beta \to \infty} \frac{F(\beta x)}{\beta} = \lim_{\beta \to \infty} \frac{2\ln(1 + e^{\beta x/2})}{\beta} = [x]_{+},
	\end{equation}
	which is uniform on any compact of $\{-\infty\} \cup \RR$, implies that
	\begin{equation}
		\lim_{\beta \to \infty} \rho_{\beta*}(A^{\textup{M}},B^{\textup{M}},C^{\textup{M}})
		=
		(A^{\textup{K}},B^{\textup{K}},C^{\textup{K}})
	\end{equation}
	uniformly on any compact of $\RR_{+}^3$. Thus, the $\beta \to \infty$ limit of $\Omega_{\Sigma}^{\textup{M}}(\sigma^{\beta})$ -- which must be the constant function $1$ -- converges to the function $\Xi_{\Sigma}^{\textup{K}}(\sp(\sigma))$ satisfying the geometric recursion for the initial data $(A^{\textup{K}},B^{\textup{K}},C^{\textup{K}})$.
\end{rem}

\begin{proof}[Proof of Theorem~\ref{thm:rescale:GR}]
	It is enough to prove the result for connected surfaces, and we proceed by induction on the Euler characteristic. Throughout the proof, we denote $\GG = \sp(\sigma)$. We recall that for any boundary component $b$, we have $\ell_{\sigma^{\beta}}(b) = \beta\ell_{\GG}(b) = \beta\ell_{\sigma}(b)$. For a pair of pants $P$, we have
	\[
		\Omega_{P;\beta}(\sigma^{\beta})
		=
		A_{\beta}(\vec{\ell}_{\sigma^{\beta}}(\partial P)) = A_{\beta}(\beta\vec{\ell}_{\sigma}(\partial P)) = A_{\beta}(\beta\vec{\ell}_{\GG}(\partial P)),
	\]
	which by assumption converges uniformly on any compact to $\hat{A}(\vec{\ell}_{\GG}(\partial P)) = \widehat{\Xi}_{P}(\GG)$, and is bounded by a constant independent of $\beta$ on any $\epsilon$-thick part of $\mathcal{T}_{P}^{\textup{comb}}$.

	\smallskip

	For a torus with one boundary component $T$, we have by the convention in Section~\ref{subsec:D:from:C}
	\begin{equation}\label{eqn:conv:torus}
		\Omega_{T;\beta}(\sigma^{\beta})
		=
		\sum_{\gamma \in S_{T}^{\circ}} C_{\beta}\big(\beta\ell_{\GG}(\partial T),\ell_{\sigma^{\beta}}(\gamma),\ell_{\sigma^{\beta}}(\gamma)\big).
	\end{equation}
	Let $K$ be a compact subset of  $\mathcal{T}_{\Sigma}^{\textup{comb}}$ and $\epsilon$ a lower bound of the systole on $K$. Let $\gamma \in S_{T}^{\circ}$ and $\GG \in K$. Since the function $\ell_{\GG}(\gamma)$ is bounded from below and from above on $K$, Proposition~\ref{prop:rescale:length} implies that $\beta^{-1}\ell_{\sigma^{\beta}}(\gamma)$ converges to $\ell_{\GG}(\gamma)$ uniformly for $\sigma \in \sp^{-1}(K)$. By uniform convergence of $\rho_{\beta*}C_{\beta}$ on that compact, we deduce that
	\[
		\lim_{\beta \to \infty} C_{\beta}\big(\beta\ell_{\GG}(\partial T),\ell_{\sigma^{\beta}}(\gamma),\ell_{\sigma^{\beta}}(\gamma)\big) = \hat{C}\big(\ell_{\GG}(T),\ell_{\GG}(\gamma),\ell_{\GG}(\gamma)\big)
	\]
	uniformly for $\sigma \in K$. Next, we would like to bound each term in \eqref{eqn:conv:torus} by a summable (over $\gamma$) quantity depending only on $K$ and not on $\beta$. If this holds, we can conclude by the (Banach-valued) dominated convergence theorem that $\Omega_{T;\beta}(\sigma^{\beta})$ converges to $\widehat{\Xi}_{T}(\GG)$ when $\beta \to \infty$, uniformly for $\sigma \in K$.

	\smallskip

	To prove the bound, we notice that by uniform admissibility we have for any $s > 0$
	\[
	\begin{aligned}
		\big|C_{\beta}\big(\beta\ell_{\GG}(\partial T),\ell_{\sigma^{\beta}}(\gamma),\ell_{\sigma^{\beta}}(\gamma)\big)\big|
		& \leq
		\frac{M_{\epsilon,s}}{\big(1 + [2 \beta^{-1} \ell_{\sigma^{\beta}}(\gamma) - \ell_{\GG}(\partial T)]_{+}\big)^{s}} \leq
		\frac{M_{\epsilon,s}}{\big(1 + [2\ell_{\GG}(\gamma) - \ell_{\GG}(\partial T)]_{+}\big)^{s}},
	\end{aligned} 
	\]
	where the second inequality used the lower bound $\beta \ell_{\GG}(\gamma) \le \ell_{\sigma^{\beta}}(\gamma)$ of Proposition~\ref{lem:length:bound}. Since the number of small pairs of pants is bounded by the number of oriented edges, there are at most $6$ curves $\gamma \in S_{T}^{\circ}$ for which the denominator is equal to $1$. Hence
	\[
		\sum_{\gamma \in S_{T}^{\circ}} \frac{1}{\big(1 + [2\ell_{\GG}(\gamma) - \ell_{\GG}(\partial T)]_{+}\big)^{s}}
		\leq 6 + \sum_{L > \ell_{\GG}(\partial T)/2} \frac{\#\Set{\gamma \in S_{T}^{\circ}\,\,|\,\,\ell_{\GG}(\gamma) \leq L + 1}}{\big(1 + 2L - \ell_{\GG}(\partial T)\big)^s}.
	\]
	Invoking Proposition~\ref{prop:poly:growth}, the numerator is bounded by $m_{\epsilon}(L + 1)^{2}$ for some constant $m_{\epsilon}$ depending on $\epsilon$ only. Therefore
	\[
		\big|\Omega_{T;\beta}(\sigma^{\beta})\big| \leq M_{\epsilon,s}\bigg(6 + m_{\epsilon} \sum_{L \geq 0} \frac{(1 + \ell_{\GG}(\partial T)/2 + L)^{2}}{\big(1 + 2L\big)^{s}}\bigg) .
	\]
	By choosing  $s = 4$ we find that $|\Omega_{T;\beta}(\sigma^{\beta})|$ is bounded by a polynomial of degree $2$ in $\ell_{\GG}(\partial T)$, whose coefficients are independent of $\beta$ but may depend on $\epsilon$. This proves the theorem for $\Sigma = T$.	
	
	\smallskip

	Now let $\Sigma$ be a connected bordered surface with $\chi_{\Sigma} < -1$, and assume the thesis for all $\Sigma'$ such that $\chi_{\Sigma'} > \chi_{\Sigma}$. Let $K$ be a compact subset of $\mathcal{T}_{\Sigma}^{\textup{comb}}$ and $\epsilon$ a lower bound on the systole on $K$. The geometric recursion gives
	\[
		\Omega_{\Sigma;\beta}(\sigma^{\beta})
		=
		\sum_{[P] \in \mathcal{P}_{\Sigma}} X_{P;\beta}(\vec{\ell}_{\sigma^{\beta}}(\partial P))\,\Omega_{\Sigma - P;\beta}(\sigma^{\beta}|_{\Sigma - P}).
	\]
	For each $[P] \in \mathcal{P}_{\Sigma}$, we can repeat the previous arguments to show that $X_{P}(\vec{\ell}_{\sigma^{\beta}}(\partial P))$ is converging to $\hat{X}_{P;\beta}(\vec{\ell}_{\GG}(\partial P))$ uniformly for $\GG \in K$. If the initial data were the one of Mirzakhani \eqref{eqn:Mirz:init:data}, the factor $\Omega_{\Sigma - P;\beta}$ would be the constant function $1$. Then, we could finish the proof very similarly to the case of $\Sigma = T$, and this would prove a second time the combinatorial Mirzakhani--McShane identities (Theorem~\ref{thm:Mirz:McS:comb}). However, in the general case, we only know by induction hypothesis that $\Omega_{\Sigma - P;\beta}(\tilde{\sigma}^{\beta})$ converges to $\widehat{\Xi}_{\Sigma - P}(\tilde{\GG})$ uniformly for $\tilde{\GG} = \sp(\tilde{\sigma})$ in any compact of $\mathcal{T}_{\Sigma - P}^{\textup{comb}}$. To employ this information, we shall compare the two hyperbolic structures $\sigma^{\beta}|_{\Sigma - P}$ and $\tilde{\sigma}^{\beta} \coloneqq \sp^{-1}(\beta \GG|_{\Sigma - P})$, by means of their lengths functions. Since the combinatorial cutting does not decrease the systole, Proposition~\ref{lem:length:bound} on $\Sigma$ and $\Sigma - P$ respectively imply for $\beta \geq \beta_{\epsilon}$ and $\gamma \in S_{\Sigma - P}^{\circ}$
	\[
	\begin{split}
		\ell_{\GG}(\gamma)
		& \leq
			\frac{\ell_{\sigma^{\beta}}(\gamma)}{\beta} \leq \ell_{\GG}(\gamma)\bigg(1 + \frac{\kappa_{\epsilon}}{\beta}\bigg), \\
		\ell_{\GG|_{\Sigma - P}}(\gamma)
		& \leq
			\frac{\ell_{\tilde{\sigma}^{\beta}}(\gamma)}{\beta} \leq \ell_{\GG|_{\Sigma - P}}(\gamma)\bigg(1 + \frac{\kappa_{\epsilon}}{\beta}\bigg),
	\end{split}
	\]
	where we can use for $\kappa_{\epsilon},\beta_{\epsilon}$ the constants provided by Proposition~\ref{lem:length:bound} for $\Sigma$. The combinatorial and hyperbolic cutting procedures are such that
	\[
		\forall \gamma \in S_{\Sigma - P}^{\circ}
		\qquad
		\ell_{\GG|_{\Sigma - P}}(\gamma) = \ell_{\GG}(\gamma)
		\quad \text{and} \quad
		\ell_{\sigma^{\beta}|_{\Sigma - P}}(\gamma) = \ell_{\sigma^{\beta}}(\gamma).
	\]
	Therefore
	\[
		\forall \gamma \in S_{\Sigma - P}^{\circ}
		\qquad
		\bigg(1 + \frac{\kappa_{\epsilon}}{\beta}\bigg)^{-1} \ell_{\tilde{\sigma}^{\beta}}(\gamma) \leq \ell_{\sigma^{\beta}|_{\Sigma - P}}(\gamma) \leq \ell_{\tilde{\sigma}^{\beta}}(\gamma)\bigg(1 + \frac{\kappa_{\epsilon}}{\beta}\bigg).
	\]
	Further, for any component of $\partial P \cap \Sigma^{\circ}$ and $\de\Sigma \cap \de(\Sigma - P)$, we have exactly $\ell_{\sigma^{\beta}|_{\Sigma - P}}(\gamma) = \ell_{\tilde{\sigma}^{\beta}}(\gamma)$, so that the above bounds extend to all $\gamma \in S_{\Sigma - P}$. From our description of the topology of combinatorial Teichm\"uller spaces (Theorem~\ref{thm:length:functional}), we deduce that $\beta^{-1}\sp(\sigma^{\beta}|_{\Sigma - P})$ remains in a compact of $\mathcal{T}_{\Sigma - P}^{\textup{comb}}$ independent of $\beta$ when $\GG \in K$, and that it converges to $\GG|_{\Sigma - P}$ uniformly for $\GG \in K$. By the induction hypothesis, we have
	\[
		\lim_{\beta \to \infty} \Omega_{\Sigma - P;\beta}(\sigma^{\beta}|_{\Sigma - P}) = \widehat{\Xi}_{\Sigma - P}(\GG|_{\Sigma - P})
	\]
	and the convergence is uniform for $\GG \in K$. Supplemented with a summable (over $[P] \in \mathcal{P}_{\Sigma}$) bound whose derivation is similar to the case $\Sigma = T$ and therefore omitted, this proves the theorem for $\Sigma$, and thus in full generality by induction. 
\end{proof}

\subsection{Topological recursion in the flow}
\label{subsec:TR:flow}

By Theorem~\ref{thm:Mon:Do}, we know that the Weil--Petersson measure flows in a pointwise sense, in the $\beta \to \infty$ limit ,to the Kontsevich measure. This implies that the Weil--Petersson volumes flow (by rescaling the boundary lengths by $\beta$) to the Kontsevich volumes. More generally, we are able to prove an analogue of Theorem~\ref{thm:rescale:GR} after integration over the moduli spaces.

\begin{thm} \label{thm:rescale:TR}
	Let $(A_{\beta},B_{\beta},C_{\beta})_{\beta \geq 1}$ be initial data such that $\rho_{\beta*}(A_{\beta},B_{\beta},C_{\beta})$ is uniformly strongly admissible and converges uniformly on any subset of the form $(0,M]^3 \subset \RR_{+}^{3}$ to a limit $(\hat{A},\hat{B},\hat{C})$. Then, $(\hat{A},\hat{B},\hat{C})$ is strongly admissible and
	\begin{equation}
		\lim_{\beta \to \infty} \frac{V\Omega_{g,n;\beta}(\beta L)}{\beta^{6g - 6 + 2n}}
		=
		V\widehat{\Xi}_{g,n}(L),
	\end{equation}
	with uniform convergence for $L$ in any subset of the form $(0,M]^{n} \subset \RR_{+}^n$.
\end{thm}

\begin{proof}
	We prove the result by induction on $2g - 2 + n > 0$. For $(g,n) = (0,3)$, we have $V\Omega_{0,3;\beta}(\beta L) = A_{\beta}(\beta L)$, which converges uniformly on any compact of $\RR^{3}$ intersected with $\RR_{+}^{3}$ to $\hat{A}(L) = V\widehat{\Xi}_{0,3}(L)$. For $(g,n) = (1,1)$ we have
	\[
		\frac{V\Omega_{1,1;\beta}(\beta L_{1})}{\beta^2}
		=
		\int_{\RR_{+}}  \, C_{\beta}(\beta L_{1},\ell,\ell) \frac{\ell}{\beta} \frac{\dd\ell}{\beta}
		=
		\int_{\RR_{+}} \,  \, C_{\beta}(\beta L_{1},\beta\ell,\beta\ell) \ell\dd\ell .
	\]
	Note $C_{\beta}(\beta L_{1},\beta\ell,\beta\ell)$ converges uniformly on any $(0,M]^2$ to $\hat{C}(L_{1},\ell,\ell)$. Strong uniform admissibility means we can bound the integrand by an integrable function independent of $\beta$. Moreover, the uniformity of the convergence around zero implies that we can exchange the integral and the limit, so that
	\[
		\lim_{\beta\to\infty} \frac{V\Omega_{1,1;\beta}(\beta L_{1})}{\beta^2}
		=
		\lim_{\beta\to\infty} \int_{\RR_{+}} \, C_{\beta}(\beta L_{1},\beta\ell,\beta\ell) \ell  \dd\ell 
		=
		\int_{\RR_{+}}  \, \hat{C}(L_{1},\ell,\ell)\ell  \dd\ell
		=
		V\widehat{\Xi}_{1,1}(L_{1}),
	\]
	and the convergence is uniform on any $(0,M]$. This proves the two base cases.

	\smallskip

	The general argument follows along the same lines as the $(g,n)=(1,1)$ case. Assume the result for $(g',n')$ such that $2g' - 2 + n' < 2g - 2 + n$. The topological recursion for $V\Omega_{g,n;\beta}(\beta L)$ yields
	\[
		\begin{split}
		\frac{V\Omega_{g,n;\beta}(\beta L_1,\ldots,\beta L_n)}{\beta^{6g-6+2n}}
		=
		\sum_{m = 2}^n \int_{\mathbb{R}_{+}} \!\! \, B_{\beta}(\beta L_1,\beta L_m,\ell) \frac{V\Omega_{g,n - 1;\beta}(\ell,\beta L_2,\ldots,\widehat{\beta L_m},\ldots,\beta L_n)}{\beta^{6g-6+2(n-1)}}      				\frac{\ell}{\beta} \frac{\dd\ell}{\beta}
		& \\
		+ \frac{1}{2} \int_{\mathbb{R}_{+}^2} \!\! \, C_{\beta}(\beta L_1,\ell,\ell')
		\bigg(
			\frac{V\Omega_{g - 1,n + 1;\beta}(\ell,\ell',\beta L_2,\ldots,\beta L_n)}{\beta^{6(g-1)-6+2(n+1)}} 
		& \\
		+ \sum_{\substack{J \sqcup J' = \{L_2,\ldots,L_m\} \\ h + h ' = g}}
		 	\frac{V\Omega_{h,1 + \#J;\beta}(\ell,\beta J)}{\beta^{6h-6+2(1+\#J)}}
		 	\frac{V\Omega_{h',1 + \#J';\beta}(\ell',\beta J')}{\beta^{6h'-6+2(1+\#J')}} 
		 \bigg) \frac{\ell \ell'}{\beta^2}  \frac{\dd\ell \, \dd\ell'}{\beta^2} . &
		\end{split}
	\]
	Now rescaling the integration variables, we get
	\[
	\begin{split}
		\frac{V\Omega_{g,n;\beta}(\beta L_1,\ldots,\beta L_n)}{\beta^{6g-6+2n}}
		=
		\sum_{m = 2}^n \int_{\mathbb{R}_{+}} \!\!\!\!  \, B_{\beta}(\beta L_1,\beta L_m,\beta\ell) \frac{V\Omega_{g,n - 1;\beta}(\beta\ell,\beta L_2,\ldots,\widehat{\beta L_m},\ldots,\beta L_n)}{\beta^{6g-6+2(n-1)}}\ell  \dd \ell 
		& \\
		+ \frac{1}{2} \int_{\mathbb{R}_{+}^2} \!\!\!\! C_{\beta}(\beta L_1,\beta\ell,\beta\ell')
		\bigg(
			\frac{V\Omega_{g - 1,n + 1;\beta}(\beta\ell,\beta\ell',\beta L_2,\ldots,\beta L_n)}{\beta^{6(g-1)-6+2(n+1)}}
		& \\
		 	+ \sum_{\substack{J \sqcup J' = \{L_2,\ldots,L_m\} \\ h + h ' = g}} \!\!\! \frac{V\Omega_{h,1 + \#J;\beta}(\beta\ell,\beta J)}{\beta^{6h-6+2(1+\#J)}} \frac{V\Omega_{h',1 + \#J';\beta}(\beta\ell',\beta J')}{\beta^{6h'-6+2(1+\#J')}}
		\bigg) \ell\ell'  \dd \ell \, \dd \ell'. &
	\end{split}
	\]
	By our induction assumption, all of the $\beta^{-(6g'-6+2n')}V\Omega_{g',n';\beta}$'s appearing on the right-hand side of the equation converge uniformly on any $(0,M]^{n}$. For $X_{\beta} \in \set{B_{\beta},C_{\beta}}$, we assumed that $X_{\beta}(\beta L)$ converges to $\hat{X}(L)$ uniformly on any $(0,M]^3$ and moreover that $X_{\beta}(\beta L)$ is uniformly admissible. Therefore using the strong uniform admissibility we can bound the integrals around infinity by integrable functions independently of $\beta$ and that can be uniformly chosen on compact sets of $L$. The uniformity in compacts around zero then implies that we can interchange the integral and the limit, which again by our induction assumption reproduces the topological recursion for $V\hat{\Xi}_{g,n}(L)$ uniformly on any $(0,M]^n$.
\end{proof}

\begin{cor}
	If $\rho_{\beta*}(A,B,C)$ is uniformly strongly admissible and converges uniformly on any compact of $\RR^{3}$ intersected with $\RR_{+}^{3}$ to a limit $(\hat{A},\hat{B},\hat{C})$, then $(\hat{A},\hat{B},\hat{C})$ is strongly admissible and
	\begin{equation}
		\lim_{\beta \to \infty} \frac{V\Omega_{g,n}(\beta L)}{\beta^{6g - 6 + 2n}}
		=
		V\widehat{\Xi}_{g,n}(L).
	\end{equation}
	In particular, if $V\Omega_{g,n}(L)$ is a polynomial, then $V\Omega_{g,n}(L)$ is of degree at most $6g-6+2n$ and $V\widehat{\Xi}_{g,n}(L)$ is the homogeneous component of degree $6g - 6 + 2n$.
\end{cor}

\newpage
\section{Enumerative geometry of multicurves}
\label{sec:enumeration:curves}

In \cite{Mir08growth}, Mirzakhani defines a function on the hyperbolic moduli space with cusps that counts the number of multicurves with length bounded by some parameter $t \in \RR_{+}$. Moreover, she calculates the integral of this function over the moduli space as a sum over stable graphs. She then considers the behaviour of this integral as the parameter tends to infinity, and in \cite{Mir08earth} she shows that the limiting value gives the Masur--Veech volume of the moduli space of quadratic differentials.

\medskip

In this section, after explaining how to parametrise multicurves via embedded ribbon graphs, we briefly recall Mirzakhani's work and then give the natural extension to the combinatorial setting. We then describe twisted geometric recursion amplitudes, and by considering a special family of such functions we illustrate how to recover additive statistics of multicurves by taking Laplace transform. This gives an interpretation of the Masur--Veech polynomials in \cite{ABCDGLW19} as a Laplace transform of averages of additive statistics. Moreover, this suggest a new family of polynomials, calculated by topological recursion, that calculate averages of additive statistics of hyperbolic length.

\subsection{Parametrising multicurves}
\label{subsec:param:mult:curv}

Consider an open cell $\mathfrak{Z}_{\Sigma,G}$ on the combinatorial Teichm\"uller space $\mathcal{T}_{\Sigma}^{\textup{comb}}$, that is a trivalent ribbon graph $G$ together with a marking $f\colon \Sigma \to |G|$. Assigning to each edge $e \in E_{G}$ the number of times a non-backtracking representative on $G$ of $c \in M_{\Sigma}$ passes through it, we obtain a map
\begin{equation}
	\mathfrak{m}_{\mathfrak{Z}_{\Sigma,G}} \colon M_{\Sigma} \longrightarrow \NN^{E_{G}}.
\end{equation}
We show that in fact this map gives a parametrisation of $M_{\Sigma}$. Let us first introduce some notation to describe its image.

\begin{defn}
	Given a trivalent ribbon graph $G$, a \emph{corner} is an ordered triple $\Delta = (e,e',e'')$ where $e,e',e''$ are edges incident to a vertex in the cyclic order. Equivalently, a corner consists of a vertex $v$ together with the choice of an incident edge $e$. We say that a corner belongs to a face $\mathfrak{f} \in F_{G}$ if $e'$ and $e''$ are edges around that face. We denote by $\mathfrak{C}(\mathfrak{f})$ the set of corners belonging to $\mathfrak{f}$, and by $\mathfrak{C}_{G}$ the set of all corners. If we have an assignment of real numbers $(x_e)_{e \in E_{G}}$ and $\Delta = (e,e',e'')$ is a corner, we set $x_{\Delta} = x_{e''} + x_{e'} - x_{e}$.
\end{defn}

\begin{lem}\label{lem:parametrisation:multicurves}
	The map $\mathfrak{m}_{\mathfrak{Z}_{\Sigma,G}}$ is a bijection between $M_{\Sigma}$ and the set
	\begin{equation}
	\begin{split}
		Z_{G} & =
			\Set{ m \in \NN^{E_{G}} |
				{\begin{array}{@{}c@{}}
					\forall \Delta \in \mathfrak{C}_{G},
					\;
					m_{\Delta} \in 2\NN \\
					\forall \mathfrak{f} \in F_{G},
					\;
					\exists \Delta \in \mathfrak{C}(\mathfrak{f}) \text{ s.t. } m_{\Delta} = 0
				\end{array}}
			}.
	\end{split}
	\end{equation}
\end{lem}

\begin{figure}
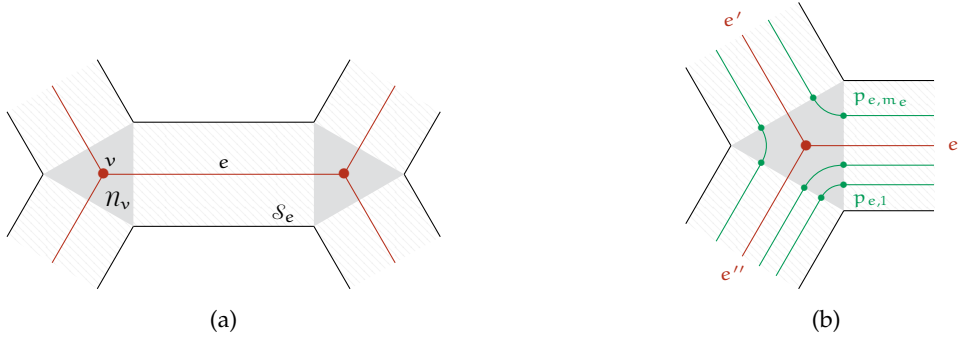

	\centering
	\begin{subfigure}[t]{.48\textwidth}
		\centering

		\caption{}
		\label{fig:decomposition:SN:a}
	\end{subfigure}
	\begin{subfigure}[t]{.48\textwidth}
		\centering
		%
		\caption{}
		\label{fig:decomposition:SN:b}
	\end{subfigure}
	\caption{Decomposition of $\Sigma$ into strips $\mathscr{S}_e$ around edges and triangular neighbourhoods $\mathscr{N}_v$ around vertices.}
	\label{fig:decomposition:SN}
\end{figure}

\begin{proof}
	As $\mathfrak{m}_{\mathfrak{Z}_{\Sigma,G}}$ is additive under union, it is enough to prove that for any simple closed curve $\gamma \in S_{\Sigma}^{\circ}$, $\mathfrak{m}_{\mathfrak{Z}_{\Sigma,G}}(\gamma) \in Z_{G}$ and there is a unique multicurve corresponding to each $m \in Z_{G}$.

	\smallskip

	For the first part, we decompose the geometric realisation $|G|$ into strips $\mathscr{S}_e$ for each edge $e$ and small triangular neighbourhoods $\mathscr{N}_v$ of each vertex the vertices as in Figure~\ref{fig:decomposition:SN:a}, and pullback this structure to $\Sigma$ via $f$. If $\gamma \in S_{\Sigma}^{\circ}$ is a simple closed curve in $\Sigma$, we can isotope $\gamma$ to a non-backtracking simple representative $\gamma$ that has $m_{e}$ parallel paths in the strip corresponding to $e \in E_{G}$. At each vertex $v$, it is possible to draw pairwise non-intersecting arcs connecting inside $\mathscr{N}_{v}$ the endpoints of $\gamma$ in $\de \mathscr{N}_{v}$ in a non-backtracking way if and only if $m_{\Delta} \in 2\NN$ for each corner $\Delta$. When these conditions hold, there is in fact a unique way (up to isotopy relative to $\de \mathscr{N}_{v}$) to draw such arcs as in Figure~\ref{fig:decomposition:SN:b}. Namely, we can label the points in $\de \mathscr{N}_{v} \cap \gamma \cap \mathscr{S}_{e}$ by $p_{e,1},\ldots,p_{e,m_{e}}$ following the cyclic order around $v$. Then, $p_{e,i}$ must be connected to
	\begin{itemize}
		\item $p_{e'',m_{e''} + 1 - i}$ for $1 \le i \leq \tfrac{1}{2} m_{\Delta'}$
		\item $p_{e',\frac{1}{2} m_{\Delta'} + 1 - i}$ for $\tfrac{1}{2} m_{\Delta'} < i \le \tfrac{1}{2}(m_{\Delta'} + m_{\Delta''}) = m_{e}$.
	\end{itemize}
	This proves that $\mathfrak{m}_{\mathfrak{Z}_{\Sigma,G}}$ is injective on $S_{\Sigma}^{\circ}$ and its image is included in $\set{ m \in \NN^{E_{G}} | m_{\Delta} \in 2\NN}$.

	\medskip

	Let $(e_{i})_{i \in \ZZ/N\ZZ}$ be the sequence of edges around a face $\mathfrak{f}$, and $\Delta_i$ be the corner containing both $e_i$ and $e_{i + 1}$. Then, the $\tfrac{1}{2} m_{\Delta_i}$ arcs in $\gamma \cap \mathscr{S}_{e_i}$ which are closest to $\mathfrak{f}$ are connected to the $\tfrac{1}{2}m_{\Delta_i}$ arcs in $\gamma \cap \mathscr{S}_{e_{i + 1}}$ which are the closest to $\mathfrak{f}$. In particular, the $\tfrac{1}{2} \min_{i} \set{ m_{\Delta_i} }$ arcs which are (in each strip around $\mathfrak{f}$) the closest to $\mathfrak{f}$ are connected and form loops, which are homotopic to the boundary component of $\Sigma$ that $\mathfrak{f}$ represents. By definition of $S_{\Sigma}^{\circ}$, we must have $\min_{i}\set{ m_{\Delta_i} } = 0$. This proves that $\mathfrak{m}_{\mathfrak{Z}_{\Sigma,G}}(S_{\Sigma}^{\circ}) \subset Z_{G}$ and the first part of the claim.

	\medskip

	Conversely, if we are given $m \in Z_{G}$, we draw $m_{e}$ parallel arcs in $\mathscr{S}_{e}$, and connect them inside each $\mathscr{N}_{v}$ in the unique non-intersecting and non-backtracking  (as explained above) way. We obtain a collection of simple closed curves, none of them being homotopic to a boundary component of $\Sigma$.
\end{proof}

As a consequence, we obtain a polynomial growth of the number of essential simple closed curve of bounded combinatorial length.

\begin{prop}\label{prop:poly:growth}
	Let $\Sigma$ be a connected bordered surface of type $(g,n)$ and $\epsilon > 0$. For any $\GG \in \mathcal{T}_{\Sigma}^{\textup{comb}}$ such that $\sys_{\GG} \geq \epsilon$, and any $t > 0$, we have 
	\begin{equation}\label{eqn:polynomial:growth}
		\#\Set{ \gamma \in S_{\Sigma}^{\circ} | \phantom{\big|} \ell_{\GG}(\gamma) \le t } \le m_{\epsilon} \, t^{6g-6+2n},
	\end{equation}
	for some positive constant $m_{\epsilon}$ depending on $\epsilon$ and the topology $(g,n)$ of the surface.
\end{prop}

\begin{proof}
	Let $\gamma \in S_{\Sigma}^{\circ}$ and denote by $E_{\GG}(\gamma)$ the number (with multiplicity) of edges $\gamma$ passes through. We recall from the proof of Proposition~\ref{lem:length:bound} that for $\sys_{\GG} \geq \epsilon$ we have
	\[
		E_{\GG}(\gamma) \leq \frac{(6g - 6 + 3n)^2}{\epsilon}\,\ell_{\GG}(\gamma).
	\]
	If the underlying graph $G$ is trivalent, from this bound and Lemma~\ref{lem:parametrisation:multicurves} we deduce
	\[
		\#\Set{ \gamma \in S_{\Sigma}^{\circ} | \phantom{\big|} \ell_{\GG}(\gamma) \le t }
		\le
		\#\Set{x \in Z_{G} | \sum_{e \in E_{G}} x_{e} \leq \frac{(6g - 6 + 3n)^2}{\epsilon} \, t}
	\]
	and the claim follows, as the latter is the set of integer points with bounded $L^1$-norm in a polytope of dimension $6g - 6 + 2n$. If $G$ is not trivalent, we resolve it in an arbitrary way into a trivalent ribbon graph $G'$ and the same argument works with $G'$.
\end{proof}

\begin{rem}
	A second proof of Proposition~\ref{prop:poly:growth} can be obtained by invoking that a similar result holds on $\mathcal{T}_{\Sigma}$ with respect to hyperbolic lengths -- \cite{Riv01,Mir08growth} for punctured surfaces, extended to bordered surfaces in \cite[Section 7.2]{ABO17} -- and flowing it to the combinatorial setting via Proposition~\ref{lem:length:bound}.
\end{rem}

\subsection{Unit balls in measured foliations: the hyperbolic case}
\label{subsec:unit:balls:hyp}

Let $\Sigma$ be a connected bordered surface of type $(g,n)$. In Section~\ref{subsec:param:mult:curv}, we considered the space ${\rm MF}_{\Sigma}^{\star}$ pa\-ra\-me\-tri\-sing measured foliations on $\Sigma$ up to Whitehead equivalence, allowing all types of boundary behaviours. We consider now the subset ${\rm MF}_{\Sigma} \subset {\rm MF}_{\Sigma}^{\star}$ that contains only those measured foliations whose boundary is made up of a union of singular leaves (that is, points at the boundary are of parallel type, \emph{cf.} Figure~\ref{fig:points:foliation:c} and ~\ref{fig:points:foliation:f}). In other words, on each boundary component there is at least one singularity, and the singularities on the boundary are connected by singular leaves. For convenience, we include in ${\rm MF}_{\Sigma}$ the empty foliation. Notice that ${\rm MF}_{\Sigma}$ is disjoint from the image of $\mathcal{T}_{\Sigma}^{\textup{comb}}$ in ${\rm MF}_{\Sigma}^{\star}$. 

\medskip

The space ${\rm MF}_{\Sigma}$ has dimension $6g - 6 + 2n$ and admits a canonical piecewise linear structure. It also admits an integral structure given by $M_{\Sigma}$, the set of multicurves in which components are not allowed to be homotopic to boundary components of $\Sigma$. We denote by $\mu_{{\rm Th}}$ the Thurston measure on ${\rm MF}_{\Sigma}$ associated to this piecewise linear integral structure\footnote{On ${\rm MF}_{\Sigma}$, the normalisation by counting lattice points differs by a power of $2$ from the normalisation obtained by taking the top power of the Thurston symplectic form.}. For an open set $A \subset {\rm MF}_{\Sigma}$, we have
\begin{equation}
	\mu_{{\rm Th}}(A) = \lim_{t \to +\infty} \frac{\#(t \cdot A \cap M_{\Sigma})}{t^{6g - 6 + 2n}}.
\end{equation}
The Thurston measure allows one to define a counting of multicurves associated to any length function on ${\rm MF}_{\Sigma}$. In fact, ${\rm MF}_{\Sigma}$ is a completion of the set of formal $\QQ_+$-linear combinations of simple closed curves. 

\medskip

In the hyperbolic context, for $\sigma \in \mathcal{T}_{\Sigma}$ the hyperbolic length function $\ell_{\sigma} \colon M_{\Sigma} \rightarrow \RR_{+}$ has a unique continuous extension $\ell_{\sigma} \colon {\rm MF}_{\Sigma} \rightarrow \RR_{\geq 0}$ that is compatible with the piecewise linear structure. Mirzakhani \cite{Mir08growth} introduces a function measuring the volume of the unit ball in ${\rm MF}_{\Sigma}$ with respect to $\ell_{\sigma}$.

\begin{defn}\label{defn:Mirz:function}
	The function $\mathscr{B}_{\Sigma} \colon \mathcal{T}_{\Sigma} \rightarrow [0,+\infty]$ is defined by
	\begin{equation}
		\mathscr{B}_{\Sigma}(\sigma) = \mu_{\textup{Th}}\bigl(
			\Set{ \mathcal{F} \in {\rm MF}_{\Sigma} | \ell_{\sigma}(\mathcal{F}) \le 1 }
		\bigr).
	\end{equation}
	It is manifestly $\Mod_{\Sigma}$-invariant, hence descends to a function $\mathscr{B}_{g,n}$ on $\mathcal{M}_{g,n}$ for $\Sigma$ of type $(g,n)$.
\end{defn}

By definition of the Thurston measure, this function describes the asymptotic growth of the number of multicurves of length $\leq t$, when $t \rightarrow \infty$, namely
\begin{equation}
	\mathscr{B}_{\Sigma}(\sigma) = \lim_{t \to +\infty} \frac{\#\Set{ c \in M_{\Sigma} | \ell_{\sigma}(c) \le t }}{t^{6g-6+2n}}.
\end{equation}
The main properties of $\mathscr{B}_{\Sigma}$ established by Mirzakhani for punctured surfaces -- \emph{i.e.} on $\mathcal{M}_{g,n}(0)$ -- can easily be generalised to the case of bordered surfaces.

\begin{thm}\label{thm:Mirz:properties} \cite[Proposition~3.2, Theorem~3.3]{Mir08growth}
	The function $\mathscr{B}_{\Sigma}$ takes values in $\RR_{+}$, is continuous on $\mathcal{T}_{\Sigma}$, and $\mathscr{B}_{g,n}$ is integrable on $\mathcal{M}_{g,n}(L)$ with respect to $\mu_{\textup{WP}}$.
	\hfill $\blacksquare$
\end{thm}

We remark that finiteness of $\mathscr{B}_{\Sigma}$ comes from the hyperbolic analogue -- known since \cite{Riv01} for punctured surfaces and extended to bordered surfaces in \cite[Theorem~7.2]{ABO17} -- of Proposition~\ref{prop:poly:growth}.

\medskip

More generally, the same properties hold for the counting function
\begin{equation}
	\mathcal{N}_{\Sigma}(\sigma;t) = \#\Set{ c \in M_{\Sigma} | \ell_{\sigma}(c) \le t },
\end{equation}
and Mirzakhani computed in \cite[Section~5]{Mir08growth} the integral of $\mathcal{N}_{g,n}(X;t)$ over $\mathcal{M}_{g,n}(L)$ as a sum over stable graphs. Here we present an equivalent version of her statement, in a form that suits better our next purpose.

\begin{lem}\label{lem:Mirz:count}
	The integral of $\mathcal{N}_{g,n}(X;t)$ over $\mathcal{M}_{g,n}(L)$ is a polynomial in $t^2,L_1^2,\dots,L_n^2$. Furthermore, its Laplace transform in the cutoff variable $t$ is given by
	\begin{equation}\label{eqn:int:bgn:t}
	\begin{split}
		& \int_{\RR_{+}} \left( \int_{\mathcal{M}_{g,n}(L)} \mathcal{N}_{g,n}(X;t) \, \dd\mu_{\textup{WP}}(X) \right) e^{-st} \dd t = \\
		& =
		\frac{1}{s} \sum_{\Gamma \in \mathcal{G}_{g,n}} \frac{1}{\#\Aut(\Gamma)}
			\int_{\RR_{+}^{E_{\Gamma}}}
			\prod_{v \in V_{\Gamma}}
				V_{g(v),n(v)}^{\textup{WP}}\bigl(
					(\ell_e)_{e \in E_{v}},(L_{\lambda})_{\lambda \in \Lambda_{v}}
				\bigr)
			\prod_{e \in E_{\Gamma}} \frac{\ell_e \, \dd\ell_e}{e^{s\ell_e} - 1}
	\end{split}
	\end{equation}
	and it is a polynomial of degree $6g - 5 + 2n$ in $s^{-1}$.
	\hfill $\blacksquare$
\end{lem}

As a consequence, using the final value theorem for the Laplace transform, we have an expression of the integral of $\mathscr{B}_{g,n}$ over $\mathcal{M}_{g,n}(L)$ as a sum over stable graphs (see also \cite{ABCDGLW19,DGZZ19}).

\begin{thm}\label{thm:Mirz:asymptotic:count} \cite[Theorem~5.3]{Mir08growth}
	The integral of $\mathscr{B}_{g,n}$ is computed as
	\begin{equation}\label{eqn:Bhyp:sum:stable:graphs}
	\begin{split}
		& (6g - 6 + 2n)! \int_{\mathcal{M}_{g,n}(L)} \mathscr{B}_{g,n} \, \dd\mu_{\textup{WP}} = \\
		& =
		 \sum_{\Gamma \in \mathcal{G}_{g,n}} \frac{1}{\#\Aut(\Gamma)}
			\int_{\RR_{+}^{E_{\Gamma}}}
			\prod_{v \in V_{\Gamma}}
				V_{g(v),n(v)}^{\textup{K}}\bigl(
					(\ell_e)_{e \in E_{v}},(0_{\lambda})_{\lambda \in \Lambda_{v}}
				\bigr)
			\prod_{e \in E_{\Gamma}} \frac{\ell_e \, \dd\ell_e}{e^{\ell_e} - 1}.
	\end{split}
	\end{equation}
	\hfill $\blacksquare$
\end{thm}

Notice the presence of Kontsevich volumes as vertices weights. This comes from the fact that the top coefficient of $s^{-1}$ in equation \eqref{eqn:int:bgn:t} is given by extraction of the largest possible power of the $\ell_{e}$, and that the top degree of $V_{g,n}^{\textup{WP}}$ is equal to $V_{g,n}^{\textup{K}}$. We return in \S~\ref{subsubsec:connection:comb:counting} to the fact that the right-hand side of Equation~\eqref{eqn:int:bgn:t} can be calculated by the topological recursion.

\subsection{Unit balls in measured foliations: the combinatorial case}
\label{subsec:unit:balls:comb}

Following the steps of the previous paragraph, we can also study the volume of the combinatorial unit balls in ${\rm MF}_{\Sigma}$. For $\GG \in \mathcal{T}_{\Sigma}^{\textup{comb}}$, we have a combinatorial length function $\ell_{\GG} \colon {\rm MF}_{\Sigma} \rightarrow \RR_{\geq 0}$, induced by the one on the larger space of measured foliations ${\rm MF}_{\Sigma}^{\star}$ and already used in Section~\ref{subsec:length:top} to define the combinatorial length of curves.  

\begin{defn}\label{defn:Mirz:comb:function}
	The function $\mathscr{B}_{\Sigma}^{\textup{comb}} \colon \mathcal{T}_{\Sigma}^{\textup{comb}} \rightarrow [0,+\infty]$ is defined by
	\begin{equation}
		\mathscr{B}_{\Sigma}^{\textup{comb}}(\GG) = \mu_{\textup{Th}}\bigl(
			\Set{ \mathcal{F} \in {\rm MF}_{\Sigma} | \ell_{\GG}(\mathcal{F}) \le 1 }
		\bigr).
	\end{equation}
\end{defn}

By construction, this function measures the asymptotic growth of the number of multicurves with respect to combinatorial lengths
\begin{equation}
	\mathscr{B}_{\Sigma}^{\textup{comb}}(\GG) = \lim_{t \to +\infty} \frac{\#\Set{ c \in M_{\Sigma} | \ell_{\GG}(c) \le t }}{t^{6g-6+2n}}.
\end{equation}
In fact, it can be obtained as a limit of the hyperbolic setting via the rescaling flow.

\begin{lem}\label{lem:B:hyp:comb}
	For any connected bordered surface $\Sigma$ of type $(g,n)$, we have
	\begin{equation}
		\lim_{\beta \rightarrow \infty} \beta^{6g - 6 + 2n} \, \mathscr{B}_{\Sigma}(\sigma^{\beta})
		=
		\mathscr{B}_{\Sigma}^{\textup{comb}}(\sp(\sigma))
	\end{equation}
	and the limit is uniform for $\sp(\sigma)$ in any thick part of $\mathcal{T}_{\Sigma}^{\textup{comb}}$.
\end{lem}

\begin{proof}
	Let $\epsilon > 0$ and $\GG \in \mathcal{T}_{\Sigma}^{\textup{comb}}$ such that $\sys_{\GG} \geq \epsilon$. We denote $\sigma = \sp^{-1}(\GG)$ and $\sigma^{\beta} = \sp^{-1}(\beta\GG)$. Since the length of a multicurve is the sum of lengths of its connected components, Proposition~\ref{lem:length:bound} implies that for any $c \in M_{\Sigma}$, we have $\beta \ell_{\GG}(c) \leq \ell_{\sigma^{\beta}}(c) \leq (\beta + \kappa_{\epsilon})\ell_{\GG}(c)$. Therefore
	\[
		\Set{c \in M_{\Sigma} | \ell_{\GG}(c) \leq \frac{t}{\beta + \kappa_{\epsilon}}}
		\subseteq
		\Set{c \in M_{\Sigma} | \ell_{\sigma^{\beta}}(c) \leq t \vphantom{\frac{t}{\beta}}}
		\subseteq
		\Set{c \in M_{\Sigma} | \ell_{\sigma^{\beta}}(c) \leq \frac{t}{\beta}}
	\]
	and thus
	\[
		\frac{\mathscr{B}^{\textup{comb}}_{\Sigma}(\GG)}{(\beta + \kappa_{\epsilon})^{6g - 6 + 2n}}
		\leq
		\mathscr{B}_{\Sigma}(\sigma^{\beta})
		\leq
		\frac{\mathscr{B}^{\textup{comb}}_{\Sigma}(\GG)}{\beta^{6g - 6 + 2n}}.
	\]
	Multiplying by $\beta^{6g - 6 + 2n}$ and taking the limit $\beta \rightarrow \infty$ yields the result.
\end{proof}

We shall now establish the combinatorial analogues of Theorem~\ref{thm:Mirz:properties}, Lemma~\ref{lem:Mirz:count} and Theorem~\ref{thm:Mirz:asymptotic:count}.

\medskip

We first notice that $\mathscr{B}_{\Sigma}^{\textup{comb}}$ takes values in $\RR_{+}$, as a consequence of Proposition~\ref{prop:poly:growth}, and is continuous since the length function $\ell_{\ast} \colon \mathcal{T}_{\Sigma}^{\textup{comb}} \to \RR_{+}^{S_{\Sigma}}$ is continuous. The same holds for the counting function
\begin{equation}
	\mathcal{N}_{\Sigma}^{\textup{comb}}(\GG;t) = \#\Set{c \in M_{\Sigma} | \ell_{\GG}(c) \le t}.
\end{equation}
The integrability result is more delicate. We derive it from the following two statements: the integration formulae of Lemma~\ref{lem:comb:Mirz:count} and Theorem~\ref{thm:comb:Mirz:asymptotic:count}.

\begin{lem}\label{lem:comb:Mirz:count}
	The integral of $\mathcal{N}_{g,n}^{\textup{comb}}(\bm{G};t)$ over $\mathcal{M}_{g,n}^{\textup{comb}}(L)$ is a polynomial in $t^2,L_1^2,\dots,L_n^2$. Furthermore, its Laplace transform in the cut-off variable $t$ is given by
	\begin{equation}\label{eqn:Ncomb:sum:stable:graphs}
	\begin{split}
		& \int_{\RR_{+}} \left( \int_{\mathcal{M}_{g,n}^{\textup{comb}}(L)} \mathcal{N}_{g,n}^{\textup{comb}}(\bm{G};t) \, \dd\mu_{\textup{K}}(\bm{G}) \right) e^{-st} \dd t \\
		& =
		\frac{1}{s} \sum_{\Gamma \in \mathcal{G}_{g,n}} \frac{1}{\#\Aut(\Gamma)}
			\int_{\RR_{+}^{E_{\Gamma}}}
			\prod_{v \in V_{\Gamma}}
				V_{g(v),n(v)}^{\textup{K}}\bigl(
					(\ell_e)_{e \in E_{v}},(L_{\lambda})_{\lambda \in \Lambda_{v}}
				\bigr)
			\prod_{e \in E_{\Gamma}} \frac{\ell_e \, \dd\ell_e}{e^{s\ell_e} - 1}.
	\end{split}
	\end{equation}
\end{lem}

\begin{proof}
	The proof is a variation on the ideas of \cite[Section~5]{Mir08growth}, adapted to our notation and the combinatorial setting. Consider the class of a multicurve on $\Sigma$ under the action of the mapping class group. It is determined by the data $(\Gamma,a)$ of a stable graph $\Gamma \in \mathcal{G}_{g,n}$, together with an integral tuple $a \in \ZZ_{+}^{E_{\Gamma}}$ giving the multiplicity of each component. We then consider the frequency of the multicurve in the mapping class orbit $(\Gamma,a)$ with combinatorial length bounded by $t$:
	\[
		s_{\Sigma}^{\textup{comb}}(\GG;\Gamma,a,t) = \#\Set{ c \in (\Gamma,a) | \ell_{\GG}(c) \le t}.
	\]
It is a nonnegative mapping class group invariant function of $\GG$. Thus, we can apply the integration formula of Proposition~\ref{prop:integration} to the function $s_{g,n}^{\textup{comb}}$ induced on the moduli space $\mathcal{M}_{g,n}^{\textup{comb}}(L)$ (see also Remarks~\ref{rem:integration:lemma:automorphism} and \ref{rem:integration:lemma:nonnegative}). With the notation of Section~\ref{subsec:integration}, consider the function $f \colon \RR_+^n \times \RR_{+}^{E_{\Gamma}} \to \RR$ given by $\Xi_{\Gamma} \equiv 1$ and
	\[
		f(L,\ell) =
		\begin{cases}
			1 & \text{if } \sum_{e \in E_{\Gamma}} a_e \ell_e \le t, \\
			0 & \text{otherwise}.
		\end{cases}
	\]
Then we get $\Xi_{\Sigma}^{f,\Gamma}(\GG) = \mathcal{N}_{\Sigma}^{\textup{comb}}(\GG;t)$, and in particular
	\[
		\int_{\mathcal{M}_{g,n}^{\textup{comb}}(L)} s_{g,n}^{\textup{comb}}(\bm{G};\Gamma,a,t) \, \dd\mu_{\textup{K}}(\bm{G})
		=
		\frac{1}{\#\Aut(\Gamma,a)} \int_{\braket{a,l} \le t}
			\prod_{v \in V_{\Gamma}}
				V_{g(v),n(v)}^{\textup{K}}\bigl(
					(\ell_e)_{e \in E_{v}},(L_{\lambda})_{\lambda \in \Lambda_{v}}
				\bigr)
			\prod_{e \in E_{\Gamma}} \ell_e \, \dd\ell_e,
	\]
	where $\braket{a,\ell} = \sum_{e \in E_{\Gamma}} a_e \ell_e$. It is easy to see that such integrals are polynomial in $t^2, L_1^2,\dots,L_n^2$. Furthermore, notice that for any polynomial $p(\ell)$ in the variables $(\ell_e^2)_{e \in E_{\Gamma}}$, we have
	\[
		\int_{\RR_+} \biggl(
			\int_{\braket{a,\ell} \le t} p(\ell) \prod_{e \in E_{\Gamma}} \ell_e \, \dd\ell_e
		\biggr) e^{-st} \dd t
		=
		\frac{1}{s} \int_{\RR_{+}^{E_{\Gamma}}} p(\ell) \prod_{e \in E_{\Gamma}} e^{-s a_e \ell_e} \ell_e \, \dd\ell_e.
	\]
	From this, together with the relation 
	\[
		\mathcal{N}_{g,n}^{\textup{comb}}(\bm{G};t) = \sum_{\Gamma \in \mathcal{G}_{g,n}} \sum_{a \in \ZZ_{+}^{E_{\Gamma}}} \frac{\#\Aut(\Gamma,a)}{\#\Aut(\Gamma)} \, s_{g,n}^{\textup{comb}}(\bm{G};\Gamma,a,t)
	\]
	and Tonelli's theorem, we find
	\[
	\begin{split}
		& \int_{\RR_{+}} \biggl( \int_{\mathcal{M}_{g,n}^{\textup{comb}}(L)} \mathcal{N}_{g,n}^{\textup{comb}}(\bm{G};t) \, \dd\mu_{\textup{K}}(\bm{G}) \biggr) e^{-st} \dd t = \\
		& =
		\sum_{\Gamma \in \mathcal{G}_{g,n}}
			\frac{1}{\#\Aut(\Gamma)}
			\sum_{a \in \ZZ_{+}^{E_{\Gamma}}}
				\int_{\RR_{+}} \biggl( \int_{\braket{a,l} \le t}
				\prod_{v \in V_{\Gamma}}
				V_{g(v),n(v)}^{\textup{K}}\bigl(
					(\ell_e)_{e \in E_{v}},(L_{\lambda})_{\lambda \in \Lambda_{v}}
				\bigr)
				\prod_{e \in E_{\Gamma}} \ell_e \, \dd\ell_e
			\biggr) e^{-st} \dd t \\
		& =
		\frac{1}{s} \sum_{\Gamma \in \mathcal{G}_{g,n}} \frac{1}{\#\Aut(\Gamma)}
			\sum_{a \in \ZZ_{+}^{E_{\Gamma}}}
			\int_{\RR_{+}^{E_{\Gamma}}}
			\prod_{v \in V_{\Gamma}}
				V_{g(v),n(v)}^{\textup{K}}\bigl(
					(\ell_e)_{e \in E_{v}},(L_{\lambda})_{\lambda \in \Lambda_{v}}
				\bigr)
			\prod_{e \in E_{\Gamma}} e^{-s a_e \ell_e} \ell_e \, \dd\ell_e \\
		& =
		\frac{1}{s} \sum_{\Gamma \in \mathcal{G}_{g,n}} \frac{1}{\#\Aut(\Gamma)}
			\int_{\RR_{+}^{E_{\Gamma}}}
			\prod_{v \in V_{\Gamma}}
				V_{g(v),n(v)}^{\textup{K}}\bigl(
					(\ell_e)_{e \in E_{v}},(L_{\lambda})_{\lambda \in \Lambda_{v}}
				\bigr)
			\prod_{e \in E_{\Gamma}} \frac{\ell_e \, \dd\ell_e}{e^{s\ell_e} - 1}.
	\end{split}
	\]
\end{proof}

\begin{thm}\label{thm:comb:Mirz:asymptotic:count}
	The integral of $\mathscr{B}_{g,n}^{\textup{comb}}$ is computed as
	\begin{equation}\label{eqn:Bcomb:sum:stable:graphs}
	\begin{split}
		(6g - 6 + 2n)! \int_{\mathcal{M}_{g,n}^{\textup{comb}}(L)} \mathscr{B}_{g,n}^{\textup{comb}} \, \dd\mu_{\textup{K}}
		=
		\sum_{\Gamma \in \mathcal{G}_{g,n}} \frac{1}{\#\Aut(\Gamma)}
			\int_{\RR_{+}^{E_{\Gamma}}}
			\prod_{v \in V_{\Gamma}}
				V_{g(v),n(v)}^{\textup{K}}\bigl(
					(\ell_e)_{e \in E_{v}},(0_{\lambda})_{\lambda \in \Lambda_{v}}
				\bigr)
			\prod_{e \in E_{\Gamma}} \frac{\ell_e \, \dd\ell_e}{e^{\ell_e} - 1}.
	\end{split}
	\end{equation}
\end{thm}

\begin{proof}
	From Lemma~\ref{lem:comb:Mirz:count}, with the change of variable $\ell_e \mapsto \ell_e/s$, we have
	\[
	\begin{split}
		& \int_{\RR_{+}} \biggl( \int_{\mathcal{M}_{g,n}^{\textup{comb}}(L)} \mathcal{N}_{g,n}^{\textup{comb}}(\bm{G};t) \, \dd\mu_{\textup{K}}(\bm{G}) \biggr) e^{-st} \dd t = \\
		& =
		\frac{1}{s^{6g-5+2n}} \biggl(
		\sum_{\Gamma \in \mathcal{G}_{g,n}} \frac{1}{\#\Aut(\Gamma)}
			\int_{\RR_{+}^{E_{\Gamma}}}
			\prod_{v \in V_{\Gamma}}
				V_{g(v),n(v)}^{\textup{K}}\bigl(
					(\ell_e)_{e \in E_{v}},(0_{\lambda})_{\lambda \in \Lambda_{v}}
				\big)
			\prod_{e \in E_{\Gamma}} \frac{\ell_e \, \dd\ell_e}{e^{\ell_e} - 1}
		+
		O(s) \biggr).
	\end{split}
	\]
	Here we used the relation $\#E_{\Gamma} + \sum_{v \in V_{\Gamma}} (3g(v) - 3 + n(v)) = 3g - 3 + n$. From properties of the Laplace transform, we find
	\[
	\begin{split}
		& \int_{\RR_{+}} \biggl( \int_{\mathcal{M}_{g,n}^{\textup{comb}}(L)} \frac{\mathcal{N}_{g,n}^{\textup{comb}}(\bm{G};t)}{t^{6g-6+2n}} \, \dd\mu_{\textup{K}}(\bm{G}) \biggr) e^{-st} \dd t = \\
		& =
		\frac{1}{s} \biggl( \frac{1}{(6g-6+2n)!} \sum_{\Gamma \in \mathcal{G}_{g,n}} \frac{1}{\#\Aut(\Gamma)}
			\int_{\RR_{+}^{E_{\Gamma}}}
			\prod_{v \in V_{\Gamma}}
				V_{g(v),n(v)}^{\textup{K}}\bigl(
					(\ell_e)_{e \in E_{v}},(0_{\lambda})_{\lambda \in \Lambda_{v}}
				\bigr)
			\prod_{e \in E_{\Gamma}} \frac{\ell_e \, \dd\ell_e}{e^{\ell_e} - 1}
		+ O(s)
		\biggr),
	\end{split}
	\]
	and from the final value theorem for the Laplace transform, we have the thesis.
\end{proof}

Note that Lemma~\ref{lem:comb:Mirz:count} and Theorem~\ref{thm:comb:Mirz:asymptotic:count} are valid as equalities in $[0,+\infty]$, since we only used Tonelli's theorem. Since $V^{\textup{K}}$ are polynomials, the right-hand sides are actually finite, hence the integrability.

\begin{cor}
	The functions $\mathcal{N}_{g,n}^{\textup{comb}}$ and $\mathscr{B}_{g,n}^{\textup{comb}}$ are integrable on $\mathcal{M}_{g,n}(L)$ with respect to $\mu_{\textup{K}}$.
	\hfill $\blacksquare$
\end{cor}

By comparing Equations~\eqref{eqn:Bhyp:sum:stable:graphs} and \eqref{eqn:Bcomb:sum:stable:graphs}, we also deduce that the integrals of the $\mathscr{B}$-functions are the same in the combinatorial and in the hyperbolic world and independently of boundary lengths.

\begin{cor}
	The following equality holds
	\begin{equation}
		\int_{\mathcal{M}_{g,n}(L)} \mathscr{B}_{g,n} \, \dd\mu_{\textup{WP}}
		=
		\int_{\mathcal{M}_{g,n}^{\textup{comb}}(L)} \mathscr{B}_{g,n}^{\textup{comb}} \, \dd\mu_{\textup{K}}.
	\end{equation}
	Furthermore, both quantities are independent of $L \in \RR_{+}^n$.
	\hfill $\blacksquare$
\end{cor}

\subsection{Geometric recursion for combinatorial length statistics of multicurves}
\label{subsec:GR:comb:length:statistics}

We can rephrase the results of the previous paragraph by saying that we were able to compute combinatorial length statistics of multicurves with respect to the Heaviside function $\theta(t - \ell)$. Indeed the combinatorial Mirzakhani's counting function $\mathcal{N}_{\Sigma}^{\textup{comb}}$ can be rewritten as
\begin{equation}
	\mathcal{N}_{\Sigma}^{\textup{comb}}(\GG;t)
	=
	\sum_{c \in M_{\Sigma}} \theta\bigl( t - \ell_{\GG}(c) \bigr).
\end{equation}
In this paragraph, we see the overlap between twisted GR amplitudes and additive statistics of lengths of multicurves. Furthermore, we connect such amplitudes to the counting $\mathcal{N}_{\Sigma}^{\textup{comb}}$ via Laplace transform.

\subsubsection{Geometric recursion and statistics of multicurves}

An immediate consequence of Section~\ref{subsec:twisted:GR} is a GR formula for the combinatorial length statistic of multicurves with exponentially decaying weight: Take $\Xi_{\Sigma} = \Xi_{\Sigma}^{\textup{K}} \equiv 1$ to be the Kontsevich amplitudes, and $F(\ell;s) = e^{-s\ell}$ for multicurves, or equivalently $f(\ell;s) = \frac{1}{e^{s\ell} - 1}$ for primitive multicurves (\emph{cf.} Remark~\ref{rem:pure:statistics:multicurves}). This choice is of particular interest, because for every $\GG \in \mathcal{T}_{\Sigma}^{\textup{comb}}$ and every multicurve $c \in M_{\Sigma}$,
\[
	\prod_{\gamma \in \pi_{0}(c)} e^{-s\ell_{\GG}(\gamma)} = e^{-s\ell_{\GG}(c)}.
\]
Hence, the specialisation of Theorem~\ref{thm:GR:statistics} to such choices gives the following result.

\begin{cor}
	The family of combinatorial GR initial data
	\begin{equation}\label{eqn:GR:initial:data:comb:statistics:exp}
	\begin{aligned}
		A(L_{1},L_{2},L_{3};s) & =  1, \\
		B(L_1,L_2,\ell;s) & = B^{\textup{K}}(L_1,L_2,\ell) + \frac{1}{e^{s\ell} - 1}, \\
		C(L_1,\ell,\ell';s) & = C^{\textup{K}}(L_1,\ell,\ell')
			+ \frac{B^{\textup{K}}(L_1,\ell,\ell')}{e^{s\ell} - 1} + \frac{B^{\textup{K}}(L_1,\ell',\ell)}{e^{s\ell'} - 1}
			+ \frac{1}{(e^{s\ell} - 1)(e^{s\ell'} - 1)}, \\
		D_T(\GG;s) & = \sum_{c \in M_{T}} e^{-s\ell_{\GG}(c)},
	\end{aligned}
	\end{equation}
	are strongly admissible, and the associated GR amplitudes $\Xi_{\Sigma}(\GG;s)$ satisfy
	\begin{equation}\label{eqn:GR:comb:statistics:exp}
		\Xi_{\Sigma}(\GG;s) = \sum_{c \in M_{\Sigma}} e^{-s\ell_{\GG}(c)}.
	\end{equation}
	Furthermore, the TR amplitudes $V\Xi_{g,n}(L;s)$ satisfy
	\begin{equation}\label{eqn:TR:comb:statistics:exp}
		V\Xi_{g,n}(L;s)
		=
		\sum_{\Gamma \in \mathcal{G}_{g,n}} \frac{1}{\#\Aut(\Gamma)}
			\int_{\RR_{+}^{E_{\Gamma}}}
			\prod_{v \in V_{\Gamma}} V_{g(v),n(v)}^{\textup{K}}\big((\ell_e)_{e \in E(v)},(L_{\lambda})_{\lambda \in \Lambda(v)}\big)
			\prod_{e \in E_{\Gamma}} \frac{\ell_e \dd\ell_e}{e^{s\ell_e} - 1}.
	\end{equation}
	\hfill $\blacksquare$
\end{cor}

\subsubsection{Connection with the combinatorial Mirzakhani's counting function}
\label{subsubsec:connection:comb:counting}

Taking the Laplace transform in the cut-off parameter $t$ of the counting function $\mathcal{N}_{\Sigma}^{\textup{comb}}(\GG;t)$ and its average $\int_{\mathcal{M}_{g,n}^{\textup{comb}}(L)} \mathcal{N}_{g,n}^{\textup{comb}}(\bm{G};t) \dd\mu_{\textup{K}}(\bm{G})$, and comparing them with Equations~\eqref{eqn:GR:comb:statistics:exp} and \eqref{eqn:TR:comb:statistics:exp}, we get the following corollary, which retrieves Lemma~\ref{lem:comb:Mirz:count}.

\begin{thm}
	The Laplace transform of the counting functions
	\begin{equation}
		s \int_{\RR_{+}} \mathcal{N}_{\Sigma}^{\textup{comb}}(\GG;t) \, e^{-st} \dd t
	\end{equation}
	coincides with the GR amplitudes $\Xi_{\Sigma}(\GG;s)$ computed from the initial data \eqref{eqn:GR:initial:data:comb:statistics:exp}. Further, their averages over the combinatorial moduli spaces
	\begin{equation}\label{eqn:Laplace:average:multicurve:count}
		V\Xi_{g,n}(L;s) = s \int_{\RR_{+}} \left(
			\int_{\mathcal{M}_{g,n}^{\textup{comb}}(L)} \mathcal{N}_{g,n}^{\textup{comb}}(\bm{G};t) \, \dd\mu_{\textup{K}}(\bm{G})
		\right) e^{-st} \dd t
	\end{equation}
	are computed by TR, and by the sum over stable graphs \eqref{eqn:Ncomb:sum:stable:graphs}.
\end{thm}

As a consequence, we find an efficient method to compute the quantities $\int_{\mathcal{M}_{g,n}^{\textup{comb}}(L)} \mathcal{N}_{g,n}^{\textup{comb}}(\bm{G};t) \dd\mu_{\textup{K}}(\bm{G})$, that are the average number of multicurves with combinatorial length bounded by $t$. Namely, one can recursively compute the TR amplitudes $V\Xi_{g,n}(L;s)$, then take an inverse Laplace transform. A list of such averages for low $2g - 2 + n$ can be find in Table~\ref{tab:multicurve:poly}.

\medskip 

In \cite{ABCDGLW19} the authors introduced some polynomials $V^{\textup{MV}}_{g,n}(L)$, indexed by $(g,n)$, and called them Masur--Veech polynomials. They are computed by topological recursion, and their constant terms are the Masur--Veech volumes of the principal strata of the moduli spaces of quadratic differentials (\emph{cf.} next section). Here we propose to upgrade the definition, and consider the $1$-parameter family of polynomials
\begin{equation}\label{eqn:one:parameter:MV:poly}
	V_{g,n}^{\textup{MV}}(L;s)
	=
	\sum_{\Gamma \in \mathcal{G}_{g,n}} \frac{1}{\#\Aut(\Gamma)}
		\int_{\RR_{+}^{E_{\Gamma}}}
		\prod_{v \in V_{\Gamma}} V_{g(v),n(v)}^{\textup{K}}\big((\ell_e)_{e \in E(v)},(L_{\lambda})_{\lambda \in \Lambda(v)}\big)
		\prod_{e \in E_{\Gamma}} \frac{\ell_e \dd\ell_e}{e^{s\ell_e} - 1},
\end{equation}
which satisfy $V^{\textup{MV}}_{g,n}(L;s) = s^{-6g+6-2n}V^{\textup{MV}}_{g,n}(sL)$. The above argument shows that such polynomials coincide with \eqref{eqn:Laplace:average:multicurve:count}, providing therefore a more direct geometric interpretation for them.

\begin{rem}\label{rem:hyp:GR:multicurves:counting}
	As mentioned previously, the hyperbolic quantities $\mathcal{N}_{g,n}(\sigma;t)$ and their averages can be calculated by GR/TR. Namely, the family of hyperbolic GR initial data
	\begin{equation}\label{eqn:GR:initial:data:hyp:statistics:exp}
	\begin{aligned}
		A(L_{1},L_{2},L_{3};s) & =  1, \\
		B(L_1,L_2,\ell;s) & = B^{\textup{M}}(L_1,L_2,\ell) + \frac{1}{e^{s\ell} - 1}, \\
		C(L_1,\ell,\ell';s) & = C^{\textup{M}}(L_1,\ell,\ell')
			+  \frac{B^{\textup{M}}(L_1,\ell,\ell')}{e^{s\ell} - 1} + \frac{B^{\textup{M}}(L_1,\ell',\ell)}{e^{s\ell'} - 1}
			+ \frac{1}{(e^{s\ell} - 1)(e^{s\ell'} - 1)}, \\
		D_T(\sigma;s) & = \sum_{c \in M_{T}} e^{-s\ell_{\sigma}(c)},
	\end{aligned}
	\end{equation}
	are strongly admissible, and the associated GR amplitudes $\Omega_{\Sigma}(\sigma;s)$ satisfy
	\begin{equation}\label{eqn:GR:hyp:statistics:exp}
		\Omega_{\Sigma}(\sigma;s) 
		= \sum_{c \in M_{\Sigma}} e^{-s\ell_{\sigma}(c)}
		= s \int_{\RR_{+}} \mathcal{N}_{\Sigma}(\sigma;t) \, e^{-st} \dd t.
	\end{equation}
	Furthermore, the TR amplitudes $V\Omega_{g,n}(L;s)$ satisfy
	\begin{equation}\label{eqn:TR:hyp:statistics:exp}
	\begin{split}
		V\Omega_{g,n}(L;s)
		& =
		\sum_{\Gamma \in \mathcal{G}_{g,n}} \frac{1}{\#\Aut(\Gamma)}
			\int_{\RR_{+}^{E_{\Gamma}}}
			\prod_{v \in V_{\Gamma}} V_{g(v),n(v)}^{\textup{WP}}\big((\ell_e)_{e \in E(v)},(L_{\lambda})_{\lambda \in \Lambda(v)}\big)
			\prod_{e \in E_{\Gamma}} \frac{\ell_e \dd\ell_e}{e^{s\ell_e} - 1} \\
		& =
		s \int_{\RR_{+}} \left(
			\int_{\mathcal{M}_{g,n}(L)} \mathcal{N}_{g,n}(\sigma;t) \, \dd\mu_{\textup{WP}}(\sigma)
		\right) e^{-st} \dd t.
	\end{split}
	\end{equation}
	The quantities $V\Omega_{g,n}(L;s)$ should be thought of as the hyperbolic analogue of the $1$-parameter family of Masur--Veech polynomials of \eqref{eqn:one:parameter:MV:poly}. Notice that, as the Weil--Petersson volumes are not homogeneous, the dependence on $s$ is more complicated in this case.
\end{rem}

\subsection{Combinatorial formulae for Masur--Veech volumes}
\label{subsec:MV:formulae}

If $\Sigma$ is a smooth surface of genus $g$ with $n > 0$ labelled \emph{punctures}, we denote $\mathfrak{T}_{\Sigma}$ (resp. $\mathfrak{M}_{g,n}$) the corresponding Teichm\"uller space (resp. moduli space). On the bundle $Q\mathfrak{T}_{\Sigma} \rightarrow \mathfrak{T}_{\Sigma}$ of meromorphic quadratic differentials with simple poles at the punctures, Masur \cite{Mas82} and Veech \cite{Vee82} introduced a piecewise linear structure via period coordinates, and an associated measure $\mu_{\textup{MV}}$ which is $\Mod_{\Sigma}$-invariant. This induces a measure on the total space of the bundle $Q^1\mathfrak{M}_{g,n} \rightarrow \mathfrak{M}_{g,n}$ of quadratic differentials of unit area, by taking for a measurable set $Y \subseteq Q^1\mathfrak{M}_{g,n}$
\begin{equation}
	\mu_{\textup{MV}}^{1}(Y) = (12g - 12 + 4n) \, \mu_{\textup{MV}}\big(\Set{ tq | t \in (0,\tfrac{1}{2}) \, \text{ and } \, q \in Y}\big).
\end{equation}
The total mass of this measure is finite, and is called \emph{Masur--Veech volume} (for the open stratum of the moduli space of quadratic differentials):
\begin{equation}\label{eqn:MV:volume:def}
	MV_{g,n} = \mu_{\textup{MV}}^{1}(Q^1\mathfrak{M}_{g,n}).
\end{equation}
Many works have been recently devoted to the understanding of $MV_{g,n}$ and closely related quantities, motivated (and related to) dynamical questions on the space of measured foliations (and on the Teichm\"uller space). It is sometimes convenient to change normalisations and divide $MV_{g,n}$ by the constants
\begin{equation}
	a_{g,n} = 2^{4g - 2 + n}(4g - 4 + n)!\cdot (6g - 6 + 2n) \qquad {\rm or}\qquad a'_{g,n} = \frac{2^{4g - 2 + n}(4g - 4 + n)!}{(6g - 7 + 2n)!}.
\end{equation}
Our present work links the combinatorial geometry to $MV_{g,n}$, adding on previously known approaches. They can be summarised in the following way.
\begin{prop}\label{prop:MV:equivalence}
	The following quantities are all equal:
	\begin{itemize}
		\item[\rm{(i)}]
			$MV_{g,n}$ defined by \eqref{eqn:MV:volume:def}.
		\item[\rm{(ii)}]
			$a'_{g,n}$ times the sum over stable graphs in Theorem~\ref{thm:Mirz:asymptotic:count}.
		\item[\rm{({iii}\textsubscript{0})}]
			$a_{g,n} \int_{\mathfrak{M}_{g,n}} \mathscr{B}_{\Sigma} \dd\mu_{\textup{WP}}$.
		\item[\rm{(iii)}]
			$a_{g,n}  \int_{\mathcal{M}_{g,n}(L)} \mathscr{B}_{\Sigma} \dd\mu_{\textup{WP}}$, for any $L \in \RR_{+}^n$.
		\item[\rm{(iv)}]
			$a'_{g,n}$ times the constant term of the hyperbolic TR amplitudes results by twisting the Kontsevich initial data by $f(\ell) = \frac{1}{e^{\ell} - 1}$.
		\item[\rm{(iv\textsubscript{+})}]
			$a'_{g,n} \lim\limits_{s \to 0} s^{6g-5+2n} V\Omega_{g,n}(L;s)$, where $V\Omega_{g,n}(L;s)$ are the $1$-parameter family of hyperbolic TR amplitudes for the initial data \eqref{eqn:GR:initial:data:hyp:statistics:exp}.
		\item[\rm{(v)}]
			$a_{g,n}  \int_{\mathcal{M}_{g,n}^{\textup{comb}}(L)} \mathscr{B}_{\Sigma}^{\textup{comb}} \dd\mu_{\textup{K}}$ for any $L \in \RR_{+}^n$.
		\item[\rm{(vi)}]
			$a'_{g,n}$ times the constant term of the combinatorial TR amplitudes results by twisting the Kontsevich initial data by $f(\ell) = \frac{1}{e^{\ell} - 1}$.
		\item[\rm{(vi\textsubscript{+})}]
			$a'_{g,n} \lim\limits_{s \to 0} s^{6g-5+2n} V\Xi_{g,n}(L;s)$, where $V\Xi_{g,n}(L;s)$ are the $1$-parameter family of combinatorial TR amplitudes for the initial data \eqref{eqn:GR:initial:data:comb:statistics:exp}.
		\end{itemize}
\end{prop}
In \cite{CMSBGL} one can find two more formulae for $MV_{g,n}$, proved starting from (i): the first one involves intersection theory on $\overline{\mathfrak{M}}_{g,n}$, and the second is another topological recursion distinct from (iv\textsubscript{+}) and (vi\textsubscript{+}). As they require very different methods (involving algebraic geometry of moduli spaces of curves), we do not discuss them.

\medskip

There are many different paths and methods to establish all the equalities in Proposition~\ref{prop:MV:equivalence}. To clarify the picture, we list what currently is known: \\[.5ex]
\noindent
	(i)~=~(iii\textsubscript{0}). Comes from the work of Bonahon \cite{Bonahon93} and Mirzakhani \cite{Mir08earth}, see \textit{e.g.} \cite{DGZZ19} or \cite[Lemma 3.1]{ABCDGLW19}. \\[.5ex]
\noindent
	(ii)~=~(iii\textsubscript{0}). Proved by Mirzakhani in \cite[Theorem~5.3]{Mir08growth}, studying the Laplace transform of the counting function as reported here in Section~\ref{subsec:unit:balls:hyp}. Proved differently in \cite{DGZZ19} by first relating (iii\textsubscript{0}) with the asymptotic number of square-tiled surfaces, and combinatorial techniques to evaluate this number. \\[.5ex]
\noindent
	(iii\textsubscript{0})~=~(iii). Proved in \cite[Corollary 3.6]{ABCDGLW19}. \\[.5ex]
\noindent (ii)~=~(iii)~=~(iv). Proved in \cite{ABCDGLW19} via the study of (asymptotics of integrals of) statistics of hyperbolic lengths of multicurves and another Laplace transform method. \\[.5ex]
\noindent
	(iii)~=~(iv\textsubscript{+}). Follows from Mirzakhani's work, together with the specialisation of \cite[Theorem~10.1]{ABO17} to the Mirzakhani GR initial data twisted by $f(\ell;s)~=~\frac{1}{e^{s\ell} - 1}$ and the Laplace transform method. \\[.5ex]
\noindent
	(ii)~=~(v)~=~(vi)~=~(vi\textsubscript{+}) are the results of the previous section. Notice that (vi)~=~(vi\textsubscript{+}) is just a consequence of the polynomiality of the $1$-parameter family of Masur--Veech polynomials.

\medskip

One can also give an interpretation of the area Siegel--Veech constant via derivative statistics of combinatorial lengths by adapting \cite[Section 4]{ABCDGLW19}. 

\medskip

There are however some missing paths that would be desirable to complete in the future. 

\medskip

A first question consists in proving (iii)~=~(v) via the rescaling flow (Section~\ref{sec:hyp:comb}). We can certainly write
\begin{equation}
	\int_{\mathcal{M}_{g,n}(\beta L)}
		\mathscr{B}_{g,n} \, \dd\mu_{\textup{WP}}
	=
	\int_{\mathcal{M}_{g,n}^{\textup{comb}}(L)}
		\beta^{6g - 6 + 2n}(\mathscr{B}_{g,n} \circ R_{\beta}) \, J_{\beta} \, \dd\mu_{\textup{K}},
	\qquad
	J_{\beta} = \frac{1}{\beta^{6g - 6 + 2n}} \, \frac{R_{\beta}^*\dd\mu_{\textup{WP}}}{\dd \mu_{\textup{K}}}
\end{equation}
We know by Theorem~\ref{thm:Mon:Do} due to Mondello that $J_{\beta}$ converges pointwise to $1$, and by Lemma~\ref{lem:B:hyp:comb} that $\beta^{6g - 6 + 2n}(\mathscr{B}_{g,n} \circ R_{\beta})$ converges to $\mathscr{B}_{g,n}^{\textup{comb}}$ uniformly on compacts. This is however not sufficient, as we would need an effective and integrable enough bound independent of $\beta$ to conclude that (iii)~=~(v) by dominated convergence, and such a bound is not currently available. Note that one cannot hope for a bound by constant, because $\mathscr{B}_{\Sigma}^{\textup{comb}}(\GG)$ can diverge when $\sys_{\GG} \rightarrow 0$. Getting effective and uniform bounds for the Jacobian $J_{\beta}$ over $\mathcal{T}_{\Sigma}^{\textup{comb}}$ is a question of broader interest: it would allow to study the behaviour for large boundary lengths $L$ of the integral on $\mathcal{M}_{g,n}(L)$ of a larger class of functions.

\medskip

So far the combinatorial expressions (v), (vi) and (vi\textsubscript{+}) are only linked to (i) or (iii\textsubscript{0}) indirectly, via the sum over stable graphs (ii). We believe it would be interested to find a direct geometric proof of the equality (i) = (v) in a similar way that (i) = (iii\textsubscript{0}) was proved in \cite{Mir08earth}, \emph{i.e.} give a direct relation between the Masur--Veech measure on the moduli space of quadratic differentials and the combinatorial enumerative geometry of measured foliations. If such a link could be made precise, one may hope that the role of hyperbolic geometry for the study of flat surfaces could be substituted by combinatorial geometry, which may bring some useful simplifications.

\clearpage
\appendix
\section{Topology on the combinatorial Teichm\"uller spaces}
\label{app:topology:comb}

The combinatorial Teichm\"uller spaces have been endowed with a natural topology associated to a dual construction: the proper arc complex. The main references for this discussion are \cite{Luo07,Mon09}.

\begin{defn}
	Fix a connected bordered surface $\Sigma$. Define the \emph{arc complex} $\mathscr{A}_{\Sigma}$ to be the simplicial complex whose vertices are isotopy classes $\alpha$ of arcs in $\Sigma$ with endpoints on the boundary $\de\Sigma$ which are homotopically non-trivial relatively to $\de \Sigma$. A simplex in $\mathscr{A}_{\Sigma}$ is a collection $\bm{\alpha} = (\alpha_1,\dots,\alpha_k)$ of distinct vertices such that the arcs $\alpha_i$ admit representatives which do not intersect. The non-proper subcomplex $\mathscr{A}_{\Sigma}^{\infty}$ of $\mathscr{A}_{\Sigma}$ consists of those simplices $\bm{\alpha}$ such that one connected component of $\Sigma \setminus \bigcup_{i=1}^k \alpha_i$ is not simply connected. The simplices in $\mathscr{A}_{\Sigma} \setminus \mathscr{A}_{\Sigma}^{\infty}$ are called proper.
\end{defn}

Consider the geometric realisation spaces $|\mathscr{A}_{\Sigma}|$ and $|\mathscr{A}_{\Sigma}^{\infty}|$. The geometric realisation of a simplicial complex comes with two natural topologies: the coherent topology, namely the finest topology that makes the realisation of all simplicial maps continuous, and the metric topology, for which every $k$-simplex is isometric to the standard simplex $\Delta^k \subset \RR^{k+1}$ and every attachment map is a local isometry. The metric topology on $|\mathscr{A}_{\Sigma}|$ is coarser than the coherent one, but they agree where the complex is locally finite, and this is the case for $|\mathscr{A}_{\Sigma}| \setminus |\mathscr{A}_{\Sigma}^{\infty}|$.

\begin{figure}[ht]
\centering

	\caption{The simplicial structure of $|\mathscr{A}_{P}|$ and $|\mathscr{A}_{P}^{\infty}|$ and the set $|\mathscr{A}_{P}| \setminus |\mathscr{A}_{P}^{\infty}|$ associated to a pair of pants $P$. The latter is in natural bijection with the slice $\set{\GG \in \mathcal{T}_P^{\textup{comb}} | \ell_{\GG}(\de_1 P) + \ell_{\GG}(\de_2 P) + \ell_{\GG}(\de_3 P) = 1}$.}
\end{figure}

\medskip

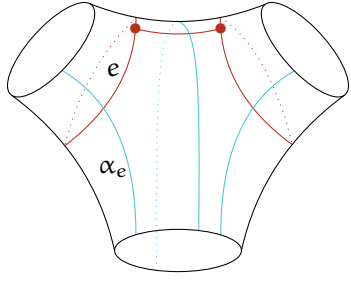
\begin{figure}[t]
	\centering
	\begin{tikzpicture}[x=1pt,y=1pt,scale=.5]
		\draw [BrickRed] (138.437, 647.363)  .. controls (174.145, 674.6563) and (192.079, 707.0133) .. (192.239, 744.434);
		\draw [BrickRed, dotted] (192.239, 744.434)  .. controls (165.1863, 717.6947) and (147.2523, 685.3377) .. (138.437, 647.363);
		\draw [BrickRed] (255.958, 744.512)  .. controls (256.1287, 706.6633) and (274.2667, 674.544) .. (310.372, 648.154);
		\draw [BrickRed, dotted] (255.963, 744.43)  .. controls (295.831, 719.864) and (310.372, 648.154) .. (310.372, 648.154);
		\node [BrickRed]  at (191.9284, 735.8591) {$\bullet$};
		\node [BrickRed]  at (256.3274, 735.8938) {$\bullet$};
		\draw [BrickRed] (191.928, 735.812)  .. controls (213.3093, 730.604) and (234.776, 730.6313) .. (256.328, 735.894);
		\draw [SkyBlue] (136.2739, 704.2748)  .. controls (173.64, 682.5228) and (192.2649, 641.0878) .. (192.1485, 579.9698);
		\draw [SkyBlue] (311.9937, 704.0063)  .. controls (274.7567, 682.7512) and (256.2126, 641.3548) .. (256.3613, 579.8169);
		\draw [SkyBlue] (239.6786, 583.1224)  .. controls (240.6412, 688.1337) and (240.0562, 736.2139) .. (224.3392, 741.3805);
		\draw [SkyBlue, dotted] (207.8515, 552.9327)  .. controls (207.3401, 687.9037) and (208.3974, 735.865) .. (223.9771, 741.3811);

		\draw (128, 720) circle[radius=1, cm={32,32,-8,16,(0,0)}];
		\draw (320, 720) circle[radius=1, cm={-32,32,8,16,(0,0)}];
		\draw (224, 568) ellipse[x radius=48, y radius=16];
		\draw (98.2606, 685.437)  .. controls (139.4202, 655.1457) and (165.3333, 616) .. (176, 568);
		\draw (272, 568)  .. controls (282.6667, 616) and (308.6583, 655.1993) .. (349.975, 685.598);
		\draw (156.827, 755.056)  .. controls (201.609, 737.0187) and (246.0413, 736.824) .. (290.124, 754.472);

		\node at (175, 705) {$e$};
		\node at (175, 630) {$\alpha_e$};

	\end{tikzpicture}
	\caption{Example of duality between a combinatorial structure (red) and an arc system (light blue).
	\label{fig:edge:arc:duality}}
\end{figure}

Consider then the topological space $(|\mathscr{A}_{\Sigma}| \setminus |\mathscr{A}_{\Sigma}^{\infty}|) \times \RR_{+}$, called \emph{metrised arc complex}, whose points are of the form $x = \sum_{i=1}^k \ell_i \, \alpha_i$, where $\ell_i > 0$ and $(\alpha_1,\dots,\alpha_k)$ is a proper simplex. There is a bijection between the spaces $\mathcal{T}_{\Sigma}^{\textup{comb}}$ and $(|\mathscr{A}_{\Sigma}| \setminus |\mathscr{A}_{\Sigma}^{\infty}|) \times \RR_{+}$. Given a combinatorial structure $\GG$, we define for each edge $e$ the dual arc $\alpha_{e}$ connecting the two (possibly equal) boundary components on the two sides of $e$ (see Figure~\ref{fig:edge:arc:duality}). Thus, we define the map $\mathcal{T}_{\Sigma}^{\textup{comb}} \to (|\mathscr{A}_{\Sigma}| \setminus |\mathscr{A}_{\Sigma}^{\infty}|) \times \RR_{+}$ by setting $\GG \mapsto \sum_{e \in E_{\GG}} \ell_{\GG}(e) \, \alpha_e$. The map is clearly invertible and we topologise the combinatorial Teichm\"uller space $\mathcal{T}_{\Sigma}^{\textup{comb}}$ by pulling back the topology of the proper arc complex. Further, through the spine construction, it is possible to prove the following

\begin{thm}\cite{Luo07,Mon09}
	There is an homeomorphism
	\begin{equation}
		\sp \colon \mathcal{T}_{\Sigma} \longrightarrow (|\mathscr{A}_{\Sigma}| \setminus |\mathscr{A}_{\Sigma}^{\infty}|) \times \RR_{+},
	\end{equation}
	equivariant under the action of the mapping class group.
	\hfill $\blacksquare$
\end{thm}

\section{Examples of cutting, gluing, and combinatorial Fenchel--Nielsen coordinates}
\label{app:cut:gluing}
\textsc{Cutting and gluing.}
In order to make the presentation of the cutting and gluing algorithms clearer, we present some neat examples of these procedures. The first two (Figures~\ref{fig:cut:glue:04:1} and \ref{fig:cut:glue:04:2}) cover the case of a sphere $\Sigma$ with four boundary components, the last two (Figures~\ref{fig:cut:glue:11:1} and \ref{fig:cut:glue:11:2}) cover the case of a torus $T$ with one boundary component.

\medskip

In all examples, the cutting algorithm is presented in lexicographic order, \emph{i.e.} the images are labelled by (a), (b), (c), \emph{et cetera}, while the gluing algorithm is presented with the same images, but in reversed lexicographic order. The combinatorial structure is depicted in red, the associated measured foliation in blue (only the singular leaves are reported). The cutting curve $\gamma$ is depicted in green, and it coincide with the curve obtain in the gluing algorithm after identification of two boundary components $\gamma_- \sim \gamma_+$. The component $\gamma_-$ is always located on the left side of the figure, while $\gamma_+$ is located on the right side. The identification points $p_{\pm} \in \gamma_{\pm}$ are depicted in grey, and no twist is performed (that is, $p_-$ is identified with $p_+$). Further, the letters $a,b,c,\dots$ are referring to edge lengths of the cutting process, while the letters $r,s,t,\dots$ are referring to edge lengths of the gluing process.

\medskip

\textsc{Combinatorial Fenchel--Nielsen coordinates.}
We present two computations of combinatorial Fen\-chel--Niel\-sen coordinates on a torus $T$ with one boundary component, relative to two different cells (Figure~\ref{fig:FN:1} and Figure~\ref{fig:FN:2}). Furthermore, an illustration of Penner's formulae (Proposition~\ref{prop:Penner:formulae}) is presented.

\medskip

In all examples, the combinatorial structure is depicted in red, the associated measured foliation in blue (only the singular leaves are reported), the pants decompositions $\mathscr{P},\mathscr{P}'$ in green and the collection $\mathscr{S},\mathscr{S}'$ in yellow.

\clearpage
\begin{minipage}[t]{.46\textwidth}
	\textsc{Cutting.}
	Consider $\GG \in \mathcal{T}_{\Sigma}^{\textup{comb}}$ and $\gamma$ a cutting curve as in Figure~\ref{fig:cut:glue:04:1:a}. We have
	\[
		\ell_{\GG}(\gamma) = a+b+e+f.
	\]
	After cutting, we obtain two pairs of pants $P_{\pm}$ and two combinatorial structures $\GG_{\pm} \in \mathcal{T}_{P_\pm}^{\textup{comb}}$.
\end{minipage}
\hfill
\begin{minipage}[t]{.46\textwidth}
	\textsc{Gluing.}
	Consider two combinatorial structures $\GG_{\pm} \in \mathcal{T}_{P_{\pm}}^{\textup{comb}}$ and two boundary components $\gamma_{\pm}$ of $P_{\pm}$ of the same length as in Figure~\ref{fig:cut:glue:04:1:f}. We have
	\[
		\ell_{\GG_-}(\gamma_-) = r + s,
		\qquad
		\ell_{\GG_+}(\gamma_+) = u + v.
	\]
	After gluing (with a choice of $p_{\pm}$), we obtain a sphere $\Sigma$ with four boundary components and a combinatorial structure $\GG \in \mathcal{T}_{\Sigma}^{\textup{comb}}$.
\end{minipage}
\vspace*{\fill}
\begin{figure}[h!]
	\centering
	\begin{subfigure}[t]{.48\textwidth}
	\centering

	\caption{}
	\label{fig:cut:glue:04:1:a}
	\end{subfigure}
	\begin{subfigure}[t]{.48\textwidth}
	\centering
	%
	\caption{}
	\label{fig:cut:glue:04:1:b}
	\end{subfigure}
	\begin{subfigure}[t]{.48\textwidth}
	\centering
	%
	\caption{}
	\label{fig:cut:glue:04:1:c}
	\end{subfigure}
	\begin{subfigure}[t]{.48\textwidth}
	\centering
	%
	\caption{}
	\label{fig:cut:glue:04:1:d}
	\end{subfigure}
	\begin{subfigure}[t]{.48\textwidth}
	\centering
	%
	\caption{}
	\label{fig:cut:glue:04:1:e}
	\end{subfigure}
	\begin{subfigure}[t]{.48\textwidth}
	\centering
	%
	\caption{}
	\label{fig:cut:glue:04:1:f}
	\end{subfigure}
	\caption{}
	\label{fig:cut:glue:04:1}
\end{figure}
\vfill
~

\clearpage
\afterpage{
\begin{landscape}
\thispagestyle{empty}
\begin{minipage}[t]{.62\textheight}
	\textsc{Cutting.}
	Consider $\GG \in \mathcal{T}_{\Sigma}^{\textup{comb}}$ and $\gamma$ a cutting curve as in Figure~\ref{fig:cut:glue:04:2:a}. We have
	\[
		\ell_{\GG}(\gamma) = 2a + b + 2c + 2d + e.
	\]
	After cutting, we obtain two pairs of pants $P_{\pm}$ and two combinatorial structures $\GG_{\pm} \in \mathcal{T}_{P_\pm}^{\textup{comb}}$. From Figure~\ref{fig:cut:glue:04:2:c} to Figure~\ref{fig:cut:glue:04:2:d}, a Whitehead move is performed.
\end{minipage}
\hfill
\begin{minipage}[t]{.62\textheight}
	\textsc{Gluing.}
	Consider two combinatorial structures $\GG_{\pm} \in \mathcal{T}_{P_{\pm}}^{\textup{comb}}$ and two boundary components $\gamma_{\pm}$ of $P_{\pm}$ of the same length as in Figure~\ref{fig:cut:glue:04:2:g}. We have
	\[
		\ell_{\GG_-}(\gamma_-) = r + s + 2t,
		\qquad
		\ell_{\GG_+}(\gamma_+) = u + v + 2w.
	\]
	After gluing (with a choice of $p_{\pm}$), we obtain a sphere $\Sigma$ with four boundary components and a combinatorial structure $\GG \in \mathcal{T}_{\Sigma}^{\textup{comb}}$. From Figure~\ref{fig:cut:glue:04:2:d} to Figure~\ref{fig:cut:glue:04:2:c}, a Whitehead move is performed.
\end{minipage}
\begin{figure}[h]
	\centering
	\begin{subfigure}[t]{.325\textwidth}
	\centering

	\caption{}
	\label{fig:cut:glue:04:2:a}
	\end{subfigure}
	\begin{subfigure}[t]{.325\textwidth}
	\centering
	%
	\caption{}
	\label{fig:cut:glue:04:2:b}
	\end{subfigure}
	\begin{subfigure}[t]{.325\textwidth}
	\centering
	%
	\caption{}
	\label{fig:cut:glue:04:2:c}
	\end{subfigure}
	\begin{subfigure}[t]{.325\textwidth}
	\centering
	%
	\caption{}
	\label{fig:cut:glue:04:2:d}
	\end{subfigure}
	\begin{subfigure}[t]{.38\textwidth}
	\centering
	%
	\caption{}
	\label{fig:cut:glue:04:2:e}
	\end{subfigure}
	\begin{subfigure}[t]{.38\textwidth}
	\centering
	%
	\caption{}
	\label{fig:cut:glue:04:2:f}
	\end{subfigure}
	\begin{subfigure}[t]{.45\textwidth}
	\centering
	%
	\caption{}
	\label{fig:cut:glue:04:2:g}
	\end{subfigure}
	\caption{}
	\label{fig:cut:glue:04:2}
\end{figure}
\vfill
\end{landscape}
}

\clearpage
\begin{minipage}[t]{.46\textwidth}
	\textsc{Cutting.}
	Consider $\GG \in \mathcal{T}_{\Sigma}^{\textup{comb}}$ and $\gamma$ a cutting curve as in Figure~\ref{fig:cut:glue:04:1:a}. We have
	\[
		\ell_{\GG}(\gamma) = a + c.
	\]
	After cutting, we obtain a pairs of pants $P$ and a combinatorial structure $\GG' \in \mathcal{T}_{P}^{\textup{comb}}$.
\end{minipage}
\hfill
\begin{minipage}[t]{.46\textwidth}
	\textsc{Gluing.}
	Consider a combinatorial structure $\GG' \in \mathcal{T}_{P}^{\textup{comb}}$ and two boundary components $\gamma_{\pm}$ of $P$ of the same length as in Figure~\ref{fig:cut:glue:04:1:f}. We have
	\[
		\ell_{\GG'}(\gamma_-) = r,
		\qquad
		\ell_{\GG'}(\gamma_+) = s.
	\]
	After gluing (with a choice of $p_{\pm}$), we obtain a torus $T$ with one boundary component and a combinatorial structure $\GG \in \mathcal{T}_{T}^{\textup{comb}}$.
\end{minipage}
\vspace*{\fill}
\begin{figure}[h!]
	\centering
	\begin{subfigure}[t]{.48\textwidth}
	\centering

	\caption{}
	\label{fig:cut:glue:11:1:a}
	\end{subfigure}
	\begin{subfigure}[t]{.48\textwidth}
	\centering
	%
	\caption{}
	\label{fig:cut:glue:11:1:b}
	\end{subfigure}
	\begin{subfigure}[t]{.48\textwidth}
	\centering
	%
	\caption{}
	\label{fig:cut:glue:11:1:c}
	\end{subfigure}
	\begin{subfigure}[t]{.48\textwidth}
	\centering
	%
	\caption{}
	\label{fig:cut:glue:11:1:d}
	\end{subfigure}
	\begin{subfigure}[t]{.48\textwidth}
	\centering
	%
	\caption{}
	\label{fig:cut:glue:11:1:e}
	\end{subfigure}
	\begin{subfigure}[t]{.48\textwidth}
	\centering
	%
	\caption{}
	\label{fig:cut:glue:11:1:f}
	\end{subfigure}
	\caption{}
	\label{fig:cut:glue:11:1}
\end{figure}
\vfill
~

\clearpage
\begin{minipage}[t]{.46\textwidth}
	\textsc{Cutting.}
	Consider $\GG \in \mathcal{T}_{\Sigma}^{\textup{comb}}$ and $\gamma$ a cutting curve as in Figure~\ref{fig:cut:glue:04:2:a}. We have
	\[
		\ell_{\GG}(\gamma) = 2a + b + c.
	\]
	After cutting, we obtain a pairs of pants $P$ and a combinatorial structure $\GG' \in \mathcal{T}_{P}^{\textup{comb}}$.
\end{minipage}
\hfill
\begin{minipage}[t]{.46\textwidth}
	\textsc{Gluing.}
	Consider a combinatorial structure $\GG' \in \mathcal{T}_{P}^{\textup{comb}}$ and two boundary components $\gamma_{\pm}$ of $P$ of the same length as in Figure~\ref{fig:cut:glue:04:2:f}. We have
	\[
		\ell_{\GG'}(\gamma_-) = r + s,
		\qquad
		\ell_{\GG'}(\gamma_+) = r + t.
	\]
	After gluing (with a choice of $p_{\pm}$), we obtain a torus $T$ with one boundary component and a combinatorial structure $\GG \in \mathcal{T}_{T}^{\textup{comb}}$.
\end{minipage}
\vspace*{\fill}
\begin{figure}[h!]
	\centering
	\begin{subfigure}[t]{.48\textwidth}
	\centering

	\caption{}
	\label{fig:cut:glue:11:2:a}
	\end{subfigure}
	\begin{subfigure}[t]{.48\textwidth}
	\centering
	%
	\caption{}
	\label{fig:cut:glue:11:2:b}
	\end{subfigure}
	\begin{subfigure}[t]{.48\textwidth}
	\centering
	%
	\caption{}
	\label{fig:cut:glue:11:2:c}
	\end{subfigure}
	\begin{subfigure}[t]{.48\textwidth}
	\centering
	%
	\caption{}
	\label{fig:cut:glue:11:2:d}
	\end{subfigure}
	\begin{subfigure}[t]{.48\textwidth}
	\centering
	%
	\caption{}
	\label{fig:cut:glue:11:2:e}
	\end{subfigure}
	\begin{subfigure}[t]{.48\textwidth}
	\centering
	%
	\caption{}
	\label{fig:cut:glue:11:2:f}
	\end{subfigure}
	\caption{}
	\label{fig:cut:glue:11:2}
\end{figure}
\vfill
~

\clearpage
\textsc{Combinatorial Fenchel--Nielsen coordinates.}
Consider the combinatorial structure $\GG \in \mathcal{T}_{T}^{\textup{comb}}$ on a torus $T$ with one boundary component as in Figure~\ref{fig:FN:1}. We have
\[
	L = 2a + 2b + 2c.
\]
Further, the combinatorial Fenchel--Nielsen coordinates with respect to the seamed pants decomposition $(\mathscr{P},\mathscr{S}) = (\gamma,\beta)$ of Figure~\ref{fig:FN:1:a} are computed with the help of Figure~\ref{fig:FN:1:b} and are given by 
\[
	\ell = a + c,
	\qquad\qquad
	\tau = c, 
\]
while the combinatorial Fenchel--Nielsen coordinates with respect to the seamed pants decomposition $(\mathscr{P}',\mathscr{S}') = (\gamma',\beta')$ of Figure~\ref{fig:FN:1:c} are computed with the help of Figure~\ref{fig:FN:1:d} and are given by 
\[
	\ell' = b + c,
	\qquad\qquad
	\tau' = -c.
\]
This is in accordance with Penner's formulae:
\[
\begin{split}
	\ell' & = |\tau| + \left[\tfrac{L - 2\ell}{2} \right]_{+} \\
	& = |c| + b \\
	& = b + c,
\end{split}
	\qquad\qquad
\begin{split}
	\tau' & = - \sgn(\tau)
	\left|
		\ell - \left[ \tfrac{L - 2\ell'(\ell,\tau)}{2} \right]_{+}
	\right| \\
	& = - \big|(a + c) - a \big| \\
	& = -c.
\end{split}
\]
\vspace*{\fill}
\begin{figure}[h]
	\centering
	\begin{subfigure}[t]{.46\textwidth}
	\centering

	\caption{}
	\label{fig:FN:1:a}
	\end{subfigure}
	\begin{subfigure}[t]{.46\textwidth}
	\centering
	%
	\caption{}
	\label{fig:FN:1:b}
	\end{subfigure}
	\begin{subfigure}[t]{.46\textwidth}
	\centering
	%
	\caption{}
	\label{fig:FN:1:c}
	\end{subfigure}
	\begin{subfigure}[t]{.46\textwidth}
	\centering
	%
	\caption{}
	\label{fig:FN:1:d}
	\end{subfigure}
	\caption{}
	\label{fig:FN:1}
\end{figure}
\vfill
~

\clearpage
\textsc{Combinatorial Fenchel--Nielsen coordinates.}
Consider the combinatorial structure $\GG \in \mathcal{T}_{T}^{\textup{comb}}$ on a torus $T$ with one boundary component as in Figure~\ref{fig:FN:2}. We have
\[
	L = 2a + 2b + 2c.
\]
Further, the combinatorial Fenchel--Nielsen coordinates with respect to the seamed pants decomposition $(\mathscr{P},\mathscr{S}) = (\gamma,\beta)$ of Figure~\ref{fig:FN:2:a} are computed with the help of Figure~\ref{fig:FN:2:b} and are given by 
\[
	\ell = 2a + b + c,
	\qquad\qquad
	\tau = a + b, 
\]
while the combinatorial Fenchel--Nielsen coordinates with respect to the seamed pants decomposition $(\mathscr{P}',\mathscr{S}') = (\gamma',\beta')$ of Figure~\ref{fig:FN:2:c} are computed with the help of Figure~\ref{fig:FN:2:d} and are given by 
\[
	\ell' = a + b,
	\qquad\qquad
	\tau' = - 2a - b.
\]
This is in accordance with Penner's formulae:
\[
\begin{split}
	\ell' & = |\tau| + \left[\tfrac{L - 2\ell}{2} \right]_{+} \\
	& = |a + b| + 0 \\
	& = a + b,
\end{split}
	\qquad\qquad
\begin{split}
	\tau' & = - \sgn(\tau)
		\left|
			\ell - \left[ \tfrac{L - 2\ell'(\ell,\tau)}{2} \right]_{+}
		\right| \\
	& = - \big| (2a + b + c) - c \big| \\
	& = - 2a - b.
\end{split}
\]
\vspace*{\fill}
\begin{figure}[h]
	\centering
	\begin{subfigure}[t]{.46\textwidth}
	\centering

	\caption{}
	\label{fig:FN:2:a}
	\end{subfigure}
	\begin{subfigure}[t]{.46\textwidth}
	\centering
	%
	\caption{}
	\label{fig:FN:2:b}
	\end{subfigure}
	\begin{subfigure}[t]{.46\textwidth}
	\centering
	%
	\caption{}
	\label{fig:FN:2:c}
	\end{subfigure}
	\begin{subfigure}[t]{.46\textwidth}
	\centering
	%
	\caption{}
	\label{fig:FN:2:d}
	\end{subfigure}
	\caption{}
	\label{fig:FN:2}
\end{figure}
\vfill
~

\clearpage

\section{Integral points in the moduli space and factors of \texorpdfstring{$2$}{2}}
\label{app:integral:points:factor:2}

Let $G$ be a ribbon graph of type $(g,n)$, and $S \subset E_{G}$ such that $\# S = n$ and the dual graph $G_S^*$ of $S$ (considered as a subgraph of $G$) is connected and has a single cycle, of odd length. We label the edges so that $S=\{e_1,\dots,e_n\}$ and $E_{G}=\{e_1,\dots,e_{6g-g+3n}\}$. If $\bm{G}$ is a metric structure on $G$, we call $\ell_i = \ell_{\bm{G}}(e_i)$. The multiplicity $A_{i,e} \in \{0,1,2\}$ of an edge $e \in E_G$ around the $i$-th face define an \emph{adjacency matrix} $A$ of size $n \times(6g - 6 + 3n)$. We recall the following well-known fact (see \textit{e.g.} \cite[Theorem~2.2]{DE14}).

\begin{lem} \label{lem:dav:eyn}
	The restriction $\hat{A}$ of the adjacency matrix $R$ to the first $n$ columns is invertible, and $|\det(\hat{A})| = 2$.
	\hfill $\blacksquare$
\end{lem}

It is useful for \S~\ref{subsubsec:disc:integr:TR} to characterise the integral points in $\mathcal{M}_{g,n}^{\textup{comb}}(L)$ in terms of integral points in $\mathcal{M}_{g,n}^{\textup{comb}}$.

\begin{lem}\label{lem:index:2}
	$\ell_{n+1},\dots,\ell_{6g-6+3n} \in \ZZ$ and $\sum\limits_{i=1}^n L_i$ is even if and only if $\ell_1,\dots,\ell_{6g-6+3n} \in \ZZ$.
\end{lem}

\begin{proof}
	It is sufficient to prove the implication, since the converse is obvious. The result is elementary and well-known, but as we did not identify a reference for it, we include a proof. Suppose that $\ell_{n+1},\dots,\ell_{6g-6+3n}$ are integers and $\sum_{i = 1}^n L_i =0 \pmod{2}$. By hypothesis, $G_{S}^*$ is the union of trees rooted at a cycle with $2p+3$ edges for some $p \geq 0$, and its vertices are labelled from $1$ to $n$. For $e \in S$ and $i \in \{1,\dots,n\}$, let $d_{e,i}$ be the graph distance in $G_{S}^*$ between $e^*$ (the dual of $e$) and the vertex $i$. If $e$ is adjacent to the face $i$, then $d_{e,i}=0$. There is a natural notion of descendent vertices of an edge belonging to a tree of $G_{S}^*$, given that the trees are rooted at the cycles of $G_{S}^*$. The inverse of $\hat{A}$ can now be made explicit.
	\begin{itemize}
		\item If $e^*$ does not belong to the cycle of $G_{S}^*$
		\[
			\hat{A}^{-1}_{e,i} = 
			\begin{cases}
				(-1)^{d_{e,i}}	& \text{if $i$ is a descendent of $e^*$} \\
				0 				& \text{otherwise.}
			\end{cases}
		\]
		\item If $e^*$ belongs to the cycle of $G_{S}^*$,
		\[
			\hat{A}^{-1}_{e,i} = \frac{(-1)^{d_{e,i}}}{2}.
		\]
	\end{itemize}
	If $e^*$ does not belong to the cycle of $G_S^*$, for all $i\in \{1,\dots,n\}$ we have $\hat{A}^{-1}_{e,i}\in\{-1,0,1\}$, so
	\[
		\ell_{e} = \sum_{i=1}^{n} \hat{A}^{-1}_{e,i} \left(
				L_i - \sum_{k = n+1}^{6g-6+3n} A_{i,e_k} \ell_k
			\right)
	\]
	belongs to $\ZZ$. If $e^*$ belongs to the cycle of $G_S^*$, for all $i \in \{1,\dots,n\}$ we have $2\hat{A}^{-1}_{e,i} = 1 \pmod{2}$. This implies the following:
	\[
		\begin{split}
			2 \ell_{e}
			& = 
			\sum_{i = 1}^{n} 2 \hat{A}^{-1}_{e,i} \left(
					L_i - \sum_{k = n+1}^{6g-6+3n} A_{i,e_k} \ell_k
				\right) \\
			& =
			\sum_{i = 1}^{n} L_i - \sum_{k = n+1}^{6g-6+3n} \left(
					\sum_{i=1}^{n} A_{i,e_k}
				\right) \ell_k \,\, \pmod{2} \\
			& =
			0 \pmod{2},
	\end{split}
	\]
	where the second to last line comes from the hypothesis $\sum_{i=1}^n L_i =0 \pmod{2}$ and $\sum_{i=1}^n A_{i,e} = 2$. The result of this calculation is that $\ell_e \in \ZZ$. 
\end{proof}

\newpage
\section*{Index of notations}
\label{sec:notation}
%
{\small
\renewcommand{\arraystretch}{1.3}
\begin{center}
\begin{tabulary}{\textwidth}{c|c|L}
$G$, $\bm{G}$
	& Def.~\ref{defn:MRG}
	& Ribbon graph, metric ribbon graph
	\\
$\Aut(G) \supset \Aut(\bm{G})$
	& Def.~\ref{defn:automorphism}
	& Their automorphism groups
	\\	
$|\bm{G}|$
	& Def.~\ref{defn:geom:realis:MRG}
	& Geometric realisation of a metric ribbon graph
	\\
$(\bm{G},f)$, $\GG$
	& Def.~\ref{defn:comb:struct}
	& Combinatorial marking, combinatorial structure
	\\
$S_{\Sigma}$
	& Sec.~\ref{subsubsec:notation:curves}
	& Set of homotopy classes of simple closed curves on $\Sigma$
	\\
$S_{\Sigma}^{\circ} \subset S_{\Sigma}$
	& Sec.~\ref{subsubsec:notation:curves}
	& Subset of essential simple closed curves $\Sigma$
	\\
$M_{\Sigma}$, $M_{\Sigma}'$
	& Sec.~\ref{subsubsec:notation:curves}
	& Set of multicurves of $\Sigma$, set of primitive multicurves on $\Sigma$
	\\
$\mathcal{P}_{\Sigma,m_0}$, $\mathcal{P}_{\Sigma,m_0}^{\varnothing}$, $\mathcal{P}_{\Sigma,m_0}^{m}$
	& Def.~\ref{defn:homotopy:class:pop}
	& Sets of homotopy classes of embedded pairs of pants in $\Sigma$: with $m_0$ as a boundary, with $m_0$ and an essential simple curve as boundaries, with $m_0$ and $m$ as boundaries
	\\
$\mathcal{A}_{\Sigma,m_0}$, $\mathcal{A}_{\Sigma,m_0}^{\varnothing}$, $\mathcal{A}_{\Sigma,m_0}^{m}$
	& Def.~\ref{defn:homotopy:class:arcs}
	& Set of non-trivial homotopy classes of proper embeddings of arcs: with initial point on $m_0$, with both endpoints on $m_0$, with initial point on $m_0$ and final point on $m$
	\\
$Q \colon \mathcal{A}_{\Sigma,m_0} \to \mathcal{P}_{\Sigma,m_0}$
	& Eqn.~\eqref{defn:arcs:to:pants}
	& Map associating to an arc the pair of pants obtained by going around the paths and the boundaries at its endpoints
	\\
$\Mod_{\Sigma}$, $\Mod_{\Sigma}^{\partial}$
	& Def.~\ref{defn:bord:MCG}
	& Mapping class group of $\Sigma$, pure mapping class group of $\Sigma$
	\\
$\mathcal{T}_{\Sigma}^{\textup{comb}}$, $\mathcal{T}_{\Sigma}^{\textup{comb}}(L)$
	& Def.~\ref{defn:comb:struct}
	& Combinatorial Teichm\"uller space of $\Sigma$, with fixed boundary lengths $L \in \RR_+^n$
	\\
$\mathcal{M}_{g,n}^{\textup{comb}}$, $\mathcal{M}_{g,n}^{\textup{comb}}(L)$
	& Def.~\ref{defn:Mgncomb}
	& Corresponding moduli spaces
	\\
$\mathcal{M}^{\textup{comb},\ZZ}_{g,n}$, $\mathcal{M}_{g,n}^{\textup{comb},\ZZ}(L)$
	& Def.~\ref{defn:MgncombZ}
	& Their subsets of integral points
	\\
$\mathcal{T}_{\Sigma}$, $\mathcal{T}_{\Sigma}(L)$
	& Sec.~\ref{subsubsec:ordinary:combinatorial}
	& Teichm\"uller space of $\Sigma$, with fixed boundary lengths $L \in \RR_{+}^n$
	\\
$\mathcal{M}_{g,n}$, $\mathcal{M}_{g,n}(L)$
	& Sec.~\ref{subsubsec:ordinary:combinatorial}
	& Corresponding moduli spaces of bordered Riemann surfaces
	\\
$\mathfrak{M}_{g,n}$, $\overline{\mathfrak{M}}_{g,n}$
	& Sec.~\ref{sec:comb:spaces}
	& Moduli space of genus $g$ smooth complex curves with $n$ marked points, Deligne--Mumford compactification
	\\
$\mathfrak{Z}_{\Sigma,G}$, $\mathfrak{Z}_{\Sigma,G}(L)$
	& Sec.~\ref{sec:comb:spaces}
	& Cell of $\mathcal{T}_{\Sigma}$ coming from a marked ribbon graph $(G,f)$, with fixed boundary lengths $L \in \RR_+^n$
	\\
$\sp \colon \mathcal{T}_{\Sigma} \to \mathcal{T}_{\Sigma}^{\textup{comb}}$
	& Eqn.~\eqref{defn:spine:map}
	& The spine homeomorphism
	\\
$\Phi_{\beta}$, $R_{\beta}$
	& Defn.~\ref{defn:rescaling:flow}
	& The rescaling flows
	\\
${\rm MF}_{\Sigma}^{\star}$
	& Def.~\ref{defn:measured:foliations}
	& Measured foliations on $\Sigma$ with all types of behaviour at $\partial \Sigma$
	\\
${\rm MF}_{\Sigma}$
	& Sec.~\ref{subsec:unit:balls:hyp}
	& Measured foliations on $\Sigma$ with points at the boundary of parallel type
	\\
$\mathcal{F}_{\ast} \colon \mathcal{T}_{\Sigma}^{\textup{comb}} \to {\rm MF}_{\Sigma}^{\star}$
	& Eqn.~\eqref{eqn:comb:Teich:to:MF}
	& Inclusion of the combinatorial Teichm\"uller space into the space of measured foliations
	\\
$\ell_{\ast} \colon \mathcal{T}_{\Sigma}^{\textup{comb}} \to \RR_{+}^{S_{\Sigma}}$
	& Thm.~\ref{thm:length:functional}
	& Length spectrum map
	\\
$\omega_{\textup{K}}$, $\mu_{\textup{K}}$
	& Eqn.~\eqref{eqn:Kon:form}, \eqref{eqn:Kon:volume:form}
	& Kontsevich $2$-form and the associated measure
	\\ 
$\omega_{\textup{WP}}$, $\mu_{\textup{WP}}$
	& Sec.~\ref{subsubsec:integr:TR}
	& Weil--Petersson symplectic form and the associated measure
	\\ 
$\Xi_{g,n}$
	& Sec.~\ref{subsubsec:GRcomb}
	& Functions on $\mathcal{M}_{g,n}^{\textup{comb}}$ that are output of combinatorial GR
	\\
$\Omega_{g,n}$
	& Sec.~\ref{subsubsec:GR}
	& Functions on $\mathcal{M}_{g,n}$ that are output of hyperbolic GR
	\\
$V\Xi_{g,n}(L)$, $N\Xi_{g,n}(L)$
	& Eqn.~\eqref{defn:Vamplitudes}, \eqref{eqn:lattice:pnt:count}
	& The integrals $\int_{\mathcal{M}_{g,n}^{\textup{comb}}(L)} \Xi_{g,n} \, \dd\mu_{\textup{K}}$, the sums $\sum_{\bm{G} \in \mathcal{M}_{g,n}^{\textup{comb},\ZZ}(L)} \frac{\Xi_{g,n}(\bm{G})}{\#\Aut(\bm{G})}$
	\\
$V\Omega_{g,n}(L)$
	& Eqn.~\eqref{defn:Vamplitudes}
	& The integrals $\int_{\mathcal{M}_{g,n}(L)} \Omega_{g,n} \, \dd\mu_{\textup{WP}}$
	\\
$\mathscr{B}_{\Sigma}$, $\mathscr{B}_{\Sigma}^{\textup{comb}}$
	& Def.~\ref{defn:Mirz:function}, \ref{defn:Mirz:comb:function}
	& Mirzakhani function in the hyperbolic and combinatorial setting
	\\
$\mu_{\textup{MV}}, MV_{g,n}$
	& Sec.~\ref{subsec:MV:formulae}
	& The Masur--Veech measure, the Masur--Veech volumes
\end{tabulary}
\end{center}
}

\newpage


\providecommand{\bysame}{\leavevmode\hbox to3em{\hrulefill}\thinspace}
\providecommand{\MR}{\relax\ifhmode\unskip\space\fi MR }
\providecommand{\MRhref}[2]{%
  \href{http://www.ams.org/mathscinet-getitem?mr=#1}{#2}
}
\providecommand{\href}[2]{#2}

\end{document}